\documentclass[11pt]{article}

\usepackage{amsmath,amssymb,amsfonts,enumitem}
\usepackage{pstricks,pst-node}

\usepackage{xr}
\numberwithin{equation}{section}
\newtheorem{thm}{Theorem}[section] 
\newtheorem{prp}[thm]{Proposition}
\newtheorem{lmm}[thm]{Lemma}   
\newtheorem{crl}[thm]{Corollary} 
\newtheorem{dfn}[thm]{Definition} 

\newcounter{stat}

\renewcommand{\frak}{\mathfrak}
\renewcommand{\Bbb}{\mathbb}

\def\C{\mathbb C}
\def\cC{\mathcal C}
\def\D{\frak D}
\def\cD{\mathcal D}
\def\E{\mathbb E}
\def\F{\mathfrak F}
\def\cI{\mathcal I}
\def\I{\mathfrak i}
\def\J{\mathcal J}
\def\M{\mathfrak M}
\def\cM{\mathcal M}
\def\O{\mathcal O}
\def\bP{\mathbb P}
\def\Q{\mathbb Q}
\def\R{\mathbb R}
\def\T{\mathcal T}
\def\U{\mathcal U}
\def\X{\mathfrak X}
\def\Z{\mathbb Z}

\def\e_ref#1{(\ref{#1})}
\def\under#1{\underline{#1}}
\def\ov#1{\overline{#1}}
\def\wt#1{\widetilde{#1}}
\def\ti#1{\tilde{#1}}

\def\lan{\langle}
\def\ran{\rangle}
\def\lr#1{\lan#1\ran}
\def\blr#1{\big\lan#1\big\ran}
\def\llrr#1{\lan\!\lan#1\ran\!\ran}
\def\bllrr#1{\big\lan\!\big\lan#1\big\ran\!\big\ran}
\def\lra{\longrightarrow}
\def\Lra{\Longrightarrow}
\def\Llra{\Longleftrightarrow}

\def\al{\alpha}
\def\be{\beta}
\def\de{\delta}
\def\ep{\epsilon}
\def\ga{\gamma}
\def\io{\iota}
\def\ka{\kappa}
\def\na{\nabla}
\def\om{\omega}
\def\si{\sigma}
\def\th{\theta}
\def\ups{\upsilon}
\def\ze{\zeta}
\def\ve{\varepsilon}
\def\vph{\varphi}
 \def\vsi{\varsigma}
\def\vt{\vartheta}

\def\Ga{\Gamma}
\def\La{\Lambda}
\def\Om{\Omega}
\def\Si{\Sigma}
\def\Th{\Theta}

\def\Aut{\textnormal{Aut}}
\def\const{\textnormal{const}}
\def\codim{\textnormal{codim}}
\def\ev{\textnormal{ev}}
\def\GW{\textnormal{GW}}
\def\id{\textnormal{id}}
\def\ind{\textnormal{ind}~\!}
\def\End{\textnormal{End}}
\def\Hol{\textnormal{Hol}}
\def\Im{\textnormal{Im}~\!}
\def\reg{\textnormal{reg}}
\def\simp{\textnormal{simp}}

\def\P{\Bbb{P}^n}
\def\PP{\Bbb{P}^2}
\def\PPP{\Bbb{P}^3}

\def\i{\infty}
\def\eset{\emptyset}

\begin{document}

\title{A Sharp Compactness Theorem for\linebreak
Genus-One Pseudo-Holomorphic Maps}

\author{Aleksey Zinger\thanks{Partially supported by an NSF Postdoctoral Fellowship}}

\date{\today}
\maketitle

\begin{abstract}
\noindent
For each compact almost Kahler manifold $(X,\om,J)$ and an element $A$ of $H_2(X;\Z)$,
we describe a natural closed subspace $\ov\M_{1,k}^0(X,A;J)$ of 
the moduli space $\ov\M_{1,k}(X,A;J)$ of stable $J$-holomorphic genus-one maps
such that $\ov\M_{1,k}^0(X,A;J)$ contains all stable maps with smooth domains.
If $(\P,\om,J_0)$ is the standard complex projective space,
$\ov\M_{1,k}^0(\P,A;J_0)$ is an irreducible component of~$\ov\M_{1,k}(\P,A;J_0)$. 
We also show that if an almost complex structure~$J$ on~$\P$ is sufficiently close
to~$J_0$, the structure of the space $\ov\M_{1,k}^0(\P,A;J)$
is similar to that of $\ov\M_{1,k}^0(\P,A;J_0)$.
This paper's compactness and structure theorems lead to new invariants
for some symplectic manifolds, which are generalized to arbitrary symplectic 
manifolds in a separate paper.
Relatedly, the smaller moduli space $\ov\M_{1,k}^0(X,A;J)$ is useful
for computing the genus-one Gromov-Witten invariants, 
which arise from the larger moduli space $\ov\M_{1,k}(X,A;J)$.
\end{abstract}

\thispagestyle{empty}

\tableofcontents

\section{Introduction}
\label{intro_sec}

\subsection{Background and Motivation}
\label{back_subs}

\noindent
Gromov-Witten invariants of symplectic manifolds have been a subject of 
much research in the past decade, 
as they play a prominent role in both symplectic topology and 
theoretical physics.
In~order to define GW-invariants of a compact symplectic manifold~$(X,\om)$, 
one fixes an almost complex structure~$J$ on~$X$,
which is compatible with~$\om$ or at least is tamed by~$\om$.
For each class $A$ in $H_2(X;\Z)$ and a pair~$(g,k)$ of nonnegative integers,
let $\ov\M_{g,k}(X,A;J)$ be the moduli space of (equivalence classes of) 
stable $J$-holomorphic maps from genus-$g$ Riemann surfaces with $k$ marked points in 
the homology class~$A$.
The expected, or virtual, dimension of this moduli space is given~by
$$\dim_{g,k}(X,A)
\equiv\dim^{vir}\,\ov\M_{g,k}(X,A;J)
=2\big( \lr{c_1(TX),A}+(n\!-\!3)(1\!-\!g)+k\big),$$
if the real dimension of $X$ is $2n$.
While in general $\ov\M_{g,k}(X,A;J)$ is not a smooth manifold, or even a variety,
of dimension $\dim_{g,k}(X,A)$,
it is shown in \cite{FuOn}, \cite{LT}, and in the algebraic case in~\cite{BeFa}
that $\ov\M_{g,k}(X,A;J)$ determines a rational homology class of
dimension~$\dim_{g,k}(X,A)$.
In turn, this virtual fundamental class of $\ov\M_{g,k}(X,A;J)$
is used to define GW-invariants of~$(X,\om)$.\\

\noindent
We denote by $\M_{g,k}^0(X,A;J)$ the subspace of $\ov\M_{g,k}(X,A;J)$
consisting of the stable maps $[{\cal C},u]$ such that the domain ${\cal C}$
is a smooth Riemann surface.
If $(\P,\om;J_0)$ is the $n$-dimensional complex projective space with 
the standard Kahler structure and $\ell$ is the homology class of 
a complex line in~$\P$,
$$\M_{g,k}^0(\P,d)\equiv\M_{g,k}^0(\P,d\ell;J_0)$$
is in fact a smooth orbifold of dimension~$\dim_{g,k}(\P,d\ell)$,
at least for $d\!\ge\!2g\!-\!1$.
In addition, from the point of view of algebraic geometry,
$\ov\M_{0,k}(\P,d)$ is an irreducible algebraic orbivariety of 
dimension $\dim_{0,k}(\P,d\ell)/2$.
From the point of view of symplectic topology, 
$\ov\M_{0,k}(\P,d)$ is a compact topological orbifold
stratified by smooth orbifolds of even dimensions and 
$\M_{0,k}^0(\P,d)$ is the main stratum of~$\ov\M_{0,k}(\P,d)$.
In particular, $\M_{0,k}^0(\P,d)$ is a dense open subset 
of~$\ov\M_{0,k}(\P,d)$.\\

\begin{figure}
\begin{pspicture}(-1.1,-1.8)(10,1.25)
\psset{unit=.4cm}
\rput{45}(0,-4){\psellipse(5,-1.5)(2.5,1.5)
\psarc[linewidth=.05](5,-3.3){2}{60}{120}\psarc[linewidth=.05](5,0.3){2}{240}{300}
\pscircle[fillstyle=solid,fillcolor=gray](5,-4){1}\pscircle*(5,-3){.2}
\pscircle[fillstyle=solid,fillcolor=gray](6.83,.65){1}\pscircle*(6.44,-.28){.2}
\pscircle[fillstyle=solid,fillcolor=gray](3.17,.65){1}\pscircle*(3.56,-.28){.2}}
\rput(.2,-.9){$d_1$}\rput(3.1,2.3){$d_2$}\rput(7.8,-2.5){$d_3$}
\psarc(15,-1){3}{-60}{60}\psline(17,-1)(22,-1)\psline(16.8,-2)(21,-3)\psline(16.8,0)(21,1)
\rput(15.2,-3.5){$(1,0)$}\rput(22.4,1){$(0,d_1)$}
\rput(23.4,-1){$(0,d_2)$}\rput(22.4,-3){$(0,d_3)$}
\rput(33,-1){\begin{tabular}{l}$d_1\!+\!d_2\!+\!d_3\!=\!d$\\
$d_1,d_2,d_3\!>\!0$\end{tabular}}
\end{pspicture}
\caption{The Domain of an Element of $\M_{1,k}^3(\P,d)$}
\label{m3_fig}
\end{figure}
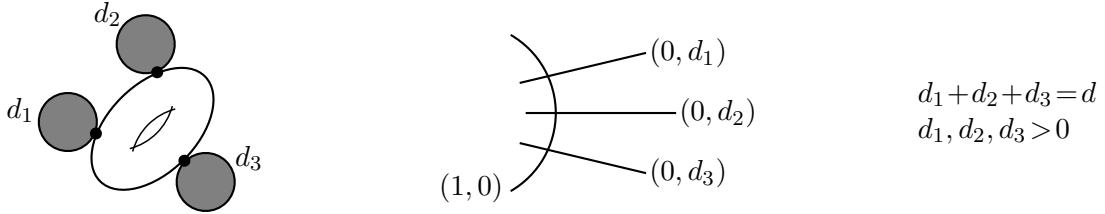

\noindent
If $g\!\ge\!1$, none of these additional properties holds even for~$(\P,\om,J_0)$.
For example, the moduli space $\ov\M_{1,k}(\P,d)$ has many irreducible
components of various dimensions.
One of these components contains $\M_{1,k}^0(\P,d)$;
we denote this component by~$\ov\M_{1,k}^0(\P,d)$.
In other words, $\ov\M_{1,k}^0(\P,d)$ is the closure 
of $\M_{1,k}^0(\P,d)$ in $\ov\M_{1,k}(\P,d)$.
The remaining components of $\ov\M_{1,k}(\P,d)$
can be described as follows.
If $m$ is a positive integer, let $\M_{1,k}^m(\P,d)$ be the subset
of $\ov\M_{1,k}(\P,d)$ consisting of the stable maps $[{\cal C},u]$
such that ${\cal C}$ is an elliptic curve~$E$ with $m$~rational 
components attached directly to~$E$, $u|_E$ is constant, and
the restriction of~$u$ to each rational component is non-constant.
Figure~\ref{m3_fig} shows the domain of an element of $\M_{1,k}^3(\P,d)$,
from the points of view of symplectic topology and of algebraic geometry.
In the first diagram, each shaded disc represents a sphere;
the integer next to each rational component ${\cal C}_i$ indicates the degree 
of~$u|_{{\cal C}_i}$.
In the second diagram, the components of ${\cal C}$ are represented by curves,
and the pair of integers next to each component ${\cal C}_i$ shows 
the genus of ${\cal C}_i$ and the degree of~$u|_{{\cal C}_i}$.
We denote by $\ov\M_{1,k}^m(\P,d)$ the closure of $\M_{1,k}^m(\P,d)$
in $\ov\M_{1,k}(\P,d)$.
The space $\ov\M_{1,k}^m(\P,d)$ has a number of irreducible components.
These components are indexed by the splittings of the degree~$d$
into $m$ positive integers and by the distributions of the $k$~marked points
between the $m\!+\!1$ components of the domain.
However, all of these components are algebraic orbivarieties of dimension,
both expected and actual,
\begin{equation*}\begin{split}
\dim_{1,k}^m(\P,d\ell)\equiv \dim\ov\M_{1,k}^m(\P,d)
&= 2\big( d(n\!+\!1)+k+n-m\big)\\
&=\dim_{1,k}(\P,d\ell)+2(n\!-\!m).
\end{split}\end{equation*}
In particular, $\M_{1,k}^0(\P,d)$ is not dense in $\ov\M_{1,k}(\P,d)$.
From the point of view of symplectic topology, 
$\ov\M_{1,k}(\P,d)$ is a union of compact topological orbifolds
and is stratified by smooth orbifolds of even dimensions.
However, $\ov\M_{1,k}(\P,d)$ contains several main strata,
and some of them are of dimension larger than~$\dim_{1,k}(\P,d\ell)$.\\

\noindent
The above example shows that $\ov\M_{1,k}^0(\P,d)$
is a true compactification of the moduli space $\M_{1,k}^0(\P,d)$,
while $\ov\M_{1,k}(\P,d)$ is simply a compact space 
containing~$\M_{1,k}^0(\P,d)$, albeit one with a nice obstruction theory.
One can view $\ov\M_{1,k}(\P,d)$ as 
a {\it geometric}-genus compactification of $\M_{1,k}^0(\P,d)$
and its subspace $\ov\M_{1,k}^0(\P,d)$ as 
an {\it arithmetic}-genus compactification.
Since the beginning of the Gromov-Witten theory,
it has been believed, or at least considered feasible, 
that an analogue of $\ov\M_{1,k}^0(\P,d)$ can be defined for 
every compact almost Kahler manifold $(X,\om,J)$, positive genus~$g$, 
and nonzero homology class $A\!\in\!H_2(X;\Z)$.
In this paper, we show that this is indeed the case if $g\!=\!1$.\\

\noindent
We describe an analogue $\ov\M_{1,k}^0(X,A;J)$ of the subspace 
$\ov\M_{1,k}^0(\P,d)$ of $\ov\M_{1,k}(\P,d)$ for every compact 
almost Kahler manifold $(X,\om,J)$ 
and homology class $A\!\in\!H_2(X;\Z)$ as the subset of elements of $\ov\M_{1,k}(X,A;J)$
that satisfy one of two conditions.
By Theorem~\ref{comp_thm}, $\ov\M_{1,k}^0(X,A;J)$ 
is a closed subspace of $\ov\M_{1,k}(X,A;J)$ and thus is compact.
This compactification of $\M_{1,k}^0(X,A;J)$ satisfies 
the following desirable properties:
\begin{enumerate}[label=(P\arabic*)]
\item\label{subman_item}
{\it naturality with respect to embeddings:}
if $(Y,\om,J)$ is a compact submanifold of $(X,\om,J)$, then
$$\ov\M_{1,k}^0(Y,A;J) = \ov\M_{1,k}^0(X,A;J)\cap \ov\M_{1,k}(Y,A;J)
\subset\ov\M_{1,k}(X,A;J);$$
\item\label{forget_item}
{\it naturality with respect to forgetful maps}:
if $k\!\ge\!1$, the pre-image of $\ov\M_{1,k-1}^0(X,A;J)$ under the forgetful map
$$\ov\M_{1,k}(X,A;J)\lra \ov\M_{1,k-1}(X,A;J)$$
is $\ov\M_{1,k}^0(X,A;J)$;
\item\label{sharp_item} 
{\it sharpness for regular $(X,\om,J)$:}
if $J$ satisfies the regularity conditions of Definition~\ref{g1reg_dfn},
then $\ov\M_{1,k}^0(X,A;J)$ is the closure of $\M_{1,k}^0(X,A;J)$
in $\ov\M_{1,k}(X,A;J)$.
\end{enumerate}
By~\ref{subman_item} and~\ref{forget_item}, $\ov\M_{1,k}^0(X,A;J)$, like 
$\ov\M_{1,k}(X,A;J)$, is a natural compactification of $\M_{1,k}^0(X,A;J)$.
By~\ref{sharp_item}, $\ov\M_{1,k}^0(X,A;J)$, in contrast to $\ov\M_{1,k}(X,A;J)$, 
is a {\it sharp} compactification of $\M_{1,k}^0(X,A;J)$,
subject to the naturality conditions~\ref{subman_item} and~\ref{forget_item}.
The first two properties of $\ov\M_{1,k}^0(X,A;J)$ are immediate
from Definition~\ref{degen_dfn}.
The last property is part of Corollary~\ref{str_crl}.
It is well-known that the regularity conditions of Definition~\ref{g1reg_dfn}
are satisfied by the standard complex structure~$J_0$ on~$\P$, 
and thus the definition of $\ov\M_{1,k}^0(\P,d\ell;J_0)$ 
given in Subsection~\ref{comp_subs} agrees with 
the description of $\ov\M_{1,k}^0(\P,d)$ given above.\\

\noindent
Theorem~\ref{str_thm} describes, under the regularity conditions of 
Definition~\ref{degen_dfn}, a neighborhood of every ``interesting" stratum 
of $\ov\M_{1,k}^0(X,A;J)$, i.e.~a stratum consisting of genus-one maps 
that are constant on the principal component.
In addition to implying~\ref{sharp_item}, Theorem~\ref{str_thm} shows that 
$\ov\M_{1,k}^0(X,A;J)$ carries a rational fundamental class 
and can be used to define Gromov-Witten style intersection numbers via pseudocycles,
as in Chapter~5 of \cite{McSa} or Section~1 of~\cite{RT1},
whenever $J$ is regular; see Subsection~\ref{geomcons_subs} below.
As the regularity requirements of Definition~\ref{degen_dfn} are open
conditions on the space of $\om$-tame almost complex structures~$J$
by Theorem~\ref{reg_thm},
Theorem~\ref{str_thm} also implies that the general topological structure 
of $\ov\M_{1,k}^0(X,A;J)$ remains unchanged under small changes in $J$ 
near a regular~$J_0$.\\

\noindent
The results of this paper have already found a variety of applications:
\begin{enumerate}[label=(A\arabic*)]
\item\label{redinv_item} $\ov\M_{1,k}^0(X,A;J)$ gives rise to 
new, {\tt reduced}, genus-one GW-invariants of arbitrary symplectic 
manifolds~(\cite{Z6});
\item\label{HP_item} in contrast to the standard genus-one GW-invariants,
the reduced invariants of a complete intersection and the ambient space
are related as geometrically expected (\cite{LZ},\cite{Z7});
\item\label{desing_item} Theorem~\ref{str_thm} is used in \cite{VaZ}
to construct a natural desingularization of  $\ov\M_{1,k}^0(\P,d)$
and thus a natural smooth compactification of the Hilbert scheme
of smooth genus-one curves in~$\bP^n$ for $n\!\ge\!3$;
\item\label{BCOV_item} \ref{redinv_item}-\ref{desing_item}
are used in~\cite{Z8} to finally confirm the 1993 prediction of~\cite{BCOV}
for genus-one GW-invariants of a quintic threefold;
\item\label{CI_item} \ref{redinv_item}-\ref{desing_item},
along with~\cite{Z9}, have made it possible to compute 
(standard) genus-one GW-invariants  of arbitrary complete intersections.\\
\end{enumerate}

\begin{figure}
\begin{pspicture}(-1.1,-2)(10,1.25)
\psset{unit=.4cm}
\psellipse[fillstyle=solid,fillcolor=gray](6,2)(1.2,.6)
\psarc[linewidth=.06](4.38,2){.42}{270}{90}
\psarc[linewidth=.06](7.62,2){.42}{90}{270}
\psarc[linewidth=.06](4.38,1.16){.42}{90}{235}
\psarc[linewidth=.06](7.62,1.16){.42}{315}{90}
\psline[linewidth=.06](4.08,.86)(5.8,-.86)
\psline[linewidth=.06](7.92,.86)(6.2,-.86)
\psarc[linewidth=.06](5.32,-1.34){.68}{0}{45}
\psarc[linewidth=.06](6.68,-1.34){.68}{135}{180}
\psarc[linewidth=.06](5,-1.34){1}{270}{0}
\psarc[linewidth=.06](7,-1.34){1}{180}{270}
\psarc[linewidth=.06](5,-2.84){.5}{90}{270}
\psarc[linewidth=.06](7,-2.84){.5}{270}{90}
\psarc[linewidth=.06](5,-4.34){1}{270}{90}
\psarc[linewidth=.06](7,-4.34){1}{90}{270}
\psline[linewidth=.06](5,-5.34)(4.38,-5.34)
\psline[linewidth=.06](7,-5.34)(7.62,-5.34)
\psarc[linewidth=.06](4.38,.92){1.5}{90}{180}
\psarc[linewidth=.06](7.62,.92){1.5}{0}{90}
\psline[linewidth=.06](2.88,.92)(2.88,-3.84)
\psline[linewidth=.06](9.12,.92)(9.12,-3.84)
\psarc[linewidth=.06](4.38,-3.84){1.5}{180}{270}
\psarc[linewidth=.06](7.62,-3.84){1.5}{270}{0}
\pscircle*(4.8,2){.2}\pscircle*(7.2,2){.2}
\pscircle*(6,-1.34){.2}\pscircle*(6,-4.34){.2}
\rput(4,-1.3){${\cal C}_1$}\rput(8,-1.3){${\cal C}_2$}
\rput(2.4,-3){$0$}\rput(9.6,-3){$0$}
\rput(5.3,3.2){${\cal C}_3$}\rput(6.6,3.1){$d$}
\psline[linewidth=.06](26,2)(30,2)
\psarc[linewidth=.06](28,2){1.41}{160}{225}
\psarc[linewidth=.06](28,2){1.41}{315}{20}
\psline(27,1)(29,-1)\psline(29,1)(27,-1)
\psarc[linewidth=.06](28,-2){1.41}{135}{225}
\psarc[linewidth=.06](28,-2){1.41}{315}{45}
\psline(27,-3)(28.5,-4.5)\psline(29,-3)(27.5,-4.5)
\rput(26.3,-.7){${\cal C}_1$}\rput(29.7,-.7){${\cal C}_2$}
\pscircle*(26.59,2){.2}\pscircle*(29.41,2){.2}
\pscircle*(28,0){.2}\pscircle*(28,-4){.2}
\rput(25.6,-3.2){$(0,0)$}\rput(30.4,-3.2){$(0,0)$}
\rput(25.4,1.9){${\cal C}_3$}\rput(31.2,2){$(0,d)$}
\rput(17,-1){``tacnode"}
\pnode(14.6,-1.1){A1}\pnode(19.3,-1.1){A2}
\pnode(6,0.2){B1}\pnode(28,0.7){B2}
\ncarc[linewidth=.05,nodesep=0,arcangleA=10,arcangleB=10,ncurv=.7]{-}{A1}{B1}
\ncarc[linewidth=.05,nodesep=0,arcangleA=-10,arcangleB=-10,ncurv=.8]{-}{A2}{B2}
\pnode(4.9,1.6){B1a}\pnode(7.1,1.6){B1b}
\ncarc[linewidth=.05,nodesep=0,arcangleA=10,arcangleB=10,ncurv=.7]{->}{B1}{B1a}
\ncarc[linewidth=.05,nodesep=0,arcangleA=-10,arcangleB=-10,ncurv=.7]{->}{B1}{B1b}
\pnode(26.8,1.8){B2a}\pnode(29.2,1.8){B2b}
\ncarc[linewidth=.05,nodesep=0,arcangleA=10,arcangleB=10,ncurv=.7]{->}{B2}{B2a}
\ncarc[linewidth=.05,nodesep=0,arcangleA=-10,arcangleB=-10,ncurv=.7]{->}{B2}{B2b}
\end{pspicture}
\caption{A Condition on Limits in Genus Two}
\label{g2_fig}
\end{figure}
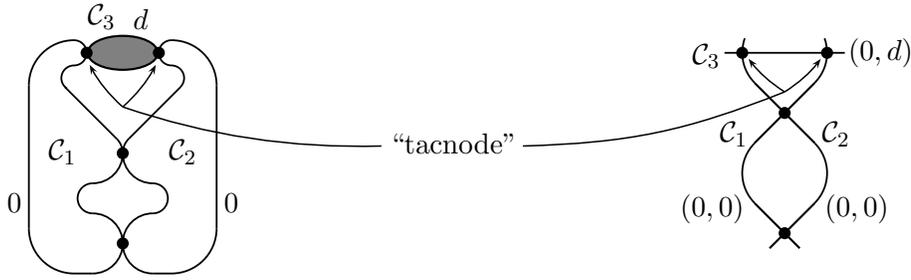

\noindent
If it is possible to define subspaces $\ov\M_{g,k}^0(X,A;J)$ of 
$\ov\M_{g,k}(X,A;J)$ analogous to the space                     \linebreak
$\ov\M_{1,k}^0(X,A;J)$ for $g\!\ge\!2$,
their description is likely to be  more complicated.
The space                                                         
$\ov\M_{1,k}^0(X,A;J)$ contains all stable maps
$[{\cal C},u]$ in $\ov\M_{1,k}(X,A;J)$ such that the restriction
of $u$ to the principal component~$\cC_P$ is nonconstant 
or such that $u|_{\cC_P}$ is constant and the restrictions
to the rational components satisfy a certain fairly simple degeneracy condition;
see Definition~\ref{degen_dfn}.
Thus, in the genus-one case the elements in $\ov\M_{1,k}(X,A;J)$
are split into two classes, according to their restriction to the principal component.
In the genus-two case, these classes would need to be split further.
For example, suppose the domain of an element 
$[{\cal C},u]$ of $\ov\M_{2,k}(\P,d)$ consists of three rational curves,
${\cal C}_1$, ${\cal C}_2$, and ${\cal C}_3$, such that ${\cal C}_1$ and ${\cal C}_2$
share two nodes and ${\cal C}_3$ has a node in common with ${\cal C}_1$ and ${\cal C}_2$;
see Figure~\ref{g2_fig}.
If $u|_{{\cal C}_1}$ and $u|_{{\cal C}_2}$ are constant,
$[{\cal C},u]$ lies in the closure of 
$\M_{2,k}^0(\P,d)$ in $\ov\M_{2,k}(\P,d)$ if and only if
the branches of the curve $u({\cal C})\!=\!u({\cal C}_3)$ corresponding
to the two nodes of~${\cal C}_3$ form a generalized tacnode,
i.e.~either one of them is a cusp or the two branches have the same tangent line;
see~\cite{Z2} for the $n\!=\!2$ case.\\

\noindent
The author would like to thank J.~Li for suggesting the problem of
computing the genus-one GW-invariants of a quintic threefold,
which led to the present paper.
The author first learned of the arithmetic/geometric-genus compactification
terminology in the context of stable maps from G.~Tian a number of years ago.

\subsection{Compactness Theorem}
\label{comp_subs}

\noindent
In this subsection, we describe the subspace $\ov\M_{1,k}^0(X,A;J)$
of $\ov\M_{1,k}(X,A;J)$; it is a closed subspace by Theorem~\ref{comp_thm}.
We specify what we mean by a regular almost structure~$J$
in Definitions~\ref{g0reg_dfn} and~\ref{g1reg_dfn}.
If $J$ is genus-one $A$-regular, the moduli space $\ov\M_{1,k}^0(X,A;J)$
has a regular structure, which is described by  Theorem~\ref{str_thm}.
Since the rather detailed statement of this theorem is
notationally involved, we postpone stating it until after 
we introduce additional notation in Subsections~\ref{notation0_subs}
and~\ref{notation1_subs}.
In this subsection, we instead state Corollary~\ref{str_crl},
which describes the two most important consequences of Theorem~\ref{str_thm}.\\

\noindent
An element $[{\cal C},u]$ of $\ov\M_{1,k}(X,A;J)$ is the equivalence class
of a pair consisting of a prestable genus-one complex curve~${\cal C}$
and a $J$-holomorphic map $u\!:{\cal C}\!\lra\!X$.
The prestable curve~${\cal C}$ is a union of the principal curve~$\cC_P$,
which is either a smooth torus or a circle of spheres, and
trees of rational bubble components, which together will be denoted by~$\cC_B$.
Let
$$\M_{1,k}^{\{0\}}(X,A;J)= \big\{ [{\cal C},u]\!\in\!\ov\M_{1,k}(X,A;J)\!:
u|_{\cC_P}~\hbox{is not constant}\big\} \supset\M_{1,k}^0(X,A;J).$$
The space $\M_{1,k}^{\{0\}}(X,A;J)$ will be a subset 
of the moduli space $\ov\M_{1,k}^{0}(X,A;J)$.\\

\noindent
Every bubble component $\cC_i\!\subset\!\cC_B$ is a sphere and has a distinguished 
singular point, which will be called the {\tt attaching node of~$\cC_i$}.
This is the node of $\cC_i$ that lies either on~$\cC_P$
or on a bubble $\cC_h$ that lies between $\cC_i$ and~$\cC_P$.
For example, if $\cC$ is as shown in Figure~\ref{chi_fig},
the attaching node of $\cC_{h_3}$ is the node $\cC_{h_3}$ shares with the torus.
Since ${\cal C}_i$ is a sphere, we can represent every element of
$\ov\M_{1,k}(X,A;J)$ by a pair $({\cal C},u)$ such that
the attaching node of every bubble component $\cC_i\!\subset\!\cC_B$
is the south pole, or the point $\i\!=\!(0,0,-1)$, of $S^2\!\subset\!\Bbb{R}^3$.
Let $e_{\i}\!=\!(1,0,0)$ be a nonzero tangent vector to $S^2$ at the south pole.
Then the vector  
$${\cal D}_i({\cal C},u)\equiv  
d\big\{u|_{{\cal C}_i}\big\}\big|_{\i}e_{\i} \in T_{u|_{{\cal C}_i}(\i)}X$$
describes the differential of the $J$-holomorphic map $u|_{{\cal C}_i}$
at the attaching node.
While this element of $T_{u|_{{\cal C}_i}(\i)}X$ depends on the choice of
a representative for an element of $\ov\M_{1,k}(X,A;J)$,
the linear subspace $\Bbb{C}\!\cdot\!{\cal D}_i({\cal C},u)$
of $T_{u|_{{\cal C}_i}(\i)}X$ is determined by the equivalence class~$[{\cal C},u]$.
If $u|_{{\cal C}_i}$ is not constant,
the branch of the rational $J$-holomorphic curve $u({\cal C}_i)\!\subset\!X$
corresponding to the attaching node of ${\cal C}_i$ has a cusp if
and only if ${\cal D}_i({\cal C},u)\!=\!0$.
If ${\cal D}_i({\cal C},u)\!\neq\!0$, 
$\Bbb{C}\!\cdot\!{\cal D}_i({\cal C},u)$ is the line tangent to
the  branch of $u({\cal C}_i)\!\subset\!X$
corresponding to the attaching node of~${\cal C}_i$.\\ 

\noindent
Suppose $[{\cal C},u]\!\in\!\ov\M_{1,k}(X,A;J)\!-\!\M_{1,k}^{\{0\}}(X,A;J)$,
i.e.~$u|_{\cC_P}$ is constant.
In such a case, we will call the bubble sphere $\cC_i\!\subset\!\cC_B$
{\tt first-level $({\cal C},u)$-effective} if $u|_{\cC_i}$ is not constant,
but $u|_{\cC_h}$ is constant for every bubble component 
$\cC_h\!\subset\!\cC_B$ that lies between $\cC_i$ and~$\cC_P$.
We denote by $\chi({\cal C},u)$ the set of  first-level 
$({\cal C},u)$-effective bubbles; see Figure~\ref{chi_fig}.
In this figure, as in Figures~\ref{m3_fig} and~\ref{g2_fig},
we show the domain~${\cal C}$ of the stable map $({\cal C},u)$
and shade the components of the domain on which the map~$u$ is not constant.
Note that $u$ maps the attaching nodes of all elements of $\chi({\cal C},u)$
to the same point in~$X$.

\begin{figure}
\begin{pspicture}(-1.1,-1.8)(10,1.25)
\psset{unit=.4cm}
\psellipse(8,-1.5)(1.5,2.5)
\psarc[linewidth=.05](6.2,-1.5){2}{-30}{30}\psarc[linewidth=.05](9.8,-1.5){2}{150}{210}
\pscircle[fillstyle=solid,fillcolor=gray](5.5,-1.5){1}\pscircle*(6.5,-1.5){.2}
\pscircle[fillstyle=solid,fillcolor=gray](3.5,-1.5){1}\pscircle*(4.5,-1.5){.2}
\pscircle(10.5,-1.5){1}\pscircle*(9.5,-1.5){.2}
\pscircle[fillstyle=solid,fillcolor=gray](11.91,-.09){1}\pscircle*(11.21,-.79){.2}
\pscircle[fillstyle=solid,fillcolor=gray](11.91,-2.91){1}\pscircle*(11.21,-2.21){.2}
\rput(5.5,0){$h_1$}\rput(3.5,0){$h_2$}\rput(10.3,0){$h_3$}
\rput(13.5,0.1){$h_4$}\rput(13.5,-2.9){$h_5$}
\rput(8,-5){\small ``tacnode"}
\pnode(8,-5){A1}\pnode(6.5,-1.5){B1}
\ncarc[nodesep=.35,arcangleA=-25,arcangleB=-15,ncurv=1]{->}{A1}{B1}
\pnode(8,-4.65){A2}\pnode(10.3,-1.5){B2}
\ncarc[nodesep=0,arcangleA=40,arcangleB=30,ncurv=1]{-}{A2}{B2}
\pnode(11,-.95){B2a}\pnode(11.02,-2.02){B2b}
\ncarc[nodesep=0,arcangleA=0,arcangleB=10,ncurv=1]{->}{B2}{B2a}
\ncarc[nodesep=0,arcangleA=0,arcangleB=10,ncurv=1]{->}{B2}{B2b}
\rput(25,-1.5){$\chi({\cal C},u)\!=\!\{h_1,h_4,h_5\}$}
\end{pspicture}
\caption{An Illustration of Definition~\ref{degen_dfn}} 
\label{chi_fig}
\end{figure}

\begin{dfn}
\label{degen_dfn}
If $(X,\om,J)$ is a compact almost Kahler manifold, $A\!\in\!H_2(X;\Z)^*$, and 
$k\!\in\!\bar{\Z}^+$, the {\tt main component} of the space 
$\ov\M_{1,k}(X,A;J)$ is the subset $\ov\M_{1,k}^0(X,A;J)$ 
consisting of the elements $[{\cal C},u]$ of $\ov\M_{1,k}(X,A;J)$
such~that\\
${}\quad$ (a) $u|_{\cC_P}$ is not constant, or\\
${}\quad$ (b) $u|_{\cC_P}$ is constant and 
$\dim_{\Bbb{C}}\text{Span}_{(\Bbb{C},J)}\{{\cal D}_i({\cal C},u)\!:
i\!\in\!\chi({\cal C},u)\}<|\chi({\cal C},u)|$.
\end{dfn}

\noindent
We now clarify this definition.
We denote $H_2(X;\Z)\!-\!\{0\}$ by $H_2(X;\Z)^*$ and
the set of nonnegative integers by~$\bar{\Z}^+$.
We call a triple $(X,\om,J)$ an {\tt almost Kahler manifold}
if $\om$ is a symplectic form on $X$ and $J$ is an almost complex structure on $X$,
which is tamed by~$\om$, i.e.
$$\om(v,Jv)>0\qquad\forall~v\in TX\!-\!X.$$
Definition~\ref{degen_dfn} actually involves only the almost complex structure $J$,
but one typically considers the moduli spaces $\ov\M_{g,k}(X,A;J)$ 
only for $\om$-tamed almost  complex structures $J$, for some  symplectic form~$\om$;
otherwise, $\ov\M_{g,k}(X,A;J)$ may not be compact.
An element 
$$[{\cal C},u]\in\ov\M_{1,k}(X,A;J)\!-\!\M_{1,k}^{\{0\}}(X,A;J)$$
belongs to $\ov\M_{1,k}^0(X,A;J)$
if and only if the branches of $u({\cal C})$ corresponding to the attaching nodes
of the first-level effective bubbles of $[{\cal C},u]$ form 
a {\it generalized tacnode}.
In the case of Figure~\ref{chi_fig}, this means that  either\\
${}\quad$ (a) for some $i\!\in\!\{h_1,h_4,h_5\}$, 
the branch of $u|_{{\cal C}_i}$ at the attaching node of ${\cal C}_i$ has a cusp, or\\
${}\quad$ (b) for all $i\!\in\!\{h_1,h_4,h_5\}$,
the branch of $u|_{{\cal C}_i}$ at the attaching node of ${\cal C}_i$ is smooth, but the\\
${}\qquad~~$ dimension of the span of the three lines tangent to these
branches is less than three.

\begin{thm}
\label{comp_thm}
If $(X,\om)$ is a compact symplectic manifold, $\under{J}\!\equiv\!(J_t)_{t\in[0,1]}$
is a $C^1$-continuous family of $\om$-tamed almost complex structures on $X$,
$A\!\in\!H_2(X;\Z)^*$, and $k\!\in\!\bar{\Z}^+$, then the moduli space
$$\ov\M_{1,k}^0(X,A;\under{J})\equiv\!
\bigcup_{t\in[0,1]}\!\!\ov\M_{1,k}^0(X,A;J_t)$$
is compact.
\end{thm}

\noindent
If $(X,J)$ is an algebraic variety, the claim of Theorem~\ref{comp_thm},
with $J_t\!=\!J$ constant,
is an immediate consequence of well-known results in algebraic geometry.
In fact, Lemma~2.4.1 in~\cite{Va} can be used to generalize the statement
of Theorem~\ref{comp_thm} to higher genera, provided $(X,J)$ is 
a complex algebraic surface.\\

\noindent
If $J_t\!=\!J$ is constant and genus-one $A$-regular in the sense of 
Definition~\ref{g1reg_dfn} below, Theorem~\ref{comp_thm} follows immediately from
the first statement of Theorem~\ref{str_thm}.
If $J_t$ is genus-one $A$-regular for all~$t$, but not necessarily constant,
Theorem~\ref{comp_thm} follows from the Gromov Compactness Theorem and Corollary~\ref{reg1_crl3}.
In Section~\ref{compthm_sec}, we combine the main ingredients of the proof of 
Theorem~\ref{str_thm} with the local setting of~\cite{LT} to obtain Theorem~\ref{comp_thm}
with $J_t\!=\!J$ constant for an arbitrary almost Kahler manifold.
The proof for a general family $\under{J}$ is similar
and is described in detail, in an even more general case,
in Section~5 of~\cite{Z6}.\\

\noindent
If $u\!:\cC\!\lra\!X$ is a smooth map from a Riemann surface and $A\!\in\!H_2(X;\Z)$, 
we write
$$u\le_{\om}A \qquad\hbox{if}\qquad
u_*[\cC]=A \quad\hbox{or}\quad
\lr{\om,u_*[\cC]}\!<\!\lr{\om,A}.$$

\begin{dfn}
\label{g0reg_dfn}
Suppose $(X,\om,J)$ is a compact almost Kahler manifold and $A\!\in\!H_2(X;\Z)$.
The almost complex structure $J$ is {\tt genus-zero $A$-regular} if
for every $J$-holomorphic map $u\!:\Bbb{P}^1\!\lra\!X$
such that $u\!\le_{\om}\!A$,\\
${}\quad$ (a) the linearization $D_{J,u}$ of 
the $\bar{\partial}_J$-operator at $u$ is surjective;\\
${}\quad$ (b) for all $z\!\in\!\Bbb{P}^1$, the map
$\D_{J,u}^z\!\!:\ker D_{J,u}\!\lra\!T_{u(z)}X$, $\D_{J,u}^z(\xi)\!=\!\xi(z)$, is onto.
\end{dfn}

\begin{dfn}
\label{g1reg_dfn}
Suppose $(X,\om,J)$ is a compact almost Kahler manifold and $A\!\in\!H_2(X;\Z)$.
The almost complex structure $J$ is {\tt genus-one $A$-regular} if\\
${}\quad$ (a) $J$ is genus-zero $A$-regular;\\
${}\quad$ (b) for every nonconstant $J$-holomorphic map $u\!:\Bbb{P}^1\!\lra\!X$
such that $u\!\le_{\om}\!A$,\\
${}\qquad$ (b-i) for all $z\!\in\!\Bbb{P}^1$ and $v\!\in\!T_z\Bbb{P}^1\!-\!\{0\}$, the map
$\D_{J,u}^{z,v}\!\!:\ker\D_{J,u}^z\!\lra\!T_{u(z)}X$, 
$\D_{J,u}^{z,v}(\xi)\!=\!\na_v\xi$,\\
${}\qquad\qquad~$ is onto;\\
${}\qquad$ (b-ii) for all $z\!\in\!\Bbb{P}^1$ and $z'\!\in\!\Bbb{P}^1\!-\!\{z\}$, the map
$\D_{J,u}^{z,z'}\!\!:\ker\D_{J,u}^z\!\lra\!T_{u(z')}X$, $\D_{J,u}^{z,z'}(\xi)\!=\!\xi(z')$, 
${}\qquad\qquad~$ is onto.\\
${}\quad$ (c) for every smooth genus-one Riemann surface $\Si$ and
every non-constant $J$-holomorphic map\\
${}\qquad~~$ $u\!:\Si\!\lra\!X$ such that $u\!\le_{\om}\!A$,
the linearization $D_{J,u}$ of the $\bar{\partial}_J$-operator
at $u$ is surjective.
\end{dfn}

\noindent
In (b-i) of Definition~\ref{g1reg_dfn},
$\na_v\xi$ denotes the covariant derivative of~$\xi$ along~$v$,
with respect to a connection $\na$ in~$TX$.
Since $\xi(z)\!=\!0$, the value of $\na_v\xi$ is in fact independent 
of the choice of~$\na$.
If $J$ is an integrable complex structure, 
the surjectivity statements of~(a) and~(b)
in Definition~\ref{g0reg_dfn} and of~(c) in Definition~\ref{g1reg_dfn}
can be written~as
$$H^1(\Bbb{P}^1;u^*TX)=\{0\},\qquad
H^1\big(\Bbb{P}^1;u^*TX\!\otimes\!{\cal O}_{\Bbb{P}^1}(-1)\big)=\{0\},
\quad\hbox{and}\quad H^1(\Si;u^*TX)=\{0\},$$
respectively. In the integrable case, the two surjectivity statements of~(b)
in Definition~\ref{g1reg_dfn} are equivalent and can be written~as
$$H^1\big(\Bbb{P}^1;u^*TX\!\otimes\!{\cal O}_{\Bbb{P}^1}(-2)\big)=\{0\}.$$
It is well-known that the standard complex structure~$J_0$ on $\P$
is genus-one $d\ell$-regular for every $d\!\in\!\Bbb{Z}$;
see Corollaries~6.3 and~6.5 in~\cite{Z3}, for example.\\

\noindent
If $J$ is a genus-zero $A$-regular almost complex structure on~$X$, 
the structure of the moduli space $\ov\M_{0,k}(X,A;J)$ is regular
for every $k\!\in\!\bar\Z^+$.
In other words, $\ov\M_{0,k}(X,A;J)$ is stratified by smooth oriented orbifolds 
of even dimensions and the neighborhood of each stratum has the expected form.
One of the results of this paper is that
if $J$ is genus-one $A$-regular, 
the structure of the moduli space $\ov\M_{1,k}^0(X,A;J)$ is regular
for every $k\!\in\!\bar\Z^+$;
see Theorem~\ref{str_thm} and Subsection~\ref{reg1_subs1}.
In particular, we have

\begin{crl}[of Theorem~\ref{str_thm}]
\label{str_crl}
Suppose $(X,\om,J)$ is a compact almost Kahler manifold, \linebreak
$A\!\in\!H_2(X;\Z)^*$, and $k\!\in\!\bar{\Z}^+$.
If $J$ is genus-one $A$-regular,
then the closure of $\M_{1,k}^0(X,A;J)$ in $\ov\M_{1,k}(X,A;J)$
is $\ov\M_{1,k}^0(X,A;J)$.
Furthermore, $\ov\M_{1,k}^0(X,A;J)$ has the general
topological structure of a unidimensional algebraic orbivariety and thus carries
a rational fundamental class.
\end{crl}

\noindent
The first statement of Corollary~\ref{str_crl} follows from
the first claim of Theorem~\ref{str_thm}, along with standard gluing arguments
such as in Chapter~5 of~\cite{McSa}; see also Subsection~\ref{reg1_subs1}.
The middle statement of Corollary~\ref{str_crl} summarizes Theorem~\ref{str_thm},
while the last one is obtained at the end of Subsection~\ref{str_subs}.\\

\noindent
We will also show that the genus-zero and genus-one regularity properties 
are well-behaved under small perturbations:

\begin{thm}
\label{reg_thm}
Suppose $(X,\om,J)$ is a compact almost Kahler manifold and $A\!\in\!H_2(X;\Z)^*$.
If $g\!=\!0,1$ and the almost complex structure $J$ is genus-$g$ $A$-regular, 
then there exists $\de_J(A)\!\in\!\Bbb{R}^+$ with the property that 
if $\tilde{J}$ is an almost complex structure
on $X$ such that $\|\tilde{J}\!-\!J\|_{C^1}\!\le\!\de_J(A)$,
then $\tilde{J}$ is genus-$g$ $A$-regular.
Furthermore, if $J$ is genus-one $A$-regular, 
$k\!\in\!\bar{\Z}^+$, and $\under{J}\!=\!(J_t)_{t\in[0,1]}$ 
is a continuous family of almost complex structures on~$X$ such that
$J_0\!=\!J$ and $\|J_t\!-\!J\|_{C^1}\!\le\!\de_J(A)$ for all  $t\!\in\![0,1]$, 
then the moduli space $\ov\M_{1,k}^0(X,A;\under{J})$ has the general topological 
structure of a unidimensional algebraic orbivariety with boundary and
$$\partial\ov\M_{1,k}^0(X,A;\under{J})
=\ov\M_{1,k}^0(X,A;J_1)-\ov\M_{1,k}^0(X,A;J_0).$$\\
\end{thm}

\noindent
The norms $\|\ti{J}\!-\!J\|_{C^1}$ and $\|J_t\!-\!J\|_{C^1}$ are 
computed using a fixed connection in the vector bundle~$TX$,
e.g.~the Levi-Civita connection of the metric on $X$
determined by~$(\om,J)$.
The regularity claims of Theorem~\ref{reg_thm}
follow from the compactness of the moduli spaces $\ov\M_{g,k}(X,A;J)$ 
and Corollaries~\ref{reg0_crl1}, \ref{reg0_crl2}, \ref{reg0_crl3},
\ref{reg1_crl1}, and \ref{reg1_crl2}.
The final claim of Theorem~\ref{reg_thm} follows from 
a family version of Theorem~\ref{str_thm}.
It can in fact be used to show that under the assumptions of Theorem~\ref{reg_thm}
$$\ov\M_{1,k}^0(X,A;\under{J})\approx [0,1]\times \ov\M_{1,k}^0(X,A;J_0).$$
The conclusion, as stated, can be obtained with weaker regularity assumptions
on~$\under{J}$.\\

\noindent
The key ingredients in the proofs of Theorems~\ref{reg_thm} and~\ref{str_thm} are 
the gluing constructions of~\cite{Z4}, adapted to the present situation, 
and the power series expansions of Theorem~2.8 and Subsection~4.1 in~\cite{Z3},
applied via a technical result of~\cite{FlHSa}.
The power series of~\cite{Z3} give estimates on the behavior of derivatives of
holomorphic maps under gluing and on the obstructions to smoothing
holomorphic maps from singular domains.
The technical result of~\cite{FlHSa} shows that locally a $J$-holomorphic map
is very close to a holomorphic one.

\subsection{Some Geometric Implications}
\label{geomcons_subs}

\noindent
Theorem~\ref{comp_thm} implies that under certain assumptions on $A$ and $J$
the number of genus-one degree-$A$ $J$-holomorphic {\it curves} that pass through 
a collection of cycles in $X$ of the appropriate total codimension is finite.
Furthermore, each such curve is isolated to first order, as explained below.
Throughout this subsection, we assume that the dimension of $X$ is $2n\!\ge\!4$.\\

\noindent
A {\tt simple} $J$-holomorphic map into $X$ is a $J$-holomorphic map $u\!:\Si\!\lra\!X$ 
such that $u$ is one-to-one outside of finitely many points of~$\Si$ and
the irreducible components of $\Si$ on which $u$ is constant.
A {\tt genus-$g$ degree-$A$ $J$-holomorphic curve} $\cC$ is the image $u(\Si)$
of an element $[\Si,u]$ of $\ov\M_{g,k}(X,A;J)$ such that $u$ is simple and
the total genus of the components on which $u$ is not constant is~$g$.
Let $\cM_g(X,A;J)$ be the space of all genus-$g$ degree-$A$ 
$J$-holomorphic curves in~$X$. 
The expected dimension of this space is~$\dim_{g,0}(X,A)$.\\

\noindent
A $J$-holomorphic curve $\cC\subset\!X$ will be called {\tt regular} if
the operator $D_{J,u}$ is surjective for a (or equivalently, every)
stable-map parametrization $u\!:\Si\!\lra\!\cC$ of $\cC$ as above.
We will call a regular curve $\cC\subset\!X$ {\tt essentially embedded}
if $\cC$ is an irreducible curve that has no singularities if $n\!\ge\!3$
and its only singularities are simple nodes if $n\!=\!2$.
In other words, if $u\!:\Si\!\lra\!\cC$ is a parametrization of $\cC$ with $k\!=\!0$, 
then $\Si$ is a smooth Riemann surface of genus~$g$.
Furthermore, if $n\!\ge\!3$, $u$ is an embedding.
If $n\!=\!2$, then
$$\dim_{\C}\hbox{Span}_{\C}\big\{\Im du|_z\!: z\!\in\!u^{-1}(q)\big\}
=\big|u^{-1}(q)\big| \qquad\forall q\!\in\!X.$$
In particular, $u$ is an immersion.\\

\noindent
Let $\mu\!=\!(\mu_1,\ldots,\mu_k)$ be a $k$-tuple of cycles in $X$ of total
(real) codimension $\dim_{g,0}(X,A)\!+\!2k$, i.e.
$$\codim\,\mu\equiv\sum_{l=1}^{l=k}\codim\,\mu_l
=\dim_{g,0}(X,A)\!+\!2k=\dim_{g,k}(X,A).$$
We denote by $\cM_g(X,A;J,\mu)$ the set of genus-$g$ degree-$A$ $J$-holomorphic curves
that pass through every cycle $\mu_1,\ldots,\mu_k$, i.e.:
\begin{gather*}
\cM_g(X,A;J,\mu)=\big\{\cC\!\in\!\cM_g(X,A;J)\!: 
\cC\!\cap\!\mu_l\!\neq\!\eset~\forall l\!\in\![k]\big\},\\
\hbox{where}\qquad [k]=\big\{1,\ldots,k\}.
\end{gather*}
We will call an element $\cC$ of $\cM_g(X,A;J,\mu)$ {\tt isolated to first order} 
if for every parametrization
$$u\!:\Si\lra\cC\subset X$$
of $\cC$, where $\Si$ is a curve with $k$ marked points $y_1,\ldots,y_k$
such that $u(y_l)\!\in\!\mu_l$ for all $l\!\in\![k]$,
\begin{equation*}\begin{split}
\xi(z)\in \Im du|_z \quad
\forall~&z\!\in\!\Si~\hbox{s.t.}~du|_z\neq0
~~\hbox{and}\\
&\xi\!\in\!\ker D_{J,u}~\hbox{s.t.}~
\xi(y_l)\!\in\!T_{u(y_l)}\mu_l\!+\!\Im du|_{y_l}~\forall l\!\in\![k].
\end{split}\end{equation*}
If $\cM_g(X,A;J)$ is a smooth manifold with the expected tangent bundle and 
the constraints $\mu_1,\ldots,\mu_k$ are in general position, then
$\cM_g(X,A;J,\mu)$ is a discreet set consisting of elements isolated to first order.
Below we describe some circumstances under which this set is also finite.\\

\noindent
We recall that $A\!\in\!H_2(X;\Z)$ is called {\tt spherical} if
$$A=f_*[S^2]\in H_2(X;\Z)$$
for some smooth map $f\!:S^2\!\lra\!X$.
A symplectic manifold $(X,\om)$ is {\tt weakly monotone}
if for every spherical homology class~$A$ such that $\om(A)\!>\!0$,
either 
$$\lr{c_1(TX),A}\ge0  \qquad\hbox{or}\qquad \lr{c_1(TX),A}\le 2\!-\!n,$$
where $2n\!=\!\dim X$, as before.
In particular, all symplectic manifolds of (real) dimensions $2$, $4$, and~$6$
are weakly monotone.
So are all complex projective spaces, which are in fact monotone;
see Chapter~5 in~\cite{McSa} for a definition.\\

\noindent
Finally, if $(X,\om)$ is a symplectic manifold, we denote by $\J(X,\om)$ 
the space of all almost complex structures on $X$
tamed by~$\om$, endowed with the $C^1$-topology.

\begin{prp}
\label{curve_prp1}
Suppose $(X,\om,J)$ is a compact almost Kahler manifold, $A\!\in\!H_2(X;\Z)^*$, 
$g\!=\!0,1$, and $J$ is genus-$g$ $A$-regular.
If $\mu$ is a $k$-tuple of cycles in $X$ of total codimension $\dim_{g,k}(X,A)$
in general position, then $\cM_g(X,A;J,\mu)$ is a finite set and
every element in $\cM_g(X,A;J,\mu)$ is irreducible, regular, and isolated to first order.
\end{prp}

\begin{prp}
\label{curve_prp2}
If $(X,\om)$ is a compact weakly monotone symplectic manifold and $A\!\in\!H_2(X;\Z)^*$,
there exists a dense open subset $\J_{\reg}(X,\om;A)$ of $\J(X,\om)$
with the following properties.
If $J\!\in\!\J_{\reg}(X,\om;A)$, $g\!=\!0,1$, and $\mu$ is a $k$-tuple 
of cycles in $X$ of total codimension $\dim_{g,k}(X,A)$
in general position, then every element in $\cM_g(X,A;J,\mu)$ is 
essentially embedded and isolated to first order.
If in addition $\blr{c_1(TX),A}\!\neq\!0$, then $\cM_g(X,A;J,\mu)$
is a finite set and its signed cardinality is independent of 
the choice of $J\!\in\!\J_{\reg}(X,\om;A)$.
\end{prp}

\noindent
Most of the genus-zero statements of these two propositions are well-known
results; see Chapters 5-7 of~\cite{McSa}, for example.
The signed cardinality of the set $\cM_0(X,A;J,\mu)$ is the corresponding
Gromov-Witten invariant, $\GW_{0,k}(A;\mu)$, of~$(X,\om)$.\\

\noindent
The remaining statements are obtained from minor extensions of some results in~\cite{McSa},
along with Theorem~\ref{comp_thm} in the genus-one case.
For each $l\!\in\![k]$, let
$$\ev_l\!:\ov\M_{g,k}(X,A;J)\lra X, \qquad
\big[\Si,y_1,\ldots,y_k;u\big]\lra u(y_l),$$
be the evaluation map for the $l$th marked point.\\

\noindent
Suppose $(X,\om,J)$ and $A$ are as in the statement of Proposition~\ref{curve_prp1}
and $J$ is genus-one $A$-regular.
By~(c) of Definition~\ref{g1reg_dfn}, the moduli space $\M_{1,k}^0(X,A;J)$
is a smooth orbifold with the expected tangent bundle.
Thus, if $\mu$ is a tuple of constraints as in the $g\!=\!1$ case of 
Proposition~\ref{curve_prp1}, then
$$\M_{1,k}^0(X,A;J,\mu) \equiv
\{b\!\in\!\M_{1,k}^0(X,A;J)\!: 
\ev_l(b)\!\in\!\mu_l~\forall l\!\in\![k]\}$$
is a zero-dimensional oriented submanifold.
If $\{b_r\}$ is a sequence of distinct elements in $\M_{1,k}^0(X,A;J,\mu)$,
by Theorem~\ref{comp_thm} a subsequence of $\{b_r\}$ must converge
to an element 
\begin{equation*}\begin{split}
b\in\ov\M_{1,k}^0(X,A;J,\mu) 
\equiv& \{b\!\in\!\ov\M_{1,k}^0(X,A;J)\!: 
\ev_l(b)\!\in\!\mu_l~\forall l\!\in\![k]\} \\
&\subset \ov\M_{1,k}(X,A;J).
\end{split}\end{equation*}
Since all elements of $\M_{1,k}^0(X,A;J,\mu)$ are isolated,
$$b\in\ov\M_{1,k}^0(X,A;J,\mu) -\M_{1,k}^0(X,A;J).$$
On the other hand, by the regularity assumptions of Definition~\ref{g1reg_dfn},
$$\partial\ov\M_{1,k}^0(X,A;J)\equiv\ov\M_{1,k}^0(X,A;J)-\M_{1,k}^0(X,A;J)$$
is a union of strata of dimensions smaller than $\dim_{1,k}(X,A)$.
Thus, if $\mu$ is a tuple of cycles of total codimension $\dim_{1,k}(X,A)$
in general position, then
$$\ov\M_{1,k}^0(X,A;J,\mu) -\M_{1,k}^0(X,A;J) =\eset.$$
It follows that $\M_{1,k}^0(X,A;J,\mu)$ is a finite set, and so is
its subset $\cM_1(X,A;J,\mu)$.\\

\noindent
We next move to Proposition~\ref{curve_prp2}. 
For any $B\!\in\!H_2(X;\Z)^*$, almost complex structure $J$ on~$X$, 
and $g\!\in\!\Z^+$, let
\begin{alignat*}{2}
&\M_{g,k}^{\simp}(X,B;J)&\subset~&\M_{g,k}^0(X,B;J), \\
&\M_{g,k}^{\{\simp\}}(X,B;J)&\subset~&\M_{g,k}^{\{0\}}(X,B;J), \qquad\hbox{and}\\
&\ov\M_{g,k}^{\simp}(X,B;J)&\subset~&\ov\M_{g,k}(X,B;J)
\end{alignat*}
denote the subspaces of simple maps.
By Chapter~3 in~\cite{McSa}, for a generic choice of~$J$, 
$D_{J,u}$ is surjective for every element $[\Si,u]$ of $\M_{g,k}^{\simp}(X,A;J)$.
Thus, as before, 
$$\cM_1(X,A;J,\mu) \approx  \M_{1,k}^{\simp}(X,A;J,\mu) \equiv
\M_{1,k}^0(X,A;J,\mu)\cap \M_{1,k}^{\simp}(X,A;J)$$
is a zero-dimensional oriented manifold, if $\mu$ is as in Proposition~\ref{curve_prp2}.
On the other hand, by the same argument as in Chapter~6 of~\cite{McSa}, 
the evaluation map 
$$\ev_{k+1}\!\times\!\ev_{k+2}\!:\M_{1,k+2}^{\simp}(X,A;J)
\lra X\!\times\!X$$
is transverse to the diagonal, for a generic choice of $J$.
Let 
$$L_{k+1}\lra \M_{1,k}^{\simp}(X,A;J)$$
be the universal tangent line bundle for the last marked point, i.e.
$$L_{k+1}|_{[\Si,u]}=T_{y_{k+1}}\Si \qquad \forall\, [\Si,u]\in\M_{1,k+1}^{\simp}(X,A;J).$$
By a small modification of the proof of Lemma~6.1.2 in~\cite{McSa}, 
the bundle section 
$$du|_{y_{k+1}}\!:\M_{1,k}^{\simp}(X,A;J) \lra L_{k+1}^*\!\otimes\!\ev_{k+1}^*TX,
\qquad [\Si,u]\lra du|_{y_{k+1}},$$
is transverse to the zero set, for a generic choice of $J$.
The key part of this modification is to view the relevant first-order equation 
as an elliptic operator acting on the space of smooth sections of the vector bundle
$n\O_{\bP^1}(1)$ over~$S^2$.
The last two transversality properties imply that for a generic element 
$[\Si,u]$ of $\M_{1,k}^{\{\simp\}}(X,A;J)$ its image $u(\Si)$ is essentially embedded.
This concludes the proof of the first statement of Proposition~\ref{curve_prp2}.\\

\noindent
By Chapters~3 and~6 in~\cite{McSa}, for a generic choice of~$J$, 
$D_{J,u}$ is surjective for every element $[\Si,u]$ of $\M_{g,k}^{\{\simp\}}(X,B;J)$.
In particular, $\M_{g,k}^{\{\simp\}}(X,B;J)$ is a finite union of smooth orbifolds
of the expected dimension.
Thus, if $\mu$ is a tuple of constraints as in the statement of Proposition~\ref{curve_prp2},
$$\big\{b\!\in\!\M_{1,k}^{\{\simp\}}(X,A;J)\!:
\ev_l(b)\!\in\!\mu_l~\forall l\!\in\![k]\big\} -\M_{1,k}^0(X,A;J)=\eset.$$
Furthermore, if $\lr{c_1(TX),A}\!\neq\!0$, every element of 
$$\ov\M_{1,k}(X,A;J,\mu) \equiv \{b\!\in\!\ov\M_{1,k}(X,A;J)\!: 
\ev_l(b)\!\in\!\mu_l~\forall l\!\in\![k]\}$$
is simple. This can be seen by considering the dimension of the image of 
the multiply covered elements of $\ov\M_{1,k}(X,A;J)$ under the evaluation map
$\ev_1\!\times\!\ldots\!\times\!\ev_k$.
This is done by passing to moduli spaces of maps consisting of simple elements;
see Chapter~5 of~\cite{McSa}.
The argument requires two separate dimension counts for multiply covered maps:
one for the elements in $\M_{1,k}^{\{0\}}(X,A;J)$
and the other for those in its complement in $\ov\M_{1,k}(X,A;J)$.
In addition to the assumption \hbox{$\lr{c_1(TX),A}\!\neq\!0$},
the weakly monotone condition on $(X,\om)$ enters directly into both dimension computations.
Finally, by the same modification of the proof of Lemma~6.1.2 in~\cite{McSa} 
as described above, but applied to tuples of genus-zero maps instead of genus-one maps,
$$\ov\M_{1,k}^{\simp}(X,A;J)\cap
\big(\ov\M_{1,k}^0(X,A;J)\!-\!\M_{1,k}^{\{0\}}(X,A;J)\big)$$
is a finite union of smooth orbifolds of dimensions less than 
$\dim_{1,k}(X,A)$. We conclude~that
$$\cM_1(X,A;J,\mu) \approx \M_{1,k}^{\simp}(X,A;J,\mu)
=\ov\M_{1,k}^0(X,A;J,\mu)$$
is a compact zero-dimensional manifold.
By a cobordism argument as in Chapter~7 of~\cite{McSa},
the signed cardinality of $\cM_1(X,A;J,\mu)$ is independent
of a generic choice of~$J$.\\

\noindent
The signed cardinality $\GW_{1,k}^0(A;\mu)$ of the set $\cM_{1,k}(X,A;J,\mu)$ 
is an integer-valued invariant of the symplectic manifold~$(X,\om)$.
The difference between this invariant for an arbitrary symplectic manifold
(when the invariant may not be an integer) and the standard 
genus-one GW-invariant is specified by Proposition~\ref{g1comp2-bdcontr_prp} in~\cite{Z6} and explicitly given by Theorems~\ref{g1diff-main_thm1} and~\ref{g1diff-main_thm2} in~\cite{Z9}.

\section{Preliminaries}
\label{prelim_sec}

\subsection{Notation: Genus-Zero Maps}
\label{notation0_subs}

\noindent
We now describe our notation for bubble maps from genus-zero Riemann surfaces,
for the spaces of such bubble maps that form
the standard stratifications of the moduli spaces of stable maps,
and for important vector bundles over them.
In general, the moduli spaces of stable maps can stratified by the dual graph.
However, in the present situation, it is more convenient to make use
of {\it linearly ordered sets}:

\begin{dfn}
\label{index_set_dfn1}
(1) A finite nonempty partially ordered set $I$ is a {\tt linearly ordered set} if 
for all \hbox{$i_1,i_2,h\!\in\!I$} such that $i_1,i_2\!<\!h$, 
either $i_1\!\le\!i_2$ \hbox{or $i_2\!\le\!i_1$.}\\
(2) A linearly ordered set $I$ is a {\tt rooted tree} if
$I$ has a unique minimal element, 
i.e.~there exists \hbox{$\hat{0}\!\in\!I$} such that $\hat{0}\!\le\!i$ 
for {all $i\!\in\!I$}.
\end{dfn}

\noindent
We use rooted trees to stratify the moduli space $\ov\M_{0,\{0\}\sqcup M}(X,A;J)$
of genus-zero stable holomorphic maps with marked points indexed by
the set $\{0\}\!\sqcup\!M$, where $M$ is a finite set.\\

\noindent
If $I$ is a linearly ordered set, let $\hat{I}$ be 
the subset of the non-minimal elements of~$I$.
For every $h\!\in\!\hat{I}$,  denote by $\io_h\!\in\!I$
the largest element of $I$ which is smaller than~$h$, i.e.
$\io_h\!=\!\max\big\{i\!\in\!I:i\!<\!h\big\}$.\\

\noindent
We identify $\C$ with $S^2\!-\!\{\i\}$ via 
the stereographic projection mapping the origin in $\C$ 
to the north pole, or the point $(0,0,1)$, in~$S^2$.
A {\tt genus-zero $X$-valued bubble map with $M$-marked points} is a~tuple
$$b=\big(M,I;x,(j,y),u\big),$$
where $I$ is a rooted tree, and
\begin{equation}\label{stablemap_e1}
x\!:\hat{I}\!\lra\!\C\!=\!S^2\!-\!\{\i\},\quad  j\!:M\!\lra\!I,\quad
y\!:M\!\lra\!\C,        \hbox{~~~and~~~} 
u\!:I\!\lra\!C^{\i}(S^2;X)
\end{equation}
are maps such that $u_h(\i)\!=\!u_{\io_h}(x_h)$ for all $h\!\in\!\hat{I}$.
We associate such a tuple with Riemann surface
\begin{equation}\label{stablemap_e2}
\Si_b=
\Big(\bigsqcup_{i\in I}\Si_{b,i}\Big)\Big/\!\sim,
\hbox{~~where}\qquad \Si_{b,i}=\{i\}\!\times\!S^2
\quad\hbox{and}\quad
(h,\i)\sim (\io_h,x_h)
~~\forall h\!\in\!\hat{I},
\end{equation}
with marked points 
$$y_l(b)\!\equiv\!(j_l,y_l)\in\Si_{b,j_l}  \qquad\hbox{and}\qquad 
y_0(b)\!\equiv\!(\hat{0},\i)\in\Si_{b,\hat{0}}$$
and with the continuous map $u_b\!:\Si_b\!\lra\!X$,
given by $u_b|_{\Si_{b,i}}\!=\!u_i$ for \hbox{all $i\!\in\!I$}.
The general structure of bubble maps is described
by tuples ${\cal T}\!=\!(M,I;j,\under{A})$, where 
$$A_i=\big\{u_b|_{\Si_{b,i}}\big\}_*[S^2] \qquad\forall i\!\in\!I.$$ 
We call such tuples {\it bubble types}.
Let ${\cal U}_{\cal T}(X;J)$ denote the subset of $\ov\M_{0,\{0\}\sqcup M}(X,A;J)$ 
consisting of stable maps $[{\cal C},u]$ such that
$$[{\cal C};u]=
\big[(\Si_b,(\hat{0},\i),(j_l,y_l)_{l\in M});u_b\big],$$
for some bubble map $b$ of type ${\cal T}$ as above,
where $\hat{0}$ is the minimal element of~$I$; 
see Section~2 in~\cite{Z4} for details. 
For $l\!\in\!\{0\}\!\sqcup\!M$, let 
$$\ev_l\!:{\cal U}_{\cal T}(X;J)\lra X$$ 
be the evaluation map corresponding to the marked point $y_l$.\\

\noindent
We denote the bundle of gluing parameters, or of smoothings at the nodes, 
over ${\cal U}_{\cal T}(X;J)$ by~${\cal FT}$.
This orbi-bundle has the form 
$${\cal FT}=\Big(\!
\bigoplus_{h\in\hat{I}}L_{h,0}\!\otimes\!L_{h,1}\Big)\big/\Aut(\T),$$
for certain line orbi-bundles $L_{h,0}$ and $L_{h,1}$.
These line bundles are the line bundles associated to certain $S^1$-principal bundles.
More precisely, there exists a subspace ${\cal U}_{\cal T}^{(0)}(X;J)$
of the space ${\cal H}_{\cal T}(X;J)$ of $J$-holomorphic maps into~$X$
of type~${\cal T}$, not of equivalence classes of such maps, such~that
$${\cal U}_{\T}(X;J)={\cal U}_{\cal T}^{(0)}(X;J)\big/
\Aut(\T)\!\propto\!(S^1)^I.$$
The line bundles $L_{h,0}$ and $L_{h,1}$ arise from this quotient;
see Subsection~2.5 in~\cite{Z4}.
In particular,
$${\cal FT}=\tilde{\cal F}{\cal T}\big/\Aut(\T)\!\propto\!(S^1)^I,
\qquad\hbox{where}\qquad
\tilde{\cal F}{\cal T}=
{\cal U}_{\cal T}^{(0)}(X;J)\!\times\!\Bbb{C}^{\hat{I}}\!\lra\!{\cal U}_{\cal T}^{(0)}(X;J).$$
We denote by ${\cal FT}^{\eset}$ and $\tilde{\cal F}{\cal T}^{\eset}$ the subsets of 
${\cal FT}$ and $\tilde{\cal F}{\cal T}$, respectively,
consisting of the elements with all components nonzero.\\

\noindent
The subset ${\cal U}_{\cal T}^{(0)}(X;J)$ of ${\cal H}_{\cal T}(X;J)$ is described by 
the conditions (B1) and (B2) in Subsection~2.5 of~\cite{Z4}. 
It is the preimage of the point $(0,1/2)^I$ in $(\Bbb{C}\!\times\!\Bbb{R})^I$
under the continuous map
$$\Psi_{\cal T}\!\equiv\!(\Psi_{{\cal T},i})_{i\in I}\!:
{\cal H}_{\cal T}(X;J)\lra(\Bbb{C}\!\times\!\Bbb{R})^I$$
defined in the proof of Proposition~3.3 in~\cite{Z4}.
The statements of the conditions (B1) and (B2) and the definition of
the map $\Psi_{\cal T}$ require a choice of a $J$-compatible metric~$g_X$.
It can be assumed~that
$$\int_{\Bbb{P}^1}|du|^2_{g_X}\ge 1$$
for every non-constant $J$-holomorphic maps $u\!:\Bbb{P}^1\!\lra\!X$.
Such a metric $g_X$ will be fixed once and for all.
If the almost complex structure~$J$ is genus-zero $A$-regular, 
where $A\!=\!\sum_{i\in I}\!A_i$,
the space ${\cal H}_{\cal T}(X;J)$  is a smooth manifold of the expected dimension;
see Chapter~3 in~\cite{McSa}.
In such a case, the map $\Psi_{\cal T}$ is smooth and transversal to every point
$(0,r_i)_{i\in I}$ such that $|r_i\!-\!\frac{1}{2}|\!\le\!\frac{1}{4}$ for all $i\!\in\!I$;
see the proof of Proposition~3.3 in~\cite{Z4}.
Let
\begin{gather}\label{chidfn0_e}
\chi({\cal T})=\big\{i\!\in\!I:A_i\!\neq\!0;~
A_h\!=\!0~\forall h\!<\!i\big\};\\
\tilde{\cal U}_{\cal T}^{(0)}(X;J)=\Psi_{\cal T}^{-1}\big(
\big\{(0,r_i)_{i\in I}\!\in\!(\Bbb{C}\!\times\!\Bbb{R})^I\!:
r_i\!=\!\frac{1}{2}~\forall i\!\in\!I\!-\!\chi({\cal T}),~
r_i\!\in\!\big(\frac{1}{4},\frac{3}{4}\big)~\forall i\!\in\!\chi({\cal T})\big\},\notag\\
\wt{\cal FT}=
\tilde{\cal U}_{\cal T}^{(0)}(X;J)\!\times\!\Bbb{C}^{\hat{I}}
\!\lra\!\tilde{\cal U}_{\cal T}^{(0)}(X;J).  \notag
\end{gather}
As before, we denote by $\wt{\cal FT}^{\eset}$ the subset of $\wt{\cal FT}$
consisting of the elements with all components nonzero.

\subsection{Notation: Genus-One Maps}
\label{notation1_subs}

\noindent
We next set up analogous notation for maps from genus-one Riemann surfaces.
In this case, we also need to specify the structure of the principal component.
Thus, we index the strata of the moduli space $\ov\M_{1,M}(X,A;J)$
by {\it enhanced linearly ordered sets}:

\begin{dfn}
\label{index_set_dfn2}
An {\tt enhanced linearly ordered set} is a pair $(I,\aleph)$,
where $I$ is a linearly ordered set, $\aleph$ is a subset of $I_0\!\times\!I_0$,
and $I_0$ is the subset of minimal elements of~$I$,
such that if $|I_0|\!>\!1$, 
$$\aleph=\big\{(i_1,i_2),(i_2,i_3),\ldots,(i_{n-1},i_n),(i_n,i_1)\big\}$$
for some bijection $i\!:\{1,\ldots,n\}\!\lra\!I_0$.
\end{dfn}

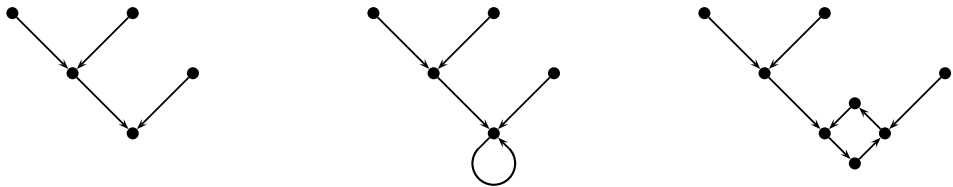
\begin{figure}
\begin{pspicture}(-1.1,-2)(10,1)
\psset{unit=.4cm}
\pscircle*(6,-3){.2}
\pscircle*(4,-1){.2}\psline[linewidth=.06]{->}(4.14,-1.14)(5.86,-2.86)
\pscircle*(8,-1){.2}\psline[linewidth=.06]{->}(7.86,-1.14)(6.14,-2.86)
\pscircle*(2,1){.2}\psline[linewidth=.06]{->}(2.14,.86)(3.86,-.86)
\pscircle*(6,1){.2}\psline[linewidth=.06]{->}(5.86,.86)(4.14,-.86)
% 2nd picture starts here
\pscircle*(18,-3){.2}\psline[linewidth=.06](17.86,-3.14)(17.5,-3.5)
\psarc(18,-4){.71}{135}{45}\psline[linewidth=.06]{->}(18.5,-3.5)(18.14,-3.14)
\pscircle*(16,-1){.2}\psline[linewidth=.06]{->}(16.14,-1.14)(17.86,-2.86)
\pscircle*(20,-1){.2}\psline[linewidth=.06]{->}(19.86,-1.14)(18.14,-2.86)
\pscircle*(14,1){.2}\psline[linewidth=.06]{->}(14.14,.86)(15.86,-.86)
\pscircle*(18,1){.2}\psline[linewidth=.06]{->}(17.86,.86)(16.14,-.86)
% 3rd picture starts here
\pscircle*(30,-2){.2}\pscircle*(30,-4){.2}\pscircle*(29,-3){.2}\pscircle*(31,-3){.2}
\psline[linewidth=.06]{->}(29.86,-2.14)(29.14,-2.86)
\psline[linewidth=.06]{->}(29.14,-3.14)(29.86,-3.86)
\psline[linewidth=.06]{->}(30.14,-3.86)(30.86,-3.14)
\psline[linewidth=.06]{->}(30.86,-2.86)(30.14,-2.14)
\pscircle*(27,-1){.2}\psline[linewidth=.06]{->}(27.14,-1.14)(28.86,-2.86)
\pscircle*(33,-1){.2}\psline[linewidth=.06]{->}(32.86,-1.14)(31.14,-2.86)
\pscircle*(25,1){.2}\psline[linewidth=.06]{->}(25.14,.86)(26.86,-.86)
\pscircle*(29,1){.2}\psline[linewidth=.06]{->}(28.86,.86)(27.14,-.86)
\end{pspicture}
\caption{Some Enhanced Linearly Ordered Sets}
\label{index_set_fig}
\end{figure}

\noindent
An enhanced linearly ordered set can be represented by an oriented connected graph.
In Figure~\ref{index_set_fig}, the dots denote the elements of~$I$.
The arrows outside the loop, if there are any, 
specify the partial ordering of the linearly ordered set~$I$.
In fact, every directed edge outside of the loop
connects a non-minimal element $h$ of $I$ with~$\io_h$.
Inside of the loop, there is a directed edge from $i_1$ to $i_2$
if and only if $(i_1,i_2)\!\in\!\aleph$.\\

\noindent
The subset $\aleph$ of $I_0\!\times\!I_0$ will be used to describe
the structure of the principal curve of the domain of stable maps in 
a stratum of the moduli space~$\ov\M_{1,M}(X,A;J)$.
If $\aleph\!=\!\eset$, and thus $|I_0|\!=\!1$,
the corresponding principal curve $\Si_{\aleph}$ 
is a smooth torus, with some complex structure.
If $\aleph\!\neq\!\eset$, the principal components form a circle of spheres:
$$\Si_{\aleph}=\Big(\bigsqcup_{i\in I_0}\{i\}\!\times\!S^2\Big)\Big/\sim,
\qquad\hbox{where}\qquad
(i_1,\i)\sim(i_2,0)~~\hbox{if}~~(i_1,i_2)\!\in\!\aleph.$$
A {\tt genus-one $X$-valued bubble map with $M$-marked points} is a tuple
$$b=\big(M,I,\aleph;S,x,(j,y),u\big),$$
where $S$ is a smooth Riemann surface of genus one if $\aleph\!=\!\eset$
and the circle of spheres $\Si_{\aleph}$ otherwise.
The objects $x$, $j$, $y$, $u$, and $(\Si_b,u_b)$ are as in 
\e_ref{stablemap_e1} and \e_ref{stablemap_e2}, except 
the sphere $\Si_{b,\hat{0}}$ is replaced by the genus-one curve $\Si_{b,\aleph}\!\equiv\!S$.
Furthermore, if $\aleph\!=\!\eset$, and thus $I_0\!=\!\{\hat{0}\}$ is a single-element set,
$u_{\hat{0}}\!\in\!C^{\i}(S;X)$ and $y_l\!\in\!S$ if $j_l\!=\!\hat{0}$.
In the genus-one case, the general structure of bubble maps is encoded by
the tuples of the form ${\cal T}\!=\!(M,I,\aleph;j,\under{A})$.
Similarly to the genus-zero case, we denote by ${\cal U}_{\cal T}(X;J)$
the subset of $\ov\M_{1,M}(X,A;J)$ 
consisting of stable maps $[{\cal C},u]$ such that
$$[{\cal C};u]=\big[(\Si_b,(j_l,y_l)_{l\in M});u_b\big],$$
for some bubble map $b$ of type ${\cal T}$ as above.\\

\noindent
If $\T\!=\!(M,I,\aleph;j,\under{A})$ is a bubble type as above, let
\begin{equation}\label{g1gendecomp_e1}\begin{split}
&I_1=\big\{h\!\in\!\hat{I}\!:\io_h\!\in\!I_0\big\},\qquad
M_0=\big\{l\!\in\!M\!:j_l\!\in\!I_0\!\big\},    \qquad\hbox{and}\\
&\qquad{\cal T}_0=\big(M_0\!\sqcup\!I_1,I_0,\aleph;j|_{M_0}\!\sqcup\!\io|_{I_1},
\under{A}|_{I_0}\big),
\end{split}\end{equation}
where $I_0$ is the subset of minimal elements of $I$.
For each $h\!\in\!I_1$, we put
\begin{equation}\label{g1gendecomp_e2}
I_h=\big\{i\!\in\!I\!:h\!\le\!i\big\}, \qquad
M_h=\big\{l\!\in\!M\!:j_l\!\in\!I_h\!\big\},   \quad\hbox{and}\quad
{\cal T}_h=\big(M_h,I_h;j|_{M_h},\under{A}|_{I_h}\big).
\end{equation}
We have a natural isomorphism
\begin{equation}\label{g1gendecomp_e3}\begin{split}
{\cal U}_{\cal T}(X;J)\approx \Big(\big\{
\big(b_0,(b_h)_{h\in I_1}\big)\!\in\!{\cal U}_{{\cal T}_0}(X;J)\!\times\!\!
\prod_{h\in I_1}\!\!{\cal U}_{{\cal T}_h}(X;J)\!:\qquad&\\
\ev_0(b_h)\!=\!\ev_{\io_h}(b_0)~\forall h\!\in\!I_1&\big\}\Big)
\big/\Aut^*(\T),
\end{split}\end{equation}
where the group $\Aut^*(\T)$ is defined by
$$\Aut^*(\T)=\Aut(\T)/\{g\!\in\!\Aut(\T)\!:g\cdot h\!=\!h~\forall h\!\in\!I_1\}.$$
This decomposition is illustrated in Figure~\ref{g1gendecomp_fig}.
In this figure, we represent an entire stratum of bubble maps
by the domain of the stable maps in that stratum.
The right-hand side of Figure~\ref{g1gendecomp_fig} 
represents the subset of the cartesian product of the three spaces
of bubble maps, corresponding to the three drawings,
on which the appropriate evaluation maps agree pairwise,
as indicated by the dotted lines and defined in~\e_ref{g1gendecomp_e3}.\\

\begin{figure}
\begin{pspicture}(-1.1,-2)(10,1.25)
\psset{unit=.4cm}
\psellipse(8,-1.5)(1.5,2.5)
\psarc[linewidth=.05](6.2,-1.5){2}{-30}{30}\psarc[linewidth=.05](9.8,-1.5){2}{150}{210}
\pscircle[fillstyle=solid,fillcolor=gray](5.5,-1.5){1}\pscircle*(6.5,-1.5){.2}
\pscircle[fillstyle=solid,fillcolor=gray](3.5,-1.5){1}\pscircle*(4.5,-1.5){.2}
\pscircle(10.5,-1.5){1}\pscircle*(9.5,-1.5){.2}
\pscircle[fillstyle=solid,fillcolor=gray](11.91,-.09){1}\pscircle*(11.21,-.79){.2}
\pscircle[fillstyle=solid,fillcolor=gray](11.91,-2.91){1}\pscircle*(11.21,-2.21){.2}
\rput(5.5,0){$h_1$}\rput(3.5,0){$h_2$}\rput(10.3,0){$h_3$}
\rput(13.5,0.1){$h_4$}\rput(13.5,-2.9){$h_5$}
\pscircle*(8,1){.2}\pscircle*(8,-4){.2}
\rput(7.9,1.7){$y_1$}\rput(7.9,-4.8){$y_2$}
\rput(17.5,-1.5){$\approx$}
\psellipse(23,-1.5)(2.5,1.5)
\pscircle*(20.5,-1.5){.2}\pscircle*(25.5,-1.5){.2}
\rput(19.8,-1.6){$y_1$}\rput(26.2,-1.6){$y_2$}
\psarc[linewidth=.05](23,.3){2}{240}{300}\psarc[linewidth=.05](23,-3.3){2}{60}{120}
\psline[linewidth=.06,linestyle=dotted](28,1)(23,0)
\psline[linewidth=.06,linestyle=dotted](28,-3)(23,-3)
\pscircle*(23,0){.2}\pscircle*(23,-3){.2}
\pscircle[fillstyle=solid,fillcolor=gray](29,1){1}\pscircle*(28,1){.2}
\pscircle[fillstyle=solid,fillcolor=gray](31,1){1}\pscircle*(30,1){.2}
\pscircle(29,-3){1}\pscircle*(28,-3){.2}
\pscircle[fillstyle=solid,fillcolor=gray](30.41,-1.59){1}\pscircle*(29.71,-2.29){.2}
\pscircle[fillstyle=solid,fillcolor=gray](30.41,-4.41){1}\pscircle*(29.71,-3.71){.2}
\rput(29,2.5){$h_1$}\rput(31,2.5){$h_2$}\rput(28.8,-1.5){$h_3$}
\rput(32,-1.4){$h_4$}\rput(32,-4.1){$h_5$}
\end{pspicture}
\caption{An Example of the Decomposition~\e_ref{g1gendecomp_e3}}
\label{g1gendecomp_fig}
\end{figure}

\noindent
Let ${\cal FT}\!\lra\!{\cal U}_{\cal T}(X;J)$
be the bundle of gluing parameters, or of smoothings at the nodes.
This orbi-bundle has the form 
$${\cal FT}=\Big(\!\bigoplus_{(h,i)\in\aleph}\!\!\!L_{h,0}\!\otimes\!L_{i,1}\oplus
\bigoplus_{h\in\hat{I}}L_{h,0}\!\otimes\!L_{h,1}\Big)\big/\Aut(\T),$$
for certain line orbi-bundles $L_{h,0}$ and $L_{h,1}$.
Similarly to the genus-zero case,
\begin{gather}
\label{g1gendecomp_e4a}
{\cal U}_{\cal T}(X;J)={\cal U}_{\T}^{(0)}(X;J)\big/
\Aut(\T)\!\propto\!(S^1)^{\hat{I}},\qquad\hbox{where}\\
\label{g1gendecomp_e4b}
{\cal U}_{\cal T}^{(0)}(X;J)=
\big\{\big(b_0,(b_h)_{h\in I_1}\big)\!\in\!{\cal U}_{\T_0}(X;J)\!\times\!\!
\prod_{h\in I_1}\!\!{\cal U}_{\T_h}^{(0)}(X;J)\!:
\ev_0(b_h)\!=\!\ev_{\io_h}(b_0)~\forall h\!\in\!I_1\big\}
\end{gather}
and ${\cal U}_{\T_h}^{(0)}(X;J)$ is the subspace of the moduli
space of holomorphic maps from genus-zero curves as in Subsection~\ref{notation0_subs}.
The line bundles $L_{h,0}$ and $L_{h,1}$ arise from the quotient~\e_ref{g1gendecomp_e4a}.
More precisely,
\begin{gather*}
{\cal FT}=\tilde{\cal F}{\cal T}\big/\Aut(\T)\!\propto\!(S^1)^{\hat{I}},
\qquad\hbox{where}\qquad
\tilde{\cal F}{\cal T}=\tilde{\cal F}_{\aleph}{\cal T}\!\oplus\!
\tilde{\cal F}_0{\cal T}\!\oplus\!\tilde{\cal F}_1{\cal T}
\!\lra\!{\cal U}_{\cal T}^{(0)}(X;J),\\
\tilde{\cal F}_0{\cal T}=\bigoplus_{h\in I_1}\tilde{\cal F}_h{\cal T}, \qquad
\tilde{\cal F}_1{\cal T}=
{\cal U}_{\cal T}^{(0)}(X;J)\!\times\!\Bbb{C}^{\hat{I}-I_1},  
\end{gather*}
and $\tilde{\cal F}_h{\cal T}$ and $\tilde{\cal F}_{\aleph}{\cal T}$
are the pullbacks by the projection map
$$\pi_P\!:{\cal U}_{\cal T}^{(0)}(X;J)\lra{\cal U}_{{\cal T}_0}(X;J)$$
of the universal tangent line $L_h{\cal T}_0$ at the $h$th marked point 
and of the bundle~${\cal FT}_0$ of gluing parameters.
In other words, if ${\frak U}_{{\cal T}_0}\!\lra\!{\cal U}_{{\cal T}_0}(X;J)$
is the semi-universal family, i.e.~the fiber at $b_0\!\in\!{\cal U}_{{\cal T}_0}(X;J)$
is the Riemann surface $\Si_{b_0}\!=\!\Si_{b_0,\aleph}$,
$L_h{\cal T}_0$ is the vertical tangent space at the point $x_h(b_0)$ of $\Si_{b_0}$.\\

\noindent
{\it Remark 1:} The above description is slightly inaccurate.
In order to insure the existence of the space ${\frak U}_{{\cal T}_0}$,
with the fibers as described, we need to replace 
the space ${\cal U}_{{\cal T}_0}(X;J)$ by a finite cover, 
analogous to the one used in~\cite{RT2}.
However, correcting this inaccuracy would complicate the notation used even further,
but would have no effect on the analysis, and thus we ignore~it.\\

\noindent
{\it Remark 2:} The rank of the bundle ${\cal FT}_0$ is $|\aleph|$, 
the number of nodes in the domain of every element of ${\cal U}_{{\cal T}_0}(X;J)$.
If $\aleph\!\neq\!\eset$, ${\cal FT}_0$ can be written as 
the quotient of the trivial bundle of rank $|\aleph|$,
over a space ${\cal U}_{{\cal T}_0}^{(0)}(X;J)$,  by 
an $\Aut(\T_0)\!\propto\!(S^1)^{\aleph}$-action
in a manner similar to the previous subsection and to Section~2 of~\cite{Z4}.
With the above identifications, the singular points of every rational component~$\Si_i$
of $\Si_{\aleph}$ are the points $0$ and~$\i$ in~$S^2$. 
Thus, the equivalence class of the restriction of a stable map in ${\cal U}_{{\cal T}_0}(X;J)$
to $\Si_i$, with its nodes, has a $\Bbb{C}^*$-family of representatives.
This family is cut down to an $S^1$-family by restricting to
the subset defined by the condition~(B2) of Subsection~2.5
in~\cite{Z4}.
This is the preimage of $1/2$ under the last, real-valued, component of the function
$\Psi_{{\cal T},i^*}$ defined in (2) of the proof of Proposition~3.3 in~\cite{Z4}.\\

\noindent
If ${\cal T}\!=\!(M,I,\aleph;j,\under{A})$ is a bubble type such that
$A_i\!=\!0$ for all minimal elements $i$ of~$I$, i.e.
$${\cal U}_{\T}(X;J)\subset  \ov\M_{1,M}(X,A;J)-\M_{1,M}^{\{0\}}(X,A;J),
\qquad\hbox{where}\quad A=\sum_{i\in I}\!A_i,$$ 
it is again useful to define a thickening of the set ${\cal U}_{\T}^{(0)}(X;J)$.
Thus, we put
\begin{equation}\label{g1gendecomp_e7}
\ti{\cal U}_{\T}^{(0)}(X;J)=
\big\{\big(b_0,(b_h)_{h\in I_1}\big)\!\in\!{\cal U}_{\T_0}(X;J)\!\times\!\!
\prod_{h\in I_1}\!\!\ti{\cal U}_{\T_h}^{(0)}(X;J)\!:
\ev_0(b_h)\!=\!\ev_{\io_h}(b_0)~\forall h\!\in\!I_1\big\},
\end{equation}
where the space $\ti{\cal U}_{\T_h}^{(0)}(X;J)$ is as in 
Subsection~\ref{notation0_subs}.
Let
\begin{gather*}
\wt{\cal FT}=\wt{{\cal F}_{\aleph}{\cal T}}\!\oplus\!\wt{{\cal F}_0{\cal T}}
\!\oplus\!\wt{{\cal F}_1{\cal T}} \!\lra\!\tilde{\cal U}_{\cal T}^{(0)}(X;J),
\qquad\hbox{where}\\
\wt{{\cal F}_{\aleph}{\cal T}}=\ti\pi_P^*{\cal FT}_0,\quad
\wt{{\cal F}_0{\cal T}}=\bigoplus_{h\in I_1}\wt{{\cal F}_h{\cal T}}, \quad
\wt{{\cal F}_h{\cal T}}=\ti\pi_P^*L_h{\cal T}_0,\quad
\wt{{\cal F}_1{\cal T}}=
\tilde{\cal U}_{\cal T}^{(0)}(X;J)\!\times\!\Bbb{C}^{\hat{I}-I_1},  
\end{gather*}
and $\ti\pi_P\!:\tilde{\cal U}_{\cal T}^{(0)}(X;J)\!\lra\!
{\cal U}_{{\cal T}_0}(X;J)$ is the projection map.
As before, we denote by $\wt{\cal FT}^{\eset}$ the subset of $\wt{\cal FT}$
consisting of the elements with all components nonzero.\\

\noindent
Suppose ${\cal T}\!=\!(M,I,\aleph;j,\under{A})$ is a bubble type as 
in the previous paragraph.
Since every holomorphic map in the zero homology class is constant,
the decomposition \e_ref{g1gendecomp_e3} is equivalent to
\begin{equation}\label{g1decomp_e1}\begin{split}
{\cal U}_{\cal T}(X;J)&\approx
\Big({\cal U}_{{\cal T}_0}(pt)\times
{\cal U}_{\bar{\cal T}}(X;J)\Big)\big/\Aut^*(\T)\\
&\subset\Big(\ov\cM_{1,k_0}\times
{\cal U}_{\bar{\cal T}}(X;J)\Big)\big/\Aut^*(\T),
\end{split}\end{equation}
where $k_0\!=\!|M_0|\!+\!|I_1|$, $\ov\cM_{1,k_0}$
is the moduli space of genus-one curves with $k_0$ marked points, and
$${\cal U}_{\bar{\cal T}}(X;J)=
\big\{(b_h)_{h\in I_1}\!\in\!\prod_{h\in I_1}\!\!{\cal U}_{{\cal T}_h}(X;J)\!:
\ev_0(b_{h_1})\!=\!\ev_0(b_{h_2})~\forall h_1,h_2\!\in\!I_1\big\}.$$
Similarly, \e_ref{g1gendecomp_e7} is equivalent~to
\begin{gather}\label{g1decomp_e3a}
\ti{\cal U}_{\cal T}^{(0)}(X;J)\approx
{\cal U}_{{\cal T}_0}(pt)\times \tilde{\cal U}_{\bar{\cal T}}^{(0)}(X;J)
\subset\ov\cM_{1,k_0}\times \tilde{\cal U}_{\bar{\cal T}}^{(0)}(X;J),
\qquad\hbox{where}\\
\label{g1decomp_e3b}
\tilde{\cal U}_{\bar{\cal T}}^{(0)}(X;J)=
\big\{(b_h)_{h\in I_1}\!\in\!\prod_{h\in I_1}\!\!\tilde{\cal U}_{{\cal T}_h}^{(0)}(X;J)\!:
\ev_0(b_{h_1})\!=\!\ev_0(b_{h_2})~\forall h_1,h_2\!\in\!I_1\big\}.
\end{gather}
We denote by
\begin{gather*}
\pi_P\!:{\cal U}_{\cal T}(X;J)\lra\ov\cM_{1,k_0}/\Aut^*(\T), \qquad
\ti\pi_P\!:\ti{\cal U}_{\cal T}^{(0)}(X;J)\lra\ov\cM_{1,k_0}, \\
\hbox{and}\qquad 
\ev_P\!:\U_{\T}(X;J),\, \ti\U_{\T}^{(0)}(X;J)\lra X
\end{gather*}
the projections onto the first component in the decompositions~\e_ref{g1decomp_e1}
and~\e_ref{g1decomp_e3a}
and the map sending each element $[{\cal C},u]$ of ${\cal U}_{\cal T}(X;J)$,
or $({\cal C},u)$ of $\tilde{\cal U}_{\cal T}^{(0)}(X;J)$,
to the image of the principal component $\cC_P$ of ${\cal C}$, 
i.e.~the point $u(\cC_P)$ in~$X$.\\

\noindent
Let $\E\!\lra\!\ov\cM_{1,k_0}$ denote the Hodge line bundle,
i.e.~the line bundle of holomorphic differentials.
For each $i\!\in\!\chi({\cal T})$, we define the bundle map 
$${\cal D}_{J,i}\!:\wt{{\cal F}_{h(i)}{\cal T}}\lra 
\ti\pi_P^*\Bbb{E}^*\otimes \ev_P^*TX,
\qquad\hbox{where}\qquad
h(i)\!=\!\min\{h\!\in\!\hat{I}\!:h\!\le\!i\}\in I_1,$$
over $\tilde{\cal U}_{\cal T}^{(0)}(X;J)$ by
\begin{gather*}
\big\{{\cal D}_{J,i}(\tilde{\ups})\big\}(\psi)=
\psi_{x_{h(i)}(b)}(\tilde{v})\cdot_J{\cal D}_ib \in T_{\ev_P(b)}X
\qquad\hbox{if}\\
\psi\!\in\!\ti\pi_P^*\Bbb{E}, \quad
\tilde{\ups}\!=\!(b,\tilde{v})\!\in\!\wt{{\cal F}_{h(i)}{\cal T}},
\quad
~b\!\in\!\tilde{\cal U}_{\cal T}^{(0)}(X;J),
\end{gather*}
and $x_{h(i)}(b)\!\in\!\Si_{b,\aleph}$ is the node joining the bubble $\Si_{b,h(i)}$
of $b$ to the principal component $\Si_{b,\aleph}$ of~$\Si_b$.
For each $\ups\!\in\!\wt{\cal FT}$, we~put
\begin{gather*}
\rho(\ups)\!=\!\big(b,(\rho_i(v))_{i\in\chi({\cal T})}\big)
\!\in\!\!\bigoplus_{i\in\chi({\cal T})}\!\!\!\wt{{\cal F}_{h(i)}{\cal T}},
\quad\hbox{where}\quad
\rho_i(v)\!=\!\!\prod_{h\in\hat{I},h\le i}\!\!\!\!\!\!v_h\in\wt{{\cal F}_{h(i)}{\cal T}},
\qquad\hbox{if}\\
\ups\!=\!\big(b;v_{\aleph},(v_i)_{i\in\hat{I}}\big),~~
b\!\in\!\tilde{\cal U}_{\cal T}^{(0)}(X;J),~~
(b,v_{\aleph})\!\in\!\wt{{\cal F}_{\aleph}{\cal T}},~~
(b,v_h)\!\in\!\wt{{\cal F}_h{\cal T}}~\hbox{if}~h\!\in\!I_1,~~
v_i\!\in\!\Bbb{C}~\hbox{if}~i\!\in\!\hat{I}\!-\!I_1.
\end{gather*}
These definitions are illustrated in Figure~\ref{str_fig}.
While the restrictions of these bundle maps to 
${\cal U}_{\cal T}^{(0)}(X;J)\!\subset\!\tilde{\cal U}_{\cal T}^{(0)}(X;J)$
do not necessarily descend to the vector bundle ${\cal FT}$
over ${\cal U}_{\cal T}(X;J)$, the~map
\begin{gather*}
{\cal D}_{\cal T}\!:{\cal FT}\lra\pi_P^*\Bbb{E}^*\otimes\ev_P^*TX
\big/\Aut^*(\T), \qquad
{\cal D}_{\cal T}(\ups)=\sum_{i\in\chi({\cal T})}\!\!\!{\cal D}_{J,i}\rho_i(\ups),
\end{gather*}
is well-defined.\\

\noindent
Finally, if ${\cal T}$ is any bubble type, for genus-zero or genus-one maps,
and $K$ is a subset of ${\cal U}_{\cal T}(X;J)$,
we denote by $K^{(0)}$ and $\tilde{K}^{(0)}$
the preimages of $K$ under the quotient projection maps
$${\cal U}_{\cal T}^{(0)}(X;J)\!\lra\!{\cal U}_{\cal T}(X;J)
\qquad\hbox{and}\qquad
\tilde{\cal U}_{\cal T}^{(0)}(X;J)\!\lra\!{\cal U}_{\cal T}(X;J),$$
respectively.
All vector orbi-bundles we encounter will be assumed to be normed.
Some will come with natural norms; for others, we choose a norm,
sometimes implicitly, once and for~all.
If \hbox{$\pi_{\frak F}\!:{\frak F}\!\lra\!{\frak X}$} is a normed vector bundle
and $\de\!:{\frak X}\!\lra\!\Bbb{R}$ is any function, possibly constant,
let
$${\frak F}_{\de}=\big\{\ups\!\in\!{\frak F}\!:
|\ups|\!<\!\de(\pi_{\frak F}(\ups))\big\}.$$
If $\Om$ is any subset of ${\frak F}$, we take 
$\Om_{\de}\!=\!\Om\cap{\frak F}_{\de}$.

\subsection{Boundary Structure Theorem}
\label{str_subs}

\noindent
In this subsection, we formulate Theorem~\ref{str_thm},
which states that an element
$$b\in\ov\M_{1,k}(X,A;J)-\M_{1,k}^{\{0\}}(X,A;J)$$
lies in the stable-map closure of the space $\M_{1,k}^0(X,A;J)$
of genus-one $J$-holomorphic maps from smooth domains
if and only if $b$ lies in $\ov\M_{1,k}^0(X,A;J)$,
provided the almost complex structure~$J$ is sufficiently regular.
In addition, Theorem~\ref{str_thm} describes a neighborhood of  every 
stratum~of 
$$\ov\M_{1,k}^0(X,A;J)-\M_{1,k}^{\{0\}}(X,A;J)$$
in $\ov\M_{1,k}^0(X,A;J)$.
If $k\!\in\!\bar{\Z}^+$, we denote by $[k]$ the set $\{1,\ldots,k\}$.

\begin{thm}
\label{str_thm}
Suppose $(X,\om,J)$ is a compact almost Kahler manifold and $A\!\in\!H_2(X;\Z)^*$.
If the regularity conditions~(a) and~(b-i) of Definition~\ref{g1reg_dfn} are satisfied
and ${\cal T}\!=\!([k],I,\aleph;j,\under{A})$ is a bubble type such that
$\sum_{i\in I}\!A_i\!=\!A$ and $A_i\!=\!0$ for all minimal elements $i$ of~$I$,
then the intersection of the closure of $\M_{1,k}^0(X,A;J)$ in 
$\ov\M_{1,k}(X,A;J)$  with ${\cal U}_{\cal T}(X;J)$
is the~set
$${\cal U}_{\T;1}(X;J)\equiv
\big\{[b]\!\in\!{\cal U}_{\cal T}(X;J)\!:
\dim_{\Bbb{C}}\hbox{Span}_{(\Bbb{C},J)}
\{{\cal D}_ib\!:i\!\in\!\chi({\cal T})\}
\!<\!|\chi({\cal T})|\big\}.$$
Furthermore, the space
$${\cal F}^1{\cal T}^{\eset}\equiv
\big\{[\ups]\!=\![b,v]\!\in\!{\cal FT}^{\eset}\!:
{\cal D}_{\cal T}(\ups)\!=\!0\big\}$$
is a smooth oriented suborbifold of ${\cal FT}$.
Finally, there exist $\de\!\in\!C({\cal U}_{\cal T}(X;J);\R^+)$,
an open neighborhood $U_{\cal T}$ of ${\cal U}_{\cal T}(X;J)$
in $\X_{1,k}(X,A)$, and
an orientation-preserving diffeomorphism
$$\phi\!:{\cal F}^1{\cal T}^{\eset}_{\de}\lra
\M_{1,k}^0(X,A;J)\cap U_{\cal T},$$
which extends to a homeomorphism 
$$\phi\!:{\cal F}^1{\cal T}_{\de}\lra\ov\M_{1,k}^0(X,A;J)\cap U_{\cal T},$$
where ${\cal F}^1{\cal T}$ is the closure of ${\cal F}^1{\cal T}^{\eset}$
in ${\cal FT}$.
\end{thm}

\noindent
We now clarify the statement of Theorem~\ref{str_thm}
and illustrate it using Figure~\ref{str_fig}.
As before, the shaded discs represent the components of the domain
on which every stable map~$[b]$ in ${\cal U}_{\cal T}(X;J)$ is non-constant.
A stable map 
$$[{\cal C},u]\in{\cal U}_{\cal T}(X;J)
\subset \ov\M_{1,k}(X,A;J)-\M_{1,k}^{\{0\}}(X,A;J),$$
is in the stable-map closure of $\M_{1,k}^0(X,A;J)$ if and only if
$[{\cal C},u]$ satisfies condition (b) of Definition~\ref{degen_dfn}.\\

\begin{figure}
\begin{pspicture}(-1.1,-1.8)(10,1.25)
\psset{unit=.4cm}
\psellipse(5,-1.5)(1.5,2.5)
\psarc[linewidth=.05](3.2,-1.5){2}{-30}{30}\psarc[linewidth=.05](6.8,-1.5){2}{150}{210}
\pscircle[fillstyle=solid,fillcolor=gray](2.5,-1.5){1}\pscircle*(3.5,-1.5){.2}
\pscircle[fillstyle=solid,fillcolor=gray](.5,-1.5){1}\pscircle*(1.5,-1.5){.2}
\pscircle(7.5,-1.5){1}\pscircle*(6.5,-1.5){.2}
\pscircle[fillstyle=solid,fillcolor=gray](8.91,-.09){1}\pscircle*(8.21,-.79){.2}
\pscircle[fillstyle=solid,fillcolor=gray](8.91,-2.91){1}\pscircle*(8.21,-2.21){.2}
\rput(2.5,0){$h_1$}\rput(.5,0){$h_2$}\rput(7.3,0){$h_3$}
\rput(10.5,0.1){$h_4$}\rput(10.5,-2.9){$h_5$}
\rput(5,-5){\small ``tacnode"}
\pnode(5,-5){A1}\pnode(3.5,-1.5){B1}
\ncarc[nodesep=.35,arcangleA=-25,arcangleB=-15,ncurv=1]{->}{A1}{B1}
\pnode(5,-4.65){A2}\pnode(7.3,-1.5){B2}
\ncarc[nodesep=0,arcangleA=40,arcangleB=30,ncurv=1]{-}{A2}{B2}
\pnode(8,-.95){B2a}\pnode(8.02,-2.02){B2b}
\ncarc[nodesep=0,arcangleA=0,arcangleB=10,ncurv=1]{->}{B2}{B2a}
\ncarc[nodesep=0,arcangleA=0,arcangleB=10,ncurv=1]{->}{B2}{B2b}
\rput(25,-1.5){\begin{small}\begin{tabular}{l}
$\chi(\T)\!=\!\{h_1,h_4,h_5\}$,~
$\rho(\ups)\!=\!(\ups_{h_1},\ups_{h_3}\ups_{h_4},\ups_{h_3}\ups_{h_5})$\\
\\
${\cal F}^1\T^{\eset}=\big\{[b;v_{h_1},v_{h_2},v_{h_3},v_{h_4},v_{h_5}]\!:
v_{h_2},v_{h_4},v_{h_5}\!\in\!\Bbb{C}^*$;\\
${}\qquad\qquad\qquad v_{h_1}\!\in\!T_{x_{h_1}}\Si_P\!-\!\{0\},~
v_{h_3}\!\in\!T_{x_{h_3}}\Si_P\!-\!\{0\}$;\\
${}\qquad\qquad\qquad v_{h_1}\cD_{J,h_1}b\!+\!v_{h_3}v_{h_4}\cD_{J,h_4}b\!+\!
v_{h_3}v_{h_5}\cD_{J,h_5}b\!=\!0\big\}$
\end{tabular}\end{small}}
\end{pspicture}
\caption{An Illustration of Theorem~\ref{str_thm}} 
\label{str_fig}
\end{figure}

\noindent
Standard arguments show that the regularity condition (a) of Definition~\ref{g1reg_dfn} 
implies that the space ${\cal U}_{\cal T}^{(0)}(X;J)$ is a smooth manifold,
while ${\cal U}_{\cal T}(X;J)$ is a smooth orbifold; 
see Chapter~3 in~\cite{McSa}, for example.
Thus, the total space of the bundle ${\cal FT}^{\eset}$ is also a smooth orbifold.
The second claim of Theorem~\ref{str_thm} is immediate from
the transversality of the bundle~map
$${\cal D}_{\cal T}\!:{\cal FT}^{\eset}\lra\pi_P^*\Bbb{E}^*\otimes\ev_P^*TX
\big/\Aut^*(\T),  \qquad
{\cal D}_{\cal T}(\ups)=\!\sum_{i\in\chi({\cal T})}\!\!\!{\cal D}_{J,i}\rho_i(\ups),$$
to the zero set.
In turn, this transversality property is an immediate consequence of 
the regularity conditions~(a) and~(b-i) of Definition~\ref{g1reg_dfn}.\\

\noindent
The middle claim of Theorem~\ref{str_thm} is needed 
to make sense of the remaining statement.
This final claim, proved in Section~\ref{strthm_sec}, 
describes a normal neighborhood of ${\cal U}_{\T;1}(X;J)$
in~$\ov\M_{1,k}^0(X,A;J)$ and implies the first statement of  Theorem~\ref{str_thm}.\\

\noindent
{\it Remark:}  The regularity assumptions on $J$ used in Theorem~\ref{str_thm}
do not guarantee that the entire space $\M_{1,k}^0(X,A;J)$ is smooth.
However, the proof of Theorem~\ref{reg_thm} implies that
$\M_{1,k}^0(X,A;J)$ is smooth near each stratum ${\cal U}_{{\cal T};1}(X;J)$
of $\ov\M_{1,k}^0(X,A;J)$.
This can be seen from the $\tilde{J}\!=\!J$ case of Corollary~\ref{reg1_crl2}
and standard Implicit Function Theorem arguments such as those in Chapter~3 of~\cite{McSa}.\\

\noindent
{\it Proof of Corollary~\ref{str_crl}:}
It remains to construct a fundamental class for $\ov\M_{1,k}^0(X,A;J)$.
Theorem~\ref{str_thm} describes a neighborhood in $\ov\M_{1,k}^0(X,A;J)$
of every stratum $\ov\M_{1,k}^0(X,A;J)\cap{\cal U}_{\cal T}(X;J)$
for a bubble type ${\cal T}\!=\!(M,I,\aleph;j,\under{A})$ such that
$A_i\!=\!0$ for all minimal elements $i\!\in\!I$.
If ${\cal T}$ is a bubble type such that $A_i\!\neq\!0$ for some minimal element $i\!\in\!I$,
a neighborhood~of 
$$\ov\M_{1,k}^0(X,A;J)\cap{\cal U}_{\cal T}(X;J)={\cal U}_{\cal T}(X;J)$$
in $\ov\M_{1,k}^0(X,A;J)$ is homeomorphic to a neighborhood of ${\cal U}_{\cal T}(X;J)$
in the corresponding bundle of gluing parameters~${\cal FT}$,
as can be seen from Subsection~\ref{reg1_subs1} and the continuity arguments
of Subsections~3.9 and~4.1 in~\cite{Z4}.
It follows that there exist arbitrary small neighborhoods $U$~of 
$$\partial\ov\M_{1,k}^0(X,A;J)\equiv \ov\M_{1,k}^0(X,A;J)-\M_{1,k}^0(X,A;J)$$
such that
$$H_l\big(U;\Q)=\{0\}\qquad\forall~k\ge 2(\lr{c_1(TX),A}\!+\!k)-1.$$
Since the moduli space $\M_{1,k}^0(X,A;J)$ is a smooth oriented orbifold,
$$\dim_{\R}\M_{1,k}^0(X,A;J)=2(\lr{c_1(TX),A}\!+\!k),$$
and the complement of $U$ in $\ov\M_{1,k}^0(X,A;J)$ is compact,
$\M_{1,k}^0(X,A;J)$ determines a~class
\begin{equation*}\begin{split}
\big[\ov\M_{1,k}^0(X,A;J)\big]&\in 
H_{2(\lr{c_1(TX),A}+k)}\big(\ov\M_{1,k}^0(X,A;J),U;\Q)\\
&\qquad\qquad \approx
H_{2(\lr{c_1(TX),A}+k)}\big(\ov\M_{1,k}^0(X,A;J);\Q),
\end{split}\end{equation*}
as claimed.
The isomorphism between the two homology groups is induced by inclusion.

\section{A Genus-Zero Gluing Procedure}
\label{gluing0_sec}

\subsection{The Genus-Zero Regularity Properties}
\label{reg0_subs1}

\noindent
In this subsection, we prove the $g\!=\!0$ case of the first claim of Theorem~\ref{reg_thm}.
It follows from Corollary~\ref{reg0_crl1} and the compactness of
the moduli space $\ov\M_{0,1}(X,A;J)$.
Corollary~\ref{reg0_crl1} is obtained by a rather straightforward 
argument via the analytic part of~\cite{LT}.
Throughout this subsection, we assume that $J$ is a genus-zero $A$-regular 
almost complex structure on~$X$.\\

\noindent
In order to prove Theorem~\ref{reg_thm}, 
we need to describe smooth maps $u\!:\Bbb{P}^1\!\lra\!X$, with one or two marked points, 
that lie close to each stratum ${\cal U}_{\cal T}(X;J)$ of $\ov\M_{0,1}(X,A;J)$ 
and~of
$$\ov\M_{0,2}(X,A;J)\approx \ov\M_{0,\{0\}\sqcup\{1\}}(X,A;J).$$
We denote by $\X_{0,M}(X,A)$ the space of equivalence classes
of all smooth maps into~$X$ from genus-zero Riemann surfaces 
with marked points indexed by the set $\{0\}\!\sqcup\!M$
in the homology class~$A$ and by
$\X_{0,M}^0(X,A)$ the subset of $\X_{0,M}(X,A)$
consisting of the maps with smooth domains, i.e.~$\Bbb{P}^1$ in this case.\\

\noindent
Let ${\cal T}\!=\!(M,I;j,\under{A})$ be a bubble type such that 
$\sum_{i\in I}\!A_i\!=\!A$,
i.e.~${\cal U}_{\cal T}(X;J)$ is a stratum of the moduli space
$\ov\M_{0,\{0\}\sqcup M}(X,A;J)$.
We will proceed as in Subsections~3.3 and~3.6 of~\cite{Z4}.
Subsections~2.1 and~2.3 in~\cite{Z3} describe a special case of 
the same construction in circumstances very similar to the present situation.\\

\noindent
For each sufficiently small element $\ups\!=\!(b,v)$ of $\wt{\cal FT}^{\eset}$,
where $b\!=\!(\Si_b,u_b)$ is an element of $\tilde{\cal U}_{\cal T}^{(0)}(X;J)$, let 
$$q_{\ups}\!:\Si_{\ups}\lra\Si_b$$
be the basic gluing map constructed in Subsection~2.2 of~\cite{Z4}.
In this case, $\Si_{\ups}$ is the projective line $\Bbb{P}^1$ with $|M|\!+\!1$ marked points.
Let
$$b(\ups)=\big(\Si_{\ups},u_{\ups}\big),
\qquad\hbox{where}\qquad      u_{\ups}=u_b\circ q_{\ups},$$
be the approximately holomorphic map corresponding to~$\ups$.
The primary marked point $y_0(\ups)$ of $\Si_{\ups}$ is 
the point~$\i$ of $\Si_{\ups}\!\approx\!S^2$.\\

\noindent
Let $\na^J$ be the $J$-linear connection induced by the Levi-Civita connection of
the metric~$g_X$.
Since the linearization~$D_{J,b}$ of the $\bar{\partial}_J$-operator at~$b$ is surjective
by Definition~\ref{g0reg_dfn},
if $\ups\!\in\!\wt{\cal FT}^{\eset}$ is sufficiently small, the linearization
$$D_{J,\ups}\!:\Ga(\ups)\!\equiv\!L^p_1(\Si_{\ups};u_{\ups}^*TX)
\lra \Ga^{0,1}(\ups;J)\!\equiv\!
L^p(\Si_{\ups};\La^{0,1}_{J,j}T^*\Si_{\ups}\!\otimes\!u_{\ups}^*TX)$$
of the $\bar{\partial}_J$-operator at~$b(\ups)$,
defined via the connection $\na^J$, is also surjective.
In particular, we can obtain a decomposition
\begin{equation}\label{gadecomp_e}
\Ga(\ups)=\Ga_-(\ups)\oplus\Ga_+(\ups)
\end{equation}
such that the linear operator 
$D_{J,\ups}\!:\Ga_+(\ups)\!\lra\!\Ga^{0,1}(\ups;J)$ is an isomorphism,
while 
$$\Ga_-(\ups)=\big\{\xi\circ q_{\ups}\!:
\xi\!\in\!\Ga_-(b)\!\equiv\!\ker D_{J,b}\big\}.$$
For the purposes of this subsection,
the space $\Ga_+(\ups)$ can be taken to be the $L^2$-orthogonal complement
of $\Ga_-(\ups)$, but for use in later subsections it is more convenient to~take
\begin{equation}\label{gaplusdfn_e}
\Ga_+(\ups)=\big\{\ze\!\in\!\Ga(\ups)\!:\ze(\hat{0},\i)\!=\!0;~
\llrr{\ze,\xi}_{\ups,2}\!=\!0~\forall\xi\!\in\!\Ga_-(\ups)
~\hbox{s.t.}~\xi(\hat{0},\i)\!=\!0\big\},
\end{equation}
where $(\hat{0},\i)$ is the primary marked point, i.e.~the south pole 
of the sphere $\Si_{\ups}\!\approx\!S^2$.
This choice of $\Ga_+(\ups)$ is permissible by Definition~\ref{g0reg_dfn}.
The $L^2$-inner product on~$\Ga(\ups)$ used in~\e_ref{gaplusdfn_e}
is defined via the metric~$g_X$ on~$X$ and 
the metric~$g_{\ups}$ on~$\Si_{\ups}$ induced by the pregluing construction.
The Banach spaces $\Ga(\ups)$ and $\Ga^{0,1}(\ups;J)$ carry the norms 
$\|\cdot\|_{\ups,p,1}$ and $\|\cdot\|_{\ups,p}$, respectively,
which are also defined by the pregluing construction.
Throughout this paper, $p$ denotes a real number greater than two.
The norms $\|\cdot\|_{\ups,p,1}$ and $\|\cdot\|_{\ups,p}$
are equivalent to the ones used in~\cite{LT}.
In particular, the norms of $D_{J,\ups}$ and of the inverse of its restriction
to $\Ga_+(\ups)$ have fiberwise uniform upper bounds, 
i.e.~dependent only on $[b]\!\in\!{\cal U}_{\cal T}(X;J)$,
and not on~$v\!\in\!\Bbb{C}^{*\hat{I}}$.
We denote by 
$$\pi_{\ups;-}\!:\Ga(\ups)\lra\Ga_-(\ups) \qquad\hbox{and}\qquad
\pi_{\ups;+}\!:\Ga(\ups)\lra\Ga_+(\ups)$$
the projection maps corresponding to the decomposition~\e_ref{gadecomp_e}.
The relevant facts concerning the objects described in this paragraph
are summarized in Lemma~\ref{reg0_lmm1}:

\begin{lmm}
\label{reg0_lmm1}
Suppose $(X,\om,J)$ is a compact almost Kahler manifold and $A\!\in\!H_2(X;\Bbb{Z})$.
If $J$ is a genus-zero $A$-regular almost complex structure and 
${\cal T}\!=\!(M,I;j,\under{A})$ is a bubble type such that $A\!=\!\sum_{i\in I}\!A_i$,
there exist $\de,C\!\in\!C({\cal U}_{\cal T}(X;J);\Bbb{R}^+)$ and 
an open neighborhood $U_{\cal T}$ of ${\cal U}_{\cal T}(X;J)$ in $\X_{0,M}(X,A)$
with the following properties:\\
(1) for all $\ups\!=\!(b,v)\!\in\!\wt{\cal FT}_{\de}^{\eset}$,
\begin{gather*}
\|\pi_{\ups;-}\xi\|_{\ups,p,1}\le C(b)\|\xi\|_{\ups,p,1}
\quad\forall  \xi\!\in\!\Ga(\ups), \qquad
\|D_{J,\ups}\xi\|_{\ups,p}\le C(b)|\ups|^{1/p}\|\xi\|_{\ups,p,1} 
\quad\forall\xi\!\in\!\Ga_-(\ups),\\
\hbox{and}\qquad
C(b)^{-1}\|\xi\|_{\ups,p,1}\le \|D_{J,\ups}\xi\|_{\ups,p} \le C(b)\|\xi\|_{\ups,p,1}
\quad\forall  \xi\!\in\!\Ga_+(\ups);
\end{gather*}
(2) for every $[\tilde{b}]\!\in\!\X_{0,M}^0(X,A)\cap U_{\cal T}$, 
there exist $\ups\!\in\!\tilde{\cal F}{\cal T}_{\de}^{\eset}$ and 
$\ze\!\in\!\Ga_+(\ups)$ such that $\|\ze\|_{\ups,p,1}\!<\!\de(b)$
and $[\exp_{b(\ups)}\!\ze]\!=\![\tilde{b}]$.
\end{lmm}

\noindent
The first two bounds in~(1) follow immediately from 
the definition of the spaces~$\Ga_-(\ups)$.
The third estimate can be deduced from the facts~that
\begin{gather}\label{reg0_lmm1e1}
\|\xi\|_{\ups,p,1}\le C(b)
\big(\|D_{J,\ups}\xi\|_{\ups,p}\!+\!\|\xi\|_{\ups,p}\big)
\quad\hbox{and}\quad
\|\xi\|_{C^0}\le C(b)\|\xi\|_{\ups,p,1} 
\qquad\forall\xi\!\in\!\Ga(\ups),\\
\hbox{and}\qquad \lim_{\ups\lra b}\Ga_-(\ups)=\Ga_-(b)
\quad\hbox{if}\quad 
b\!=\!(\Si_b,u_b)\in\tilde{\cal U}_{\cal T}^{(0)}(X;J);\notag
\end{gather}
see Subsection~3.5 in~\cite{Z4}.
In (2) of Lemma~\ref{reg0_lmm1}, $\exp_{b(\ups)}\!\ze$ denotes the stable map 
that has the same domain and marked points as the map~$b(\ups)$,
but the map into~$X$ is $\exp_{u_{\ups}}\!\ze$,
where $\exp$ is the exponential map of the connection~$\na^J$.
The final claim of Lemma~\ref{reg0_lmm1} also follows from 
the above properties of $\Ga_-(\ups)$,
along with the uniformly smooth dependence of the spaces~$\Ga_-(\ups)$ on~$\ups$;
see Section~4 of~\cite{Z4}.
In fact, for each $[\tilde{b}]$ in $U_{\cal T}\cap\X_{0,M}^0(X,A)$,
the corresponding pair $(\ups,\ze)$ is unique, up to 
the action of the group $\Aut(\T)\!\propto\!(S^1)^I$.

\begin{crl}
\label{reg0_crl1}
If $(X,\om,J)$, $A$, and ${\cal T}$ are as in Lemma~\ref{reg0_lmm1} and $M\!=\!\eset$, 
for every precompact open subset $K$ of ${\cal U}_{\cal T}(X;J)$,
there exist $\de_K,C_K\!\in\!\Bbb{R}^+$
and an open neighborhood $U_K\!\subset\!U_{\cal T}$ of $K$ in $\X_{0,\eset}(X,A)$
with the following properties:\\
(1) requirements (1) and (2) of Lemma~\ref{reg0_lmm1} are satisfied;\\
(2) if $\tilde{J}$ is an almost complex structure on $X$ 
s.t.~$\|\tilde{J}\!-\!J\|_{C^1}\!<\!\de_K$ and 
$[\tilde{b}]\!\in\!U_K\!\cap\!\X_{0,1}^0(X,A)$, there exists 
a smooth map $\tilde{u}\!:\Bbb{P}^1\!\lra\!X$ such that
$[\tilde{b}]\!=\![\Bbb{P}^1,\tilde{u}]$ and,
for a choice of linearization of $\bar{\partial}_{\tilde{J}}$ at~$\tilde{u}$,
the operators $D_{\tilde{J},\tilde{u}}$ and $\D_{\tilde{J},\tilde{u}}^{\i}$ are surjective.
\end{crl}

\noindent
{\it Remark:} If the map $\tilde{u}$ is $\tilde{J}$-holomorphic, 
i.e.~$\bar{\partial}_{\tilde{J}}$ vanishes at~$\tilde{u}$, 
there is only one linearization of $\bar{\partial}_{\tilde{J}}$ at~$\tilde{u}$, 
though there are different ways of writing it explicitly.
In the proof of this corollary, 
whether or not $\tilde{u}$ is a $\tilde{J}$-holomorphic map, 
$D_{\tilde{J},\tilde{u}}$ denotes the linearization of 
$\bar{\partial}_{\tilde{J}}$ at~$\tilde{u}$ with respect to the connection~$\na^J$;
see Chapter~3 in~\cite{McSa}.\\

\noindent
{\it Proof:} (1) By (2) of Lemma~\ref{reg0_lmm1}, it is sufficient to check 
the surjectivity claims for every smooth map $\ti{u}\!=\!\exp_{u_{\ups}}\!\ze$,
where $\ups\!=\!(b,v)\!\in\!\ti{\cal F}{\cal T}_{\de_K}^{\eset}\big|_{K^{(0)}}$
and $\|\ze\|_{\ups,p,1}\!<\!\de_K$. 
If 
$$\xi\in\Ga(\tilde{u})\equiv L^p_1(\Si_{\ups};\tilde{u}^*TX),$$
we define $\tilde{\xi}\!\in\!\Ga(\ups)$ by 
$$\tilde{\xi}(z)=\Pi_{\ze(z)}^{-1}\xi(z)
\qquad\forall~ z\!\in\!\Bbb{P}^1,$$
where $\Pi_{\ze(z)}$ is the parallel transport along the geodesic
$t\!\lra\!\exp_{u_{\ups}(z)}\!t\ze(z)$ with respect to~$\na^J$.
By~\e_ref{reg0_lmm1e1},
\begin{equation}\label{reg0_crl1e1}\begin{split}
\|D_{J,\ups}\ti{\xi}\|_{\ups,p} 
&\le  \|D_{\ti{J},\ti{u}}\xi\|_{\ups,p}+
C_K\big(\|\ti{J}\!-\!J\|_{C^1}\!+\!\|\ze\|_{\ups,p,1}\big) \|\xi\|_{\ups,p,1}\\
&=\|D_{\ti{J},\ti{u}}\xi\|_{\ups,p}+
C_K\big(\|\ti{J}\!-\!J\|_{C^1}\!+\!\|\ze\|_{\ups,p,1}\big)\|\ti{\xi}\|_{\ups,p,1}
\qquad\forall~\xi\!\in\!\Ga(\ti{u}).
\end{split}
\end{equation}
Thus, by (1) of Lemma~\ref{reg0_lmm1},
\begin{equation}\label{reg0_crl1e3}\begin{split}
\|\pi_{\ups;+}\ti{\xi}\|_{\ups,p,1}
&\le C_K'\big(\|\ti{J}\!-\!J\|_{C^1}\!+\!\|\ze\|_{\ups,p,1}\big)\|\ti{\xi}\|_{\ups,p,1}\\
\Lra\qquad  \|\ti{\xi}\|_{\ups,p,1}&\le C_K\|\pi_{\ups;-}\ti{\xi}\|_{\ups,p,1}
\end{split}
\qquad\forall~ \xi\!\in\!\ker D_{\ti{J},\ti{u}},
\end{equation}
if $\de_K$ is sufficiently small. 
By \e_ref{reg0_crl1e3} and (a) of Definition~\ref{g0reg_dfn},
$$\dim\ker D_{\ti{J},\ti{u}}\le \dim\Ga_-(\ups)
=\ind D_{J,b}=\ind D_{\ti{J},\ti{u}}.$$
In particular, the operator $D_{\ti{J},\ti{u}}$ is surjective.\\
(2) The surjectivity of the map $\D_{\ti{J},\ti{u}}^{\i}$ is proved similarly.
Let 
\begin{alignat*}{1}
&\ti{\pi}_{\ups;-}\!:\Ga_-(\ups)\lra
\ti{\Ga}_-(\ups)\!\equiv\!
\{\xi\!\in\!\Ga_-(\ups)\!:\xi(\i)\!=\!0\}\approx\ker\D_{J,u_b}^{\i}
\quad\hbox{and}\\
&\ti{\pi}_{\ups;+}\!:\Ga_-(\ups)\lra
\ti{\Ga}_+(\ups)\!\equiv\!\{\xi\!\in\!\Ga_-(\ups)\!:
\llrr{\xi,\xi'}_{\ups,2}\!=\!0~\forall\xi'\!\in\!\ti{\Ga}_-(\ups)\}
\end{alignat*}
be the $L^2$-orthogonal projections onto $\ti{\Ga}_-(\ups)$
and its orthogonal complement in $\Ga_-(\ups)$.
Then,
\begin{equation}\label{reg0_crl1e5}
\|\xi\|_{\ups,p,1}\le C_K|\xi(\i)|
\qquad\forall~\xi\!\in\!\ti{\Ga}_+(\ups),
\end{equation}
since the analogous bound holds for the map 
$$\D_{J,b}^{\i}\!:\ker D_{J,u}\lra T_{u(\i)}X,$$
by Definition~\ref{g0reg_dfn}.
Combining \e_ref{reg0_lmm1e1}-\e_ref{reg0_crl1e5}, we obtain
\begin{gather*}\begin{split}
\|\tilde{\pi}_{\ups;+}\pi_{\ups;-}\tilde{\xi}\|_{\ups,p,1}
&\le C_K|\pi_{\ups;-}\tilde{\xi}(\i)|
\le C_K'\big(|\tilde{\xi}(\i)|\!+\!|\pi_{\ups;+}\tilde{\xi}(\i)|\big)\\
&\le C_K''\big(|\xi(\i)|+(\|\tilde{J}\!-\!J\|_{C^1}\!+\!\|\ze\|_{\ups,p,1})
\|\tilde{\xi}\|_{\ups,p,1}\big)\qquad\forall~ \xi\!\in\!\ker D_{\tilde{J},\tilde{u}}
\end{split}\\
\Lra\qquad  \|\tilde{\xi}\|_{\ups,p,1}\le C_K
\|\tilde{\pi}_{\ups;-}\pi_{\ups;-}\tilde{\xi}\|_{\ups,p,1}
\qquad\qquad\forall~\xi\!\in\!\ker\D_{\tilde{J},\tilde{u}}^{\i},
\end{gather*}
if $\de_K$ is sufficiently small. Thus,
\begin{equation}\label{reg0_crl1e11}
\dim\ker \D_{\ti{J},\ti{u}}^{\i}\le \dim\tilde{\Ga}_-(\ups)
=\ind\D_{J,b}^{\i}=\ind\D_{\tilde{J},\tilde{u}}^{\i},
\end{equation}
and the operator $\D_{\tilde{J},\tilde{u}}^{\i}$ is surjective.

\subsection{Some Power-Series Expansions}
\label{powerseries_subs}

\noindent
In Subsection~2.5 of~\cite{Z3} we describe the behavior of all derivatives of 
rational $J_0$-holomorphic maps into~$\P$ near each stratum $\U_{\T}(\P;J_0)$
by making use of special properties of the standard complex structure~$J_0$ on~$\P$.
In this subsection, we obtain analogous estimates for modified derivatives
of $J$-holomorphic maps into $X$ for an arbitrary genus-zero $A$-regular
almost complex structure~$J$; see Lemma~\ref{reg0_lmm4}.
We use these estimates a number of times in the rest of the paper.\\

\noindent
If $b\!=\!(\Si_b,u_b)$ is as element of $\tilde{\cal U}_{\cal T}^{(0)}(X;J)$ as
in the previous subsection,
the tangent bundle $T_b\tilde{\cal U}_{\cal T}^{(0)}(X;J)$ of 
$\tilde{\cal U}_{\cal T}^{(0)}(X;J)$ at~$b$ consists of the pairs $(\under{w},\xi)$, 
where $\xi\!\in\!\ker D_{J,b}$ and $\under{w}\!\in\!\Bbb{C}^{\hat{I}}$ encode 
the change in~$u_b$ and in the position of the nodes on~$\Si_b$, respectively,
that satisfy a certain balancing condition; see Subsection~2.5 in~\cite{Z4}.
We denote by $\ti{T}_b\ti{\cal U}_{\T}^{(0)}(X;J)$ 
the subspace of the tuples $(0,\xi)$ of $T_b\tilde{\cal U}_{\cal T}^{(0)}(X;J)$.
In particular,
$$\tilde{T}_b\tilde{\cal U}_{\cal T}^{(0)}(X;J)\subset
\Ga_-(b)\!\equiv\!\big\{(\xi_h)_{h\in I}\!\in\!
\bigoplus_{h\in I}\ker D_{J,u_{b,h}}\!\!:
\xi_h(\i)\!=\!\xi_{\io_h}(x_h(b))~\forall h\!\in\!\hat{I}\big\},$$
where $u_{b,h}\!=\!u_b|_{\Si_{b,h}}$.
If $i\!\in\!\chi({\cal T})$, where $\chi({\cal T})$ is as in \e_ref{chidfn0_e},
the image of the projection map
\begin{gather*}
\big\{(\xi_h)_{h\in I}\!\in\!\tilde{T}_b\tilde{\cal U}_{\cal T}^{(0)}(X;J)\!:
\xi_i(\hat{0},\i)\!=\!0\big\}\lra 
\big\{\ze\!\in\!\ker D_{J,u_{b,i}}\!:\ze(\hat{0},\i)\!=\!0\big\},\\
\xi\!=\!(\xi_h)_{h\in I}\lra\xi|_{\Si_{b,i}}\!=\!\xi_i,
\end{gather*}
has real codimension two. Its complement corresponds to 
the infinitesimal translations in $\Bbb{C}\!\subset\!\Si_{b,i}$.
Thus, if $J$ satisfies the regularity conditions (a) and (b-i) of Definition~\ref{g1reg_dfn}
and $\hat{0}$ is the minimal element of~$I$,
for all $i\!\in\!\chi({\cal T})$ the~map
$$\D_{J,u;i}\!: \tilde{\Ga}_-(b)\!\equiv\!\big\{
\xi\!\in\!\tilde{T}_b\tilde{\cal U}_{\cal T}^{(0)}(X;J)\!: \xi(\hat{0},\i)\!=\!0\big\} 
\lra  T_{\ev_0(b)}X, \qquad
\D_{J,u;i}(\xi)=\na_{e_{\i}}^J\xi_i,$$
is surjective.

\begin{lmm}
\label{reg0_lmm2}
If $(X,\om,J)$, $A$, and ${\cal T}$ are as in Lemma~\ref{reg0_lmm1},
for every precompact open subset $K$ of ${\cal U}_{\cal T}(X;J)$,
there exist $\de_K,\ep_K,C_K\!\in\!\Bbb{R}^+$ and
an open neighborhood $U_K\!\subset\!U_{\cal T}$ of $K$ in ${\frak X}_{0,M}(X,A)$
with the following properties:\\
(1) the requirements (1) and (2) of Corollary~\ref{reg0_crl1} are satisfied;\\
(2) if $\tilde{J}$ is an almost complex structure on $X$ 
s.t.~$\|\tilde{J}\!-\!J\|_{C^1}\!\le\!\de_K$,\\
${}\quad$ (2a) for all 
$\ups\!=\!(b,v)\!\in\!\wt{\cal FT}_{\de_K}^{\eset}\big|_{\tilde{K}^{(0)}}$, 
the equation
$$\bar{\partial}_{\tilde{J}}\exp_{u_{\ups}}\!\ze=0,\qquad
\ze\!\in\!\Ga_+(\ups),~\|\ze\|_{\ups,p,1}\!<\!\ep_K,$$
${}\qquad~~~$ has a unique solution $\ze_{\tilde{J},\ups}$;\\
${}\quad$ (2b) the map 
$\tilde{\phi}_{\tilde{J}}\!:\wt{\cal FT}_{\de_K}^{\eset}|_{\tilde{K}^{(0)}}
\!\lra\!\M_{0,\{0\}\sqcup M}^0(X,A;\tilde{J})\cap U_K$,
$\ups\!\lra\![\exp_{b(\ups)}\ze_{\tilde{J},\ups}]$, is smooth;\\
${}\quad$ (2c) for all $\ups\!=\!(b,v)\!\in\!\wt{\cal FT}_{\de_K}^{\eset}|_{\tilde{K}^{(0)}}$,
$\ev_0(\ti\phi_{\ti{J}}(\ups))\!=\!\ev_0(b)$;\\
${}\quad$ (2d) for all 
$\ups\!=\!(b,v)\!\in\!\wt{\cal FT}_{\de_K}^{\eset}\big|_{\tilde{K}^{(0)}}$,
\begin{equation}\label{reg0_lmm2e}
\big\|\ze_{\tilde{J},\ups}\big\|_{\ups,p,1},
\big\|\na^{\cal T}\ze_{\tilde{J},\ups}\big\|_{\ups,p,1}
\le C_K\big(\|\tilde{J}\!-\!J\|_{C^1}\!+\!|\ups|^{1/p}\big),
\end{equation}
where $\na^{\cal T}\ze_{\tilde{J},\ups}$ denotes the differential of the bundle map 
$\ups\!\lra\!\ze_{\tilde{J},\ups}$ along $\tilde{\Ga}_-(b)$ with respect to a connection 
in the bundle $\Ga((\cdot,v))$ over~$\tilde{\cal U}_{\cal T}^{(0)}(X;J)$.
\end{lmm}

\noindent
{\it Remark:}
Let $\vph\!:\tilde{\Ga}_-(\cdot)\!\lra\!\tilde{\cal U}_{\cal T}^{(0)}(X;J)$ 
be a smooth map such that
$$d\vph|_{(b,0)}\!:\tilde{\Ga}_-(b)
\lra T_b\tilde{\cal U}_{\cal T}^{(0)}(X;J)$$
is the inclusion map for all $b\!\in\!\tilde{\cal U}_{\cal T}^{(0)}(X;J)$ and 
$$\Si_{\vph(\vsi)}=\Si_b   \quad\hbox{and}\quad
\ev_0(\vph(\vsi))=\ev_0(b)
\qquad\forall~b\!\in\!\ti{\cal U}_{\cal T}^{(0)}(X;J),~
\vsi\!\in\!\tilde{\Ga}_-(b).$$
Let
$$\ti\vph\!:
\pi_{\tilde{\Ga}_-(\cdot)}^*\Ga((\cdot,v))\!=\!
\pi_{\tilde{\Ga}_-(\cdot)}^*\Ga(\Si_{(\cdot,v)};u_{(\cdot,v)}^*TX)
\lra\Ga((\cdot,v))$$
be a lift of $\vph$ to a vector-bundle homomorphism that restricts to 
the identity over~$\tilde{\cal U}_{\cal T}^{(0)}(X;J)$.
For example, we can take $\tilde{\vph}$ to be given by
\begin{gather*}
\big\{\tilde{\vph}(\vsi;\xi)\big\}(z)=\Pi_{\ze(q_{(b,v)}(z))}\xi(z)
\qquad\forall~z\!\in\!\Si_{(b,v)}\!=\!\Si_{(\vph(\vsi),v)}\\
\hbox{if}\qquad \vsi\!\in\!\tilde{\Ga}_-(b)
\quad\hbox{and}\quad
\vph(\vsi)=(\Si_b,\exp_{u_b}\ze(z)),
\end{gather*}
where $\Pi_{\ze(q_{(b,v)}(z))}\xi(z)$ is the parallel transport of $\xi(z)$
along the geodesic 
$$\ga_{\ze(z)}\!:[0,1]\lra X,\qquad
\tau\lra\exp_{u_{(b,v)}(z)}\tau\ze(q_{(b,v)}(z)).$$
We can then define 
$\na^{\cal T}\ze_{\tilde{J},(\cdot,v)}\!:\tilde{\Ga}_-(\cdot)\!\lra\!\Ga((\cdot,v))$ by
$$\big\{\na^{\cal T}_{\vsi}\ze_{\tilde{J},\cdot}\big\}\big|_{(b,v)}
=\lim_{t\lra0}
  \frac{\ze_{\tilde{J},(\vph(t\vsi),v)}-
   \tilde{\vph}(\ze_{\tilde{J},(b,v)})}{t}
 \in \Ga((b,v))$$ 
for $X\!\in\!\tilde{\Ga}_-(b)$.
Finally, a choice of metric on $\tilde{\cal U}_{\cal T}^{(0)}(X;J)$ 
determines $\|\na^{\cal T}\ze_{\tilde{J},\ups}\|_{\ups,p,1}$.\\

\noindent
Claim (2a) of Lemma~\ref{reg0_lmm2} and the first bound in~(2d)
follow immediately from~(1) of Lemma~\ref{reg0_lmm1} and~\e_ref{reg0_lmm1e1}
via a quadratic expansion of the $\bar{\partial}_{\tilde{J}}$-operator at~$u_{\ups}$
and the Contraction Principle; see Subsection~3.6 in~\cite{Z4}.
Claim~(2c) is a consequence of~\e_ref{gaplusdfn_e}.
The smoothness of the map $\tilde{\phi}_{\tilde{J}}$ follows from the smooth dependence of
solutions of the equation in~(2a) on the parameters.
The second bound in~(2d) is obtained from the uniform behavior of these parameters;
see Subsection~3.4 in~\cite{Z4}.\\

\noindent
If $b\!\in\!\tilde{\cal U}_{\cal T}^{(0)}(X;J)$ is as before, 
the domain~$\Si_b$ of~$b$ has the~form
\begin{equation}\label{singpts_e}
\Si_b=\Big(\bigsqcup_{i\in I}\{i\}\!\times\!S^2\Big)\Big/\sim,
\qquad\hbox{where}\qquad
(h,\i)\sim(\io_h,x_h(b))~~\forall h\!\in\!\hat{I},
\end{equation}
and $x_h(b)\!\in\!S^2\!-\!\{\i\}$.
The basic gluing map $q_{\ups}\!:\Si_{\ups}\!\lra\!\Si_b$ 
used in this paper is a homeomorphism outside of at most $|\hat{I}|$ circles in~$\Si_{\ups}$
and is holomorphic outside of the annuli 
$$\ti{A}_{\ups,h}^{\pm}\equiv q_{\ups}^{-1}\big(A_{b,h}^{\pm}(|\ups_h|)\big),$$ 
with $h\!\in\!\hat{I}$,  where
\begin{equation*}\begin{split}
&A_{b,h}^-(\de)= A_h^-(\de)
=\big\{(h,z)\!\in\!\{h\}\!\times\!S^2\!:|z|\!\ge\!\de^{-1/2}/2\big\},\\
&A_{b,h}^+(\de)
=\big\{(\io_h,z)\!\in\!\{\io_h\}\!\times\!S^2\!:|z\!-\!x_h(b)|\!\le\!2\de^{1/2}\},
\end{split}\end{equation*}
for any $\de\!\in\!\Bbb{R}^+$.
For each $h\!\in\!\hat{I}$, $\tilde{A}_{\ups,h}^+\!\cup\tilde{A}_{\ups,h}^-$ is 
the thin neck of $\Si_{\ups}$ corresponding to smoothing the node of~$\Si_b$
joining the spheres $\Si_{b,\io_h}$ and~$\Si_{b,h}$.
If $\de\!\in\!\Bbb{R}^+$, let
\begin{gather*}
\partial^-A_{b,h}^-(\de)=
\big\{(h,z)\!\in\!\{h\}\!\times\!S^2\!:|z|\!=\!\de^{-1/2}/2\big\},\quad
\Si_b^0(\de)=\bigg(\Si_b-\bigcup_{i\in\chi(\T)}\bigcup_{h\ge i}\!\Si_{b,h}\bigg)
\cup\bigcup_{i\in\chi({\cal T})}\!\!\!\!A_{b,i}^-(\de);\\
\tilde{A}_{\ups,h}^{\pm}(\de)=q_{\ups}^{-1}(A_{b,h}^{\pm}(\de))\subset\Si_{\ups},
\quad
\partial^-\tilde{A}_{\ups,h}^-(\de)=q_{\ups}^{-1}\big(\partial^-A_{b,h}^-(\de)\big),
\quad  \Si_{\ups}^0(\de)=q_{\ups}^{-1}\big(\Si_b^0(\de)\big).
\end{gather*}
In the case of Figure~\ref{reg0lmm4_fig}, 
$\Si_{\ups}^0(\de)$ consists of the two non-shaded components, with the node joining
them turned into a thin neck, the three thin necks corresponding to the nodes 
attaching the bubbles $h_1$, $h_4$, and~$h_5$,
and small annuli extending from each of these three necks into
the interior of the corresponding bubble, provided $|\ups|\!<\!\de$.
If $\ups\!=\!(b,v)$, with $v\!\in\!\Bbb{C}^{*\hat{I}}$,
the complement of~$\Si_{\ups}^0(\de)$ is the union of~$|\chi({\cal T})|$
disks that support nearly all of the map~$u_{\ups}\!\equiv\!u_b\!\circ\!q_{\ups}$.

\begin{lmm}
\label{reg0_lmm3}
If $(X,\om,J)$, $A$, and ${\cal T}$ are as in Lemma~\ref{reg0_lmm1},
for every precompact open subset $K$ of ${\cal U}_{\cal T}(X;J)$,
there exist $\de_K,\ep_K,C_K\!\in\!\Bbb{R}^+$ and
an open neighborhood $U_K\!\subset\!U_{\cal T}$ of $K$ in ${\frak X}_{0,M}(X,A)$
with the following properties:\\
(1) the requirements (1) and (2) of Corollary~\ref{reg0_crl1} are satisfied;\\
(2) if $\tilde{J}$ is an almost complex structure on $X$ 
s.t.~$\|\tilde{J}\!-\!J\|_{C^1}\!\le\!\de_K$, there exist a smooth map
$$\tilde{\phi}_{\tilde{J}}\!:\wt{\cal FT}_{\de_K}^{\eset}|_{\tilde{K}^{(0)}}\lra
\M_{0,\{0\}\sqcup M}^0(X,A;\tilde{J})\cap U_K$$
such that the requirements (2a)-(2d) of Lemma~\ref{reg0_lmm2} are satisfied.
Furthermore, for every $b\!\in\!\tilde{K}^{(0)}$
and $\ups\!=\!(b,v)\!\in\!\wt{\cal FT}_{\de_K}^{\eset}|_{\tilde{K}^{(0)}}$,
there exist
\begin{alignat*}{2}
&\Phi_b\in L^p_1\big(\Si_b^0(\de_K);\End(\ev_0^*TX)\big), &\qquad
&\vt_b\in\Hol_J\big(\Si_b^0(\de_K);\ev_0^*TX\big),\\
&\Phi_{\ti{J},\ups}\in L^p_1\big(\Si_{\ups}^0(\de_K);\End(\ev_0^*TX)\big), &\qquad
&\vt_{\ti{J},\ups}\in\Hol_J\big(\Si_{\ups}^0(\de_K);\ev_0^*TX\big),
\end{alignat*}
such that\\
${}\quad$ (2a) the maps $b\!\lra\!(\Phi_b,\vt_b)$ and 
$\ups\!\lra\!(\Phi_{\tilde{J},\ups},\vt_{\tilde{J},\ups})$ are smooth;\\
${}\quad$ (2b) for all $b\!\in\!\tilde{K}^{(0)}$,
\begin{gather*}
\exp_{\ev_0(b)}\big(\Phi_b(z)\vt_b(z)\big)=u_b(z) \qquad\forall z\in\Si_b^0(\de_K);\\
\Phi_b|_{\Si_b^0(0)}=Id, \qquad\hbox{and}\qquad
\big\|\Phi_b\!-\!Id\big\|_{b,p,1}, \,
\big\|\na^{\cal T}(\Phi_b\!-\!Id)\big\|_{b,p,1}\le\frac{1}{2};
\end{gather*}
${}\quad$ (2c)  for all 
$\ups\!=\!(b,v)\!\in\!\wt{\cal FT}_{\de_K}^{\eset}\big|_{\tilde{K}^{(0)}}$,
\begin{gather*}
\exp_{\ev_0(b)}\big(\Phi_{\tilde{J},\ups}(z)\vt_{\tilde{J},\ups}(z)\big)
=\tilde{u}_{\ups}(z)\qquad\forall z\in\!\Si_{\ups}^0(\de_K)\\
\hbox{and}\qquad
\big\|\Phi_{\tilde{J},\ups}\!-\!\Phi_b\!\circ\!q_{\ups}\big\|_{\ups,p,1},
\big\|\na^{\cal T}(\Phi_{\tilde{J},\ups}\!-\!\Phi_b\!\circ\!q_{\ups})\big\|_{\ups,p,1}
\le  C_K\big(\|\tilde{J}\!-\!J\|_{C^1}\!+\!|\ups|^{1/p}\big),
\end{gather*}
if $\tilde{u}_{\ups}\!=\!\exp_{u_{\ups}}\!\ze_{\tilde{J},\ups}$.
\end{lmm}

\noindent
In the statement of this lemma, $\Hol_J(\Si_b^0(\de_K);\ev_0^*TX)$
and $\Hol_J(\Si_{\ups}^0(\de_K);\ev_0^*TX)$
denote the spaces of holomorphic maps from $\Si_b^0(\de_K)$
and  $\Si_{\ups}^0(\de_K)$ into the complex vector space $(T_{\ev_0(b)}X,J)$.
In brief, the substance of Lemma~\ref{reg0_lmm3} is that 
a $J$-holomorphic map can be well approximated by
a holomorphic map on a neighborhood of the primary marked point,
or any other point, of the domain.
Due to Lemma~\ref{reg0_lmm2}, Lemma~\ref{reg0_lmm3} is essentially 
a parametrized version of Theorem~2.2 in~\cite{FlHSa}, 
and only a couple of additional ingredients are needed.
The crucial fact used in~\cite{FlHSa} is that the operator
$$\Xi\!:L^p_1\big(S^2;\End_{\Bbb{C}}(\Bbb{C}^n)\big)\lra
L^p\big(S^2;\La^{0,1}T^*S^2\!\otimes\!\End_{\Bbb{C}}(\Bbb{C}^n)\big)
\oplus\End_{\Bbb{C}}(\Bbb{C}^n),         \quad
\Xi(\Th)\!=\!(\bar{\partial}\Th,\Th(0)),$$
is an isomorphism.
The map~$\Xi$ is still an isomorphism if $S^2$ is replaced by
a tree of spheres~$\Si$ and $0$ by any point on~$\Si$.
Furthermore, if $y$ is a smooth point of~$\Si$
for all sufficiently small smoothings~$\ups$ of the nodes,
the operators
$$\Xi_{\ups}\!:L^p_1\big(\Si_{\ups};\End_{\Bbb{C}}(\Bbb{C}^n)\big)\lra
L^p\big(\Si_{\ups};\La^{0,1}T^*\Si_{\ups}\!\otimes\!\End_{\Bbb{C}}(\Bbb{C}^n)\big)
\oplus\End_{\Bbb{C}}(\Bbb{C}^n),            \quad
\Xi_{\ups}(\Th)\!=\!(\bar{\partial}\Th,\Th(y)),$$
are also isomorphisms.
In addition, for some $C\!\in\!\Bbb{R}^+$ and for all sufficiently small smoothings~$\ups$,
\begin{equation}\label{reg0_lmm3e1}
C^{-1}\|\Xi_{\ups}\Th\|_{\ups,p}\le 
\|\Th\|_{\ups,p,1}\le C\|\Xi_{\ups}\Th\|_{\ups,p}
\qquad\forall~\Th\!\in\!L^p_1\big(\Si_{\ups};\End_{\Bbb{C}}(\Bbb{C}^n)\big).
\end{equation}
If all components of $\ups$ are nonzero, 
$\Si_{\ups}$ is topologically a sphere, but should instead be viewed
as a tree of spheres joined by thin necks.
As before, we denote by~$\|\cdot\|_{\ups,p,1}$ and~$\|\cdot\|_{\ups,p}$
the norms induced by the pregluing construction above.
In particular, \e_ref{reg0_lmm3e1} can be viewed as a special case 
of~(1) of Lemma~\ref{reg0_lmm1}.
We need to use the norms $\|\cdot\|_{\ups,p,1}$ and~$\|\cdot\|_{\ups,p}$,
since these are the norms used in Lemmas~\ref{reg0_lmm1} and~\ref{reg0_lmm2}.
Keeping track of all norms in the proof of Theorem~2.2 in~\cite{FlHSa},
we see that the maps $\Phi_b$, $\vt_b$, $\Phi_{\tilde{J},\ups}$,
and $\vt_{\tilde{J},\ups}$ satisfying~(2b) and the first condition in~(2c)
exist, provided that $\de_K$ is sufficiently small.
The last two estimates in~(2c) are obtained by an argument similar
to Subsection~4.1 in~\cite{Z4}.\\

\begin{figure}
\begin{pspicture}(-1.1,-1.8)(10,1.25)
\psset{unit=.4cm}
\pscircle(5,-1.5){1.5}\pscircle*(5,-3){.24}\rput(5,-3.8){$(\hat{0},\i)$}
\pscircle[fillstyle=solid,fillcolor=gray](3.23,.27){1}\pscircle*(3.94,-.44){.2}
\pscircle[fillstyle=solid,fillcolor=gray](1.23,.27){1}\pscircle*(2.23,.27){.2}
\pscircle(6.77,.27){1}\pscircle*(6.06,-.44){.2}
\pscircle[fillstyle=solid,fillcolor=gray](6.77,2.27){1}\pscircle*(6.77,1.27){.2}
\pscircle[fillstyle=solid,fillcolor=gray](8.77,.27){1}\pscircle*(7.77,.27){.2}
\rput(3.3,1.8){$h_1$}\rput(1.3,1.8){$h_2$}\rput(5.3,1){$h_3$}
\rput(8.4,2.4){$h_4$}\rput(10.4,.4){$h_5$}
\pnode(2,-3){A1}\rput(.9,-3){$x_{h_1}(b)$}\pnode(4,-.46){B1}
\ncarc[nodesep=.35,arcangleA=-30,arcangleB=-50,ncurv=.8]{->}{A1}{B1}
\pnode(7.3,-3){A2}\rput(8.4,-2.9){$x_{h_3}(b)$}\pnode(6,-.45){B2}
\ncarc[nodesep=.35,arcangleA=30,arcangleB=40,ncurv=.8]{->}{A2}{B2}
\pnode(8.8,-1.55){A3}\rput(9.9,-1.4){$x_{h_5}(b)$}\pnode(7.85,.25){B3}
\ncarc[nodesep=.35,arcangleA=45,arcangleB=110,ncurv=1]{->}{A3}{B3}
\rput(24.7,0){\begin{tabular}{l}
$\chi({\cal T})\!=\!\{h_1,h_4,h_5\}$\\
$\rho(\ups)=(\ups_{h_1},\ups_{h_3}\ups_{h_4},\ups_{h_3}\ups_{h_5})$\\
$x_{h_5}(\ups)=x_{h_3}(b)+v_{h_3}x_{h_5}(b)$\\
\\
${\cal D}_{\hat{0}}^{(1)}\tilde{b}_{\tilde{J}}(\ups)
\cong {\cal D}_{J,h_1}^{(1)}\ups_{h_1}\!+\!
{\cal D}_{J,h_4}^{(1)}\ups_{h_3}\ups_{h_4}\!+\!{\cal D}_{J,h_5}^{(1)}\ups_{h_3}\ups_{h_5}$
\end{tabular}}
\end{pspicture}
\caption{An Example of the Estimates of Lemma~\ref{reg0_lmm4}}
\label{reg0lmm4_fig}
\end{figure}

\noindent
Before we can state Lemma~\ref{reg0_lmm4}, we need to introduce additional notation.
For each $\ups\!=\!(b,v)$, where $b\!\in\!\ti\U_{\T}^{(0)}(X;J)$
and $v\!\in\!\Bbb{C}^{\hat{I}}$, and $i\!\in\!\chi({\cal T})$, let
$$\rho_i(v)=\!\prod_{\hat{0}<h\le i}\!\!\!v_h\in\Bbb{C},
\quad \rho_i(\ups)=(b,\rho_i(v));\qquad
x_i(\ups)=\!\!\sum_{\hat{0}<i'\le i}\!\Big(x_{i'}(b)
\!\!\!\prod_{\hat{0}<h<i'}\!\!\!\!v_h\Big)~\in~\Bbb{C},$$
where $x_i(b)$ is as in~\e_ref{singpts_e};
see Figure~\ref{reg0lmm4_fig}.
If $K$ and $\vt_b$ are as in Lemma~\ref{reg0_lmm3},
$b\!=\!(\Si_b,u_b)\!\in\!\tilde{K}^{(0)}$,
$v\!\in\!\Bbb{C}$, $i\!\in\!\chi({\cal T})$, and $r\!\in\!\Z^+$, we~put
$${\cal D}_i^{(r)}b=\frac{1}{r!}
\frac{d^r}{dw_i^r}\vt_{b,i}(w_i)\Big|_{w_i=0}\in T_{\ev_0(b)}X
\quad\hbox{and}\quad  {\cal D}_{J,i}^{(r)}(b,v)=v\cdot_J{\cal D}_i^{(r)}b,$$
where $\vt_{b,i}\!=\!\vt_b|_{\Si_{b,i}}$,
$\cdot_J$ is the scalar multiplication in~$(TX,J)$,
and $w_i$ is the standard holomorphic coordinate centered 
at the point $\i$ in $\Si_{b,i}\!=\!S^2$.
If $\ups\!\in\!\wt{\cal FT}_{\de_K}^{\eset}|_{\tilde{K}^{(0)}}$, 
we similarly set
$$\ti{b}_{\ti{J}}(\ups)=(\Si_{\ups},\ti{u}_{\ups}),
\qquad
{\cal D}_{\hat{0}}^{(r)}\tilde{b}_{\tilde{J}}(\ups)
=\frac{1}{r!}\frac{d^r}{dw^r}\vt_{\tilde{J},\ups}(w)\Big|_{w=0}\in T_{\ev_0(b)}X,$$
where $w$ is the standard holomorphic coordinate centered 
at the point $\i$ in $\Si_{\ups}\!\approx\!S^2$.
The value of ${\cal D}_i^{(r)}b$ depends on the choice of $\vt_b$
in Lemma~\ref{reg0_lmm3}, 
which can be uniquely prescribed by the choice of $\de_K\!\in\!\Bbb{R}^+$.
Alternatively, one can replace small positive numbers $\de_K$ dependent
on compact subsets $K$ of ${\cal U}_{\cal T}(X;J)$ by
a small function $\de\!:{\cal U}_{\cal T}(X;J)\!\lra\!\Bbb{R}^+$,
which can be used to choose a holomorphic~map 
$$\vt_b\!:\Si_b^0(\de(b))\lra(T_{\ev_0(b)}X,J)$$
for each $J$-holomorphic stable map $b\in\!\ti{\cal U}_{\T}^{(0)}(X;J)$.
Of course, the definition of $\cD_{\hat{0}}^{(r)}\ti{b}_{\ti{J}}(\ups)$ depends on
even more choices, including those involved in the gluing construction
of Lemma~\ref{reg0_lmm2}.

\begin{lmm}
\label{reg0_lmm4}
If $(X,\om,J)$, $A$, and ${\cal T}$ are as in Lemma~\ref{reg0_lmm1},
for every precompact open subset $K$ of ${\cal U}_{\cal T}(X;J)$,
there exist $\de_K,\ep_K,C_K\!\in\!\Bbb{R}^+$ and
an open neighborhood $U_K\!\subset\!U_{\cal T}$ of $K$ in ${\frak X}_{0,M}(X,A)$
with the following properties:\\
(1) the requirements (1) and (2) of Corollary~\ref{reg0_crl1} are satisfied;\\
(2) if $\tilde{J}$ is an almost complex structure on $X$ 
s.t.~$\|\tilde{J}\!-\!J\|_{C^1}\!\le\!\de_K$, there exist $\tilde{\phi}_{\tilde{J}}$,
$\Phi_b$ and $\vt_b$ for each $b\!\in\!\tilde{K}^{(0)}$,
and $\Phi_{\tilde{J},\ups}$ and $\vt_{\tilde{J},\ups}$ for each 
$\ups\!\in\!\wt{\cal FT}_{\de_K}^{\eset}|_{\tilde{K}^{(0)}}$,
such that the requirements (2a)-(2d) of Lemma~\ref{reg0_lmm2} and
(2a)-(2c) of Lemma~\ref{reg0_lmm3} are satisfied.
Furthermore,  for each $k\!\in\!\Z^+$ and $i\!\in\!\chi({\cal T})$
there exists a smooth~map
$$\ve_{\ti{J},i}^{(k)}\!:\wt{\cal FT}_{\de_K}^{\eset}\big|_{\ti{K}^{(0)}}\lra\ev_0^*TX,$$
such that:\\
${}\quad$ (2a) for all $r\!\in\!\Z^+$ and 
$\ups\!=\!(b,v)\!\in\!\wt{\cal FT}_{\de_K}^{\eset}\big|_{\tilde{K}^{(0)}}$,
$$\cD_{\hat{0}}^{(r)}\ti{b}_J(\ups)=
\sum_{k=1}^{k=r}\binom{r\!-\!1}{k\!-\!1}\sum_{i\in\chi({\cal T})}\!\!\!
x_i^{r-k}(\ups)\big\{\cD_{J,i}^{(k)}\!+\!\ve_{\ti{J},i}^{(k)}(\ups)\big\}\rho_i^k(\ups)
\in T_{\ev_0(b)}X;$$
${}\quad$ (2b) for all $k\!\in\!\Z^+$, $i\!\in\!\chi({\cal T})$, and
$\ups\!=\!(b,v)\!\in\!\wt{\cal FT}_{\de_K}^{\eset}\big|_{\tilde{K}^{(0)}}$,
$$|\ve_{\ti{J},i}^{(k)}|_{\ups},|\na^{\T}\ve_{\ti{J},i}^{(k)}|_{\ups}
\le C_K\de_K^{-k/2} \big(\|\ti{J}\!-\!J\|_{C^1}\!+\!|\ups|^{1/p}\big).$$
\end{lmm}

\noindent
{\it Proof:} (1) We apply, with some modifications, 
the argument for the analytic estimates of Theorem~2.8 in~\cite{Z3} 
to holomorphic functions $\vt_b$ and $\vt_{\ti{J},\ups}$,
instead of the functions $u_b$ and  $\tilde{u}_{\ups}$
which are $J$-holomorphic and $\tilde{J}$-holomorphic in this case.
We will use coordinates $z$ on $S^2\!-\!\{\i\}\!\approx\!\Bbb{C}$ and $w\!=\!z^{-1}$
on $S^2\!-\!\{0\}$.
Since $\vt_{\tilde{J},\ups}$ is holomorphic on $\Si_{\ups}^0(\de_K)$,
\begin{equation}\label{reg0_lmm4e1}\begin{split}
{\cal D}_{\hat{0}}^{(r)}\tilde{b}_J(\ups)
&=\frac{1}{r!}\frac{\partial^r}{\partial w^r}\vt_{\tilde{J},\ups}(w)\big|_{w=0}
=\frac{1}{2\pi\I}\oint_{\partial\Si_{\ups}^0(\de_K)}
\vt_{\tilde{J},\ups}(w)\frac{dw}{w^{r+1}}\\
&=\frac{1}{2\pi\I}\sum_{i\in\chi({\cal T})}\oint_{\partial^-\tilde{A}_{\ups,i}^-(\de_K)}
\!\!\!\vt_{\tilde{J},\ups}(w)\frac{dw}{w^{r+1}}\\
&=-\frac{1}{2\pi\I}\sum_{i\in\chi({\cal T})}\oint_{\partial^-\tilde{A}_{\ups,i}^-(\de_K)}
\!\!\!\Phi_{\tilde{J},\ups}^{-1}(z)
\big(\Phi_b\vt_b\circ q_{\ups}(z)+\tilde{\ze}_{\tilde{J},\ups}(z)\big)z^{r-1}dz,
\end{split}\end{equation}
where $\tilde{\ze}_{\tilde{J},\ups}\!\in\!C^{\i}(\Si_{\ups}^0(\de_K);T_{\ev_0(b)}X)$
is defined by
$$\exp_{\ev_0(b)}\big(\exp_{\ev_0(b)}^{-1}\!u_{\ups}(z)+\ti\ze_{\ti{J},\ups}\big)
=\ti{u}_{\ups}(z)\equiv\exp_{u_{\ups}(z)}\ze_{\ti{J},\ups}(z).$$
(2) In order to estimate each integral on the last line of~\e_ref{reg0_lmm4e1},
we expand $z^{r-1}$ around $x_i(\ups)$, the center of the circle
$\partial^-\tilde{A}_{\ups,i}^-(\de_K)$, as a polynomial in $\tilde{z}_i\!=\!z\!-\!x_i(\ups)$:
\begin{equation}\label{reg0_lmm4e3}\begin{split}
&\oint_{\partial^-\tilde{A}_{\ups,i}^-(\de_K)}\!\!\!
\Phi_{\tilde{J},\ups}^{-1}(z)
\big(\Phi_b\vt_b\circ q_{\ups}(z)+\tilde{\ze}_{\tilde{J},\ups}(z)\big)z^{r-1}dz\\
&\qquad=\sum_{k=1}^{k=r}\binom{r\!-\!1}{k\!-\!1}x_i^{r-k}(\ups)
\oint_{\partial^-\tilde{A}_{\ups,i}^-(\de_K)}\!\!\!
\Phi_{\tilde{J},\ups}^{-1}(\tilde{z}_i)
\big(\Phi_b\vt_b\circ q_{\ups}(\tilde{z}_i)+\tilde{\ze}_{\tilde{J},\ups}(\tilde{z}_i)\big)
\tilde{z}_i^{k-1}d\tilde{z}_i.
\end{split}\end{equation}
By construction, $z_i\!=\!q_{\ups}(\ti{z}_i)\!=\!\rho_i^{-1}(v)\ti{z}_i$ near
$\partial^-\ti{A}_{\ups,i}^-(\de_K)$, if $z_i$ is the standard holomorphic coordinate on 
$\Si_{\ups,i}\!-\!\{\i\}\!=\!\{i\}\!\times\!(S^2\!-\!\{\i\})$ and $|\ups|\!<\!\de_K$.
Since $\vt_{b,i}\!\equiv\!\vt_b|_{\Si_{b,i}}$ is $J_{\ev_0(b)}$-holomorphic,
\begin{equation}\label{reg0_lmm4e5}\begin{split}
\oint_{\partial^-\tilde{A}_{\ups,i}^-(\de_K)}\!\!\!
\vt_b(q_{\ups}(\tilde{z}_i))\tilde{z}_i^{k-1}d\tilde{z}_i
&=\rho_i^k(v)
\oint_{\partial^-A_i^-(\de_K)}\!\!\vt_{b,i}(z_i)\, z_i^{k-1}dz_i\\
&=-\rho_i^k(v)\oint_{\partial^-A_i^-(\de_K)}\!\!\vt_{b,i}(w_i)\, \frac{dw_i}{w_i^{k+1}}\\
&=-\rho_i^k(v)\, \frac{2\pi\I}{k!}\,
\frac{\partial^k}{\partial w_i^k}\vt_{b,i}(w_i)\big|_{w_i=0}
=-2\pi\I\,{\cal D}_{J,i}^{(k)}\rho_i^k(\ups),
\end{split}\end{equation}
where $w_i\!=\!z_i^{-1}$ is the standard holomorphic coordinate 
on $\{i\}\!\times\!(S^2\!-\!\{0\})$.
Similarly,
\begin{equation}\label{reg0_lmm4e7a}\begin{split}
\oint_{\partial^-\tilde{A}_{\ups,i}^-(\de_K)}\!\!\!\big(
\big\{\Phi_{\tilde{J},\ups}^{-1}(\tilde{z}_i)\Phi_b(q_{\ups}(\tilde{z}_i))\!-\!Id\big\}
\vt_b(q_{\ups}(\tilde{z}_i))
+\Phi_{\tilde{J},\ups}^{-1}(\tilde{z}_i)\tilde{\ze}_{\tilde{J},\ups}(\tilde{z}_i)\big)
\ti{z}_i^{k-1}d\tilde{z}_i\qquad&\\
=-2\pi\I\,\ve_{\ti{J},i}^{(k)}(\ups)\rho_i^k(\ups)&,
\end{split}\end{equation}
where
\begin{equation}\label{reg0_lmm4e7b}
\ve_{\ti{J},i}^{(k)}(\ups) =\frac{1}{2\pi{\frak i}}
\oint_{\partial^-A_i^-(\de_K)}\!\!\big(
\big\{\Phi_{\tilde{J},\ups}^{-1}(w_i)\Phi_b(q_{\ups}(w_i))\!-\!Id\big\}
\vt_b(q_{\ups}(w_i))
+\Phi_{\tilde{J},\ups}^{-1}(w_i)\tilde{\ze}_{\tilde{J},\ups}(w_i)\big)
\frac{dw_i}{w_i^{k+1}}.
\end{equation}
The expansion in (2a) of the lemma follows immediately from 
\e_ref{reg0_lmm4e1}-\e_ref{reg0_lmm4e7a}.
By definition, $\partial^-A_i^-(\de_K)$ is a circle of
radius $2\de_K^{1/2}$ around the south pole in the sphere~$\Si_{b,i}$.
Thus, part~(2b) of the lemma follows from~\e_ref{reg0_lmm4e7b},
along with~\e_ref{reg0_lmm2e} and the bounds in~(2b) and~(2c) of Lemma~\ref{reg0_lmm3}.

\subsection{The Genus-One Regularity Properties (b-i) and (b-ii)}
\label{reg0_subs2}

\noindent
In this subsection, we show that if $J$ is an almost complex structure on $X$
that satisfies the regularity conditions~(a) and (b-i) of Definition~\ref{g1reg_dfn}, 
then so does every almost complex structure $\tilde{J}$ on $X$ which is
sufficiently close to~$J$.
This claim follows from Corollary~\ref{reg0_crl2} 
and the compactness of the moduli space $\ov\M_{0,1}(X,A;J)$.
We also show that  if $J$ is an almost complex structure on $X$
that satisfies the regularity conditions~(a), (b-i), and~(b-ii) of Definition~\ref{g1reg_dfn}, 
then so does every nearby almost complex structure $\tilde{J}$ on~$X$.
This conclusion is immediate from Corollary~\ref{reg0_crl3} 
and the compactness of the moduli space $\ov\M_{0,2}(X,A;J)$.\\

\noindent
If $\T$ is a bubble type as in Corollary~\ref{reg0_crl1} with $A_{\hat{0}}\!=\!0$,
where $\hat{0}$ is the minimal element of~$I$, 
the analogue of~\e_ref{reg0_crl1e5} does not hold for 
the map $\D_{J,b}^{\i,v}$ for any fixed nonzero vector~$v$ tangent to~$\Bbb{P}^1$ at~$\i$.
The reason is that $\D_{J,b}^{\i,v}$ is the zero homomorphism on $\ker\D_{J,b}^{\i}$,
since the map~$u_b$ is constant on the component~$\Si_{b,\hat{0}}$
of the domain~$\Si_b$ of~$b$ which contains the marked point~$\i$.
In particular, $\D_{\tilde{J},\tilde{u}}^{\i,v}$ need not be surjective 
for a smooth map $\tilde{u}\!:\Bbb{P}^1\!\lra\!X$
arbitrary close to the moduli space~$\ov\M_{0,1}(X,A;J)$.
Thus, a different approach is required to understand the behavior
of the regularity condition~(b-i) of Definition~\ref{g1reg_dfn}
near~${\cal U}_{\cal T}(X;J)$.\\

\noindent
Claim~(c) of Theorem~\ref{reg_thm} can alternatively be viewed as 
a statement concerning the behavior of the first derivatives $du|_{\i}$
of $J$-holomorphic maps.
Lemma~\ref{reg0_lmm4} describes the behavior of modified 
first and higher-order derivatives of $\tilde{J}$-holomorphic maps near 
a stratum ${\cal U}_{\T}(X;J)$ with ${\cal T}$ as in the previous paragraph.
We use the estimate for the higher-order derivatives to describe the behavior
of the regularity condition~(b-ii) of Definition~\ref{g1reg_dfn}
near~${\cal U}_{\T}(X;J)$.

\begin{crl}
\label{reg0_crl2}
Suppose $(X,\om,J)$, $A\!\neq\!0$, and ${\cal T}$ are as in Lemma~\ref{reg0_lmm1} 
and $M\!=\!\eset$.
If the almost complex structure~$J$ satisfies the regularity conditions~(a)
and~(b-i) of Definition~\ref{g1reg_dfn},
for every precompact open subset $K$ of ${\cal U}_{\cal T}(X;J)$,
there exist $\de_K,C_K\!\in\!\Bbb{R}^+$
and an open neighborhood $U_K\!\subset\!U_{\cal T}$ of $K$ in $\X_{0,\eset}(X,A)$
with the following properties:\\
(1) requirements (1) and (2) of Lemma~\ref{reg0_lmm1} are satisfied;\\
(2) if $\tilde{J}$ is an almost complex structure on $X$ 
s.t.~$\|\tilde{J}\!-\!J\|_{C^1}\!<\!\de_K$ and 
$[\Bbb{P}^1,\tilde{u}]\!\in\!\M_{0,\{0\}}^0(X,A;\tilde{J})\cap U_K$,
the operators $D_{\tilde{J},\tilde{u}}$, $\D_{\tilde{J},\tilde{u}}^{\i}$,
and $\D_{\tilde{J},\tilde{u}}^{\i,e_{\i}}$ are surjective.
\end{crl}

\noindent
{\it Proof:} (1) By Corollary~\ref{reg0_crl1},  it remains to show that the operator 
$\D_{\tilde{J},\tilde{u}}^{\i,e_{\i}}$ is surjective.
If \hbox{${\cal T}\!=\!(\eset,I;,\under{A})$} with $A_{\hat{0}}\!\neq\!0$,
the argument used twice in the proof of Corollary~\ref{reg0_crl1} 
can be repeated once more to show that 
the operator $\D_{\tilde{J},\tilde{u}}^{\i,e_{\i}}$ is also surjective
for any smooth map~$\tilde{u}$ sufficiently close to~$K$.
Thus, we will assume that $A_{\hat{0}}\!=\!0$.\\
(2) Let $q$ be any point in $X$. 
By (a) of Definition~\ref{g1reg_dfn},
$$\tilde{\cal U}_{\cal T}^{(0)}(q;J)\equiv
\{b\!\in\!\tilde{\cal U}_{\cal T}^{(0)}(X;J)\!:\ev_0(b)\!=\!q\}$$
is a smooth submanifold of $\tilde{\cal U}_{\cal T}^{(0)}(X;J)$.
By Corollary~\ref{reg0_crl1}, $\M_{0,\{0\}}^0(X,A;\tilde{J})$ is a smooth orbifold,
while
$$\M_{0,\{0\}}^0(q;\tilde{J})\equiv
\{b\!\in\!\M_{0,\{0\}}^0(X,A;\tilde{J})\!:\ev_0(b)\!=\!q\}$$
is a smooth suborbifold of $\M_{0,\{0\}}^0(X,A;\tilde{J})$.
By Lemma~\ref{reg0_lmm2}, every element of $\M_{0,\{0\}}^0(q;\tilde{J})\cap U_K$
has a representative of the form $(\Bbb{P}^1,\tilde{u})$,
where
$$\tilde{u}\!=\!\exp_{u_{\ups}}\!\ze_{\tilde{J},\ups},
\qquad\ups\!=\!(b,v)\in\wt{\cal FT}_{\de_K}^{\eset}\big|_{\tilde{K}^{(0)}},$$
and $\ze_{\tilde{J},\ups}$ is as in~(2a) of Lemma~\ref{reg0_lmm2}.
By~(2c) of Lemma~\ref{reg0_lmm2},
\begin{equation}\label{reg0_crl2e1}\begin{split}
&\frac{d}{dt}\exp_{u_{(\vph(t\xi),v)}}\!\ze_{\tilde{J},(\vph(t\xi),v)}\big|_{t=0}
\in\ker \D_{\tilde{J},\tilde{u}}^{\i}\!\subset\!\ker D_{\tilde{J},\tilde{u}}\\
&\quad\forall~ 
\xi\!\in\!\tilde{T}_b\tilde{\cal U}_{\cal T}^{(0)}(X;J)\cap
T_b\tilde{\cal U}_{\cal T}^{(0)}(q;J)=\tilde{\Ga}_-(b),
\end{split}\end{equation}
where $\vph$ is the map defined in the remark following Lemma~\ref{reg0_lmm2}.
We will show that the map 
\begin{equation}\label{reg0_crl2e2}
\tilde{\Ga}_-(b)\lra T_qX, \qquad
\xi\lra \frac{d}{dt}\big\{
d\{\exp_{u_{(\vph(t\xi),v)}}\!\ze_{\ti{J},(\vph(t\xi),v)}\}|_{\i}e_{\i}\big\}\big|_{t=0},
\end{equation}
is surjective. Along with \e_ref{reg0_crl2e1}, this claim implies 
Corollary~\ref{reg0_crl2}.\\
(3) Let $\Phi_b$, $\vt_b$, $\Phi_{\tilde{J},\ups}$, and $\vt_{\tilde{J},\ups}$ 
be as in Lemma~\ref{reg0_lmm3}. 
Since for all 
$\ups'\!=\!(b',v)\!\in\!\wt{\cal FT}_{\de_K}^{\eset}\big|_{\tilde{K}^{(0)}}$,
$\Phi_{\tilde{J},\ups'}$ is an $L^p_1$-map on $\Si_{\ups'}^0(\de_K)\!\subset\!\Bbb{P}^1$,
while the $J_{\ev_0(b')}$-holomorphic map $\vt_{\tilde{J},\ups'}$ vanishes at $\i\!\in\!\Si_{\ups'}^0(\de_K)$,
\begin{equation*}\begin{split}
d\big\{\!\exp_{u_{\ups'}}\ze_{\tilde{J},\ups'}\big\}\big|_{\i}e_{\i}
&=d\{\Phi_{\tilde{J},\ups'}\vt_{\tilde{J},\ups'}\}\big|_{\i}e_{\i}\\
&=\{\Phi_{\ti{J},\ups'}(\i)\}
\big(d\vt_{\ti{J},\ups'}|_{\i}e_{\i}\big).
\end{split}\end{equation*}
Thus, by the $r\!=\!1$ case of (2a) of Lemma~\ref{reg0_lmm4},
\begin{equation}\label{reg0_crl2e3}\begin{split}
d\big\{\!\exp_{u_{\ups'}}\!\ze_{\tilde{J},\ups'}\big\}\big|_{\i}e_{\i}
&=\{\Phi_{\ti{J},\ups'}(\i)\}\, {\cal D}_{\hat{0}}^{(1)}\tilde{b}_{\tilde{J}}(\ups')\\
&=\{\Phi_{\ti{J},\ups'}(\i)\}
\sum_{i\in\chi({\T})}\!\!\!
\big\{(\cD_i^{(1)}b')\!+\!\ve_{\ti{J},i}^{(1)}(\ups')\big\}\rho_i(v').
\end{split}\end{equation}
Replacing $b'$ with $\vph(t\xi)$ in~\e_ref{reg0_crl2e3}
and differentiating at $t\!=\!0$, we obtain
\begin{equation}\label{reg0_crl2e5}\begin{split}
&\frac{d}{dt}\big\{
d\{\exp_{u_{(\vph(t\xi),v)}}\!\ze_{\ti{J},(\vph(t\xi),v)}\}|_{\i}e_{\i}\big\}\big|_{t=0}\\
&\qquad\qquad\qquad=
\big\{\na^{\cal T}_{\xi}\Phi_{\ti{J},\ups}(\i)\big\}
\sum_{i\in\chi({\cal T})}\!\!\!
\big\{(\cD_i^{(1)}b)\!+\!\ve_{\ti{J},i}^{(1)}(\ups)\big\}\rho_i(v)\\
&\qquad\qquad\qquad\qquad+\{\Phi_{\ti{J},\ups}(\i)\}
\sum_{i\in\chi({\cal T})}\!\!\!
\big\{\na^{\cal T}_{\xi}({\cal D}_i^{(1)}b)\!+\!
\na^{\T}_{\xi}\ve_{\ti{J},i}^{(1)}(\ups)\big\}\rho_i(v).
\end{split}\end{equation}
By (2b) and (2c) of Lemma~\ref{reg0_lmm3} and (2b) of Lemma~\ref{reg0_lmm4},
\begin{equation}\label{reg0_crl2e7}\begin{split}
&\big|\Phi_{\ti{J},\ups}^{-1}(\i)\!-\!Id\big|,
\big|\ve_{\ti{J},h}^{(1)}(\ups)\big|\le
C_K\big(\|\ti{J}\!-\!J\|_{C^1}\!+\!|\ups|^{1/p}\big);\\
&\big|\na^{\T}_{\xi}\Phi_{\ti{J},\ups}(\i)\big|,
\big|\na^{\T}_{\xi}\ve_{\ti{J},h}^{(1)}(\ups)\big|\le
C_K\big(\|\ti{J}\!-\!J\|_{C^1}\!+\!|\ups|^{1/p}\big)\|\xi\|_{b,p,1}.
\end{split}\end{equation}
On the other hand, since $\vt_b(\i)\!=\!0$,
\begin{equation}\label{reg0_crl2e9}\begin{split}
\na^{\cal T}_{\xi}({\cal D}_i^{(1)}b)
&=\Phi_b(\i)\,
\frac{d}{dt}\na^{\cal T}_{\xi}\big(du_{\vph(t\xi),i}|_{\i}e_{\i}\big)\big|_{t=0}
=\frac{d}{dt}\big(du_{\vph(t\xi),i}|_{\i}e_{\i}\big)\big|_{t=0} \\
&=\na_{e_{\i}}^J(\xi|_{\Si_{b,i}})\equiv \D_{J,b;i}\xi
\qquad\forall~ \xi\!\in\!\ti{\Ga}_-(b).
\end{split}\end{equation}
By (a) and~(b-i) of Definition~\ref{g1reg_dfn}, the map
$$\D_{J,b;i}\!:\ti{\Ga}_-(b)\lra T_{\ev_0(b)}X$$
is surjective for all $i\!\in\!\chi({\cal T})$;
see the paragraph preceding Lemma~\ref{reg0_lmm2}.
Since $\rho_i(\ups)\!\in\!\Bbb{C}^*$ for all
$\ups\!=\!(b,v)\!\in\!\wt{\cal FT}_{\de_K}^{\eset}\big|_{\tilde{K}^{(0)}}$ 
and $i\!\in\!\chi({\cal T})$,
it follows from \e_ref{reg0_crl2e5}-\e_ref{reg0_crl2e9} that
the map in~\e_ref{reg0_crl2e2} is also surjective,
provided~$\de_K$ is sufficiently small.\\

\noindent
{\it Remark:} At the end of the argument above, we use the fact that
$\chi(\T)\!\neq\!\eset$.
This is the case if and only if $A\!\neq\!0$.

\begin{crl}
\label{reg0_crl3}
Suppose $(X,\om,J)$, $A\!\neq\!0$, and ${\cal T}$ are as in Lemma~\ref{reg0_lmm1} 
and $M\!=\!\{1\}$.
If the almost complex structure~$J$ satisfies the regularity conditions~(a), (b-i),
and~(b-ii) of Definition~\ref{g1reg_dfn},
for every precompact open subset $K$ of ${\cal U}_{\cal T}(X;J)$,
there exist $\de_K,C_K\!\in\!\Bbb{R}^+$
and an open neighborhood $U_K\!\subset\!U_{\cal T}$ of $K$ in $\X_{0,\{1\}}(X,A)$
with the following properties:\\
(1) requirements (1) and (2) of Lemma~\ref{reg0_lmm1} are satisfied;\\
(2) if $\tilde{J}$ is an almost complex structure on $X$ 
s.t.~$\|\tilde{J}\!-\!J\|_{C^1}\!<\!\de_K$ and 
$$[\Bbb{P}^1,y_1,\tilde{u}]\in\M_{0,\{0,1\}}^0(X,A;\ti{J})\cap U_K,$$
the operators $D_{\tilde{J},\tilde{u}}$, $\D_{\tilde{J},\tilde{u}}^{\i}$,
and $\D_{\tilde{J},\tilde{u}}^{\i,y_1}$ are surjective.
\end{crl}

\noindent
{\it Proof:} (1) By Corollary~\ref{reg0_crl1},  it remains to show that the operator 
$\D_{\tilde{J},\tilde{u}}^{\i,y_1}$ is surjective.   If 
$${\cal T}=(\{1\},I;j,\under{A})$$ 
is a bubble type such that $A_i\!\neq\!0$ for some $i\!\le\!j_1$,
the argument used in the proof of Corollary~\ref{reg0_crl1} 
can be repeated once more to show that 
the operator $\D_{\tilde{J},\tilde{u}}^{\i,y_1}$ is also surjective
for any smooth map~$(y_1,\tilde{u})$, with two marked points, sufficiently close to~$K$.
Thus, we will assume that $A_i\!=\!0$ for all $i\!\le\!j_1$.
In this case, $\ev_0(b)\!=\!\ev_1(b)$ for all $b\!\in\!{\cal U}_{\T}(X;J)$,
as there are no non-ghost components between the marked points $(\hat{0},\i)$ and~$y_1$.
In the case of Figure~\ref{reg0lmm4_fig}, this means that $j_1\!\in\!\{\hat{0},h_3\}$,
i.e.~$y_1$ lies on one of the non-shaded bubbles.\\
(2) For any point $q\!\in\!X$, let $\tilde{\cal U}_{\cal T}^{(0)}(q;J)$
and $\M_{0,\{0,1\}}^0(q;\tilde{J})$ be defined as in~(2) of 
the proof of Corollary~\ref{reg0_crl2}.
By Lemma~\ref{reg0_lmm2}, every element of $\M_{0,\{0,1\}}^0(q;\tilde{J})\cap U_K$
has a representative of the form $(\Bbb{P}^1,\tilde{y}_1,\tilde{u})$,
where
$$\tilde{u}\!=\!\exp_{u_{\ups}}\!\ze_{\tilde{J},\ups},\quad
\tilde{y}_1\!=\!y_1(\ups)\!=\!\!\sum_{\hat{0}<i<j_1}\!\!\!\!x_i(b)
\!\!\!\prod_{\hat{0}<h<i}\!\!\!\!v_h
+y_1\!\!\!\prod_{\hat{0}<h\le j_1}\!\!\!\!\!v_h\in\Bbb{C},\qquad
\ups\!=\!(b,v)\in\wt{\cal FT}_{\de_K}^{\eset}\big|_{\tilde{K}^{(0)}},$$
and $\ze_{\tilde{J},\ups}$ is as in~(2a) of Lemma~\ref{reg0_lmm2}.
We will show that the map 
\begin{equation}\label{reg0_crl3e2}
\tilde{\Ga}_-(b)\lra T_qX,  \qquad
\xi\lra \frac{d}{dt}
\big\{\!\exp_q\ze_{\tilde{J},(\vph(t\xi),v)}(\tilde{y}_1)\big\}\big|_{t=0},
\end{equation}
is surjective. 
Note that $q\!=\!u_{(\vph(t\xi),v)}(\hat{0},\i)$.
Along with \e_ref{reg0_crl2e1}, this claim implies 
Corollary~\ref{reg0_crl3}.\\
(3) Let $\Phi_b$, $\vt_b$, $\Phi_{\tilde{J},\ups}$, and $\vt_{\tilde{J},\ups}$
be as in Lemma~\ref{reg0_lmm3}.
Since $\vt_{\tilde{J},\ups}$ is a $J_{\ev_0(b)}$-holomorphic map on $\Si_{\ups}^0(\de_K)$,
vanishing at~$\i$,
and $\tilde{y}_1\!\in\!\Si_{\ups}^0(\de_K)$, for all 
$\ups'\!=\!(b',v)\!\in\!\wt{\cal FT}_{\de_K}^{\eset}\big|_{\tilde{K}^{(0)}}$,
\begin{equation}\label{reg0_crl3e3}\begin{split}
\vt_{\tilde{J},\ups'}(\tilde{y}_1)
&=\sum_{r=1}^{\i}\frac{1}{r!}\frac{d^r}{dw^r}\vt_{\tilde{J},\ups'}(w)\big|_{w=0}
\cdot\tilde{y}_1^{-r}
=\sum_{r=1}^{\i}{\cal D}_{\hat{0}}^{(r)}\tilde{b}_{\tilde{J}}(\ups')
\cdot\tilde{y}_1^{-r}\\
&=\sum_{k=1}^{\i}\sum_{i\in\chi({\cal T})}\!\!\!
(\tilde{y}_1\!-\!x_i(\ups'))^{-k}
\big\{\cD_{J,i}^{(k)}\!+\!\ve_{\ti{J},i}^{(k)}(\ups')\big\}\rho_i^k(\ups')
\in T_{\ev_0(b')}X,
\end{split}\end{equation}
by (2a) of Lemma~\ref{reg0_lmm4}; see the proof of 
Lemma~4.2 in~\cite{Z5} for more details.
For each $i\!\in\!\chi({\cal T})$, we denote by $h(i)$
the largest element of $\hat{I}$ such that $h(i)\!\le\!i$
and $\io_{h(i)}\!\le\!j_1$.
We~set
$$x_{h;1}(b')=x_{h(i)}(b') \qquad\hbox{and}\qquad
y_{1;h}(b')=\begin{cases}
y_1(b'),&\hbox{if}~\io_{h(i)}\!=\!j_1;\\
x_{\tilde{h}}(b'),&\hbox{if}~
\tilde{h}\!\le\!j_1~\hbox{and}~\io_{\tilde{h}}\!=\!\io_{h(i)}.
\end{cases}$$
If $\ups'\!=\!(b',v)$, we put
\begin{gather*}
\rho_{i;1}(v')=\!\prod_{h(i)\le h\le i}\!\!\!\!v_h,\qquad
x_{i;1}(\ups')=\sum_{h(i)\le i'\le i}\!\!\Big(x_{i'}(b')
\!\!\!\prod_{h(i)\le h<i'}\!\!\!\!\!\!v_h\Big)~\in~\Bbb{C},\\
y_{1;i}(\ups')=\sum_{\io_{h(i)}<i'<j_1}\!\!\!x_{i'}(b')
\!\!\!\!\!\prod_{\io_{h(i)}<h<i'}\!\!\!\!\!\!v_h
+y_1\!\!\prod_{\io_{h(i)}<h\le j_1}\!\!\!\!\!\!\!v_h~\in~\Bbb{C}.
\end{gather*}
It is straightforward to see from the definitions that
\begin{equation}\label{reg0_crl3e5}
\big(y_1(\ups')\!-\!x_i(\ups')\big)^{-1}\rho_i(\ups')=
\big(y_{1;i}(\ups')\!-\!x_{i;1}(\ups')\big)^{-1}\rho_{i;1}(\ups').
\end{equation}
By \e_ref{reg0_crl3e3} and~\e_ref{reg0_crl3e5},
\begin{equation}\label{reg0_crl3e7}\begin{split}
&\vt_{\tilde{J},\ups'}(\tilde{y}_1)
=\sum_{i\in\chi({\cal T})}\!\!\!(y_{1;i}(b')\!-\!x_{i;1}(b'))^{-1}
\big\{({\cal D}_{J,i}^{(1)}b')\!+\!\ve_{\ti{J},i}^{(0)}(\ups')\big\}\rho_{i;1}(v)
\in T_{\ev_0(b')}X,\\
&\qquad\qquad\hbox{and}\qquad
|\ve_{\ti{J},i}^{(0)}(\ups')|,|\na^{\T}\ve_{\ti{J},i}^{(0)}(\ups')|
\le  C_K\big(\|\tilde{J}\!-\!J\|_{C^1}\!+\!|\ups'|^{1/p}\big).
\end{split}\end{equation}
Finally, for all $\xi\!\in\!\tilde{T}_b\tilde{\cal U}_{\cal T}^{(0)}(X;J)$,
\begin{equation}\label{reg0_crl3e9}
\tilde{y}_1(\vph(t\xi),v)=\tilde{y}_1(b,v), \quad
\tilde{y}_{1;i}(\vph(t\xi))=\tilde{y}_{1;i}(b),~~
x_{i;1}(\vph(t\xi))=x_{i;1}(b)~~~\forall i\!\in\!\chi({\cal T}).
\end{equation}
The surjectivity of the map in~\e_ref{reg0_crl3e2} follows from \e_ref{reg0_crl3e7} 
and \e_ref{reg0_crl3e9}, along with~(2) of Lemma~\ref{reg0_lmm3},  
by the same argument as in~(3) of the proof of Corollary~\ref{reg0_crl2}.\\

\noindent
Corollaries~\ref{reg0_crl2} and~\ref{reg0_crl3} complete 
the proof of the parts of Theorem~\ref{reg_thm} that concern
genus-zero stable maps.
However, this is a convenient point to deduce a few more conclusions from 
Lemma~\ref{reg0_lmm4}.
We use Corollary~\ref{reg0_crl4} in the next three sections.
Suppose ${\cal T}\!=\!(M,I;j,\under{A})$ is any bubble type and 
$\ups\!=\!(b,v)\!\in\!\wt{\cal FT}$ is sufficiently small.
With notation as in the proof of Corollary~\ref{reg0_crl2},
we define the homomorphism
\begin{gather*}
\ti{R}_{\ti{J},\ups}\!: \ti\Ga_-(b)
\lra\ker\D_{\ti{J},\ti{u}_{\ups}}^{\i} \subset\ti{\Ga}(\ups;\ti{J}),
\qquad\hbox{where}\quad
\ti{\Ga}(\ups;\ti{J})\!=\!L^p_1(\Si_{\ups};\ti{u}_{\ups}^*TX),~~
\ti{u}_{\ups}\!=\!\exp_{b,u_{\ups}}\!\ze_{\ti{J},\ups},\\
\hbox{by}\qquad
\ti{R}_{\ti{J},\ups}\xi=
\frac{d}{dt}\exp_{u_{(\vph(t\xi),v)}}\!\ze_{\ti{J},(\vph(t\xi),v)}\big|_{t=0}.
\end{gather*}
We denote the image of $\tilde{R}_{\tilde{J},\ups}$ by $\tilde{\Ga}_-(\ups;\tilde{J})$.
Let 
\begin{alignat*}{1}
A_{\ups,\hat{0}}^-(\de)&=
\big\{(\hat{0},z)\!\in\!\{\hat{0}\}\!\times\!S^2\!: 
|z|\!\ge\!\de^{-1/2}/2\big\}\subset\Si_{\ups};\\
\partial^-A_{\ups,\hat{0}}^-(\de)&=
\big\{(\hat{0},z)\!\in\!\{\hat{0}\}\!\times\!S^2\!: 
|z|\!=\!\de^{-1/2}/2\big\}\subset\Si_{\ups}.
\end{alignat*}

\begin{crl}
\label{reg0_crl4}
If $(X,\om,J)$, $A$, and ${\cal T}$ are as in Lemma~\ref{reg0_lmm1},
for every precompact open subset $K$ of ${\cal U}_{\cal T}(X;J)$,
there exist $\de_K,\ep_K,C_K\!\in\!\Bbb{R}^+$ and
an open neighborhood $U_K\!\subset\!U_{\cal T}$ of $K$ in ${\frak X}_{0,M}(X,A)$
with the following properties:\\
(1) the requirements (1) and (2) of Lemma~\ref{reg0_lmm1} are satisfied;\\
(2) if $\tilde{J}$ is an almost complex structure on $X$ 
s.t.~$\|\tilde{J}\!-\!J\|_{C^1}\!\le\!\de_K$, 
there exists $\ze_{\tilde{J},\ups}\!\in\!\Ga_+(\ups)$ for each
$\ups\!=\!(b,v)\!\in\!\wt{\cal FT}_{\de_K}^{\eset}\big|_{\tilde{K}^{(0)}}$
such that:\\
${}\quad$ (2a) the requirements (2a)-(2d) of Lemma~\ref{reg0_lmm2} are satisfied;\\
${}\quad$ (2b) for all 
$\ups\!=\!(b,v)\!\in\!\wt{\cal FT}_{\de_K}^{\eset}\big|_{\tilde{K}^{(0)}}$
and $\de\!\le\!4\de_K$, if $\tilde{u}_{\ups}\!=\!\exp_{u_{\ups}}\!\ze_{\tilde{J},\ups}$,
$$\|d\tilde{u}_{\ups}|_{A_{\ups,\hat{0}}^-(\de)}\big\|_{\ups,p}
\le C_K\de^{1/p}\!\!\sum_{i\in\chi({\cal T})}\!\!\!\big|\rho_i(\ups)\big|;$$
${}\quad$ (2c) for all 
$\ups\!=\!(b,v)\!\in\!\wt{\cal FT}_{\de_K}^{\eset}\big|_{\tilde{K}^{(0)}}$
and $\de\!\le\!4\de_K$, if
$\exp_{\ev_0(b)}\tilde{f}_{\ups}(z)\!=\!\tilde{u}_{\ups}(z)$
for $z\!\in\!A_{\ups,\hat{0}}^-(\de_K)$,
\begin{gather*}
|\tilde{f}_{\ups}(w)|\le C_K|w|\!\! \sum_{i\in\chi({\cal T})}\!\!\!\big|\rho_i(\ups)\big|
\qquad \forall w\!\in\!A_{\ups,\hat{0}}^-(\de_K) ;\\
\Big|\oint_{\partial^-A_{\ups,\hat{0}}^-(\de)}\!
\tilde{f}_{\ups}(w)\frac{dw}{w^2}-2\pi\I
\!\!\sum_{i\in\chi({\cal T})}\!\!\!\!{\cal D}_{J,i}\rho_i(\ups)\Big|
\le C_K\big(\|\tilde{J}\!-\!J\|_{C^1}\!+\!|\ups|^{1/p}\!+\!\de^{(p-2)/p}\big)
\!\!\sum_{i\in\chi({\cal T})}\!\!\!\big|\rho_i(\ups)\big|,
\end{gather*}
where $w$ is the standard holomorphic coordinate on the complement of $(0,0,1)$ in $S^2$;\\
${}\quad$ (2d) for all 
$\ups\!=\!(b,v)\!\in\!\wt{\cal FT}_{\de_K}^{\eset}\big|_{\tilde{K}^{(0)}}$,
$\xi\!\in\!\tilde{\Ga}_-(b)$, and $\de\!\in\!(0,4\de_K)$,
\begin{gather*}
\|\tilde{R}_{\tilde{J},\ups}\xi\|_{\ups,p,1}\le C_K\|\xi\|_{b,p,1};\\
\begin{split}
&\Big|\oint_{\partial^-A_{\ups,\hat{0}}^-(\de)}\!
\{R_{\tilde{J},\ups}\xi\}(w)\frac{dw}{w^2}-2\pi\I
\!\!\sum_{i\in\chi({\cal T})}\!\!\!\rho_i(v)\D_{J,b;i}\xi\Big|\\
&\qquad\qquad\qquad\qquad\qquad\qquad\qquad
\le C_K\big(\|\tilde{J}\!-\!J\|_{C^1}\!+\!|\ups|^{1/p}\!+\!\de^{(p-2)/p}\big)
\!\!\sum_{i\in\chi({\cal T})}\!\!\!\big|\rho_i(\ups)\big|\cdot\|\xi\|_{\ups,p,1};
\end{split}\end{gather*}
${}\quad$ (2e) for all 
$\ups\!=\!(b,v)\!\in\!\wt{\cal FT}_{\de_K}^{\eset}\big|_{\tilde{K}^{(0)}}$,
$\xi\!\in\!\tilde{\Ga}_-(\ups;\tilde{J})$, and $\de\!\le\!4\de_K$,
\begin{equation*}\begin{split}
&|\xi|_w\le C_K|w|\!\! \sum_{i\in\chi({\cal T})}\!\!\!\big|\rho_i(\ups)\big|
\cdot \|\xi\|_{\ups,p,1}  \quad\forall w\!\in\!A_{\ups,\hat{0}}^-(\de);\\
&\big\|\na\xi|_{A_{\ups,\hat{0}}^-(\de)}\big\|_{\ups,p} \le 
 C_K\de^{1/p}\!\!\sum_{i\in\chi({\cal T})}\!\!\!\big|\rho_i(\ups)\big| \cdot \|\xi\|_{\ups,p,1},
\end{split}\end{equation*}
where $w$ is as in (2c).
\end{crl}

\noindent
The first part of the proof of Corollary~\ref{reg0_crl1} shows that 
$$C_K^{-1}\|\xi\|_{b,p,1}\le \|\tilde{R}_{\tilde{J},\ups}\xi\|_{\ups,p,1}
\le C_K\|\xi\|_{b,p,1}
\qquad\forall~\xi\!\in\!\tilde{\Ga}_-(b).$$
In fact, $C_K$ can be required to be arbitrary close to one in this case.
By the proof of~\e_ref{reg0_crl3e7},
\begin{equation}\label{reg0_crl4e1}
\tilde{f}_{\ups}(w)=\Phi_{\tilde{J},\ups}(w)\,\vt_{\tilde{J},\ups}(w) 
=\Phi_{\tilde{J},\ups}(w)\cdot\!\!\sum_{i\in\chi({\cal T})}\!\!\!w
\big\{({\cal D}_{J,i}^{(1)}b)\!+\!\ve_{\ti{J},i}^{(0)}(\ups)\big\}\rho_i(v)+\ve(\ups,w)
\end{equation}
for all $w\!\in\!A_{\ups,\hat{0}}^-(\de_K)$, where
$$|\ve(\ups,w)|\le C_K|w|^2\!\!
\sum_{i\in\chi({\cal T})}\!\!\!\big|\rho_i(\ups)\big|
\quad\hbox{and}\quad
|\partial_w\ve(\ups,w)|\le C_K|w|\!\!
\sum_{i\in\chi({\cal T})}\!\!\!\big|\rho_i(\ups)\big|.$$
Both estimates in (2c) are immediate from~\e_ref{reg0_crl4e1},
(2b) and (2c) of Lemma~\ref{reg0_lmm3}, and Holder's inequality, since
$$({\cal D}_i^{(1)}b)\cdot\rho_i(v)={\cal D}_{J,i}\rho_i(\ups).$$
Differentiating \e_ref{reg0_crl4e1} with respect to~$w$ and integrating,
we obtain~(2b).
The first bound in~(2e) is obtained by
differentiating \e_ref{reg0_crl4e1} with respect to~$\xi$, 
as in~(3) of the proof of Corollary~\ref{reg0_crl2}.
The second estimate in~(2e) follows by
differentiating the resulting expression with respect to~$w$ and integrating.
Finally, in the remaining statement of~(2d), 
each element $\xi(w)$  of $T_{\tilde{u}_{\ups}(w)}X$ is 
identified with its preimage in $T_{\ev_0(b)}X$ 
under the parallel transport along the geodesics.
This estimate follows by differentiating~\e_ref{reg0_crl4e1}
with respect to~$\xi$.
Due to the first bound in~(2e), the parallel transport and the geodesics can be defined
either with respect to the $J$-compatible connection~$\na^J$ 
or with the respect to the $\ti{J}$-compatible connection
$$\na^{\ti{J}}=\frac{1}{2}\big(\na^J-\ti{J}\na^J\ti{J})$$
in the bundle $TX\!\lra\!X$.

\section{Genus-One Gluing Procedures}
\label{gluing1_sec}

\subsection{A One-Step Gluing Construction}
\label{reg1_subs1}

\noindent
Our next goal is to show that the regularity condition (c) of 
Definition~\ref{g1reg_dfn} is well-behaved under small perturbations
of the almost complex structure~$J$.
Corollaries~\ref{reg1_crl1} and~\ref{reg1_crl2},
along with the compactness of the moduli space $\ov\M_1(X,J;A)$,
show that this is indeed the case.
They conclude the proof of the $g\!=\!1$ case of
the first claim in Theorem~\ref{reg_thm}.\\

\noindent
We denote by $\X_{1,M}(X,A)$ the space of equivalence classes
of all smooth maps into~$X$  from genus-one Riemann surfaces  
with marked points indexed by the set~$M$ in the homology class~$A$   
and by $\X_{1,M}^0(X,A)$ the subset of $\X_{1,M}(X,A)$
consisting of maps with smooth domains,  i.e.~smooth tori in this case.
Similarly to the previous section, we need to describe smooth
maps $\tilde{u}\!:\Si\!\lra\!X$ in $\X_{1,\eset}(X,A)$ that lie close
to each stratum ${\cal U}_{\cal T}(X;J)$ of the moduli space 
$\ov\M_{1,\eset}(X,J;A)$.
If ${\cal U}_{\cal T}(X;J)$ is contained in $\M_{1,\eset}^{\{0\}}(X,A;J)$,
the surjectivity of the operator $D_{\tilde{J},\tilde{u}}$ 
can be shown by an argument similar to the proof of Corollary~\ref{reg0_crl1}.
This case is handled in this subsection.
We will assume that $J$ is an almost complex structure that satisfies 
the regularity conditions~(a), (b-ii), and~(c) of Definition~\ref{g1reg_dfn}.\\

\noindent
Let ${\cal T}\!=\!(M,I,\aleph;j,\under{A})$ be a bubble type such 
$\sum_{i\in I}\!A_i\!=\!A$ and $A_i\!\neq\!0$ for some minimal element of~$I$.
We proceed similarly to Subsection~\ref{reg0_subs1}.
For each sufficiently small element $\ups\!=\!(b,v)$ of 
$\wt{\cal FT}^{\eset}$, let 
$$b(\ups)=\big(\Si_{\ups},j_{\ups},u_{\ups}\big), 
\qquad\hbox{where}\qquad  u_{\ups}=u_b\circ q_{\ups},$$
be the corresponding approximately holomorphic stable map.
Since the stable map $u_b$ is not constant on the principal component,
the linearization~$D_{J,b}$ of the $\bar{\partial}_J$-operator at~$b$ is surjective
by~(a), (b-ii), and~(c) of Definition~\ref{g1reg_dfn}.
Thus, if $\ups$ is sufficiently small, the linearization
$$D_{J,\ups}\!:\Ga(\ups)\!\equiv\!L^p_1(\Si_{\ups};u_{\ups}^*TX)
\lra \Ga^{0,1}(\ups;J)\!\equiv\!
L^p(\Si_{\ups};\La^{0,1}_{J,j_{\ups}}T^*\Si_{\ups}\!\otimes\!u_{\ups}^*TX),$$
of the $\bar{\partial}_J$-operator at~$b(\ups)$ ,
defined via the connection~$\na^J$, is also surjective.
In particular, we can obtain an orthogonal decomposition
$$\Ga(\ups)=\Ga_-(\ups)\oplus\Ga_+(\ups)$$ 
such that
the linear operator $D_{J,\ups}\!:\Ga_+(\ups)\!\lra\!\Ga^{0,1}(\ups;J)$ is an isomorphism,
while $\Ga_-(\ups)$ is close to $\Ga_-(b)\!\equiv\!\ker D_{J,b}$.
The relevant facts concerning the objects described in this paragraph
are summarized in Lemma~\ref{reg1_lmm1} below.\\

\noindent
{\it Remark 1:} The focus of the pregluing construction described in~\cite{Z4}
is attaching bubble trees of spheres to a fixed Riemann surface~$\Si$. 
The present situation is of course different.
However, the main ingredient in the pregluing construction
is a smooth family of nearly holomorphic maps 
$q_{\ups}\!:\Si_{\ups}\!\lra\!\Si_b$, constructed using a metric on~$\Si$.
All other objects that appear in the above paragraph are 
essentially determined by the map~$q_{\ups}$, 
and the homeomorphism type of $\Si_b$ plays little role.
In the case $\aleph\!=\!\eset$, 
i.e.~the principal component $\Si_{b,\aleph}$ of the domain $\Si_b$
of every element $b$ of $\tilde{\cal U}_{\cal T}^{(0)}(X;J)$ is a smooth torus,
we choose a family of Kahler metrics $\{g_{b_0}\}$ on the fibers of the semi-universal
bundle ${\frak U}_{{\cal T}_0}\!\lra\!{\cal U}_{{\cal T}_0}(X;J)$;
see Subsection~\ref{notation1_subs} for notation.
If $\ups\!=\!(b,v)$ is a small element of $\wt{\cal FT}^{\eset}$
and $b_0\!=\!\ti\pi_P(b)$,
we construct the map $q_{\ups}\!:\Si_{\ups}\!\lra\!\Si_b$
as in Subsection~2.2 of~\cite{Z4},
using the metric $g_{b_0}$ on~$\Si_{b,\aleph}$.\\

\noindent
{\it Remark 2:} In the case $\aleph\!\neq\!\eset$, i.e.~the principal component 
$\Si_{b,\aleph}$ of any bubble map $b\!\in\!\tilde{\cal U}_{\cal T}^{(0)}(X;J)$
is the circle of spheres $\Si_{\aleph}$, we split the pregluing construction into two steps.
The first step will correspond to gluing at the nodes~$\aleph$ of the principal curve and 
the second to attaching the trees of spheres to the resulting elliptic curve.
The bundle of gluing parameters ${\cal FT}_0$ over ${\cal U}_{{\cal T}_0}(X;J)$
has the~form 
$${\cal FT}_0=\bigoplus_{(h,i)\in\aleph}\!\!\!L_{(h,i)}
=\bigoplus_{(h,i)\in\aleph}\!\!\!L_{(h,i);0}\!\otimes\!L_{(h,i);1}$$
for some line bundles $L_{(h,i);0},L_{(h,i);1}\!\lra\!{\cal U}_{{\cal T}_0}(X;J)$.
In addition, there exist bundle~maps
$$\phi_{(h,i);0}\!:L_{(h,i);0}\lra{\frak U}_{{\cal T}_0}
\qquad\hbox{and}\qquad 
\phi_{(h,i);1}\!:L_{(h,i);1}\lra{\frak U}_{{\cal T}_0}$$
over ${\cal U}_{{\cal T}_0}(X;J)$ such that for all $b_0\!\in\!{\cal U}_{{\cal T}_0}(X;J)$
$$\phi_{(h,i);0}|_{b_0}\!:L_{(h,i);0}|_{b_0}\lra\Si_{b_0,h}
\qquad
\phi_{(h,i);1}|_{b_0}\!:L_{(h,i);1}|_{b_0}\lra\Si_{b_0,i}$$
are biholomorphisms that take $(b_0,0)$ to the node $(h,i)$ of $\Si_{b_0}$.
Let 
$$\si_{{\cal T}_0}\!:\tilde{\frak U}_{{\cal T}_0}\lra
{\cal FT}_{0;\de},
\qquad\hbox{where}\qquad \de\!\in\!C^{\i}({\cal U}_{{\cal T}_0}(X;J);\Bbb{R}^+),$$
be a semi-universal family of deformations of the elliptic curve $\Si_{\aleph}$,
along with the marked points indexed by $M_0\!\sqcup I_1$, 
where $M_0$ and $I_1$ are the sets of marked points lying on $\Si_{\aleph}$
and  of first-level bubbles of the elements of $\tilde{\cal U}_{\cal T}^{(0)}(X;J)$, 
respectively; see Subsection~\ref{notation1_subs}.
In particular, 
$$\tilde{\frak U}_{{\cal T}_0}\big|_{{\cal U}_{{\cal T}_0}^{(0)}(X;J)}=
{\frak U}_{{\cal T}_0}.$$
A small neighborhood in $\tilde{\frak U}_{{\cal T}_0}$ 
of the section ${\cal Z}_{h,i}$ of $\si_{{\cal T}_0}$ over ${\cal U}_{{\cal T}_0}(X;J)$
corresponding to the node~$(h,i)$ of $\Si_{\aleph}$ can be identified with the~set
\begin{equation*}\begin{split}
U_{(h,i)}=\big\{(b_0,v;x,y)\!:
b_0\!\in\!{\cal U}_{{\cal T}_0}(X;J);~
(b_0,v)\!\in\!{\cal FT}_0,~(b_0,x)\!\in\!L_{(h,i);0},
(b_0,y)\!\in\!L_{(h,i);1};&\\
~|v|,|x|,|y|\!<\!\de(b);~xy\!=\!v_{(h,i)}&\big\},
\end{split}\end{equation*}
in such a way that $\si_{{\cal T}_0}(b_0,v;x,y)=\!(b_0,v)$
and the fibers of $\si_{{\cal T}_0}$ are identified holomorphically.
For each $(b_0,v)\!\in\!{\cal FT}_{0;\de}$,
we set $\Si_{(b_0,v)}\!=\!\si_{{\cal T}_0}^{-1}(b_0,v)$. 
Let
$$\tilde{\si}_{{\cal T}_0}\!:\tilde{\frak U}_{{\cal T}_0}\lra
{\frak U}_{{\cal T}_0}
=\si_{{\cal T}_0}^{-1}\big({\cal U}_{{\cal T}_0}(X;J)\big)$$
be a smooth map such that 
$\tilde{\si}_{{\cal T}_0}(\Si_{(b_0,v)})\!\subset\!\Si_{b_0}$,
$\tilde{\si}_{{\cal T}_0}|_{\Si_{(b_0,v)}}$ is holomorphic outside of
the $|\aleph|$ open sets~$U_{(h,i)}$, 
$$\tilde{\si}_{{\cal T}_0}(b_0,v;x,y)=
\begin{cases}
(b_0,0;\phi_{(h,i);0}(x)),&\hbox{if}~|x|\!\ge\!2|v_{(h,i)}|^{1/2};\\
(b_0,0;\phi_{(h,i);1}(y)),&\hbox{if}~|y|\!\ge\!2|v_{(h,i)}|^{1/2};
\end{cases}
\qquad\hbox{if}~
(b_0,v;x,y)\!\in\!U_{(h,i)},$$
and $\tilde{\si}_{{\cal T}_0}(b_0,v;x_h(b_0,v))\!=\!(b_0,0;x_h(b_0,0))$ for all 
$h\!\in\!M_0\sqcup\!I_1$
and for all $(b_0,v)\!\in\!{\cal FT}_{0;\de}$,
where $x_h(b_0,v)$ is the marked point indexed by~$h$ on $\Si_{(b_0,v)}$
and $v_{(h,i)}$ is the $L_{(h,i)}$-component of~$v$.
The last condition can be used to define the points $x_h(b_0,v)$ for $v\!\neq\!0$.
Let $q_{(b_0,v);0}$ denote the restriction of $\tilde{\si}_{{\cal T}_0}$ to~$\Si_{(b_0,v)}$.
We choose a Riemannian metric on $\tilde{\frak U}$ such that 
its restriction $g_{(b_0,v)}$ to each fiber $\Si_{(b_0,v)}$ of $\si_{{\cal T}_0}$ is Kahler.
Along the way, we have made a number of choices.
These choices will be fixed once and for all.
If $\ups\!\in\!\wt{\cal FT}$, let
$$\ups_{\aleph}\!=\!(b,v_{\aleph}),
\qquad\hbox{if}\qquad
\ups\!=\!\big(b,v_{\aleph},v_{\hat{I}}\big),~~
\big(b,v_{\aleph}\big)\!\in\!\wt{{\cal F}_{\aleph}{\cal T}},~~
\big(b,v_{\hat{I}}\big)\!\in\!
\wt{{\cal F}_0{\cal T}}\oplus\wt{{\cal F}_1{\cal T}};$$
see Subsection~\ref{notation1_subs}.
If $\ups$ is sufficiently small, 
we denote by $\Si_{\ups_{\aleph}}$ the Riemann surface 
obtained from $\Si_b$ by replacing the circle of spheres $\Si_{b,\aleph}$ with 
$|M_0|\sqcup|I_1|$ marked points, which together we denote by~$b_0$, 
by $\Si_{(b_0,v_{\aleph})}$.
Let
$$q_{\ups_{\aleph}}\!:\Si_{\ups_{\aleph}}\lra\Si_b$$
be the smooth map obtained by extending the map $q_{(b_0,v);0}$
by identity to the rational components of~$\Si_b$.
We put
\begin{gather*}
\ups^{\aleph}=
\big(b(\ups_{\aleph}),(q_{\ups_{\aleph}}^*v_h)_{h\in\hat{I}}\big),
\qquad\hbox{where}\qquad
b(\ups_{\aleph})=\big(\Si_{\ups_{\aleph}},u_{\ups_{\aleph}}\big),
\quad u_{\ups_{\aleph}}=u_b\circ q_{\ups_{\aleph}},\\
\hbox{and}\qquad
q_{\ups_{\aleph}}^*v_h= 
\begin{cases}
dq_{\ups_{\aleph}}|_{x_h(\ups_{\aleph})}^{-1}v_h
\!\in\!T_{x_h(\ups_{\aleph})}\Si_{(b_0,v_{\aleph})},&
\hbox{if}~h\!\in\!I_1;\\
v_h\!\in\!\Bbb{C},&\hbox{if}~h\!\in\!\hat{I}\!-\!I_1.
\end{cases}
\end{gather*}
Let $(\Si_{\ups},j_{\ups},g_{\ups})$ be the Riemann surface obtained by attaching 
the bubble trees of spheres to the elliptic curve $\Si_{\ups_{\aleph}}$,
using the gluing parameter~$\ups^{\aleph}$ and the metric $g_{\ups_{\aleph}}$ on 
the principal component $\Si_{(b_0,v_{\aleph})}$ of~$\Si_{\ups_{\aleph}}$,
via the procedure described in Section~2 and Subsection~3.3 of~\cite{Z4}.
We take the key basic gluing map $q_{\ups}\!:\Si_{\ups}\!\lra\!\Si_b$
to be simply the composition  $q_{\ups_{\aleph}}\!\circ\!q_{\ups^{\aleph}}$.

\begin{lmm}
\label{reg1_lmm1}
Suppose $(X,\om,J)$ is a compact almost Kahler manifold, $A\!\in\!H_2(X;\Bbb{Z})$, and
$J$ satisfies the regularity conditions (a), (b-ii), and~(c) of Definition~\ref{g1reg_dfn}.
If ${\cal T}\!=\!(M,I,\aleph;j,\under{A})$ is a bubble type such that $A\!=\!\sum_{i\in I}\!A_i$
and $A_i\!\neq\!0$ for some minimal element $i$ of~$I$,
there exist $\de,C\!\in\!C({\cal U}_{\cal T}(X;J);\Bbb{R}^+)$ and 
an open neighborhood $U_{\cal T}$ of ${\cal U}_{\cal T}(X;J)$ in $\X_{1,M}(X,A)$
with the following properties:\\
(1) for all 
$\ups\!=\!(b,v)\!\in\!\wt{\cal FT}_{\de}^{\eset}$,
\begin{gather*}
\|\pi_{\ups;-}\xi\|_{\ups,p,1}\le C(b)\|\xi\|_{\ups,p,1}
\quad\forall  \xi\!\in\!\Ga(\ups), \qquad
\|D_{J,\ups}\xi\|_{\ups,p}\le C(b)|\ups|^{1/p}\|\xi\|_{\ups,p,1} 
\quad\forall\xi\!\in\!\Ga_-(\ups),\\
\hbox{and}\qquad
C(b)^{-1}\|\xi\|_{\ups,p,1}\le \|D_{J,\ups}\xi\|_{\ups,p} \le C(b)\|\xi\|_{\ups,p,1}
\quad\forall  \xi\!\in\!\Ga_+(\ups);
\end{gather*}
(2) for every $[\tilde{b}]\!\in\!\X_{1,M}^0(X,A)\cap U_{\cal T}$, 
there exist $\ups\!=\!(b,v)\!\in\!\tilde{\cal F}{\cal T}_{\de}$ and 
$\ze\!\in\!\Ga_+(\ups)$ such that $\|\ze\|_{\ups,p,1}\!<\!\de(b)$
and $[\exp_{b(\ups)}\!\ze]\!=\![\tilde{b}]$.
\end{lmm}

\noindent
This lemma is obtained by an argument analogous to that for Lemma~\ref{reg0_lmm1}.
In particular, the bijectivity arguments in Section~4 of~\cite{Z4}, with minor modifications,
apply in the present situation.

\begin{crl}
\label{reg1_crl1}
If $(X,\om,J)$, $A$, and ${\cal T}$ are as in Lemma~\ref{reg1_lmm1} and $M\!=\!\eset$, 
for every precompact open subset $K$ of ${\cal U}_{\cal T}(X;J)$,
there exist $\de_K,C_K\!\in\!\Bbb{R}^+$
and an open neighborhood $U_K\!\subset\!U_{\cal T}$ of $K$ in $\X_{1,\eset}(X,A)$
with the following properties:\\
(1) requirements (1) and (2) of Lemma~\ref{reg1_lmm1} are satisfied;\\
(2) if $\tilde{J}$ is an almost complex structure on $X$ 
s.t.~$\|\tilde{J}\!-\!J\|_{C^1}\!<\!\de_K$ and 
$[\tilde{b}]\!\in\!U_K\cap\X_{1,\eset}^0(X,A)$, there exist
a smooth Riemann surface $\Si$ and a smooth map $\tilde{u}\!:\Si\!\lra\!X$ such that
$[\tilde{b}]\!=\![\Si,\tilde{u}]$ and 
a linearization $D_{\tilde{J},\tilde{u}}$ of $\bar{\partial}_{\tilde{J}}$ at $\tilde{u}$
is surjective.
\end{crl}

\noindent
The proof is identical to that for Corollary~\ref{reg0_crl1}.

\subsection{A Two-Step Gluing Construction}
\label{reg1_subs2}

\noindent
We prove the analogue of Corollary~\ref{reg1_crl1}
for bubble types $\T\!=\!(\eset,I,\aleph;,\under{A})$ 
such that $A_i\!=\!0$ for all minimal elements~$i$ of~$I$, i.e.
$$\U_{\T}(X;J)\subset  \ov\M_{1,\eset}(X,A;J)-\M_{1,\eset}^{\{0\}}(X,A;J),$$
in the next subsection.
In this subsection, we modify the gluing construction of~\cite{Z4} in two ways.
First, we subdivide this construction into two steps.
At the first stage, we use Lemma~\ref{reg0_lmm2} to smooth out all nodes
of the domain of a stable map that lie away from the principle component.
At the second stage, we smooth out the remaining nodes, but at this step
it may not be possible to perturb each approximately holomorphic map
into a $J$-holomorphic map.
The second modification is that the second-stage approximately holomorphic maps
are closer to being holomorphic than they would be if constructed 
as in Subsections~\ref{reg0_subs1} and~\ref{reg1_subs1} and in Subsection~3.3 of~\cite{Z4}.
This modification is motivated by the pregluing construction of Section~3 in~\cite{LT}.
The two adjustments allow us to obtain estimates on the behavior of
the operator $D_{\ti{J},\ti{u}}$ that are similar to 
the estimates of Corollary~\ref{reg0_crl2} for 
the operator~$\D_{\ti{J},\ti{u}}^{\i,e_{\i}}$.\\

\noindent
If $\T\!=\!(M,I,\aleph;j,\under{A})$ is a bubble type such that $A_i\!=\!0$
for all $i\!\in\!I_0$, let $I_h\!\subset\!I$,
for $h\!\in\!I_1$, be  as in Subsection~\ref{notation1_subs}.
We~put
$$A_h({\cal T})=\sum_{i\in I_h}A_i.$$
Let $\ti{\pi}_h\!:\ti\U_{\T}^{(0)}(X;J)\!\lra\!\ti\U_{\T_h}^{(0)}(X;J)$
be the projection corresponding to the decomposition~\e_ref{g1decomp_e3a}.\\ 
If $\ups\!=\!(b,v)\in\wt{\cal FT}$, let
\begin{gather*}
\ups_0=(b,v_{\aleph},v_0),\quad \ups_1=(b,v_1), \quad
\ups_{\{h\}}=\big(\ti{\pi}_h(b),v_{\{h\}}\big)~~\hbox{for}~~h\!\in\!I_1,
\qquad\hbox{if}\\
\ups\!=\!(b,v_{\aleph},v_0,v_1),~~~
b\!\in\!\ti{\cal U}_{\T}^{(0)}(X;J),~~
(b,v_{\aleph})\!\in\!\wt{{\cal F}_{\aleph}{\cal T}},~~
(b,v_0)\!\in\!\wt{{\cal F}_0\T},~~
v_1\!=\!\big(v_{\{h\}}\big)_{h\in I_1}\in\!\bigoplus_{h\in I_1}\C^{\hat{I}_h}.
\end{gather*}
The component $\ups_1$ of $\ups$ consists of the smoothings of the nodes 
of $\Si_b$ that lie away from the principal component.
In the case of Figure~\ref{g1gendecomp_fig} on page~\pageref{g1gendecomp_fig}, 
these are the attaching nodes of the bubbles $h_2$, $h_4$, and~$h_5$.
The bubble map $b(\ups_1)$ for $\ups_1\!\in\!\wt{{\cal F}_1{\cal T}}^{\eset}$
is of the bubble type 
$$\tilde{\cal T}\!=\!(M,I_0\cup I_1,\aleph;j,\tilde{A}),
\qquad\hbox{where}\qquad
\tilde{A}_i=\begin{cases}
0,&\hbox{if}~i\!\in\!I_0;\\
A_i({\cal T}),&\hbox{if}~i\!\in\!I_1.
\end{cases}$$
Similarly to \e_ref{g1decomp_e3a} and \e_ref{g1decomp_e3b}, we put
\begin{equation}\label{g1decomp_e5}
{\cal H}_{\ti\T}(X;J)= {\cal U}_{T_0}(pt)\times 
\big\{(b_h)_{h\in I_1}\!\in\!\prod_{h\in I_1}\!\!{\cal H}_{\ti\T_h}(X;J)\!:
\ev_0(b_{h_1})\!=\!\ev_0(b_{h_2})~\forall h_1,h_2\!\in\!I_1\big\},
\end{equation}
where ${\cal H}_{\ti\T_h}(X;J)$ is the space of all $J$-holomorphic maps
from $\Bbb{P}^1$ of type~$\tilde{\cal T}_h$.
For each $h\!\in\!I_1$, $\de\!\in\!\Bbb{R}^+$, and $\ups\!\in\!\wt{\cal FT}$
as above, let
\begin{alignat*}{1}
A_{\ups_1,h}^-(\de)&=
\big\{(h,z)\!\in\!\{h\}\!\times\!S^2\!: |z|\!\ge\!\de^{-1/2}/2\big\}\subset\Si_{\ups_1,h};\\
\partial^-A_{\ups_1,h}^-(\de)&=
\big\{(h,z)\!\in\!\{h\}\!\times\!S^2\!: |z|\!=\!\de^{-1/2}/2\big\}\subset\Si_{\ups_1,h}.\\
\end{alignat*}
Finally,  if $h\!\in\!I_1$ and $i\!\in\!\chi({\cal T})\cap I_h$, we put
$$\ti\rho_i(\ups)\equiv(b,\ti{v}_i),
\qquad\hbox{where}\qquad
\tilde{v}_i=\!\prod_{h<i'\le i}\!\!\!v_{i'}\in\Bbb{C}.$$
In the case of Figure~\ref{g1gendecomp_fig} on page~\pageref{g1gendecomp_fig},
$$\ti\rho_{h_1}(\ups)=(b,1),\qquad
\ti\rho_{h_4}(\ups)=(b,v_{h_4}),\quad\hbox{and}\quad
\ti\rho_{h_5}(\ups)=(b,v_{h_5}).$$

\begin{lmm}
\label{reg1_lmm2}
Suppose $(X,\om,J)$ is a compact almost Kahler manifold, $A\!\in\!H_2(X;\Bbb{Z})$, and
$J$ is a genus-zero $A$-regular almost complex structure.
If ${\cal T}\!=\!(M,I,\aleph;j,\under{A})$ is a bubble type such that
$\sum_{i\in I}\!A_i\!=\!A$ and $A_i\!=\!0$ for all minimal elements $i$ of~$I$,
for every precompact open subset $K$ of ${\cal U}_{\cal T}(X;J)$,
there exist $\de_K,C_K\!\in\!\Bbb{R}^+$
and an open neighborhood $U_K$ of $K$ in $\X_{1,M}(X,A)$ with 
the following property.
If $\tilde{J}$ is an almost complex structure on $X$ such that 
$\|\tilde{J}\!-\!J\|_{C^1}\!\le\!\de_K$, there exist a smooth map 
$$\tilde{\phi}_{\tilde{J};1}\!:\wt{{\cal F}_1{\cal T}}_{\de_K}^{\eset}|_{\tilde{K}^{(0)}}
\lra {\cal H}_{\tilde{\cal T}}(X;\tilde{J})$$
such that\\
(1) the image of $\Im\tilde{\phi}_{\tilde{J};1}$ under the quotient map
${\cal H}_{\tilde{\cal T}}(X;\tilde{J})\!\lra\!{\cal U}_{\tilde{\cal T}}(X;\tilde{J})$ is
${\cal U}_{\tilde{\cal T}}(X;\tilde{J})\cap U_K$;\\
(2) $\ev_P(\ti\phi_{\ti{J};1}(\ups_1))\!=\!\ev_P(b)$ for all
$\ups_1\!=\!(b,v_1)\!\in\!\wt{{\cal F}_1{\cal T}}_{\de_K}^{\eset}|_{\tilde{K}^{(0)}}$;\\
(3) if $\ups_1\!\in\!\wt{{\cal F}_1{\cal T}}^{\eset}_{\de_K}|_{\tilde{K}^{(0)}}$,
$h\!\in\!I_1$, 
$\tilde{\pi}_h(\tilde{\phi}_{\tilde{J};1}(\ups_1))\!=\!(\Bbb{P}^1,\tilde{u}_{\ups_1,h})$,
and $\de\!\le\!4\de_K$, 
$$\big\|d\tilde{u}_{\ups_1,h}|_{A_{\ups_1,h}^-(\de)}\big\|_{\ups_1,p}
\le C_K\de^{1/p}\!\!\sum_{i\in\chi({\cal T})\cap I_h}\!\!\!\!\!|\tilde{\rho}_i(\ups)|.$$\\
\end{lmm}

\noindent
The smooth map~$\tilde{\phi}_{\tilde{J};1}$ is defined~by
$$\ti\pi_h(\ti\phi_{\ti{J};1}(\ups_1))=
\big(\Bbb{P}^1,\exp_{u_{\ups_{\{h\}}}}\!\ze_{\ti{J},\ups_{\{h\}}}\big)
\quad\forall h\!\in\!I_1,$$
where $\ups_{\{h\}}$ is as above and $\ze_{\ti{J},\ups_{\{h\}}}$ 
is as in (2a) of Lemma~\ref{reg0_lmm2}.
By~(2c) of Lemma~\ref{reg0_lmm2}, 
the value of the map $\ti{\pi}_h(\ti\phi_{\ti{J};1}(\ups_1))$ at the attaching node
of the bubble~$h$ is the same for all $h\!\in\!I_1$, as needed,
and~(2) of Lemma~\ref{reg1_lmm1} is satisfied.
The bound in~(3) is simply a restatement of~(2b) of Corollary~\ref{reg0_crl4},
since
$$\big|\ti\rho_i(\ups)\big|=\big|\rho_i(\ups_{\{h\}})\big| 
\qquad\forall~h\!\in\!\chi(\T)\!\cap\!I_1.$$\\

\noindent
With notation as above, for each $\ups\!=\!(b,v)$, let 
$$\ze_{\ti{J},\ups_1}\in \Ga(\Si_{\ups_1};u_{\ups_1}^*TX)$$
be given by
\begin{equation}\label{reg1_lmm2e2}
\ze_{\ti{J},\ups_1}(z) = \begin{cases}
\ze_{\ti{J},\ups_{\{h\}}}(z), &\hbox{if}~z\!\in\!\Si_{\ups_1,h},~h\!\in\!I_1;\\
0,&\hbox{otherwise}.
\end{cases}
\end{equation}
We write
$$\ti\phi_{\ti{J};1}(\ups_1)=\ti{b}_{\ti{J}}(\ups_1)=
(\Si_{\ups_1},\ti{u}_{\ups_1})\in{\cal H}_{\ti{\cal T}}(X;\ti{J}).$$
The domain $\Si_{\ups_1}$ of the stable map $\tilde{b}_{\tilde{J}}(\ups_1)$
consists of the principal component $\Si_{\ups_1,\aleph}$,
which is either a smooth torus or a circle of spheres,
and $|I_1|$ rational bubbles~$\Si_{\ups,h}$, with $h\!\in\!I_1$,
attached directly to~$\Si_{\ups_1,\aleph}$.
The $\tilde{J}$-holomorphic map~$\tilde{u}_{\ups_1}$ is constant on $\Si_{\ups_1,\aleph}$.
Let
$$\tilde{\Ga}(\ups_1;\tilde{J})=L^p_1(\Si_{\ups_1};\tilde{u}_{\ups_1}^*TX).$$
We denote by $\Ga_-(\ups_1;\tilde{J})\!\subset\!\tilde\Ga(\ups_1;\tilde{J})$ 
the kernel of the linearization $D_{\tilde{J},\tilde{b}_{\tilde{J}}(\ups_1)}$
of the $\bar{\partial}_{\tilde{J}}$-operator at~$\tilde{b}_{\tilde{J}}(\ups_1)$.\\

\noindent
For each $h\!\in\!I_1$, $\xi\!\in\!\ti\Ga_-(\ti{\pi}_h(b))$,
and $\ti\xi\!\in\!\ti\Ga(\ups_{\{h\}};\ti{J})$ with $\ti{\xi}(h,\i)\!=\!0$,
we define $\ti{R}_b\xi\!\in\!\Ga_-(b)$
and $\ti{R}_{\ups_1}\ti{\xi}\!\in\!\ti{\Ga}(\ups_1;\ti{J})$ by
\begin{equation}\label{reg1_lmm2e2a}
\ti{R}_b\xi(z)= \begin{cases}
\xi(z),&\hbox{if}~z\!\in\!\Si_{b_h};\\
0,&\hbox{otherwise}; \end{cases}
\qquad\hbox{and}\qquad
\ti{R}_{\ups_1}\ti{\xi}(z)= \begin{cases}
\ti{\xi}(z),&\hbox{if}~z\!\in\!\Si_{\ups_1,h};\\
0,&\hbox{otherwise}; \end{cases}
\end{equation}
see the paragraphs preceding Lemma~\ref{reg0_lmm2} and Corollary~\ref{reg0_crl4}
for notation.
We put
\begin{alignat}{1}
\label{reg1_lmm2e2c}
\ti\Ga_-(b)&=\big\{\ti{R}_b\xi\!:
\xi\!\in\!\ti{\Ga}_-(\ti{\pi}_h(b)),~h\!\in\!I_1\big\}\subset \Ga_-(b);\\
\ti\Ga_-(\ups_1;\ti{J})&=\big\{\ti{R}_{\ups_1}\xi\!:
\xi\!\in\!\ti\Ga_-(\ups_{\{h\}};\ti{J}),~h\!\in\!I_1\big\} \subset\Ga_-(\ups_1;\ti{J}).\notag
\end{alignat}
Let $\ti{R}_{\ups_1,\ti{J}}\!:\ti{\Ga}_-(b)\!\lra\!\ti{\Ga}_-(\ups_1;\ti{J})$
be the homomorphism such that 
$$\ti{R}_{\ups_1,\ti{J}}\ti{R}_b\xi
=\ti{R}_{\ups_1}\ti{R}_{\ups_{\{h\}},\ti{J}}\xi
\qquad\forall~\xi\!\in\!\ti{\Ga}_-(\ti{\pi}_h(b)),~h\!\in\!I_1,$$
where $\ti{R}_{\ups_{\{h\}},\ti{J}}\!:
\ti\Ga_-(\ti{\pi}_h(b))\!\lra\!\Ga_-(\ups_{\{h\}};\ti{J})$
is the homomorphism defined just before Corollary~\ref{reg0_crl4}.\\

\noindent
If $\ups\!=\!(b,v)\!\in\!\wt{\cal FT}_{\de_K}^{\eset}|_{\tilde{K}^{(0)}}$, let
$$\ups_{\tilde{J}}=\big(\tilde{b}_{\tilde{J}}(\ups_1),v_{\aleph},v_0\big).$$
We denote by $(\Si_{\ups},j_{\ups})$ the smooth Riemann surface constructed 
as in Remark~2 of Subsection~\ref{reg1_subs1} and~by 
$$q_{\ups_0;2}\!=\!q_{\ups_{\tilde{J};\aleph}}\!\circ\!q_{\ups_{\tilde{J}}^{\aleph}}\!:
\Si_{\ups}\!=\!\Si_{(\tilde{b}_{\tilde{J}}(\ups_1),v_{\aleph},v_0)}
\lra\Si_{\ups_1}\!=\!\Si_{\tilde{b}_{\tilde{J}}(\ups_1)},$$ 
the corresponding basic gluing map.
We next construct another map 
$$\ti{q}_{\ups_0;2}\!=\!q_{\ups_{\ti{J};\aleph}}\!\circ\!
\ti{q}_{\ups_{\ti{J}}^{\aleph}}\!:\Si_{\ups}\lra\Si_{\ups_1}$$
by defining the map 
$$\ti{q}_{\ups_{\ti{J}}^{\aleph}}\!:
\Si_{\ups}\lra\Si_{(\ti{b}_{\ti{J}}(\ups_1),v_{\aleph})}.$$
By construction, $\Si_{(\ti{b}_{\ti{J}}(\ups_1),v_{\aleph})}$ is 
a smooth torus $\Si_{\ups_{\aleph},\aleph}$
with $|I_1|$ bubbles attached at the points $\{x_h(\ups_{\aleph})\}_{h\in I_1}$ 
of~$\Si_{\ups_{\aleph},\aleph}$.
For each $h\!\in\!I_1$, we identify a small neighborhood $U_h(\ups_{\aleph})$
of $x_h(\ups_{\aleph})$ in~$\Si_{\ups_{\aleph},\aleph}$ with
a neighborhood of $x_h(\ups_{\aleph})$ in $T_{x_h(\ups_{\aleph})}\Si_{\ups_{\aleph},\aleph}$,
biholomorphically and isometrically, with respect to 
the metric $g_{\ups_{\aleph}}$ on $\Si_{\ups_{\aleph},\aleph}$
of Remark~2 in Subsection~\ref{reg1_subs1}.
We assume that all of these neighborhoods are disjoint from each other
and from the $|\aleph|$ thin necks of~$\Si_{\ups_{\aleph},\aleph}$.
If $z\!\in\!U_h(\ups_{\aleph})$, we denote by $z\!-\!x_h(\ups_{\aleph})$
the corresponding element of $T_{x_h(\ups_{\aleph})}\Si_{\ups_{\aleph},\aleph}$ 
and by $|z\!-\!x_h(\ups_{\aleph})|$ its norm with respect to the metric~$g_{\ups_{\aleph}}$.
Let $\be\!:\R^+\!\lra\![0,1]$ be a smooth cutoff function such that
$$\be(r)=\begin{cases}
0,&\hbox{if}~r\!\le\!1;\\
1,&\hbox{if}~r\!\ge\!2;
\end{cases}
\qquad\hbox{and}\qquad \be'(r)\!\neq\!0~\hbox{if}~r\!\in\!(1,2).$$
For each $\ep\!\in\!\Bbb{R}^+$, we define 
$\be_{\ep}\!\in\!C^{\i}(\Bbb{R};\Bbb{R})$ by $\be_{\ep}(r)\!=\!\be(r/\sqrt{\ep})$.
If $|z\!-\!x_h(\ups_{\aleph})|\!\le\!2\sqrt{\de(b)}$, we~put
\begin{alignat*}{1}
q_{\ups_{\ti{J}}^{\aleph};h}^-(z)&=
\big(1-\be_{\de_K}(2|z\!-\!x_h(\ups_{\aleph})|)\big)
\ov{\Big(\frac{v_h}{z\!-\!x_h(\ups_{\aleph})}\Big)}\in\C;\\
q_{\ups_{\ti{J}}^{\aleph};h}^+(z)&=
x_h(\ups_{\aleph})\!+\!
\be_{\de_K}(|z\!-\!x_h(\ups_{\aleph})|)\big(z\!-\!x_h(\ups_{\aleph})\big)
\in\Si_{(\ti{b}_{\ti{J}}(\ups_1),v_{\aleph}),\aleph},
\end{alignat*}
where  $v_0\!=\!(v_h)_{h\in I_1}$.
By construction, the smooth Riemann surface $\Si_{\ups}$ is 
the main component $\Si_{(\ti{b}_{\ti{J}}(\ups_1),v_{\aleph}),\aleph}$
of $\Si_{(\ti{b}_{\ti{J}}(\ups_1),v_{\aleph})}$.
We define the map 
$\ti{q}_{\ups_{\tilde{J}}^{\aleph}}\!:
\Si_{\ups}\!\lra\!\Si_{(\ti{b}_{\ti{J}}(\ups_1),v_{\aleph})}$ by
$$\ti{q}_{\ups_{\tilde{J}}^{\aleph}}(z)=
\begin{cases}
\big(h,q_S(q_{\ups_{\tilde{J}}^{\aleph};h}^-(z))\big)
\in\Si_{(\tilde{b}_{\tilde{J}}(\ups_1),v_{\aleph}),h},
&\hbox{if}~|z\!-\!x_h(\ups_{\aleph})|/\sqrt{\de_K}\!\le\!1,~h\!\in\!I_1;\\
q_{\ups_{\tilde{J}}^{\aleph};h}^+(z)\in\Si_{(\tilde{b}_{\tilde{J}}(\ups_1),v_{\aleph}),\aleph},
&\hbox{if}~1\!\le\!|z\!-\!x_h(\ups_{\aleph})|/\sqrt{\de_K}\!\le\!2,~h\!\in\!I_1;\\
z\in\Si_{(\tilde{b}_{\tilde{J}}(\ups_1),v_{\aleph}),\aleph},
&\hbox{otherwise},
\end{cases}$$
where $q_S\!:\Bbb{C}\!\lra\!S^2$ is the standard (antiholomorphic) stereographic
projection taking the origin in $\Bbb{C}$ to the south pole in~$S^2$.
Like the map $q_{\ups_0;2}$, $\tilde{q}_{\ups_0;2}$
smooths out the nodes of the principal component and
stretches small neighborhoods of the points $x_h(\ups_{\aleph})$ 
around the $|I_1|$ bubbles. Furthermore,
\begin{equation}\label{approx_e1}
\big\|d\tilde{q}_{\ups_0;2}\big\|_{C^0}\le C(b),
\end{equation}
for some $C\!\in\!C^{\i}({\cal U}_{\cal T}(X;J);\Bbb{R}^+)$,
if the norm is taken with respect to the metrics $g_{\ups}$ on $\Si_{\ups}$
and $g_{\ups_1}$ on $\Si_{\ups_1}$, constructed via the basic gluing maps
$q_{\ups}$ and~$q_{\ups_1}$, respectively; see Subsection~3.3 in~\cite{Z4}.
The map $\ti{q}_{\ups_0;2}$ is a homeomorphism outside of $|\aleph|\!+\!|I_1|$
circles of $\Si_{\ups}$ and is biholomorphic outside of 
the $|\aleph|$ thin necks corresponding to the nodes of the principal component
of~$\Si_{\aleph}$ and the $|I_1|$ annuli 
$\ti{\cal A}_{b,h}^-\!\cup\!\ti{\cal A}_{b,h}^+$ with $h\!\in\!I_1$, where
\begin{alignat*}{1}
\ti{\cal A}_{b,h}^-&=\big\{z\!\in\!\Si_{\ups}\!:
1/2\!\le\!|z\!-\!x_h(\ups_{\aleph})|/\sqrt{\de_K}\!\le\!1\big\},\\
\ti{\cal A}_{b,h}^+&=\big\{z\!\in\!\Si_{\ups}\!:
1\!\le\!|z\!-\!x_h(\ups_{\aleph})|/\sqrt{\de_K}\!\le\!2\big\}.
\end{alignat*}
The key advantage of the map $\tilde{q}_{\ups_0;2}$ over $q_{\ups_0;2}$ 
is that
\begin{equation}\label{approx_e2}
\big\|d\ti{q}_{\ups_0;2}\big\|_{C^0(\ti{\cal A}_{b,h}^-)}
\le C(b)|v_h|\qquad\forall~h\!\in\!I_1;
\end{equation}
this bound is immediate from the definition of the norms.\\

\noindent
If $\ti\phi_{\ti{J};1}(\ups_1)\!=\!(\Si_{\ups_1},\ti{u}_{\ups_1})$ as above,
we take  
$$b_{\tilde{J}}(\ups)=(\Si_{\ups},j_{\ups},u_{\ups}),
\qquad\hbox{where}\qquad u_{\ups}=\tilde{u}_{\ups_1}\circ \tilde{q}_{\ups_0;2},$$
to be the approximately $\tilde{J}$-holomorphic map corresponding to 
the gluing parameter~$\ups$ at the present, second, stage of the gluing construction.
By~\e_ref{approx_e1} and~\e_ref{approx_e2},
\begin{equation}\label{approx_e3}
\|du_{\ups}\|_{\ups,p}\le C(b)\|d\tilde{u}_{\ups_1}\|_{\ups_1,p}
\qquad\hbox{and}\qquad
\|\bar{\partial}_{\tilde{J}}u_{\ups}\|_{\ups,p}\le  C(b)\sum_{h\in I_1}
\big\|d\ti{u}_{\ups_1}|_{{\cal A}_{\ups,h}^-}\big\|_{\ups_1,p}
|v_h|^{\frac{p-2}{p}},
\end{equation}
where 
${\cal A}_{\ups,h}^-\!=\!\ti{q}_{\ups_0;2}(\ti{\cal A}_{b,h}^-)\!\subset\!\Si_{\ups_1}$
and $\|\cdot\|_{\ups,p}$ and $\|\cdot\|_{\ups_1,p}$ are the norms corresponding
to the basic gluing maps $q_{\ups}$ and~$q_{\ups_1}$;
see Subsection~3.3 in~\cite{Z4}.
The second bound follows from the fact that the map $\tilde{u}_{\ups_1}$ is $\tilde{J}$-holomorphic 
on $\Si_{\ups_1}$ and is constant on the principal component of~$\Si_{\ups_1}$;
thus, $\bar{\partial}_{\tilde{J}}u_{\ups}$ is supported 
on the annuli $\tilde{\cal A}_{b,h}^-$, with $h\!\in\!I_1$.\\

\noindent
If $\ups\!=\!(b,v)$, we denote by
$$\Ga(\ups)=L^p_1(\Si_{\ups};u_{\ups}^*TX)
\quad\hbox{and}\quad
\Ga^{0,1}(\ups;\tilde{J})=L^p(\Si_{\ups};
\La_{\tilde{J},j_{\ups}}^{0,1}T^*\Si_{\ups}\!\otimes\!u_{\ups}^*TX)$$
the Banach completions of the corresponding spaces of smooth sections
with respect to the norms $\|\cdot\|_{\ups,p,1}$ and~$\|\cdot\|_{\ups,p}$,
induced from the basic gluing map  $q_{\ups}\!:\Si_{\ups}\!\lra\!\Si_b$, as before,
and the $\ti{J}$-compatible metric
$$g_{\tilde{J}}(\cdot,\cdot)\equiv 
\frac{1}{2}\big(g_X(J\cdot,\tilde{J}\cdot)+g_X(\tilde{J}\cdot,J\cdot)\big)$$
on~$X$.
We~put
\begin{alignat*}{1}
\Ga_-(\ups;\tilde{J})&=\big\{R_{\ups_0}\xi\!: \xi\!\in\!\Ga_-(\ups_1;\tilde{J})\big\}
\subset\Ga(\ups);\notag\\
\tilde{\Ga}_-(\ups;\tilde{J})&=\big\{R_{\ups_0}\xi\!: 
\xi\!\in\!\tilde{\Ga}_-(\ups_1;\tilde{J})\big\} \subset\Ga_-(\ups;\tilde{J}),\notag
\end{alignat*}
where $R_{\ups_0}\xi\!=\!\xi\circ\ti{q}_{\ups_0;2}$.
Let
$$R_{\ups,\ti{J}}\!=\!R_{\ups_0}\!\circ\ti{R}_{\ups_1,\ti{J}}\!:
\ti{\Ga}_-(b)\lra\ti{\Ga}_-(\ups;\ti{J}).$$
With $\ze_{\ti{J},\ups_1}$ as in~\e_ref{reg1_lmm2e2}, we set
$$\Ga_-(\ups)=\big\{ R_{\ups_0} \Pi_{\ze_{\ti{J},\ups_1}}\!(\xi\!\circ\!q_{\ups_1}) \!: 
\xi\!\in\!\Ga_-(b)\big\} \subset\Ga(\ups).$$
We denote by $\Ga_+(\ups)$ the $L^2$-orthogonal complement of $\Ga_-(\ups)$ in $\Ga(\ups)$,
as in Subsections~\ref{reg0_subs1} and~\ref{reg1_subs1}.\\

\noindent
It remains to describe the obstruction bundle, i.e.~a complement to
the image of $\Ga_+(\ups)$ under $D_{J,\ups}$, or $D_{\tilde{J},\ups}$ 
if $\tilde{J}$ is sufficiently close to~$J$.
First, we describe the Fredholm situation along~$\tilde{\cal U}_{\cal T}^{(0)}(X;J)$.
The linearization $D_{J,b}$ of the $\bar{\partial}_J$-operator along
$\tilde{\cal U}_{\cal T}^{(0)}(X;J)$ is not surjective.
From the decomposition~\e_ref{g1decomp_e3a} and
the regularity conditions of Definition~\ref{g0reg_dfn},
we see that the cokernel of $D_{J,b}$, for $b\!\in\!\tilde{\cal U}_{\cal T}^{(0)}(X;J)$,
can be identified with the vector space
$$\Ga^{0,1}_-(b;J)\approx\E_{\ti\pi_P(b)}^*\otimes_J T_{\ev_P(b)}X$$
of $(T_{\ev_P(b)}X,J)$-harmonic antilinear differentials on 
the main component $\Si_{b,\aleph}$.
In other words, if $\psi\!\in\!{\cal H}_{b,\aleph}$ is a nonzero harmonic $(0,1)$-form 
on~$\Si_{b,\aleph}$,
$$\Ga^{0,1}_-(b;J)=\big\{Y\!\cdot_J\!\psi\!: Y\!\in\!T_{\ev_P(b)}X\big\}.$$
If $\aleph\!\neq\!\eset$, i.e.~$\Si_{b,\aleph}$ is a circle of spheres,
the elements of $\Ga^{0,1}_-(b;J)$ have simple poles at the nodes of $\Si_{b,\aleph}$
with the residues adding up to zero at each node.
Let 
$$\Ga_-^{0,1}(b;\ti{J})\approx
\E_{\ti\pi_P(b)}^*\otimes_{\ti{J}}T_{\ev_P(b)}X$$
be the vector space of $(T_{\ev_P(b)}X,\tilde{J})$-harmonic differentials 
on the main component $\Si_{b,\aleph}$ of~$\Si_b$.
If $\ups\!=\!(b,v)\!\in\!\wt{\cal FT}^{\eset}_{\de_K}|_{\tilde{K}^{(0)}}$,
with notation as above, let
$$\Ga_-^{0,1}(\ups_{\aleph};\tilde{J})\approx
\E_{\Si_{\ups_{\aleph},\aleph}}^*\otimes_{\tilde{J}}T_{\ev_P(b)}X$$
be the space of $(T_{\ev_P(b)}X,\tilde{J})$-harmonic differentials on the main component $\Si_{\ups_{\aleph},\aleph}$ of~$\Si_{\ups_{\aleph}}$.
If $\aleph\!=\!\eset$, 
$\Ga_-^{0,1}(\ups_{\aleph};\tilde{J})\!\equiv\!\Ga_-^{0,1}(b;\tilde{J})$.\\

\noindent
We now construct a homomorphism 
$$R_{\ups}^{\aleph}\!: \Ga_-^{0,1}(\ups_{\aleph};\tilde{J})
\lra \Ga^{0,1}(\ups;\tilde{J}).$$
For each $h\!\in\!I_1$ and $z\!\in\!A_{\ups_1,h}^-(4\de_K)$,
we define 
$$\ti\ze_{b;\ups}(z)\in T_{\ev_P(b)}X  \qquad\hbox{by}\qquad
\exp_{\ev_P(b)}^{\tilde{J}}\tilde{\ze}_{b;\ups}(z)=\tilde{u}_{\ups_1}(z),~~
\big|\tilde{\ze}_{b;\ups}(z)\big|<r_{\tilde{J}},$$
where $\exp^{\tilde{J}}$ is the exponential map for 
the connection $\na^{\ti{J}}$ and $r_{\ti{J}}$ is its injectivity radius.
If $\eta\!\in\!\Ga_-^{0,1}(\ups_{\aleph};\ti{J})$, we define 
$R_{\ups}^{\aleph}\eta\!\in\!\Ga^{0,1}(\ups;\tilde{J})$ by
$$\{R_{\ups}^{\aleph}\eta\}_zw=\begin{cases}
0,&\hbox{if}~\sqrt{\de_K}\frac{|z\!-\!x_h(\ups_{\aleph})|}{|v_h|}
\le\frac{1}{4},~h\!\in\!I_1;\\
\be_{|v_h|^2/\de_K}\big(4|z\!-\!x_h(\ups_{\aleph})|\big)
\Pi^{\ti{J}}_{\ti\ze_{b;\ups}\ti{q}_{\ups_0;2}(z)}\eta_z(w),   &\hbox{if}~ 
\frac{1}{4}\!\le\!\sqrt{\de_K}\frac{|z\!-\!x_h(\ups_{\aleph})|}{|v_h|}
\le\frac{1}{2},~h\!\in\!I_1;\\
\eta_z(w),   &\hbox{otherwise},
\end{cases}$$
for all $z\!\in\!\Si_{\ups}$ and $w\!\in\!T_z\Si_{\ups}$,
where $\Pi^{\ti{J}}$ denotes the parallel transport of the connection $\na^{\ti{J}}$.
Let $\Ga_-^{0,1}(\ups;\ti{J})$ denote the image of 
$\Ga_-^{0,1}(\ups_{\aleph};\ti{J})$ under~$R_{\ups}^{\aleph}$.\\

\noindent
If $\eta\!\in\!\Ga_-^{0,1}(b;\ti{J})$ and 
$\ti{\eta}\!\in\!\Ga_-^{0,1}(\ups_{\aleph};\ti{J})$, 
we put
$$\|\eta\|=\sum_{h\in I_1}|\eta|_{x_h(b)}  \qquad\hbox{and}\qquad
\|\ti\eta\|=\sum_{h\in I_1}|\ti\eta|_{x_h(\ti{b}_{\ti{J}}(\ups_1),v_{\aleph})},$$
where $|\eta|_{x_h(b)}$ is the norm of $\eta|_{x_h(b)}$ with respect to
the metric $g_{\ti{J}}$ on $X$ and the metric $g_{\ti\pi_P(b)}$ on $\Si_{b,\aleph}$.
Similarly, $|\ti\eta|_{x_h(\ti{b}_{\ti{J}}(\ups_1),v_{\aleph})}$
denotes the norm of $\ti\eta|_{x_h(\ti{b}_{\ti{J}}(\ups_1),v_{\aleph})}$
with respect to $g_{\ti{J}}$ and the metric $g_{(\ti\pi_P(b),v_{\aleph})}$ on 
$\Si_{(\ti{b}_{\ti{J}}(\ups_1),v_{\aleph}),\aleph}\!=\!\Si_{(\ti\pi_P(b),v_{\aleph})}$.
If $\aleph\!\neq\eset$, we can obtain an isomorphism 
$$R_{\ups_{\aleph}}\!:
\Ga_-^{0,1}(b;\tilde{J})\lra \Ga_-^{0,1}(\ups_{\aleph};\tilde{J})$$
by requiring that
$$\{R_{\ups_{\aleph}}\eta\}_{x_h(\ups_{\aleph})}
=dq_{\ups_{\aleph}}|_{x_h(\ups_{\aleph})}^*\eta|_{x_h(b)}
\qquad\forall \eta\!\in\!\Ga_-^{0,1}(b;\tilde{J}),~h\!\in\!I_1.$$
If $\aleph\!=\eset$, we take $R_{\ups_{\aleph}}$ to be the identity map.
In either case, we denote by 
$$R_{\ups}\!:\Ga^{0,1}_-(b;\tilde{J})\lra\Ga^{0,1}_-(\ups;\tilde{J})$$
the composition $R_{\ups}^{\aleph}\!\circ\!R_{\ups_{\aleph}}$.
It is immediate from this construction that for every $q\!\in\![1,2)$,
\begin{equation}\label{obsbund_e}
\big\|R_{\ups}\eta\big\|_{\ups,q} \le C_q \|\eta\| \qquad
\forall~\eta\!\in\!\Ga_-^{0,1}(b;\tilde{J}).
\end{equation}\\

\noindent
Finally, we denote by $D_{J,\ups}$ the linearization of $\bar{\partial}_J$
at $b_{\ti{J}}(\ups)$ defined via $\na^J$ and 
by $D_{\ti{J},\ups}$ the linearization of $\bar{\partial}_{\ti{J}}$
at $b_{\ti{J}}(\ups)$ defined via~$\na^{\tilde{J}}$.
Let $D_{\ti{J},\ups}^*$ denote the formal adjoint of $D_{\ti{J},\ups}$,
defined with respect to the metrics $g_{\ups}$ on $\Si_{\ups}$
and $g_{\tilde{J}}$ on~$X$; see Chapter~3 in~\cite{McSa}.
For any $h\!\in\!I_1$ and $\de\!\in\!\Bbb{R}^+$, we take
\begin{alignat*}{1}
&A_{\ups_{\aleph},h}^+(\de)=
\big\{z\!\in\!\Si_{\ups_{\aleph},\aleph}\!: 
|z\!-\!x_h(\ups_{\aleph})|\!\le\!2\de^{1/2}\big\}
\subset\Si_{\ups_{\tilde{J};\aleph},\aleph},\\
&\tilde{A}_{\ups,h}^+(\de)\!=\!
q_{\ups_{\tilde{J}}^{\aleph}}^{-1}(A_{\ups_{\aleph},h}^+(\de)), \quad
\tilde{A}_{\ups,h}^-(\de)
\!=\!q_{\ups_0;2}^{-1}\big(A_{\ups_1,h}^-(\de)\big)
\!=\!q_{\ups}^{-1}\big(A_{b,h}^-(\de)\big)\subset \Si_{\ups},\\
&\hbox{and}\qquad\partial^-\ti{A}_{\ups,h}^-(\de)=
q_{\ups_0;2}^{-1}\big(\partial^-A_{\ups_1,h}^-(\de)\big),
\end{alignat*}
where $A_{\ups_1,h}^-(\de)$ and $\partial^-A_{\ups_1,h}^-(\de)$ 
are as in the paragraph preceding Lemma~\ref{reg1_lmm2}.
If $Y_1,Y_2\!\in\!T_qX$ for some $q\!\in\!X$, we~put
$$\lr{Y_1,Y_2}_{\tilde{J}}=
g_{\tilde{J}}(Y_1,Y_2)+\I\, g_{\tilde{J}}(Y_1,\tilde{J}Y_2)\in\Bbb{C}.$$
Similarly, if $\eta_1,\eta_2\!\in\!\Ga^{0,1}(\ups;\tilde{J})$, we put
$$\llrr{\eta_1,\eta_2}=\llrr{\eta_1,\eta_2}_{\ups,2}+
\I\,\llrr{\eta_1,\tilde{J}\eta_2}_{\ups,2}\in\Bbb{C},$$
where $\llrr{\cdot,\cdot}_{\ups,2}$ is the (real-valued) $L^2$-inner product 
on $\Ga^{0,1}(\ups;\tilde{J})$ with respect to the metric $g_{\tilde{J}}$ on~$X$.
Note that by Holder's inequality and~\e_ref{obsbund_e}
\begin{equation}\label{g1reg_lmm1e6}
\big|\bllrr{\eta',R_{\ups}\eta}_{\ups}\big|
\le C\|\eta\|\|\eta'\|_{\ups,p}
\qquad\forall~\eta\!\in\!\Ga_-^{0,1}(b;\tilde{J}),~\eta'\!\in\!\Ga^{0,1}(\ups;\tilde{J}).
\end{equation}

\begin{lmm}
\label{reg1_lmm3}
If $(X,\om,J)$, $A$, and ${\cal T}$ are as in Lemma~\ref{reg1_lmm2},
for every precompact open subset $K$ of ${\cal U}_{\cal T}(X;J)$,
there exist $\de_K,C_K\!\in\!\Bbb{R}^+$
and an open neighborhood $U_K$ of $K$ in $\X_{1,M}(X,A)$ with 
the following property.
If $\tilde{J}$ is an almost complex structure on $X$ such that 
$\|\tilde{J}\!-\!J\|_{C^1}\!\le\!\de_K$,\\
(1) the second-stage pregluing map, \hbox{$\ups\!\lra\!b_{\tilde{J}}(\ups)$},
is defined on
$\wt{\cal FT}_{\de_K}^{\eset}\big|_{\tilde{K}^{(0)}}$;\\
(2) for every $[\tilde{b}]\!\in\X_{1,M}^0(X,A)\cap U_K$, 
there exist $\ups\!\in\!\tilde{\cal F}{\cal T}_{\de_K}^{\eset}|_{K^{(0)}}$ and 
$\ze\!\in\!\Ga_+(\ups)$ such that $\|\ze\|_{\ups,p,1}\!<\!\de_K$
and $[\exp_{b_{\tilde{J}}(\ups)}^{\tilde{J}}\ze]\!=\![\tilde{b}]$;\\
(3) for all $\ups\!=\!(b,v)\!\in\!\wt{\cal FT}^{\eset}_{\de_K}|_{\tilde{K}^{(0)}}$,
\begin{gather*}
\|\bar{\partial}_{\tilde{J}}u_{\ups}\|_{\ups,p}
\le  C_K|\rho(\ups)|,
\qquad \|D_{\ti{J},\ups}\xi\|_{\ups,p}\le C_K|\ups|^{\frac{p-2}{p}}\|\xi\|_{\ups,p,1} 
\quad\forall\xi\!\in\!\Ga_-(\ups;\tilde{J}),\\
\hbox{and}\qquad
C_K^{-1}\|\xi\|_{\ups,p,1}\le \|D_{\tilde{J},\ups}\xi\|_{\ups,p} 
\le C_K\|\xi\|_{\ups,p,1}   \quad\forall\xi\!\in\!\Ga_+(\ups);
\end{gather*}
(4) for all $\ups\!=\!(b,v)\!\in\!\wt{\cal FT}^{\eset}_{\de_K}|_{\tilde{K}^{(0)}}$,
$h\!\in\!I_1$, and $\xi\!\in\!\tilde{\Ga}_-(\ups;\tilde{J})$,
$$\|D_{\ti{J},\ups}\xi\|_{\ups,p}\le C_K|\rho(\ups)| \cdot\|\xi\|_{\ups,p,1}
\qquad\hbox{and}\qquad
\big\|\xi|_{\ti{A}_{\ups,h}^+(\de_K)}\big\|_{\ups,p,1}\le C_K
|\ups|^{1/p}|\rho(\ups)| \cdot
\|\xi\|_{\ups,p,1};$$
(5) for all $\ups\!=\!(b,v)\!\in\!\wt{\cal FT}^{\eset}_{\de_K}|_{\tilde{K}^{(0)}}$,
$\xi\!\in\!\tilde{\Ga}_-(b)$, and $\eta\!\in\!\Ga_-^{0,1}(b;\tilde{J})$, 
\begin{gather*}
\|R_{\ups,\tilde{J}}\xi\|_{\ups,p,1}\le\|\xi\|_{b,p,1},\qquad
\big\|R_{\ups}\eta|_{\tilde{A}_{\ups,h}^-(4\de_K)}\big\|_{\ups,2} \le
C_K|\ups|^{1/2}\|\eta\|    \quad\forall~h\!\in\!I_1,\\
\begin{split}
&\Big|\bllrr{D_{\ti{J},\ups}R_{\ups,\ti{J}}\xi,R_{\ups}\eta}+
2\pi\I\!\!\sum_{i\in\chi({\cal T})}\!\!\!\!
\blr{\D_{J,b;i}\xi,\eta_{x_{h(i)}(b)}(\rho_i(\ups))}_b\Big|\\
&\qquad\qquad\qquad\qquad\qquad\qquad
\le C_K\big(\|\tilde{J}\!-\!J\|_{C^1}\!+\!|\ups|^{1/p}\!+\!|\ups|^{(p-2)/p}\big)
|\rho(\ups)|\cdot\|\eta\|\|\xi\|_{b,p,1};
\end{split}
\end{gather*}
(6) for all $\ups\!=\!(b,v)\!\in\!\wt{\cal FT}^{\eset}_{\de_K}|_{\tilde{K}^{(0)}}$
and $\eta\!\in\!\Ga_-^{0,1}(b;\tilde{J})$,  
\begin{equation*}\begin{split}
&\Big|\bllrr{\bar{\partial}_{\tilde{J}}u_{\ups},R_{\ups}\eta}+
2\pi{\frak i}\!\!\sum_{i\in\chi({\cal T})}\!\!\!\!
\blr{{\cal D}_{J,i}\rho_i(\ups),\eta_{x_{h(i)}(b)}}_b\Big|\\
&\qquad\qquad\qquad\qquad\qquad\qquad
\le C_K\big(\|\tilde{J}\!-\!J\|_{C^1}\!+\!|\ups|^{1/p}+\!|\ups|^{(p-2)/p}\big)
|\rho(\ups)|\cdot\|\eta\|;
\end{split}\end{equation*}
(7) for all $\ups\!=\!(b,v)\!\in\!\wt{\cal FT}^{\eset}_{\de_K}|_{\tilde{K}^{(0)}}$,
$\xi\!\in\!\Ga(\ups)$, and $\eta\!\in\!\Ga_-^{0,1}(b;\tilde{J})$, 
$$\big|\llrr{D_{\tilde{J},\ups}\xi,R_{\ups}\eta}_{\ups,2}\big|
\le C_K |\ups|^{1/2}\|\eta\|\|\xi\|_{\ups,p,1}.$$\\
\end{lmm}

\noindent
{\it Remark:} In (6) above, 
$\lr{\cdot,\cdot}_b$ denotes the combination of the inner-product defined 
before Lemma~\ref{reg1_lmm3} with a contraction. 
More precisely,
$$\lr{{\cal D}_{\tilde{J},i}(b,v),\eta_x}_b=
\ov{\psi_x(v)} \lr{{\cal D}_ib,Y}
\quad\hbox{if}\quad \eta=\psi\!\otimes\!Y \in \Bbb{E}^*\!\otimes\!\ev_P^*TX;$$
see the paragraph preceding Lemma~\ref{reg0_lmm4}.\\

\noindent
The first statement of this lemma is essentially a restatement of Lemma~\ref{reg1_lmm2}, 
in the light of the constructions following Lemma~\ref{reg1_lmm2}.
In~(2), 
$$\exp_{b_{\ti{J}}(\ups)}^{\ti{J}}\!\ze=
\big(\Si_{\ups},j_{\ups},\exp_{u_{\ups}}^{\ti{J}}\!\ze\big),
\qquad\hbox{if}\quad
b_{\ti{J}}(\ups)=(\Si_{\ups},j_{\ups},u_{\ups}).$$
The arguments of Section~4 in~\cite{Z4} can be modified, 
in a straightforward way, to show that for every $[\ti{b}]\!\in\!\X_{1,M}^0(X,A)$
sufficiently close to ${\cal U}_{\T}(X;J)$,
there exists a pair $(\ups,\ze)$ as in~(2) of Lemma~\ref{reg1_lmm3}
and this pair is unique up to the action of the group 
$\Aut(\T)\!\propto\!(S^1)^{\hat{I}}$,
i.e.~the present two-stage gluing construction retains 
the essential bijectivity property of the one-stage gluing construction in~\cite{Z4}.
The key point is that the metrics $g_{\ups}$ on $\Si_{\ups}$
and the weights used to modify the standard Sobolev norms,
as in Subsection~3.3 of~\cite{Z4},
are the same in the one-stage gluing construction and in the present case, 
while the difference between the data appearing 
in the two constructions is very small.\\

\noindent
The first bound in (3) of Lemma~\ref{reg1_lmm3} is immediate from 
the second bound in~\e_ref{approx_e3} and~(2b) of Corollary~\ref{reg0_crl4}, since 
$$\ti{q}_{\ups_0;2}(\ti{\cal A}_{b,h}^-) ={\cal A}_{\ups,h}^-
\subset A_{\ups_1,h}^-\big(|v_h|^2/\de_K\big).$$
The two bounds in~(4) follow from~(2e) of Corollary~\ref{reg0_crl4} in a similar way.
The second estimate in~(3) is obtained by the same argument as
the second bound in~\e_ref{approx_e3}.
The final claim of~(3) is a consequence of the analogous inequalities for~$D_{J,\ups}$;
see Subsection~3.5 in~\cite{Z4}.
The first inequality in~(5) is clear from the first inequality in~(2d) 
of Corollary~\ref{reg0_crl4}.
For the second one, it is enough to observe that the $L^2$-norm of a one-form 
is invariant under conformal changes of the metric on a two-dimensional domain,
while the larger radius of the annulus $\tilde{A}_{\ups,h}^-(\de)$ is $|v_h|^{1/2}$,
with respect to the metric $g_{\ups_{\aleph}}$ on~$\Si_{\ups_{\aleph},\aleph}$.\\

\noindent
For the remaining three estimates, we observe that for any $h\!\in\!I_1$,
\begin{alignat}{2}\label{dadj_e1a}
|D_{\ti{J},\ups}^*R_{\ups}\eta|_{g_{\ups},z}  &\le C_K |du_{\ups}|_{g_{\ups},z}\|\eta\|
&\qquad  &\forall z\!\in\!\ti{A}_{\ups,h}^+(\de_K),\\
\label{dadj_e1b}
|D_{\ti{J},\ups}^*R_{\ups}\eta|_{g_{\ups},w_h} 
&\le C_K|du_{\ups}|_{g_{\ups},w_h}\frac{|v_h|}{|w_h|}\|\eta\|
&\qquad &\forall w_h\!=\!\frac{v_h}{z}\!\in\!\ti{A}_{\ups,h}^-(\de_K),
 \qquad\hbox{and}\\
\label{dadj_e1c}
|D_{\ti{J},\ups}^*R_{\ups}\eta|_{g_{\ups},w_h} 
&\le C_K\big(1\!+\!|du_{\ups}|_{g_{\ups},w_h}\big)|v_h|\|\eta\|
&\qquad &\forall w_h\!=\!\frac{v_h}{z}\!\in\!
\ti{A}_{\ups,h}^-(4\de_K)\!-\!\ti{A}_{\ups,h}^-(\de_K),
\end{alignat}
where $z$ is a holomorphic coordinate on a neighborhood of $x_h(\ups_{\aleph})$
in $\Si_{\aleph}$, which is unitary with respect to 
the metric $g_{\ups_{\aleph}}$ on $\Si_{\ups_{\aleph},\aleph}$, 
and $|z|$ denotes the norm in the standard metric on $\Bbb{C}$. 
These estimates are obtained by a direct computation from 
an explicit expression for $D_{\tilde{J},\ups}^*$, 
such as the one in Chapter~3 of~\cite{McSa}, and
simple facts of Riemannian geometry, such as those in Subsection~2.1 of~\cite{Z1}.
The difference between~\e_ref{dadj_e1b} and~\e_ref{dadj_e1c}
is due to the fact that the cutoff function used in the construction
of $R_{\ups}\eta$ is constant outside of the annuli 
$\ti{A}_{\ups,h}^-(4\de_K)\!-\!\ti{A}_{\ups,h}^-(\de_K)$, with~$h\!\in\!I_1$.
An explicit computation of the contribution of this cutoff function
on $\ti{A}_{\ups,h}^-(4\de_K)\!-\!\ti{A}_{\ups,h}^-(\de_K)$
is given in Subsection~2.2 of~\cite{Z3}.
From the definition of the map $\ti{q}_{\ups_0;2}$ and the metric $g_{\ups}$,
it is easy to see that 
\begin{equation}\label{dadj_e3}
\big|d\ti{q}_{\ups_0;2}\big|_{g_{\ups},z}\le 4\frac{|v_h|}{|z|^2} 
~~~\forall z\!\in\!\tilde{A}_{\ups,h}^+(\de_K)
\quad\hbox{and}\quad
\big|d\ti{q}_{\ups_0;2}\big|_{g_{\ups},w_h}\le 4
~~~\forall w_h\!=\!\frac{v_h}{z}\!\in\!\ti{A}_{\ups,h}^-(4\de_K).
\end{equation}
By \e_ref{dadj_e1a}, the first bound in~\e_ref{dadj_e3}, 
Holder's inequality, and a change of variables, we obtain
\begin{equation}\label{dadj_e4}\begin{split}
\big\|D_{\ti{J},\ups}^*R_{\ups}\eta|_{\ti{A}_{\ups,h}^+(\de_K)}\big\|_{\ups,1}
&\le C_K\big\|d\ti{u}_{\ups_1}|_{\ti{q}_{0;2}(\ti{A}_{\ups,h}^+(\de_K)}
\big\|_{\ups_1,p}|v_h|^{\frac{p-2}{p}}\|\eta\|\\
&\le C_K'|v_h|^{\frac{p-1}{p}}\!\!\!
\sum_{i\in\chi(\T)\cap I_h}\!\!\!\!\!|\tilde{\rho}_i(\ups)|\cdot\|\eta\|,
\end{split}\end{equation}
by (2b) of Corollary~\ref{reg0_crl4},
since $\ti{q}_{0;2}(\ti{A}_{\ups,h}^+(\de_K))
\!\subset\!A_{\ups_1,h}^-(|v_h|)$.
Similarly, by \e_ref{dadj_e1b}, \e_ref{dadj_e1c},
the second bound in~\e_ref{dadj_e3}, and Holder's inequality,
\begin{equation}\label{dadj_e5}\begin{split}
\big\|D_{\ti{J},\ups}^*R_{\ups}\eta|_{\ti{A}_{\ups,h}^-(4\de_K)}\big\|_{\ups,1}
&\le C_K
\big(1\!+\!\big\|d\ti{u}_{\ups_1}|_{A_{\ups_1,h}^-(4\de_K)}\big\|_{\ups_1,p}\big)
|v_h|\cdot\|\eta\|   \\
&\le C_K'|v_h|\cdot\|\eta\|,
\end{split}\end{equation}
by (2b) of Corollary~\ref{reg0_crl4}.
Since $D_{\ti{J},\ups}^*R_{\ups}\eta$ is supported on the annuli
$\ti{A}_{\ups,h}^-(\de_K)\!\cup\!\ti{A}_{\ups,h}^+(\de_K)$,
with $h\!\in\!I_1$, by \e_ref{dadj_e4} and~\e_ref{dadj_e5},
\begin{equation}\label{dadj_e6}
\|D_{\ti{J},\ups}^*R_{\ups}\eta\|_{\ups,1}
\le C_K|\ups|^{(p-1)/p}\|\eta\|.
\end{equation}
The last inequality in~Lemma~\ref{reg1_lmm3} is immediate from~\e_ref{dadj_e6},
since $p\!>\!2$.\\

\noindent
We next prove the last estimate in (5) of Lemma~\ref{reg1_lmm3}.
By the first inequalities in (2d) and~(2e) of Corollary~\ref{reg0_crl4},
for all $\xi\!\in\!\tilde{\Ga}_-(b)$, 
\begin{equation}\label{dadj_e7}
\big|R_{\ups,\tilde{J}}\xi\big|_z\le C_K\frac{|v_h|}{|z|}
\!\sum_{i\in\chi({\cal T})\cap I_h}\!\!\!\!\!|\tilde{\rho}_i(\ups)|\cdot\|\xi\|_{b,p,1}
\qquad \forall~ z\!\in\!\tilde{A}_{\ups,h}^+(\de_K).
\end{equation}
By \e_ref{dadj_e1a}, the first bound in \e_ref{dadj_e3}, \e_ref{dadj_e7},
a change of variables, and Holder's inequality, we obtain
\begin{equation}\label{dadj_e8}\begin{split}
\big|\bllrr{R_{\ups,\tilde{J}}\xi,
D_{\tilde{J},\ups}^*R_{\ups}\eta|_{\tilde{A}_{\ups,h}^+(\de_K)}}_{\ups,2}\big| 
&\le C_K\big\|d\tilde{u}_{\ups_1}|_{A_{\ups,h}^-(|\ups_h|)}\big\|_{\ups_1,p}
|v_h|^{\frac{3}{2}-\frac{1}{p}}\|\eta\|\|\xi\|_{b,p,1}\\
&\le C_K'|v_h|^{1/2}\!\!
\sum_{i\in\chi({\cal T})\cap I_h}\!\!\!\!\!|\rho_i(\ups)|\cdot\|\eta\|\|\xi\|_{b,p,1},
\end{split}\end{equation}
by (2b) of Corollary~\ref{reg0_crl4}.
Since the map $\tilde{q}_{0;2}$ is holomorphic outside of the annuli
$\tilde{\cal A}_{b,h}^{\pm}$ with $h\!\in\!I_1$ and 
$R_{\ups}\xi$ vanishes on~$\tilde{\cal A}_{b,h}^+$,
$$\bllrr{D_{\ti{J},\ups}R_{\ups,\ti{J}}\xi,R_{\ups}\eta}_{\ups,2}=
\sum_{i\in I_1}\!\int_{\tilde{\cal A}_{b,h}^-}
\blr{D_{\tilde{J},\ups}R_{\ups}\xi,R_{\ups}\eta}dz\,d\bar{z}.$$
Since $\tilde{\cal A}_{b,h}^-\!\subset\!\tilde{A}_{\ups,h}^+(\de_K)$,
by \e_ref{dadj_e8} and integration by parts,
\begin{equation*}\begin{split}
&\Big|\bllrr{D_{\tilde{J},\ups}R_{\ups,\tilde{J}}\xi,R_{\ups}\eta}+
\!\sum_{h\in I_1}\!\!
\oint_{\partial^-\tilde{\cal A}_{b,h}^-}
\blr{R_{\ups,\tilde{J}}\xi, \Pi^{\tilde{J}}_{\tilde{\ze}_{b;\ups}\tilde{q}_{0;2}(z)}
              \eta_{x_h(b)}\partial_{\bar{z}}} \, dz\Big|\\
&\qquad\qquad\qquad\qquad\qquad\qquad\qquad\qquad\qquad
\le C_K|\ups|^{1/2}\!\!
\sum_{i\in\chi({\cal T})}\!\!\!\!|\rho_i(\ups)|\cdot\|\eta\|\|\xi\|_{b,p,1},
\end{split}\end{equation*}
Thus, by a change of variables and the definition of $R_{\ups,\tilde{J}}$,
\begin{equation}\label{dadj_e10}\begin{split}
&\Big|\bllrr{D_{\tilde{J},\ups}R_{\ups,\tilde{J}}\xi,R_{\ups}\eta}-
\!\sum_{h\in I_1}\!\!
\oint_{\partial^-A_{\ups_1,h}^-(|v_h|/\de_K)}
\blr{R_{\ups_1,\tilde{J}}\xi,
\Pi^{\tilde{J}}_{\tilde{\ze}_{b;\ups}(w_h)}\eta_{x_h(b)}v_h}\frac{dw_h}{w_h^2}\Big|\\
&\qquad\qquad\qquad\qquad\qquad\qquad\qquad\qquad\qquad\qquad
\le C_K|\ups|^{1/2}|\rho(\ups)|\cdot\|\eta\|\|\xi\|_{b,p,1},
\end{split}\end{equation}
where $w_h\!=\!v_h/z$.
The last estimate in (5) of Lemma~\ref{reg1_lmm3} is immediate from \e_ref{dadj_e10} and
the second estimate in~(2d) of Corollary~\ref{reg0_crl4}.\\

\noindent
It remains to prove part (6) of Lemma~\ref{reg1_lmm3}.
Let $\tilde{\ze}_{b;\ups}\!:A_{\ups_1,h}^-(4\de_K)\!\lra\!T_{\ev_P(b)}X$
be as above.
If $h\!\in\!I_1$ and $z\!\in\!\tilde{\cal A}_{b,h}^-$,
\begin{equation}\label{dadj_e11a}
\big|\Pi^{\ti{J}~\!-1}_{\tilde{\ze}_{b;\ups}\ti{q}_{\ups_0;2}(z)}\circ
\bar{\partial}_{\ti{J}}u_{\ups}-
\bar{\partial}_{\ti{J}_{\ev_P(b)}}(\ti\ze_{b;\ups}\!\circ\!\ti{q}_{\ups_0;2})\big|_z
\le C_X|\ti\ze_{b;\ups}\!\circ\!\ti{q}_{\ups_0;2}|_z
|d(\ti{\ze}_{b;\ups}\!\circ\!\ti{q}_{\ups_0;2})|_z;
\end{equation}
see Subsection 2.3 of~\cite{Z1}.
Thus, by integration by parts, if $\eta\!=\!Y\!\otimes\!d\bar{z}$,
\begin{equation}\label{dadj_e11}\begin{split}
&\Big|\int_{\tilde{\cal A}_{b,h}^-}\!
\blr{\bar{\partial}_{\ti{J}}u_{\ups},R_{\ups}\eta}_b-
\oint_{\partial^-\ti{\cal A}_{b,h}^-}
\!\blr{\ti\ze_{b;\ups}\ti{q}_{\ups_0;2}(z),Y}dz\Big|\\
&\qquad\qquad\qquad\qquad\qquad\qquad\qquad\qquad
\le C_X\int_{\ti{\cal A}_{b,h}^-} |\ti\ze_{b;\ups}\!\circ\!\ti{q}_{\ups_0;2}|_z
\big|d(\ti\ze_{b;\ups}\!\circ\!\ti{q}_{\ups_0;2})\big|_z \,dz\,d\bar{z}\cdot\|\eta\|,
\end{split}\end{equation}
since $\ti{\ze}_{b;\ups}$ vanishes on $\Si_{\ups_1,\aleph}$.
Since
$$\oint_{\partial^-\ti{\cal A}_{b,h}^-}
\!\blr{\ti{\ze}_{b;\ups}\ti{q}_{\ups_0;2}(z),Y}dz
=-\oint_{\partial^-A_{\ups_1,h}^-(|v_h|^2/\de_K)}
\!\blr{\ti{\ze}_{b;\ups}(w_h),\eta_{x_h(b)}(v_i)}\frac{dw_h}{w_h^2},$$
where $w_h$ is as in the two previous paragraphs,
\begin{equation}\label{dadj_e12}\begin{split}
&\Big|\oint_{\partial^-\ti{\cal A}_{b,h}^-}
\!\blr{\ti{\ze}_{b;\ups}\ti{q}_{\ups_0;2}(z),Y}dz
+2\pi\I\!\!\!\sum_{i\in I_h\cap\chi({\cal T})}\!\!\!\!\!\!\!
\blr{\cD_{J,i}\rho_i(\ups),\eta_{x_h(b)}}_b\Big|\\
&\qquad\qquad\qquad\qquad\qquad\qquad
\le C_K\big(\|\tilde{J}\!-\!J\|_{C^1}\!+\!|\ups|^{1/p}\!+\!|\ups|^{(p-2)/p}\big)
\!\!\!\sum_{i\in I_h\cap\chi({\cal T})}\!\!\!\!\!\!|\rho_i(\ups)|\cdot\|\eta\|,
\end{split}\end{equation}
by the two estimates in~(2c) of Corollary~\ref{reg0_crl4}.
On the other hand,
by Holder's inequality, change of variables, (2b) and the first estimate in~(2c)
of Corollary~\ref{reg0_crl4},
\begin{equation}\label{dadj_e14}\begin{split}
\int_{\ti{\cal A}_{b,h}^-} |\ti\ze_{b;\ups}\!\circ\!\ti{q}_{\ups_0;2}|_z
\big|d(\ti\ze_{b;\ups}\!\circ\!\ti{q}_{\ups_0;2})|_zdz\,d\bar{z}
&\le C_K\big\|\ti\ze_{b;\ups}\big\|_{C^0(\ti{q}_{\ups_0;2}(\ti{\cal A}_{b,h}^-))}
\cdot|v_h|^{\frac{p-2}{p}}
\big\|d\ti{u}_{\ups_1}|_{\ti{q}_{\ups_0;2}(\ti{\cal A}_{b,h}^-)}\big\|_{\ups_1,p}\\
&\le C_K'\!\!\! \sum_{i\in I_h\cap\chi(\T)}\!\!\!\!\!\!|\rho_i(\ups)|^2,
\end{split}\end{equation}
since $\ti{q}_{\ups_0;2}(\ti{\cal A}_{b,h}^-)\!\subset\!A_{\ups_1,h}^-(|v_h|^2/\de_K)$.
Since $\bar{\partial}_{\ti{J}}u_{\ups}$ is supported on the annuli
$\ti{\cal A}_{b,h}^-$, with $h\!\in\!I_1$,
the estimate~(6) of Lemma~\ref{reg1_lmm3} follows from \e_ref{dadj_e11}-\e_ref{dadj_e14}.

\subsection{Some Geometric Conclusions}
\label{reg1_subs3}

\noindent
We now use the two-step gluing construction of the previous subsection to
conclude the proof of Theorem~\ref{reg_thm}.

\begin{crl}
\label{reg1_crl2}
Suppose $(X,\om,J)$, $A\!\neq\!0$, and ${\cal T}$ are as in Lemma~\ref{reg1_lmm2}
and $M\!=\!\eset$.
If $J$ satisfies the regularity conditions (a) and~(b-i),
for every precompact open subset $K$ of ${\cal U}_{\cal T}(X;J)$,
there exist $\de_K,C_K\!\in\!\Bbb{R}^+$
and an open neighborhood $U_K$ of $K$ in $\X_{1,\eset}(X,A)$
with the following properties:\\
(1) all requirements of Lemma~\ref{reg1_lmm3} are satisfied;\\
(2) if $\tilde{J}$ is an almost complex structure on $X$ 
s.t.~$\|\tilde{J}\!-\!J\|_{C^1}\!<\!\de_K$ and 
$[\tilde{b}]\!\in\!\M_{1,\eset}^0(X,A;\tilde{J})$, 
the linearization $D_{\tilde{J},\tilde{b}}$ of $\bar{\partial}_{\tilde{J}}$ 
at $\tilde{b}$ is surjective.
\end{crl}

\noindent
{\it Proof:} (1) We continue with the notation preceding Lemma~\ref{reg1_lmm3}.
By Lemma~\ref{reg1_lmm3}, it can be assumed that
$$\tilde{b}=\big(\Si_{\ups},j_{\ups},\tilde{u}_{\ups}),
\quad\hbox{where}\quad
\tilde{u}_{\ups}=\exp_{u_\ups}^{\tilde{J}}\!\ze_{\tilde{J},\ups},~~
\ups\!\in\!\tilde{\cal F}{\cal T}_{\de_K}^{\eset}\big|_{K^{(0)}},~
\ze_{\tilde{J},\ups}\!\in\!\Ga_+(\ups),
~\|\ze_{\tilde{J},\ups}\|_{\ups,p,1}\!\le\!\tilde{\de}_K,$$
for some $\tilde{\de}_K\!\in\!(0,\de_K)$ to be chosen later.
Since $\bar{\partial}_{\tilde{J}}\tilde{u}_{\ups}\!=\!0$,
\begin{equation}\label{reg1_crl2e1}
\bar{\partial}_{\tilde{J}}u_{\ups}+D_{\tilde{J},\ups}\ze_{\tilde{J},\ups}+
N_{\tilde{J},\ups}\ze_{\tilde{J},\ups}=0,
\end{equation}
where $N_{\tilde{J},\ups}$ is a quadratic term. 
In particular, $N_{\tilde{J},\ups}0\!=\!0$ and
\begin{equation}\label{reg1_crl2e2}
\big\|N_{\tilde{J},\ups}\xi\!-\!N_{\tilde{J},\ups}\xi'\big\|_{\ups,p}\le 
C_K\big(\|\xi\|_{\ups,p,1}\!+\!\|\xi'\|_{\ups,p,1}\big)
\|\xi\!-\!\xi'\|_{\ups,p,1},
\end{equation}
if  $\xi,\xi'\!\in\!\Ga(\ups)$ and $\|\xi\|_{\ups,p,1},\|\xi\|_{\ups,p,1}'\!\le\!\de_K$.
By \e_ref{reg1_crl2e1}, \e_ref{reg1_crl2e2}, and (3) of Lemma~\ref{reg1_lmm3},
\begin{equation}\label{reg1_crl2e3}
\big\|\ze_{\tilde{J},\ups}\big\|_{\ups,p,1} \le C_K|\rho(\ups)|,
\end{equation}
provided $\tilde{\de}_K$ is sufficiently small.\\
(2) Since $\tilde{u}_{\ups}$ is $\tilde{J}$-holomorphic, 
all linearizations~$D_{\tilde{J},\tilde{u}}$ of $\bar{\partial}_{\tilde{J}}$ are the same.
We give an explicit expression for $D_{\tilde{J},\tilde{u}}$ and
show that the dimension of its kernel does not exceed the index of~$D_{\tilde{J},\tilde{u}}$.
For any $\xi\!\in\!\Ga(\Si_{\ups};\tilde{u}_{\ups}^*TX)$, let
$$\tilde{\xi}\!=\!\Pi_{\ze_{\tilde{J},\ups}}^{\tilde{J}~\!-1}\xi.$$
We~put
\begin{equation*}\begin{split}
\bar{\partial}_{\tilde{J}}\xi
&\equiv \Pi_{\ze_{\tilde{J},\ups}}^{\tilde{J}} \circ 
\Pi_{\tilde{\xi}+\ze_{\tilde{J},\ups}}^{\tilde{J}~\!-1} \circ 
\bar{\partial}_{\tilde{J}}\exp_{u_{\ups}}^{\tilde{J}}
\big(\tilde{\xi}\!+\!\ze_{\tilde{J},\ups}\big)\\
&=\Pi_{\ze_{\tilde{J},\ups}}^{\tilde{J}} \circ
\big(\bar{\partial}_{\tilde{J}}u_{\ups}+
D_{\tilde{J},\ups}(\tilde{\xi}\!+\!\ze_{\tilde{J},\ups})+
N_{\tilde{J},\ups}(\tilde{\xi}\!+\!\ze_{\tilde{J},\ups})\big)\\
&=\Pi_{\ze_{\tilde{J},\ups}}^{\tilde{J}} \circ 
\big(D_{\tilde{J},\ups}\tilde{\xi}+
N_{\tilde{J},\ups}(\tilde{\xi}\!+\!\ze_{\tilde{J},\ups})-
N_{\tilde{J},\ups}\ze_{\tilde{J},\ups}\big),
\end{split}\end{equation*}
by \e_ref{reg1_crl2e1}.
By \e_ref{reg1_crl2e2}, we can write
$$N_{\tilde{J},\ups}(\tilde{\xi}\!+\!\ze_{\tilde{J},\ups})-
N_{\tilde{J},\ups}\ze_{\tilde{J},\ups}
=L_{\tilde{J},\ups}\tilde{\xi}+\tilde{N}_{\tilde{J},\ups}\tilde{\xi},$$
where $\tilde{N}_{\tilde{J},\ups}$ is a quadratic term, while the linear map
$L_{\tilde{J},\ups}\!:\Ga(\ups)\!\lra\!\Ga^{0,1}(\ups;\tilde{J})$ satisfies
\begin{equation}\label{reg1_crl2e5}
\big\|L_{\ti{J},\ups}\xi\big\|_{\ups,p}
\le C_K\|\ze_{\ti{J},\ups}\|_{\ups,p,1}\|\ti\xi\|_{\ups,p,1}
\le C_K'|\rho(\ups)| \cdot \|\tilde{\xi}\|_{\ups,p,1}
\qquad\forall~\ti\xi\!\in\!\Ga(\ups),
\end{equation}
by \e_ref{reg1_crl2e3}.
We conclude that
$$D_{\tilde{J},\tilde{u}}=\Pi_{\ze_{\tilde{J},\ups}} \circ
\big( D_{\tilde{J},\ups}+L_{\tilde{J},\ups}\big)\circ \Pi_{\ze_{\tilde{J},\ups}}^{-1}.$$
Thus, it is sufficient to show that the dimension of the kernel of
$D_{\tilde{J},\ups}\!+\!L_{\tilde{J},\ups}$ does not exceed 
the index of~$D_{\tilde{J},\ups}$.\\
(3) Suppose $\xi\!\in\!\ker(D_{\tilde{J},\ups}\!+\!L_{\tilde{J},\ups})$. 
Since the dimension of $\Ga_-(\ups;\tilde{J})$ is the same as the dimension of $\Ga_-(\ups)$,
by~(3) of Lemma~\ref{reg1_lmm3}, we can write 
$$\xi=\xi_-+\xi_+, \qquad\hbox{where}\qquad
\xi_-\!\in\!\Ga_-(\ups;\tilde{J}),~~\xi_+\!\in\!\Ga_+(\ups).$$
If $\tilde{\de}_K$ is sufficiently small,
by (3) of Lemma~\ref{reg1_lmm3} and \e_ref{reg1_crl2e5},
\begin{equation}\label{reg1_crl2e7}\begin{split}
\big\|\xi_+\big\|_{\ups,p,1}
&\le C_K\big(\|D_{\tilde{J},\ups}\xi_-\|_{\ups,p}+
\|L_{\tilde{J},\ups}\xi_-\|_{\ups,p}\big)\\
&\le C_K'|\ups|^{\frac{p-2}{p}}\|\xi_-\|_{\ups,p,1}
\qquad\qquad\qquad\qquad\forall~
\xi_-\!+\!\xi_+\!\in\!\ker(D_{\tilde{J},\ups}\!+\!L_{\tilde{J},\ups}).
\end{split}
\end{equation}
Thus, the projection map, $\xi\!\lra\!\xi_-$ is an injection on 
$\ker(D_{\tilde{J},\ups}\!+\!L_{\tilde{J},\ups})$.
We denote its image by~$\Ga_-^*(\ups;\tilde{J})$.
Furthermore, by~(4) of Lemma~\ref{reg1_lmm3} and~\e_ref{reg1_crl2e5},
\begin{equation}\label{reg1_crl2e8}
\begin{aligned}
\big\|\xi_+\big\|_{\ups,p,1}
&\le C_K\big(\|D_{\ti{J},\ups}\xi_-\|_{\ups,p}+\|L_{\ti{J},\ups}\xi_-\|_{\ups,p}\big)
&\qquad \forall~&\xi_-\!\in\!\Ga_-(\ups;\ti{J}),\\
&\le C_K'|\rho(\ups)|\cdot\|\xi_-\|_{\ups,p,1}
&\qquad &\xi_-\!+\!\xi_+\!\in\!\ker(D_{\ti{J},\ups}\!+\!L_{\ti{J},\ups}).
\end{aligned}\end{equation}
(4) We now use Lemma~\ref{reg1_lmm3} to estimate the $L^2$-inner product of 
$\{D_{\tilde{J},\ups}\!+\!L_{\tilde{J},\ups}\}(\xi_-\!+\!\xi_+)$
with an element $R_{\ups}\eta$ of $\Ga_-^{0,1}(\ups;\tilde{J})$,
whenever $\xi_-\!\in\!\tilde{\Ga}_-(\ups;\tilde{J})$, $\xi_+\!\in\!\Ga_+(\ups)$, and
$\xi_-\!+\!\xi_+\!\in\!\ker\{D_{\tilde{J},\ups}\!+\!L_{\tilde{J},\ups}\}$.
By~\e_ref{g1reg_lmm1e6}, \e_ref{reg1_crl2e5}, and~\e_ref{reg1_crl2e8},
\begin{equation}\label{reg1_crl2e11}
\big|\bllrr{L_{\ti{J},\ups}\xi_+,R_{\ups}\eta}\big|
\le C_K|\rho(\ups)|^2\cdot\|\eta\|\|\xi_-\|_{\ups,p,1},
\qquad \forall~\xi_-\!\in\!\ti{\Ga}_-(\ups;\ti{J}).
\end{equation}
By \e_ref{reg1_crl2e8} and (7) of Lemma~\ref{reg1_lmm3},
\begin{equation}\label{reg1_crl2e12}
\big|\bllrr{D_{\ti{J},\ups}\xi_+,R_{\ups}\eta}\big|
\le C_K|\ups|^{1/2}|\rho(\ups)|\cdot\|\eta\|\|\xi_-\|_{\ups,p,1},
\qquad \forall~\xi_-\!\in\!\ti\Ga_-(\ups;\ti{J}).
\end{equation}
For each $h\!\in\!I_1$, by \e_ref{reg1_crl2e5} and (5) of Lemma~\ref{reg1_lmm3},
\begin{equation}\label{reg1_crl2e14}\begin{split}
\big|\bllrr{L_{\ti{J},\ups}\xi|_{\ti{A}_{\ups;h}^-(4\de_K)},R_{\ups}\eta}\big|
&\le C_K\big\|\eta|_{\ti{A}_{\ups;h}^-(4\de_K)}\big\|_{\ups,2}
\|L_{\ti{J},\ups}\xi\|_{\ups,p}\\
&\le C_K|\ups|^{1/2}|\rho(\ups)|\cdot\|\eta\|\|\xi\|_{\ups,p,1};
\qquad\forall~\xi\!\in\!\Ga(\ups).
\end{split}
\end{equation}
Since the metric $g_{\ups}$ on the annulus $\tilde{A}_{\ups;h}^+(\de_K)$
differs from the standard metric on the annulus with radii $2\sqrt{\de_K}$
and $\sqrt{|\ups_h|}$ by a factor bounded above by four and below by one-quarter,
$$\|\xi\|_{C^0}\le C_K\|\xi\|_{\ups,p,1}
\qquad\forall~ \xi\!\in\!\Ga(\ti{A}_{\ups;h}^+(\de_K);u_{\ups}^*TX);$$
see Subsection~3.1 in~\cite{Z1} and Subsection~3.3
in~\cite{Z4}.
Thus,
\begin{equation}\label{reg1_crl2e15}\begin{split}
\big\|(L_{\tilde{J},\ups}\xi)|_{\tilde{A}_{\ups;h}^+(\de_K)}\big\|_{\ups,p}
&\le C_K\|\ze_{\tilde{J},\ups}\|_{\ups,p,1}
\big\|\xi|_{\tilde{A}_{\ups;h}^+(\de_K)}\big\|_{\ups,p,1}\\
&\le C_K'\!\!\sum_{i\in\chi({\cal T})\cap I_h}\!\!\!\!\!|\rho_i(\ups)|
\cdot \big\|\xi|_{\tilde{A}_{\ups;h}^+(\de_K)}\big\|_{\ups,p,1}
\qquad\forall~\xi\!\in\!\Ga(\ups),
\end{split}
\end{equation}
since \e_ref{reg1_crl2e2} is obtained from a pointwise bound; 
see Subsection~2.4 in~\cite{Z1}.
By~\e_ref{reg1_crl2e5}, \e_ref{reg1_crl2e15}, and~(4) of Lemma~\ref{reg1_lmm3},
\begin{equation}\label{reg1_crl2e16}\begin{split}
\big|\bllrr{L_{\tilde{J},\ups}\xi_-|_{\tilde{A}_{\ups;h}^+(\de_K)},R_{\ups}\eta}\big|
&\le C_K\|\eta\|\big\|\xi_-|_{\tilde{A}_{\ups;h}^+(\de_K)}\big\|_{\ups,p,1}\\
&\le C_K|\ups|^{1/p}
\!\!\sum_{i\in\chi({\cal T})\cap I_h}\!\!\!\!\!|\rho_i(\ups)|\cdot\|\eta\|\|\xi_-\|_{\ups,p,1}
\qquad\forall~\xi_-\!\in\!\tilde{\Ga}_-(\ups;\tilde{J}).
\end{split} \end{equation}
Since the intersection of the support of $\xi_-\!\in\!\tilde{\Ga}_-(\ups;\tilde{J})$
with the support of $R_{\ups}\eta\!\in\!\Ga_-^{0,1}(\ups;\tilde{J})$ is contained
in the $|I_1|$ annuli $\tilde{A}_{\ups;h}^-(4\de_K)\!\cup\!\tilde{A}_{\ups;h}^+(\de_K)$,
by~\e_ref{reg1_crl2e14} and~\e_ref{reg1_crl2e16},
\begin{equation}\label{reg1_crl2e17}
\big|\bllrr{L_{\tilde{J},\ups}\xi_-,R_{\ups}\eta}\big|
\le C_K|\ups|^{1/p}|\rho(\ups)|\cdot\|\eta\|\|\xi_-\|_{\ups,p,1}
\quad\forall~\eta\!\in\!\Ga_-^{0,1}(b;\ti{J}),~
\xi_-\!\in\!\ti\Ga_-(\ups;\ti{J}).
\end{equation}
Finally, by \e_ref{reg1_crl2e11}, \e_ref{reg1_crl2e12}, \e_ref{reg1_crl2e17}, 
and~(5) of Lemma~\ref{reg1_lmm3}, 
\begin{equation}\label{reg1_crl2e18}\begin{split}
&\Big|\bllrr{\{D_{\ti{J},\ups}\!+\!L_{\ti{J},\ups}\}(R_{\ups,\ti{J}}\xi\!+\!\xi_+),
R_{\ups}\eta}+ 2\pi\I\!\!\sum_{i\in\chi({\cal T})}\!\!\!\!
\blr{\D_{J,b;i}\xi,\eta_{x_{h(i)}(b)}(\rho_i(\ups))}_b\Big|\\
&\qquad\qquad\qquad\qquad\qquad
\le C_K\big(\|\ti{J}\!-\!J\|_{C^1}\!+\!|\ups|^{1/p}
\!+\!|\ups|^{(p-2)/p}\big)|\rho(\ups)|\cdot\|\eta\|\|\xi\|_{b,p,1},
\end{split}\end{equation}
for all $\xi\!\in\!\ti\Ga_-(b)$ and $\xi_+\!\in\!\Ga_+(\ups)$ such that 
$R_{\ups,\tilde{J}}\xi\!+\!\xi_+\!\in\!\ker(D_{\tilde{J},\ups}\!+\!L_{\tilde{J},\ups})$.\\
(5) Let $\{\eta_r\}$ be a basis for 
$\Ga^{0,1}_-(b;\tilde{J})\!=\!\E_{\ti\pi_P(b)}^*\!\otimes_{\ti{J}}\!T_{\ev_P(b)}X$,
orthonormal with respect to the inner-product corresponding to the norm~$\|\cdot\|$.
We define the homomorphism
\begin{gather*}
\D_{\ups}\!:\Ga_-^*(\ups;\ti{J})\lra \Ga^{0,1}_-(b;\ti{J})          
\qquad\hbox{by}\qquad
\D_{\ups}\xi=\sum_r
\bllrr{\{D_{\ti{J},\ups}\!+\!L_{\ti{J},\ups}\}(\xi\!+\!\xi_+),
R_{\ups}\eta_r}\eta_r\\
\quad\hbox{if}\qquad \xi_+\!\in\!\Ga_+(\ups),~
\xi\!+\!\xi_+\!\in\!\ker(D_{\ti{J},\ups}\!+\!L_{\ti{J},\ups}).
\end{gather*}
Since the projection map 
$\ker(D_{\ti{J},\ups}\!+\!L_{\ti{J},\ups})\!\lra\!\Ga_-^*(\ups;\ti{J})$
is an isomorphism by~(3) above, the map $\D_{\ups}$ is well-defined.
By definition, $\D_{\ups}\!\equiv\!0$.
On the other hand, by \e_ref{reg1_crl2e18},
\begin{equation}\label{reg1_crl2e21}
\D_{\ups}R_{\ups,\ti{J}}\xi=
-2\pi\I\!\!\sum_{i\in\chi({\cal T})}\!\!\!
\big\{\D_{J,b;i}\!+\!\ve_i(\ups)\big\}\xi\otimes\rho_i(v)
\qquad\forall\xi\in R_{\ups,\tilde{J}}^{~-1}\Ga_-^*(\ups;\tilde{J}),
\end{equation}
where $\ve_i(\ups)\!:R_{\ups,\ti{J}}^{~-1}\Ga_-^*(\ups;\ti{J})\!\lra\!T_{\ev_P(b)}X$
is a homomorphism such that 
\begin{equation}\label{reg1_crl2e22}
|\ve_i(\ups)|\le C_K
\big(\|\tilde{J}\!-\!J\|_{C^1}\!+\!|\ups|^{1/p}\!+\!|\ups|^{(p-2)/p}\big)
\qquad\forall\ups\!\in\!\tilde{\cal F}{\cal T}^{\eset}_{\de_K}\big|_{K^{(0)}}.
\end{equation}
By (a) and (b-i) of Definition~\ref{g1reg_dfn}, the map
$$\D_{J,b;i}\!: \tilde{\Ga}_-(b)\lra T_{\ev_P(b)}X$$
is surjective for all $i\!\in\!\chi({\cal T})$; 
see the paragraph preceding Lemma~\ref{reg0_lmm2}.
Since $\rho_i(\ups)\!\neq\!0$ for all $i\!\in\!I_1$ and 
$\ups\!\in\!\tilde{\cal F}{\cal T}^{\eset}$,
it follows from \e_ref{reg1_crl2e21} and~\e_ref{reg1_crl2e22}
that if $\tilde{\de}_K$ is sufficiently small,
\begin{equation*}\begin{split}
\dim\ker D_{\tilde{J},\tilde{u}}  
&= \dim\ker(D_{\tilde{J},\ups}\!+\!L_{\tilde{J},\ups})
= \dim\Ga_-^*(\ups;\tilde{J}) = \dim\ker\D_{\ups}\\
&\le \dim\Ga_-(\ups;\tilde{J})-\dim\Ga_-^{0,1}(b;\tilde{J})
=\dim\Ga_-(b)-\dim\Ga_-^{0,1}(b;J)\\
&=\ind D_{J,b}=\ind D_{\ti{J},\ti{u}},
\end{split}\end{equation*}
as needed.\\

\noindent
Corollary~\ref{reg1_crl2} concludes the proof of the genus-one regularity
property of Theorem~\ref{reg_thm}.
Corollary~\ref{reg1_crl3} below and the Gromov compactness theorem imply that 
if $J$ is a genus-zero $A$-regular almost complex structure on $X$,
$J_r$ is a sequence of almost complex structures on $X$ 
such that $J_r\!\lra\!J$ as $r\!\lra\!\i$, and  $b_r\!\in\!\M_{1,M}^0(X,A;J_r)$,
then a subsequence of $\{b_r\}$ converges to an element of $\ov\M_{1,M}^0(X,A;J)$.

\begin{crl}
\label{reg1_crl3}
If $(X,\om,J)$, $A\!\neq\!0$, and ${\cal T}$ are as in Lemma~\ref{reg1_lmm2},
for every precompact open subset $K$ of ${\cal U}_{\T}(X;J)\!-\!{\cal U}_{{\T};1}(X;J)$,
there exist $\de_K\!\in\!\Bbb{R}^+$ and
an open neighborhood $U_K$ of $K$ in $\X_{1,M}(X,A)$ such that
$$\M_{1,M}^0(X,A;\ti{J})\cap U_K=\eset$$
if $\ti{J}$ is an almost complex structure on $X$ such that 
$\|\ti{J}\!-\!J\|_{C^1}\!<\!\de_K$.
\end{crl}

\noindent
{\it Proof:}
(1) Suppose $[\tilde{b}]\!\in\!\M_{1,M}^0(X,A;\tilde{J})\cap U_K$.\
By Lemma~\ref{reg1_lmm3}, it can be assumed that
$$\tilde{b}=\big(\Si_{\ups},j_{\ups},\tilde{u}_{\ups}),
\quad\hbox{where}\quad
\tilde{u}_{\ups}=\exp_{u_\ups}^{\tilde{J}}\!\ze_{\tilde{J},\ups},~~
\ups\!\in\!\tilde{\cal F}{\cal T}_{\de_K}^{\eset}\big|_{K^{(0)}},~
\ze_{\tilde{J},\ups}\!\in\!\Ga_+(\ups),
~\|\ze_{\tilde{J},\ups}\|_{\ups,p,1}\!\le\!\tilde{\de}_K,$$
for some $\tilde{\de}_K\!\in\!(0,\de_K)$ to be chosen later.
Since $\bar{\partial}_{\tilde{J}}\tilde{u}_{\ups}\!=\!0$,
\begin{equation}\label{reg1_crl3e1}
\bar{\partial}_{\tilde{J}}u_{\ups}+D_{\tilde{J},\ups}\ze_{\tilde{J},\ups}
+N_{\tilde{J},\ups}\ze_{\tilde{J},\ups}=0.
\end{equation}
By \e_ref{reg1_crl2e2} and (3) of Lemma~\ref{reg1_lmm3},
\begin{equation}\label{reg1_crl3e2}
\big\|N_{\tilde{J},\ups}\ze_{\tilde{J},\ups}\big\|_{\ups,p}
\le C_K\|\ze_{\tilde{J},\ups}\|_{\ups,p,1}^2      \quad\Lra\quad
\big\|\ze_{\tilde{J},\ups}\big\|_{\ups,p,1} \le C_K|\rho(\ups)|,
\end{equation}
provided $\tilde{\de}_K$ is sufficiently small.\\
(2) If $\eta\!\in\!\Ga_-^{0,1}(b;\tilde{J})$, by~\e_ref{g1reg_lmm1e6},
\e_ref{reg1_crl3e2}, and (7)~of Lemma~\ref{reg1_lmm3},
\begin{equation}\label{reg1_crl3e3}\begin{split}
\big|\bllrr{D_{\ti{J},\ups}\ze_{\ti{J},\ups},R_{\ups}\eta}\big|
&\le C_K|\ups|^{1/2}|\rho(\ups)|\cdot\|\eta\|
\qquad\hbox{and}\\
\big|\bllrr{N_{\ti{J},\ups}\ze_{\ti{J},\ups},R_{\ups}\eta}\big|
&\le C_K|\rho(\ups)|^2\cdot\|\eta\|.
\end{split}\end{equation}
By \e_ref{reg1_crl3e1}, \e_ref{reg1_crl3e3}, and (6) of Lemma~\ref{reg1_lmm3},
for all $\eta\!\in\!\Ga_-^{0,1}(b;\tilde{J})$,
\begin{equation}\label{reg1_crl3e4}
\Big|\!\!\sum_{i\in\chi({\cal T})}\!\!\!\!
\blr{{\cal D}_{J,i}\rho_i(\ups),\eta_{x_{h(i)}(b)}}_b\Big|
\le C_K\big(\|\ti{J}\!-\!J\|_{C^1}\!+\!|\ups|^{1/p}\!+\!|\ups|^{(p-2)/p}\big)
|\rho(\ups)|\cdot\|\eta\|.
\end{equation}
On the other hand, since the closure of $K$ in 
${\cal U}_{\cal T}(X;J)\!-\!{\cal U}_{{\cal T};1}(X;J)$ is compact,
\begin{equation}\label{reg1_crl3e5}
\Big|\!\!\sum_{i\in\chi({\cal T})}\!\!\!\!{\cal D}_{J,i}(b,v_i)\Big|
\ge C_K^{-1}|v|
\qquad\forall~ b\!\in\!K^{(0)},~v\!=\!(v_i)_{i\in\chi({\cal T})},
\end{equation}
for some $C_K\!\in\!\R^+$, by definition of the set 
${\cal U}_{{\cal T};1}(X;J)\!\subset\!{\cal U}_{\cal T}(X;J)$;
see Theorem~\ref{str_thm}.
Since $\Ga_-^{0,1}(b;\tilde{J})\!=\!\E^*\!\otimes\!\ev_P^*TX$,
\e_ref{reg1_crl3e4} and \e_ref{reg1_crl3e5} imply that
$$\|\tilde{J}\!-\!J\|_{C^1}+|\ups|^{1/p}+|\ups|^{(p-2)/p}\ge \tilde{C}_K^{-1},$$
as needed.

\section{Completion of Proof of Theorem~\ref{comp_thm}}
\label{compthm_sec}

\subsection{Summary and Setup} 
\label{compthm0_subs}

\noindent
In this section we sketch proofs of Propositions~\ref{cuspmap_prp}-\ref{torus_prp},
based on the arguments of Sections~\ref{gluing0_sec} and~\ref{gluing1_sec}.
Detailed proofs of generalizations of these propositions can be found in 
Section~5 of~\cite{Z6}.
These three propositions are special cases of Theorem~\ref{comp_thm},
but together they imply Theorem~\ref{comp_thm} for an arbitrary
compact almost Kahler manifold $(X,\om,J)$, $J_t\!=\!J$ constant, and $A\!\in\!H_2(X;\Z)^*$. 
They also show that a limiting curve of a sequence of $J$-holomorphic curves in $X$ of
arithmetic genus of at least one must have arithmetic genus of at least one as well,
as is the case in algebraic geometry.\\

\noindent
Suppose $\{b_r\}$ is a sequence of elements of $\ov\M_{1,M}^0(X,A;J)$
such that 
$$\lim_{r\lra\i}b_r=b\in\ov\M_{1,M}(X,A;J).$$
We need to show that $b\!\in\!\ov\M_{1,M}^0(X,A;J)$.
By Definition~\ref{degen_dfn}, it is sufficient to assume that $b$ is an
element of $\U_{\T}(X;J)$ for a bubble type
$$\T=\big(M,I,\aleph;j,\under{A})$$
such that $A_i\!=\!0$ for all minimal elements $i\!\in\!I$.\\

\noindent
We can also assume that for some bubble type
$$\T'=\big(M,I',\aleph';j',\under{A}')$$
$b_r\!\in\!\U_{\T'}(X;J)$ for all $r$.
If $A_i'\!=\!0$ for all minimal elements $i\!\in\!I'$,
the desired conclusion follows from Proposition~\ref{cuspmap_prp} below,
as it implies that the second condition in Definition~\ref{degen_dfn} 
is closed with respect to the stable map topology.
If $A_i'\!\neq\!0$ for some minimal element $i\!\in\!I'$ and $\aleph'\!\neq\!\eset$,
i.e.~the principal component of $\Si_{b_r}$ is a circle of spheres,
Proposition~\ref{node_prp} implies that $b$ satisfies the second condition in
Definition~\ref{degen_dfn}.
Finally, if $\aleph'\!=\!\eset$ and $A_i'\!\neq\!0$ for the unique minimal element
$i$ of~$I'$,  the desired conclusion follows from Proposition~\ref{torus_prp}.\\

\noindent
We note that the three propositions are applied with $b$ and $b_r$ 
that are components of the ones above.

\begin{prp}
\label{cuspmap_prp}
Suppose $(X,\om,J)$ is a compact almost Kahler manifold,
$A\!\in\!H_2(X;\Z)^*$, and $M$ is a finite set.
If $[b_r]$ is a sequence of elements in $\M_{0,\{0\}\sqcup M}^0(X,A;J)$ 
such that
$$\lim_{r\lra\i}[b_r]=[b]\in\U_{\T}(X;J)\subset\ov\M_{0,\{0\}\sqcup M}(X,A;J),$$
then either\\
${}\quad$ (a) $\dim_{\C}\hbox{Span}_{(\C,J)}\{\cD_ib\!:i\!\in\!\chi(b)\}<|\chi(b)|$, or\\
${}\quad$ (b) $\bigcap_{r=1}^{\i}\ov{\bigcup_{r'>r}\C\cdot\cD_{\hat{0}}b_{r'}}\subset
\hbox{Span}_{(\C,J)}\{\cD_ib\!:i\!\in\!\chi(b)\}$.
\end{prp}

\begin{prp}
\label{node_prp}
Suppose $(X,\om,J)$ is a compact almost Kahler manifold, 
$$n\!\in\!\Z^+, \qquad A_1,\ldots,A_n\!\in\!H_2(X;\Z)^*,$$ 
and $M_1,\ldots,M_n$ are finite sets.
If $[b_{k,r}]$ is a sequence of elements in $\M_{0,\{0,1\}\sqcup M_k}^0(X,A_k;J)$
for each $k\!\in\![n]$  such~that
\begin{gather*}
\ev_1(b_{k,r})=\ev_0(b_{k+1,r}) \quad\forall k\!\in\![n\!-\!1], 
\qquad  \ev_1(b_{n,r})=\ev_0(b_{1,r}), \qquad\hbox{and}\\
\lim_{r\lra\i}[b_{k,r}]=[b_k]\in
\U_{\T^{(k)}}(X;J)\subset \ov\M_{0,\{0,1\}\sqcup M_k}(X,A_k;J)
\quad\forall k\!\in\![n],
\end{gather*}
where each $\T^{(k)}\!\equiv\!(\{1\}\!\sqcup\!M_k,I^{(k)};j^{(k)},\under{A}^{(k)})$ 
is a bubble type such that $A_i^{(k)}\!=\!0$ for all $i\!\le\!j_1$, then 
$$\dim_{\C}\hbox{Span}_{(\C,J)}\big\{\cD_ib_k\!:i\!\in\!\chi(b_k),~k\!\in\![n]\big\}
< \sum_{k=1}^{k=n}|\chi(b_k)|.$$
\end{prp}

\begin{prp}
\label{torus_prp}
Suppose $(X,\om,J)$ is a compact almost Kahler manifold,
$A\!\in\!H_2(X;\Z)^*$, and $M$ is a finite set.
If $[b_r]$ is a sequence of elements in $\M_{1,M}^0(X,A;J)$ 
such that
$$\lim_{r\lra\i}[b_r]=[b]\in
\U_{\T}(X;J)\subset  \ov\M_{1,M}(X,A;J),$$
where $\T\!\equiv\!(M,I,\aleph;j,\under{A})$ 
is a bubble type such that $A_i\!=\!0$ for all minimal elements $i\!\in\!I$,
then $\dim_{\C}\hbox{Span}_{(\C,J)}\{\cD_ib\!:i\!\in\!\chi(b)\}<|\chi(b)|$.
\end{prp}

\noindent
We prove these three propositions in the next two subsections by combining 
the approach of Sections~\ref{gluing0_sec} and~\ref{gluing1_sec}
with some aspects of the local setting used in~\cite{LT}.
The latter makes it possible to proceed with the genus-zero gluing
construction of Subsection~\ref{reg0_subs1} and the first step of
the genus-one gluing construction of Subsection~\ref{reg1_subs2}
near a given bubble map~$b$ even if $J$ is not genus-zero regular.
The maps we encounter are not holomorphic on the entire domain,
but are holomorphic on the parts of the domain that appear in Lemma~\ref{reg0_lmm3}.
This is sufficient for the purposes of the key power series expansion in 
Lemma~\ref{reg0_lmm4}.

\subsection{Proofs of Propositions~\ref{cuspmap_prp} and~\ref{node_prp}} 
\label{comp0prp_subs}

\noindent
Let $(X,\om,J)$, $A$, $M$, $b_r$, 
$$b=(M,I;x,(j,y),u), \qquad u_i\equiv u_b|_{\Si_{b,i}},$$ 
and $\T$ be as in the statement of Proposition~\ref{cuspmap_prp}.
We put
$$I^+=\big\{i\!\in\!I\!:A_i\!\neq\!0\big\}.$$
For each $i\!\in\!I^+$, choose a finite-dimensional linear subspace 
$$\ti\Ga^{0,1}_-(b;i)\subset \Ga\big(\Si_{b,i}\!\times\!X;
\La^{0,1}_{J,j}\pi_1^*T^*\Si_{b,i}\!\otimes\!\pi_2^*TX\big)$$
such that
\begin{equation*}\begin{split}
\Ga(\Si_{b,i};\La^{0,1}_{J,j}T^*\Si_{b,i}\!\otimes\!u_i^*TX\big)
=\big\{D_{J,u_{b,i}}\xi\!:\xi\!\in\!\Ga(\Si_{b,i};u_i^*TX\big),
\, \xi(\i)\!=\!0\big\}\qquad&\\
\oplus\big\{\{\id\!\times\!u_{b,i}\}^*\eta\!:\eta\!\in\!\ti\Ga^{0,1}_-(b;i)\big\}&
\end{split}\end{equation*}
and every element of $\ti\Ga^{0,1}_-(b;i)$ vanishes on a neighborhood 
of $\i\!\in\!\Si_{b,i}$ and the nodes $x_{b,h}\!\in\!\Si_{b,i}$
with $\io_h\!=\!i$.
If $i\!\in\!I\!-\!I^0$, we denote by $\ti\Ga^{0,1}_-(b;i)$ the
zero vector space.
Let
\begin{equation*}\begin{split}
\ti\U_{\T}(X;J)= \big\{b'\!\equiv\!(M,I;x',(j,y'),u')\!:~&
 b'=\hbox{bubble map of type~}\T;\\
&\bar\partial_{J,j}u_i'\in
\{\id\!\times\!u_i'\}^*\ti\Ga^{0,1}_-(b;i)~\forall i\!\in\!I\big\}.
\end{split}\end{equation*}
By the Implicit Function Theorem, $\ti\U_{\T}(X;J)$ is a smooth manifold 
near~$b$.
Let
$$\wt{\cal FT}\equiv \ti\U_{\T}(X;J)\times\C^{\hat{I}}$$
be the bundle of smoothing parameters.\\

\noindent
Since the sequence $[b_r]$ converges to $[b]$, for all $r$ sufficiently large
there exist\\
$$b_r'\in\ti\U_{\T}(X;J), \qquad
\ups_r=(b_r',v_r)\in\wt{\cal FT}^{\eset}, \quad\hbox{and}\quad
\xi_r\in\Ga(\ups_r)\!\equiv\!\Ga(\Si_{\ups_r};u_{\ups_r}^*TX)$$
such that
\begin{gather}\label{cuspmap_e0}
\xi_r(\i)=0~~~\forall r,\qquad
\lim_{r\lra\i}b_r'=b, \qquad \lim_{r\lra\i}|v_r|=0, \qquad
\lim_{r\lra\i}\|\xi_r\|_{\ups_r,p,1}=0,\\
\hbox{and}\qquad b_r=\exp_{b(\ups_r)}\xi_r.\notag
\end{gather}
The last equality holds for a representative $b_r$ for $[b_r]$.\\

\noindent
{\it Remark:}
The existence of $b_r'$, $\ups_r$, and $\xi_r$ as above can be shown 
by an argument similar to the surjectivity argument in Section~4 of~\cite{Z4}, 
with significant simplifications.
In fact, the only facts about the bubble maps $b_r'$ we use below are that they
are constant on the degree-zero components and holomorphic on fixed neighborhoods
of the attaching nodes of the first-level effective bubbles,
i.e.~on $\Si_{b_r'}^0(\de)$ in the notation of Subsection~\ref{reg0_subs2}.
Such bubble maps $b_r'$, along with $\ups_r$ and $\xi_r$, can be constructed directly from
the maps $b_r$; see the beginning of Subsection~4.4 in~\cite{Z4}.\\

\noindent
By the same argument as in the proofs of Lemma~\ref{reg0_lmm4} and 
Corollary~\ref{reg0_crl2}, but now applied to the sequence~$(\ups_r,\xi_r)$
with sufficiently small~$\de_K$, 
\begin{equation}\label{cuspmap_e1}
\Big|\cD_{\hat{0}}b_r-
\sum_{i\in\chi(\T)}\!\!\big(\cD_ib_r'\big)\rho_i(v_r)\Big|
\le C\big(|\ups_r|^{1/p}\!+\!\|\xi_r\|_{\ups_r,p,1}\big)
\sum_{i\in\chi(\T)}\!\!\big|\rho_i(v_r)\big|.
\end{equation}
This estimate follows from equation~\e_ref{reg0_crl2e3} with 
$b'$, $\ups'$, and $\ze_{\ti{J},\ups'}$ replaced by 
$b_r'$, $\ups_r$, and $\xi_r$, respectively. 
Recall that $\Phi_{\ti{J},\ups'}(\i)\!=\!\id$ for $\ti{J}\!=\!J$.
Since $b_r'\!\lra\!b$, \e_ref{cuspmap_e1} implies that
\begin{equation}\label{cuspmap_e2}
\Big|\cD_{\hat{0}}b_r-
\sum_{i\in\chi(\T)}\!\!\big(\cD_ib\big)\rho_i(v_r)\Big|
\le C\big(|\ups_r|^{1/p}\!+\!\|\xi_r\|_{\ups_r,p,1}\big)
\sum_{i\in\chi(\T)}\!\!\big|\rho_i(v_r)\big|,
\end{equation}
where the difference is computed via a parallel transport 
of $T_{\ev_0(b_r')}X$ to $T_{\ev_0(b)}X$ with respect to the $J$-linear connection~$\na^J$.
By~\e_ref{cuspmap_e0} and \e_ref{cuspmap_e2}, $b$ must satisfy one of
the two conditions in the statement of Proposition~\ref{cuspmap_prp}.\\

\noindent
The proof of Proposition~\ref{node_prp} involves a similar extension
of Lemma~\ref{reg0_lmm4} and Corollary~\ref{reg0_crl3}.
By the assumption on the bubble types $\T^{(k)}$ made in Proposition~\ref{node_prp},
$\ev_0(b_k)\!=\!\ev_1(b_k)$ for all~$k$.
Thus,
$$\ev_1(b_k)=\ev_0(b_k)=\ev_1(b_l) \qquad\forall k,l\in[n].$$
Let $q$ denote the point $\ev_0(b_1)$.
We identify a small neighborhood of $q$ in $X$ with a small neighborhood 
of $q$ in $T_qX$ via the exponential map of the metric $g_X$ and 
the tangent space to $X$ at a point close to $q$
with $T_qX$ via the $\na^J$-parallel transport.\\

\noindent
For each pair $(k,r)$, with $r$ sufficiently large, 
let $(b_{k,r}',\ups_{k,r},\xi_{k,r})$ be an analogue of $(b_r',\ups_r,\xi_r)$ 
for~$b_{k,r}$.
As before, the key point is that the bubble maps $b_{k,r}'$ are constant on 
the degree-zero components and holomorphic on fixed neighborhoods
of the attaching nodes of the first-level effective bubbles.
Let
\begin{gather*}
\ze_{k,r}=\ev_0(b_{k,r}')\in T_qX  \qquad\hbox{and}\\
\ti\ze_{k,r}=\ev_1(b_{k,r})-\ev_0(b_{k,r})
=\ev_1(b_{k,r})-\ev_0(b_{k,r}') \in T_qX.
\end{gather*}
By the assumption on the maps $b_{k,r}$ made in the statement of Proposition~\ref{node_prp},
\begin{gather}\begin{split}
&\big|\ze_{k,r}+\ti\ze_{k,r}-\ze_{k+1,r}\big| 
\le C\big|\ze_{k,r}\big|\cdot\big|\ti\ze_{k,r}\big|
\qquad\forall k\!\in[n\!-\!1],\\
&\big|\ze_{n,r}+\ti\ze_{n,r}-\ze_{1,r}\big| 
\le C\big|\ze_{n,r}\big|\cdot\big|\ti\ze_{n,r}\big|;
\end{split}\notag\\
\label{node_e1}
\Lra\qquad 
\big|\ti\ze_{1,r}+\ldots+\ti\ze_{n,r}\big|
\le \ep_r\sum_{k=1}^{k=n}\big|\ti\ze_{k,r}\big|,
\end{gather}
for a sequence $\{\ep_r\}$ converging to $0$.
On the other hand, by the proof of~\e_ref{reg0_crl3e7},
\begin{equation}\label{node_e2}\begin{split}
&\Big|\ti\ze_{k,r}-
\sum_{i\in\chi(\T^{(k)})}\!\!(y_{1;i}(b_{k,r}')\!-\!x_{i;1}(b_{k,r}'))^{-1}
(\cD_ib_{k,r}')\rho_{i;1}(v_{k,r}) \Big| \\
&\qquad\qquad\qquad\qquad\qquad\qquad
\le C\big(|\ups_{k,r}|^{1/p}\!+\!\|\xi_{k,r}\|_{\ups_{k,r},p,1}\big)
\sum_{i\in\chi(\T^{(k)})}\!\!\big|\rho_{i;1}(v_{k,r})\big|;
\end{split}\end{equation}
see (3) of the proof of Corollary~\ref{reg0_crl3} for notation.
By \e_ref{node_e1} and \e_ref{node_e2},
\begin{equation}\label{node_e3}
\Big|\sum_{k=1}^{k=n}\sum_{i\in\chi(\T^{(k)})}\!\!(y_{1;i}(b_{k,r}')\!-\!x_{i;1}(b_{k,r}'))^{-1}
({\cal D}_ib_{k,r}')\rho_{i;1}(v_{k,r}) \Big|
\le\ti\ep_k\sum_{k=1}^{k=n}\sum_{i\in\chi(\T^{(k)})}
\!\!\big|\rho_{i;1}(v_{k,r})\big|,
\end{equation}
for a sequence $\{\ti\ep_r\}$ converging to $0$.
Since $\cD_ib_{k,r}'\!\lra\!\cD_ib_k$ as $r\!\lra\!\i$, 
\e_ref{node_e3} implies the conclusion of Proposition~\ref{node_prp}.

\subsection{Proof of Proposition~\ref{torus_prp}} 
\label{comp1prp_subs}

\noindent
Let $(X,\om,J)$, $A$, $M$, $b_r$, 
$$b=(M,I,\aleph;S,x,(j,y),u), \qquad u_i\equiv u_b|_{\Si_{b,i}},$$ 
and $\T$ be as in the statement of Proposition~\ref{torus_prp}.
Let $\T_h$ and $b_h$ for $h\!\in\!I_1$ be as in Subsection~\ref{notation1_subs}.
For each $h\!\in\!I_1$ and $i\!\in\!I_h^+$, choose a subspace
$$\ti\Ga^{0,1}_-(b;i)\!\equiv\!\ti\Ga^{0,1}_-(b_h;i)
\subset \Ga\big(\Si_{b,i}\!\times\!X;
\La^{0,1}_{J,j}\pi_1^*T^*\Si_{b,i}\!\otimes\!\pi_2^*TX\big)$$
as in the previous subsection.
If $A_i\!=\!0$, denote by $\ti\Ga^{0,1}_-(b;i)$ the zero vector space.
We define $\ti\U_{\T}(X;J)$ as at the beginning of 
Subsection~\ref{comp0prp_subs}.
Let
$$\wt{\cal FT} \lra \ti\U_{\T}(X;J)$$
be the bundle of gluing parameters.
For each $b'\!\in\!\ti\U_{\T}(X;J)$, let 
$$\ti\Ga^{0,1}_-(b') \subset \Ga\big(\Si_{b'}\!\times\!X;
\La^{0,1}_{J,j}\pi_1^*T^*\Si_{b'}\!\otimes\!\pi_2^*TX\big)$$
be the subspace obtained by extending all elements of 
$\ti\Ga^{0,1}(b';i)\!=\!\ti\Ga^{0,1}(b;i)$
by zero outside of the component $\Si_{b',i}$ of~$\Si_{b'}$.\\

\noindent
The sequence $[b_r]$ converges to $[b]$.
Thus, with notation as in Subsection~\ref{reg1_subs2}, for all $r$ sufficiently large 
there exist
\begin{gather*}
b_r'\in\ti\U_{\T}(X;J), \qquad \ups_r=(b_r',v_r)\in\wt{\cal FT}^{\eset},
\qquad \xi_{r;1}\in\Ga(\ups_{r;1}), \qquad\hbox{and}\\
\xi_{r;2}\in\Ga_+(\ups_r)\subset
\Ga\big(\ti{u}_{\ups_{r;1}}\!\circ\!\ti{q}_{\ups_{r;0};2}\big),
\qquad\hbox{where}\quad
\ti{u}_{\ups_{r;1}}=\exp_{u_{\ups_{r;1}}}\!\xi_{r;1},
\end{gather*}
such that
\begin{gather}
\label{torusmap_e0a}
\xi_{r;1}|_{\Si_{\ups_{r;1};\aleph}} =0, \quad
\bar{\partial}_J\ti{u}_{\ups_{r;1}}
\in \{q_{\ups_{r;1}}\!\times\!\ti{u}_{\ups_{r;1}}\big\}^* \ti\Ga^{0,1}_-(b_r'), 
\quad  b_r=\exp_{\ti{u}_{\ups_{r;1}}\circ\ti{q}_{\ups_{r;0};2}}\!\xi_{r;2}  
 \qquad\forall r,\\
\label{torusmap_e0b}
\lim_{r\lra\i}b_r'=b, \quad \lim_{r\lra\i}|v_r|=0, \quad
\lim_{r\lra\i}\|\xi_{r;1}\|_{\ups_{r;1},p,1}=0,\quad
\lim_{r\lra\i}\|\xi_{r;2}\|_{\ups_r,p,1}=0.
\end{gather}
We note that just as in the first step of the gluing construction in 
Subsection~\ref{reg1_subs2}, there is no obstruction to smoothing the internal
bubble nodes of the bubble map $b_r'$ subject to the second condition in~\e_ref{torusmap_e0a},
as long as $b_r'$ is sufficiently close to~$b_r$.
For defining the spaces $\Ga_+(\ups_r)$ at the second step of the gluing construction,
we take
$$\Ga_-(b')=\big\{\xi\!\in\!\Ga(b')\!: 
D_{J,b'}\xi\in\{\id\!\times\!u_{b'}\}^*\ti\Ga^{0,1}_-(b')\big\}.$$
The proof of the existence of the elements $\ups_r$, $\xi_{r;1}$, and $\xi_{r;2}$
as above is similar to the proof of the surjectivity property for 
the gluing map in Lemma~\ref{reg1_lmm3}, but simpler.\\

\noindent
Since for each $h\!\in\!I_1$ the map $\ti{u}_{\ups_{r;1}}$ is holomorphic on 
$\Si_{\ups_{r;h}}^0(\de)\!\subset\!\Si_{\ups_{r;1}}$ for $\de\!\in\!\R^+$ 
sufficiently small, 
the estimates of Corollary~\ref{reg0_crl4} apply to each map 
$\ti{u}_{\ups_{r;1}}|_{\Si_{\ups_{r;1},h}}$.
Thus, we can define an obstruction bundle $\Ga^{0,1}_-(\ups_r)$
for the second stage of the gluing construction in Subsection~\ref{reg1_subs2},
with the estimates of Lemma~\ref{reg1_lmm3} remaining valid.
The claim of Proposition~\ref{torus_prp} is then obtained by the same argument
as Corollary~\ref{reg1_crl3}, with $\ti{J}$, $u_{\ups}$, and $\ze_{\ti{J},\ups}$
replaced by $J$, $\ti{u}_{\ups_{r;1}}\!\circ\!\ti{q}_{\ups_{r;0};2}$,
and $\xi_{r;2}$, respectively.

\section{Proof of Theorem~\ref{str_thm}}
\label{strthm_sec}

\subsection{A Multi-Step Gluing Construction}
\label{g1gluing2_subs}

\noindent
The first part of the last claim of Theorem~\ref{str_thm} can be proved by showing
that a fine version of the converse to the $\ti{J}\!=\!J$ case of Corollary~\ref{reg1_crl3}
holds.
More precisely, using the two-step gluing construction of Subsection~\ref{reg1_subs2}
and the Inverse Function Theorem twice, we can construct 
an orientation-preserving diffeomorphism
$$\phi\!:{\cal F}^1{\cal T}^{\eset}_{\de} \lra
\M_{1,k}^0(X,A;J)\cap U_{\cal T}.$$
Unfortunately, one of the families of the domain spaces involved in this construction 
does not extend continuously over ${\cal F}^1\T_{\de}\!\cap\!\rho^{-1}(0)$
for a general bubble type $\T$ as in Theorem~\ref{str_thm}.
As these domain spaces are needed to apply IFT over 
${\cal F}^1\T_{\de}\!-\!{\cal F}^1{\cal T}^{\eset}$,
the above map $\phi$ cannot extend continuously over ${\cal F}^1{\cal T}_{\de}$,
except for bubble types~$\T$ such that either $|\chi(\T)|\!=\!1$ or $\chi(\T)\!=\!\hat{I}$.
In the first case, both families do extend continuously over~${\cal F}^1\T$.
In the second case, $\rho(\ups)\!=\!\ups$ for all $\ups\!\in\!{\cal FT}$
and both families extend continuously over ${\cal F}^1\T\!-\!\{0\}$.
On the other hand, as $\ups\!\lra\!0$ both perturbations approach zero.
This means that the corrections to be chosen in the domain spaces approach zero
as well and thus extend continuously over~${\cal F}^1\T_{\de}$.\\

\noindent
In this subsection, we describe a multi-step variation of the two-step gluing 
construction of Subsection~\ref{reg1_subs2}.
In the next subsection, we will use IFT multiple times to construct 
an orientation-preserving diffeomorphism
$$\phi\!:{\cal F}^1\T^{\eset}_{\de} \lra \M_{1,k}^0(X,A;J)\cap U_{\T}.$$
Some of the domain spaces involved will not extend continuously over ${\cal F}^1\T_{\de}$.
However, whenever a domain space cannot be extended to a point 
$\ups^*\!\in\!{\cal F}^1\T_{\de}$, the corresponding perturbations will approach zero 
as a sequence of elements $\ups_r\!\in\!{\cal F}^1\T_{\de}^{\eset}$ approaches~$\ups^*$.
For this reason, the above diffeomorphism~$\phi$ extends to a continuous map
$$\phi\!:{\cal F}^1\T_{\de} \lra \ov\M_{1,k}^0(X,A;J)\cap U_{\cal T}.$$
This map can shown to be a bijection by the same argument as in Subsection~4.1 of~\cite{Z4}.\\

\noindent
The multi-step gluing construction described in this subsection is suitable for 
the purposes of Subsection~\ref{reg1_subs3} and thus could have been described 
in Subsection~\ref{reg1_subs2} instead of the two-stage gluing construction.
However, describing the former in Subsection~\ref{reg1_subs3} would have 
further obscured the proofs of Corollaries~\ref{reg1_crl2} and~\ref{reg1_crl3}.
As these two corollaries appear far more central, than Theorem~\ref{str_thm},
to applications in the Gromov-Witten theory and enumerative geometry,
we have postponed describing the multi-step gluing construction until the present section.\\

\noindent
If $b\!=\!(\Si_b,u_b)$ is any genus-one bubble map such that 
$u_b|_{\Si_{b;P}}$ is constant, 
let $\Si_b^0\!\subset\!\Si_b$ be the maximum connected union of 
the irreducible components of $\Si_b$ such that $\Si_{b;P}\!\subset\!\Si_b^0$ 
and $u_b|_{\Si_b^0}$ is constant.
We~put
\begin{alignat*}{1}
\Ga_B(b) &=\big\{\xi\!\in\!\Ga(\Si_b;u_b^*TX)\!: \xi|_{\Si_b^0}\!=\!0\big\}
\qquad\hbox{and}\\
\Ga_B^{0,1}(b;J) &=\big\{\eta\!\in\!\Ga(\Si_b;\La^{0,1}_{J,j}T^*\Si_b\!\otimes\!u_b^*TX)\!:
\eta|_{\Si_b^0}\!=\!0\big\}.
\end{alignat*}
We denote by 
$$D_{J,b}^B\!:\Ga_B(b)\lra \Ga_B^{0,1}(b;J)$$
the restriction of the linearization $D_{J,b}$ of the $\bar\partial_J$-operator at~$b$
defined with respect to the connection~$\na^X$.
Let
$$\Ga_-(b)=\ker D_{J,b} \qquad\hbox{and}\qquad 
\Ga_{B;-}(b)=\ker D_{J,b}^B.$$
If $b$ is $J$-holomorphic, let $\ti\Ga_-(b)\!\subset\!\Ga_{B;-}(b)$ 
be the subspace defined in Subsection~\ref{reg1_subs2}; 
see \e_ref{reg1_lmm2e2a} and~\e_ref{reg1_lmm2e2c}.\\

\noindent
Suppose $\T\!=\!(M,I,\aleph;j,\under{A})$ is a bubble type as in Theorem~\ref{str_thm},
i.e.~$A_i\!=\!0$ for all $i\!\in\!I_0$, where $I_0\!\subset\!I$ is the subset of
minimal elements.
We~put
\begin{gather*}
\chi^0(\T)=\big\{h\!\in\!I\!: A_i\!=\!0~\forall~i\!\le\!h\big\},\qquad
\chi^-(\T)=\big\{h\!\in\!\hat{I}\!: h\!<\!i~\hbox{for some}~i\!\in\!\chi(\T)\big\}
\subset\chi^0(\T),\\
\lr{\T}=\max\big\{\big|\{h\!\in\!\hat{I}\!:h\!\le\!i\}\big|\!: i\!\in\!\chi(\T)\big\}\ge 1,
\qquad \cI_{\lr{\T}}^*=\chi(\T), \qquad \cI_{\lr{\T}}=\hat{I}-\chi(\T)-\chi^-(\T)-I_1,
\end{gather*}
where $I_1\!\subset\!I$ is as in Subsection~\ref{notation1_subs}.
For each $s\!\in\!\{0\}\!\cup\![\lr{\T}\!-\!1]$, let 
$$\cI_s=\big\{i\!\in\!\chi(\T)\!\cup\!\chi^-(\T)\!:
\big|\{h\!\in\!\hat{I}\!:h\!<\!i\}\big|\!=\!s\big\}, \qquad
\cI_s^*=\cI_s\cup\bigcup_{t=0}^{s-1}\big(\cI_t\!\cap\!\chi(\T)\big).$$
In the case of Figure~\ref{g1gendecomp_fig} on page~\pageref{g1gendecomp_fig}, 
$$\lr{\T}=2,\qquad \cI_0=\{h_1,h_3\}, \qquad \cI_1=\{h_4,h_5\}, \qquad \cI_2=\{h_2\}.$$
In general, the set $\cI_{\lr{\T}}$ could be empty, but the sets $\cI_s$ with
$s\!<\!\lr{\T}$ never are.\\

\noindent
If $b$ is a bubble map of type $\T$ as in Subsection~\ref{notation1_subs} and
$s\!\in\![\lr{\T}]$, we~put
$$\Si_b^{(s)}=\!\!\bigcup_{i\in\chi^0(\T)-\chi^-(\T)}\!\!\!\!\!\!\!\!\!\!\!\!\!\Si_{b,i}~~
\cup \bigcup_{h\in\cI_{s-1}^*}\bigcup_{i<h}\Si_{b,i} ~\subset \Si_b.$$
With notation as in Subsection~\ref{notation1_subs}, let
$$\wt{\F\T}=\! \bigoplus_{i\in\chi(\T)} \!\! \wt{{\cal F}_{h(i)}\T}
\lra\ti{\cal U}_{\T}^{(0)}(X;J).$$
If $s\!\in\![\lr{\T}]$ and $h\!\in\!\cI_{s-1}^*$, let
$$\chi_h(\T)=\big\{i\!\in\!\chi(\T)\!: h\!\le\!i\big\}, \qquad
\wt{\F_h\T}= \ti{\cal U}_{\T}^{(0)}(X;J) \times \C^{\chi_h(\T)}.$$
If in addition $\ups\!=\!(b,v)\!\in\!\wt{\cal FT}$, let
\begin{gather*}
\rho_{s;h}(\ups)=\big(b,(\rho_{h;i}(\ups))_{i\in\chi_h(\T)}\big)\in\wt{\F_h\T},
\qquad\hbox{where}\quad
\rho_{h;i}(\ups)=\!\!\prod_{h<h'\le i}\!\!\!\!\! v_{h'} \in\C;\\
\rho_s(\ups)=\big(b,(\rho_{s;h}(\ups))_{h\in\cI_{s-1}^*}\big)
\in \wt{\F_s\T}\!\equiv\!\!\bigoplus_{h\in\cI_{s-1}^*}\!\!\!\wt{\F_h\T}.
\end{gather*}
Note that $\rho(\ups)\!\in\!\wt{\F\T}$;  see Subsection~\ref{notation1_subs}.\\

\noindent
As in Subsection~\ref{reg1_subs2}, for each $\ups\!=\!(b,v)\!\in\!\wt{\cal FT}$ we put
$$\ups_0=(b,v_{\aleph},v_0).$$
Let $\ups_{\lr{0}}\!=\!\ups$.
If $s\!\in\![\lr{\T}]$, let
$$\ups_s=\big(b,(v_h)_{h\in\cI_s}\big) \qquad\hbox{and}\qquad
\ups_{\lr{s}}=\big(b,(v_h)_{h\in\cI_t,t\ge s}\big).$$
The component $\ups_{\lr{\T}}$ of $\ups$ consists of smoothings 
at the nodes of $\Si_b$ that do not lie between the principal component $\Si_{b;\aleph}$
of $\Si_b$ and the first-level effective bubbles and do not lie on~$\Si_{b;\aleph}$.
These nodes will be smoothed out at the first step of the gluing construction,
as specified by~$\ups_{\lr{\T}}$.
At the next step, we will smooth out the nodes indexed by the set $\cI_{\lr{\T}-1}$,
according the tuple of gluing parameters~$\ups_{\lr{\T}-1}$.
As in Subsection~\ref{reg1_subs2}, at the last step we will smooth out, if possible,
the nodes that lie on the principal component $\Si_{b;\aleph}$ of $\Si_b$
according to~$\ups_0$. This step will be obstructed.\\

\noindent
Suppose $\ups\!=\!(b,v)\!\in\wt{\cal FT}^{\eset}$ is a sufficiently small element.
We will inductively construct approximately $J$-holomorphic bubble maps 
$$b_s(\ups)=(\Si_{\ups_{\lr{s}}},u_{\ups,s}\big), \qquad
\forall~s=0,\ldots,\lr{\T},$$
$J$-holomorphic bubble maps 
$$\ti{b}_s(\ups)=(\Si_{\ups_{\lr{s}}},\ti{u}_{\ups,s}\big), \qquad
\forall~s=1,\ldots,\lr{\T},$$
and injective homomorphisms
$$R_{\ups,s}\!:\Ga_-(b) \lra\Ga\big(\Si_{\ups_{\lr{s}}};u_{\ups,s}^*TX\big)
\qquad\hbox{and}\qquad
\ti{R}_{\ups,s}\!:\Ga_-(b) \lra
                \Ga\big(\Si_{\ups_{\lr{s}}};\ti{u}_{\ups,s}^*TX\big),$$
such that the following properties are satisfied.
First, for all $s\!\in\![\lr{\T}]$,
$$\Si_{b_s(\ups)}^0=\Si_{\ti{b}_s(\ups)}^0=q_{\ups_{\lr{s}}}^{-1}\big(\Si_b^{(s)}\big)
\qquad\hbox{and}\qquad
u_{\ups,s}\big(\Si_{b_s(\ups)}^0\big)=\ti{u}_{\ups,s}\big(\Si_{\ti{b}_s(\ups)}^0\big)
=u_b(\Si_b^0) \!\equiv\!\ev_P(b),$$
where as before 
$$q_{\ups_{\lr{s}}}\!: \Si_{\ups_{\lr{s}}}\lra\Si_b$$
is the basic gluing map of Subsection~2.2 in~\cite{Z4}.
Second, for all $\xi\in\Ga_-(b)$
\begin{gather}\label{indminus_e0}
R_{\ups,s}\xi\big|_{\Si_{b_s(\ups)}^0}\!=\const,\quad
\ti{R}_{\ups,s}\xi\big|_{\Si_{\ti{b}_s(\ups)}^0}\!=\const,\quad
R_{\ups,s}\xi\big(\Si_{b_s(\ups)}^0\big)
=\ti{R}_{\ups,s}\xi\big(\Si_{\ti{b}_s(\ups)}^0\big)=\xi(\Si_b^0);\\
\label{indminus_e1}
C(b)^{-1}\|\xi\|_{b,p,1} \le \big\|R_{\ups,s}\xi\big\|_{\ups_{\lr{s}},p,1},
\big\|\ti{R}_{\ups,s}\xi\big\|_{\ups_{\lr{s}},p,1}   \le C(b)\|\xi\|_{b,p,1},\\
\label{indminus_e2}
\big\|D_{J,b(\ups_{\lr{s}})}R_{\ups,s}\xi\big\|_{\ups,p}, \,
\big\|D_{J,\ti{b}(\ups_{\lr{s}})}\ti{R}_{\ups,s}\xi\big\|_{\ups,p}
\le C(b)\big(|\ups|^{1/p}\!+\!|\ups|^{p-2/p}\big) \|\xi\|_{b,p,1}, 
\end{gather}
for some $C\!\in\!C({\cal U}_{\T}(X;J);\R^+)$.\\

\noindent
{\it Remark:} Similarly to Sections~\ref{gluing0_sec} and~\ref{gluing1_sec},
above and below $\|\cdot\|_{\ups_{\lr{s}},p,1}$ denotes the weighted $L^p_1$-norms 
on the spaces
$$\Ga_B\big(\Si_{\ups_{\lr{s}}};u_{\ups,s}^*TX\big)
\qquad\hbox{and}\qquad
\Ga_B\big(\Si_{\ups_{\lr{s}}};\ti{u}_{\ups,s}^*TX\big)$$
induced from the basic gluing map~$q_{\ups_{\lr{s}}}$ as in Subsection~3.3 of~\cite{Z4}.
Similarly, $\|\cdot\|_{\ups_{\lr{s}},p}$ denotes the weighted $L^p$-norms 
on the spaces
$$\Ga_B\big(\Si_{\ups_{\lr{s}}};
   \La^{0,1}_{J,j}T^*\Si_{\ups_{\lr{s}}}\!\otimes\!u_{\ups,s}^*TX\big)
\qquad\hbox{and}\qquad
\Ga_B\big(\Si_{\ups_{\lr{s}}};
   \La^{0,1}_{J,j}T^*\Si_{\ups_{\lr{s}}}\!\otimes\!\ti{u}_{\ups,s}^*TX\big).$$
We denote the corresponding completions by $\Ga_B(\ups_{\lr{s}})$, $\ti\Ga_B(\ups_{\lr{s}})$,
$\Ga_B^{0,1}(\ups_{\lr{s}};J)$, and~$\ti\Ga_B^{0,1}(\ups_{\lr{s}};J)$.\\

\noindent
For $s\!\in\!\{0\}\!\cup\![\lr{\T}]$, let $\Ga_-(\ups_{\lr{s}})$ be the image
of $R_{\ups,s}$. Similarly, if $s\!\in\![\lr{\T}]$, we denote by
$$\Ga_{B;-}(\ups_{\lr{s}})\subset\Ga_B\big(b_s(\ups)\big)
\qquad\hbox{and}\qquad
\ti\Ga_-(\ups_{\lr{s}}),\ti{\Ga}_{B;-}(\ups_{\lr{s}})\subset\Ga_B\big(\ti{b}_s(\ups)\big)$$
the image of $\Ga_{B;-}(b)$ under $R_{\ups,s}$, the image of $\ti{R}_{\ups,s}$, 
the image of $\Ga_{B;-}(b)$ under $\ti{R}_{\ups,s}$, respectively; see \e_ref{indminus_e0}.
Let $\Ga_{B;+}(\ups_{\lr{s}})$ and $\ti{\Ga}_{B;+}(\ups_{\lr{s}})$ be 
the $L^2$-orthogonal complements of 
$\Ga_{B;-}(\ups_{\lr{s}})$ and $\ti\Ga_{B;-}(\ups_{\lr{s}})$
in $\Ga_B(\ups_{\lr{s}})$ and~$\ti\Ga_B(\ups_{\lr{s}})$.
These spaces will satisfy
\begin{gather}\label{indplus_e1}
C(b)^{-1}\|\xi\|_{\ups_{\lr{s}},p,1} \le 
\big\|D_{J,b(\ups_{\lr{s}})}\xi\big\|_{\ups_{\lr{s}},p} \le 
C(b)\|\xi\|_{\ups_{\lr{s}},p,1}    \qquad\forall~\xi\in\Ga_{B;+}(\ups_{\lr{s}});\\
\label{indplus_e2}
C(b)^{-1}\|\xi\|_{\ups_{\lr{s}},p,1} \le 
\big\|D_{J,\ti{b}(\ups_{\lr{s}})}\xi\big\|_{\ups_{\lr{s}},p} \le 
C(b)\|\xi\|_{\ups_{\lr{s}},p,1}    \qquad\forall~\xi\in\ti\Ga_{B;+}(\ups_{\lr{s}}).
\end{gather}
Furthermore, 
\begin{equation}\label{indpert_e1}
\ti{u}_{\ups,s}=\exp_{u_{\ups,s}}\!\ze_{\ups,s}
\quad\hbox{for some}\quad \ze_{\ups,s}\in\Ga_{B;+}(\ups_{\lr{s}})
~~\hbox{s.t.}~~
\big\|\ze_{\ups,s}\big\|_{\ups_{\lr{s}},p,1} \le  C(b)|\ups|^{1/p}.
\end{equation}
Finally, for $\de\!\in\!C({\cal U}_{\T}(X;J);\R^+)$ sufficiently small, all maps
$$\ups\lra b_s(\ups),\, \ti{b}_s(\ups),\, \ze_{\ups,s},\, R_{\ups,\ze},\, \ti{R}_{\ups,\ze}$$
are smooth on $\wt{\cal FT}_{\de}^{\eset}$ and extend continuously 
over~$\wt{\cal FT}_{\de}$.\\

\noindent
We now describe the inductive construction referred to above.
If $\ups\!\in\!\wt{\cal FT}^{\eset}$ is as above and $b\!=\!(\Si_b,u_b)$, we put
$$u_{\ups,\lr{\T}}=u_b\circ q_{\ups_{\lr{\T}}},
\qquad R_{\ups,\lr{\T}}\xi=\xi\!\circ\!q_{\ups_{\lr{\T}}}
\quad\forall\,\xi\!\in\!\Ga_-(b).$$
The first bounds in~\e_ref{indminus_e2} and~\e_ref{indplus_e1} with $s\!=\!\lr{\T}$
hold for the same reasons as the corresponding estimates in Lemma~\ref{reg0_lmm1}.
Since the operator~$D_{J,b}^B$ is surjective, by the first bound in~\e_ref{indminus_e2}
the operator
$$D_{J,b_{\lr{\T}}(\ups)}^B\!: \Ga_{B;+}\big(\ups_{\lr{\T}}\big)
\lra  \Ga_B^{0,1}\big(\ups_{\lr{\T}};J\big)$$
is an isomorphism.
On the other hand,
by the construction of the map~$q_{\ups,\lr{\T}}$ in Subsection~2.2 in~\cite{Z4},
\begin{equation}\label{indpert_e2}
\big\|\bar\partial_Ju_{\ups,\lr{\T}}\big\|_{\ups_{\lr{\T}},p}\le C(b)|\ups|^{1/p}.
\end{equation}
Thus, by the Contraction Principle, if $\ups$ is sufficiently small,
there exists a unique small element
\begin{equation}\label{indpert_e2a}
\ze_{\ups,\lr{\T}}\in\Ga_{B;+}\big(\ups_{\lr{\T}}\big)
\qquad\hbox{s.t.}\qquad
\bar\partial_J \exp_{u_{\ups,\lr{\T}}}\!\ze_{\ups,\lr{\T}}=0.
\end{equation}
Furthermore, by~\e_ref{indpert_e2},
$$\big\|\ze_{\ups,\lr{\T}}\big\|_{\ups_{\lr{\T}},p,1} \le  C(b)|\ups|^{1/p}.$$
We thus define $\ti{b}_{\lr{\T}}(\ups)$ by the first equation in~\e_ref{indpert_e1}.\\

\noindent
If $s\!\in\![\lr{\T}\!-\!1]$, let
$$q_{\ups_s;\lr{\T}+1-s}\!: \Si_{\ups_{\lr{s}}} \lra \Si_{\ups_{\lr{s+1}}}$$
be the basic gluing map of Subsection~2.2 in~\cite{Z4} corresponding 
to the gluing parameter~$\ups_s$.
If $\ti{b}_{s+1}(\ups)$ and $\ti{R}_{\ups,s+1}$ have been defined,  we put
$$u_{\ups,s}=\ti{u}_{\ups,s+1}\circ q_{\ups_s;\lr{\T}+1-s},
\qquad 
R_{\ups,s}\xi=\ti{R}_{\ups,s+1}\xi\circ q_{\ups_s;\lr{\T}+1-s}
\quad\forall\, \xi\!\in\!\Ga_-(b).$$
The first bounds in~\e_ref{indminus_e2} and~\e_ref{indplus_e1} follow
from the second estimates in~\e_ref{indminus_e1} and~\e_ref{indminus_e2}
and from~\e_ref{indplus_e2}, with $s$ replaced by~$s\!+\!1$.
On the other hand, by the inductive construction and~\e_ref{indpert_e1},
\begin{equation}\label{indpert_e2b}\begin{split}
&\qquad \ti{u}_{\ups,s+1}=\exp_{u_{\ups_{\lr{s+1}}}}\!\ti\ze_{\ups,s+1}\\
\hbox{for some}\qquad 
&\ti\ze_{\ups,s+1}\in\Ga_B\big(\Si_{\ups_{\lr{s+1}}};u_{\ups_{\lr{s+1}}}^*TX\big)
~~\hbox{s.t.}~~
\big\|\ti\ze_{\ups,s+1}\big\|_{\ups_{\lr{s+1}},p,1} \le  C(b)|\ups|^{1/p},
\end{split}\end{equation}
where
$$u_{\ups_{\lr{s+1}}}=u_b\circ q_{\ups_{\lr{s+1}}}.$$
Thus, if $\de$ is sufficiently small, the estimate in~(2b) of Corollary~\ref{reg0_crl4} 
implies
\begin{gather}\label{indpert_e3}
\big\|d\ti{u}_{\ups,s+1}|_{A_{\ups_{\lr{s+1}},h}^-(\de)}\big\|_{\ups,p}
\le C(b)\de^{1/p}\big|\rho_{s;h}(\ups)\big| \qquad\forall\,h\in\cI_s^*,
\qquad\hbox{where}\\
A_{\ups_{\lr{s+1}},h}^-(\de)=q_{\ups_{\lr{s+1}}}^{-1}\big(
\big\{(h,z)\!\in\!\Si_{b,h}\!=\!\{h\}\!\times\!S^2\!: 
|z|\!\ge\!\de^{-1/2}/2\big\}\big). \notag
\end{gather}
It follows that
\begin{equation}\label{indpert_e4}
\big\|\bar\partial_Ju_{\ups,s}\big\|_{\ups,p}\le 
C(b)|\ups|^{1/p}\big|\rho_{s+1}(\ups)\big| \le C(b)|\ups|^{1/p}.
\end{equation}
Thus, similarly to the $s\!=\!\lr{\T}$ case above, if $\ups$ is sufficiently small,
there exists a unique small element
\begin{equation}\label{indpert_e2c}
\ze_{\ups,s}\in\Ga_{B;+}\big(\ups_{\lr{s}}\big)
\qquad\hbox{s.t.}\qquad
\bar\partial_J \exp_{u_{\ups,s}}\!\ze_{\ups,s}=0.
\end{equation}
Furthermore, by~\e_ref{indpert_e4},
$$\big\|\ze_{\ups,s}\big\|_{\ups_s,p,1} \le  C(b)|\ups|^{1/p}.$$
We again define $\ti{b}_s(\ups)$ by the first equation in~\e_ref{indpert_e1}.\\

\noindent
If $s\!\in\![\lr{\T}]$ and $\ti{b}_s(\ups)$ has been defined 
via~\e_ref{indpert_e1}, we~put
$$\ti{R}_{\ups,s}\xi= \Pi_{\ze_{\ups,s}}R_{\ups,s}\xi
\qquad\forall\,\xi\!\in\!\Ga_-(b),$$
where $\Pi_{\ze_{\ups,s}}$ is the parallel transport along the geodesics
$$\tau\lra\exp_{u_{\ups,s}}\!\tau\ze_{\ups,s}, \qquad \tau\in[0,1].$$
The bounds on $R_{\ups,s}\xi$ in~\e_ref{indminus_e0}-\e_ref{indminus_e2}
and the estimate~\e_ref{indplus_e1}, along with~\e_ref{indpert_e1},
imply the bounds on $\ti{R}_{\ups,s}\xi$ in~\e_ref{indminus_e0}-\e_ref{indminus_e2}
and the estimate~\e_ref{indplus_e2}.\\

\noindent
At the final step of this inductive construction, we put
$$u_{\ups,0}=\ti{u}_{\ups,1}\circ \ti{q}_{\ups_0;\lr{\T}+1},
\qquad 
R_{\ups,0}\xi\equiv \ti{R}_{\ups,1}\xi\circ\ti{q}_{\ups_0;\lr{\T}+1}
\quad\forall\,\xi\!\in\!\Ga_-(b)\big\},$$
where 
$$\ti{q}_{\ups_0;\lr{\T}+1}\!: \Si_{\ups}\lra\Si_{\ups_{\lr{1}}}$$
is the modified basic gluing map constructed in Subsection~\ref{reg1_subs2}
as~$\ti{q}_{\ups_0;2}$.
In order to construct this map in this case, 
we need to replace $\de_K\!\in\!\R^+$ with $\de\!\in\!C^{\i}({\cal U}_{\T}(X;J)$,
which we  view as a function on $\ti{\cal U}_{\T}^{(0)}(X;J)$ 
via the quotient projection map
$$\ti{\cal U}_{\T}^{(0)}(X;J) \lra {\cal U}_{\T}^{(0)}(X;J).$$
The homomorphism~$R_{\ups,0}$ satisfies the required properties.
Let $\Ga_+(\ups)\!\subset\!\Ga(\ups)$ be the $L^2$-orthogonal complement of~$\Ga_-(\ups)$.\\

\noindent
For each $h\!\in\!\hat{I}$ and $\de\!\in\!\bar\R^+$, let
\begin{gather*}
A_{b,h}^+(\de)=
\big\{(\io_h,z)\!\in\!\Si_{b,\io_h}\!: |z\!-\!x_h(b)|\!\le\!2\de^{1/2}\big\},\\
A_{b,h}^-(\de)= \big\{(h,z)\!\in\!\Si_{b,h}\!=\!\{h\}\!\times\!S^2\!: 
|z|\!\ge\!\de^{-1/2}/2\big\},\\
\Si_{b,h}^*(\de)=\Si_{b,h}-A_{b,h}^-(\de)-\bigcup_{\io_{h'}=h}\!\!A_{b,h'}^+(\de),\\
\ti{A}_{\ups,h}^{\pm}(\de)=q_{\ups}^{-1}\big(A_{b,h}^{\pm}(\de)\big)\subset\Si_{\ups}, 
\qquad
\ti\Si_{\ups,h}^*(\de)=q_{\ups}^{-1}\big(\Si_{b,h}^*(\de)\big).
\end{gather*}
We define the homomorphism 
$$R_{\ups}\!: \Ga^{0,1}_-(b;J)\lra \Ga^{0,1}(\ups;J)$$
similarly to Subsection~\ref{reg1_subs2}, but with two changes.
First, we replace the number $\de_K$ with the function 
$\de\!\in\!C^{\i}({\cal U}_{\T}(X;J);\R^+)$.
Second, we cut-off $R_{\ups}^{\aleph}\eta$ over the annuli 
$$\ti{A}_{\ups,i}^+(4\de(b))-\ti{A}_{\ups,i}^+(\de(b))$$ 
with $i\!\in\!\chi(\T)$, instead of $h\!\in\!I_1$;
see Subsection~2.2 in~\cite{Z3} for a version of this construction.
Let $ \Ga_-^{0,1}(\ups;J)$ be the image of~$R_{\ups}$.
We note that due to~\e_ref{indpert_e2b}, the estimates in~(3) and (6) 
of Lemma~\ref{reg1_lmm3} remain valid.
Of course, in this case $\ti{J}\!=\!J$ and  $C_K\!\in\!\R^+$ should be replaced
by an $\R^+$-valued continuous function on~${\cal U}_{\T}(X;J)$.
We summarize the key results of this construction below. 

\begin{lmm}
\label{str_lmm}
Suppose $(X,\om,J)$ is a compact almost Kahler manifold, $A\!\in\!H_2(X;\Bbb{Z})$, and
$J$ is a genus-zero $A$-regular almost complex structure.
If $\T\!=\!(M,I,\aleph;j,\under{A})$ is a bubble type such that
$\sum_{i\in I}\!A_i\!=\!A$ and $A_i\!=\!0$ for all minimal elements $i$ of~$I$,
there exist $\de,C\!\in\!C({\cal U}_{\cal T}(X;J);\Bbb{R}^+)$
and an open neighborhood $U_{\T}$ of ${\cal U}_{\T}(X;J)$ in $\X_{1,M}(X,A)$
such that\\
(1)  $b_s(\ups)$, $\ti{b}_s(\ups)$, $R_{\ups,s}$, and $\ti{R}_{\ups,s}$
as above are defined for all $\ups\!\in\!\wt{\cal FT}_{\de}^{\eset}$;\\
(2) for every $[\tilde{b}]\!\in\!\X_{1,M}^0(X,A)\cap U_K$, 
there exist $\ups\!=\!(b,v)\in\!\ti{\cal F}{\cal T}_{\de}$ and 
$\ze_{\ups,0}\!\in\!\Ga_+(\ups)$ such~that 
$$\|\ze_{\ups,0}\|_{\ups,p,1}<\de(b) \qquad\hbox{and}\qquad 
\big[\exp_{b_0(\ups)}\!\ze_{\ups,0}\big]=[\tilde{b}].$$
Furthermore, such a pair $(\ups,\ze_{\ups,0})$ is unique up to 
the $\Aut(\T)\!\propto\!(S^1)^{\hat{I}}$-action;\\
(3) for all $\ups\!=\!(b,v)\!\in\!\wt{\cal FT}^{\eset}_{\de}$,
\begin{gather*}
\|\bar{\partial}_Ju_{\ups,0}\|_{\ups,p}\le  C(b)\big|\rho(\ups)\big|,
\qquad \|D_{J,b_0(\ups)}\xi\|_{\ups,p}\le 
C(b)\big(|\ups|^{1/p}\!+\!|\ups|^{\frac{p-2}{p}}\big)\|\xi\|_{\ups,p,1} 
\quad\forall\xi\!\in\!\Ga_-(\ups),\\
\hbox{and}\qquad
C(b)^{-1}\|\xi\|_{\ups,p,1}\le \|D_{J,b_0(\ups)}\xi\|_{\ups,p} \le C(b)\|\xi\|_{\ups,p,1}
\quad\forall\xi\!\in\!\Ga_+(\ups);
\end{gather*}
(4) for all $\ups\!=\!(b,v)\!\in\!\wt{\cal FT}^{\eset}_{\de}$,
$\xi\!\in\!\Ga(\ups)$, and $\eta\!\in\!\Ga_-^{0,1}(b;J)$, 
$$\big|\llrr{D_{J,b_0(\ups)}\xi,R_{\ups}\eta}_{\ups,2}\big|
\le C(b) \big|\rho(\ups)\big|^{1/2}\|\eta\|\|\xi\|_{\ups,p,1};$$\\
(5) for all $\ups\!=\!(b,v)\!\in\!\wt{\cal FT}^{\eset}_{\de}$, 
$s\!\in\![\lr{\T}]$, and $h\!\in\!\cI_{s-1}^*$,
$$\Big|\cD_h\ti{b}_s(\ups)-
\!\!\sum_{i\in\chi_h(\T)}\!\!\!\!\cD_{J,i}\rho_{h;i}(\ups)\Big|
\le C(b)|\ups|^{1/p}|\rho_{s;h}(\ups)\big|;$$
(6) for all $\ups\!=\!(b,v)\!\in\!\wt{\cal FT}^{\eset}_{\de}$
and $\eta\!\in\!\Ga_-^{0,1}(b;J)$,  
$$\Big|\bllrr{\bar{\partial}_Ju_{\ups,0},R_{\ups}\eta}+2\pi\I
\!\!\sum_{i\in\chi({\cal T})}\!\!\!\!
\blr{{\cal D}_i\rho_i(\ups),\eta_{x_{h(i)}(b)}}_b\Big|
\le C(b)\big(|\ups|^{1/p}\!+\!|\ups|^{(p-2)/p}\big)
|\rho(\ups)\big|\cdot\|\eta\|.$$
(7) all maps
$$\ups\lra b_s(\ups),\, \ti{b}_s(\ups),\, \ze_{\ups,s},\, R_{\ups,s},\, \ti{R}_{\ups,s}$$
are smooth on $\wt{\cal FT}_{\de}^{\eset}$ and extend continuously 
over~$\wt{\cal FT}_{\de}$.
\end{lmm}

\noindent
Due to~\e_ref{indpert_e2b}, (5) is proved by the same argument as 
the $r\!=\!1$ case of the expansion~(2a) in Lemma~\ref{reg0_lmm4}.
Part~(2) of Lemma~\ref{str_lmm} holds for the same reason as 
part~(2) of Lemma~\ref{reg1_lmm3}.
Regarding part~(7) of Lemma~\ref{str_lmm}, it is immediate from the inductive construction
that each of the maps
$$\ups\lra b_s(\ups),\, \ti{b}_s(\ups),\, \, R_{\ups,s},\, \ti{R}_{\ups,s}$$
is smooth on $\wt{\cal FT}_{\de}^{\eset}$ and extends continuously 
over~$\wt{\cal FT}_{\de}$, provided this is the case for all the objects 
defined at the preceding steps of the construction.
Under these circumstances the map $\ups\!\lra\!\ze_{\ups,s}$ is also smooth,
by the smooth dependence of the solutions of~\e_ref{indpert_e2a}
and~\e_ref{indpert_e2c} on the parameters.
It extends continuously over~$\wt{\cal FT}_{\de}$ by the same argument
as in Subsection~4.1 in~\cite{Z4}.\\

\noindent
The estimate in (4) of Lemma~\ref{str_lmm} is an improvement on (7) of 
Lemma~\ref{reg1_lmm3} and is proved by a similar argument.
In this case, the support of~$D_{J,\ups}^*R_{\ups}\eta$ is contained in the union~of
the~annuli
\begin{gather*}
\ti{A}_{\ups,h}^+\big(\de(b)\big),~~~
h\!\in\!\chi^-(\T)\!\cup\!\chi(\T); \qquad
\ti{A}_{\ups,h}^-\big(\de(b)\big),~~~
h\!\in\!\chi^-(\T)\!\cup\!\chi(\T); \\
\ti\Si_{\ups,h}^*\big(\de(b)\big),~~~h\!\in\!\chi^-(\T); \qquad
\ti{A}_{\ups,h}^-\big(4\de(b)\big)\!-\!\ti{A}_{\ups,h}^-\big(\de(b)\big),
~~~h\!\in\!\chi(\T).
\end{gather*}
Similarly to the proof of (7) of Lemma~\ref{reg1_lmm3},
\begin{equation}\label{strlmm_e1}
\big|D_{J,\ups}^*R_{\ups}\eta\big|_{g_{\ups},z}
\le C_X\big|du_{\ups,0}\big|_{g_{\ups},z}|\eta|_{g_{\ups,z}}
\end{equation}
for every point $z$ of any of the annuli of the first three types above.
Thus, the estimate~\e_ref{dadj_e4} still applies to the annuli of 
the first type with $h\!\in\!I_1$.
By definition of the metric~$g_{\ups}$,
\begin{alignat}{3}\label{strlmm_e3a}
&|\eta|_{g_{\ups,z}} \le  C(b)\|\eta\|\cdot \!\!\prod_{h'\in\hat{I},h'<h}
\!\!\!\!\!\!\! v_{h'} &\qquad  &\forall~~ 
z\!\in\!\ti{A}_{\ups,h}^+\big(\de(b)\big),& ~~& h\!\in\!\chi^-(\T)\!\cup\!\chi(\T);\\
\label{strlmm_e3b}
&|\eta|_{g_{\ups,z}} \le  C(b)\|\eta\|\cdot \!\!\prod_{h'\in\hat{I},h'\le h}
\!\!\!\!\!\!\! v_{h'} &\qquad  &\forall~~ 
z\!\in\!\ti\Si_{\ups,h}^*\big(\de(b)\big), &~~& h\!\in\!\chi^-(\T);\\
\label{strlmm_e3c}
&|\eta|_{g_{\ups,z}} \le C(b)\|\eta\|\cdot|w_h|^{-1}\!\!\prod_{h'\in\hat{I},h'\le h}
\!\!\!\!\!\!\! v_{h'} &\qquad  &\forall~~ 
z\!\in\!\ti{A}_{\ups,h}^-\big(4\de(b)\big), &~~& h\!\in\!\chi^-(\T)\!\cup\!\chi(\T),
\end{alignat}
where $w_h$ is the coordinate on $\ti{A}_{\ups,h}^-(4\de(b))$ defined similarly
to $w_i$ in~(2) of the proof of Lemma~\ref{reg0_lmm4}.
On the other hand, by~\e_ref{indpert_e3} and the assumption~(a) of 
Definition~\ref{g1reg_dfn},
\begin{equation}\label{strlmm_e5}
\big\|\ze_{\ups,s}|_{\ti\Si_{\ups,i}^*(0)}\big\|_{\ups_{\lr{s}},p,1}
\le C(b)\big|\rho_{s;h}(\ups)\big| \qquad\forall~ h\!\in\!\cI_{s-1}^*,~i\!\ge\!h;
\end{equation}
the above assumption implies that the operators $D_{J,b_s(\ups)}^h$ defined in
Subsection~\ref{strthmcont_subs} below are surjective.
Note that by the inductive construction,
$$\ze_{\ups,s}|_{\ti\Si_{\ups,i}^*(0)}=0 \qquad
\forall~ h\!\in\!\cI_{s-1}^*,~i\!<\!h.$$
Combining this observation with \e_ref{strlmm_e5}, we find that
$$\big\|\ti\ze_{\ups,1}|_{\ti\Si_{\ups,h}^*(0)}\big\|_{\ups_{\lr{1}},p,1}
\le C(b)\big|\rho_{s;h}(\ups)\big| \qquad\forall~ h\!\in\!\hat{I},$$
where $\ti\ze_{\ups,1}$ is as in~\e_ref{indpert_e2b}.
Thus,
\begin{alignat}{2}\label{strlmm_e7a}
&\big\|du_{\ups,0}|_{\ti{A}_{\ups,h}^+(\de(b))}\big\|_{\ups,p,1}
\le  C(b)\prod_{h'\in\chi^-(\T)\cup\chi(\T),h'\ge h}
\!\!\!\!\!\!\!\!\!\!\!\!\!\!\!\! v_{h'} &
\qquad\qquad  &\forall~~
 h\in\chi^-(\T)\!\cup\!\chi(\T)-I_1;\\
\label{strlmm_e7b}
&\big\|du_{\ups,0}|_{\ti\Si_{\ups,h}^*(\de(b))}\big\|_{\ups,p,1}
\le  C(b)\prod_{h'\in\chi^-(\T)\cup\chi(\T),h'>h}
\!\!\!\!\!\!\!\!\!\!\!\!\!\!\!\! v_{h'} &
\qquad\qquad  &\forall~~ h\in\chi^-(\T);\\
\label{strlmm_e7c}
&\big\|du_{\ups,0}|_{\ti{A}_{\ups,h}^-(4\de(b))}\big\|_{\ups,p,1}
\le  C(b)\prod_{h'\in\chi^-(\T)\cup\chi(\T),h'>h}
\!\!\!\!\!\!\!\!\!\!\!\!\!\!\!\! v_{h'} &
\qquad\qquad  &\forall~~ h\in\chi^-(\T)\!\cup\!\chi(\T).
\end{alignat}
Combining \e_ref{strlmm_e1}-\e_ref{strlmm_e3c}, 
\e_ref{strlmm_e7a}-\e_ref{strlmm_e7c}, and Holder's inequality, 
we find that the $L^1$-norm of $D_{J,\ups}^*R_{\ups}\eta$,
with respect to the metric~$g_{\ups}$, on each of the annuli above 
is bounded by $C(b)|\rho(\ups)|\|\eta\|$.
Finally, analogously to~\e_ref{dadj_e1c}, 
$$\big|D_{J,\ups}^*R_{\ups}\eta\big|_{g_{\ups},z}
\le C(b)\big(1\!+\!\big|du_{\ups,0}\big|_{g_{\ups},z}\big)
|\eta|_{g_{\ups,z}} 
\qquad\forall~z\!\in\!\ti{A}_{\ups,h}^-\big(4\de(b)\big)\!-\!\ti{A}_{\ups,h}^-\big(\de(b)\big),
~~h\!\in\!\chi(\T).$$
Thus, by \e_ref{strlmm_e3c} and~\e_ref{strlmm_e7c} the $L^1$-norm of $D_{J,\ups}^*R_{\ups}\eta$ 
on such annuli is also bounded by $C(b)|\rho(\ups)|\|\eta\|$.\\

\noindent
{\it Remark:} The exponent $1/2$ in (4) of Lemma~\ref{str_lmm} is due to the exponent
$(p\!-\!1)/p$ in~\e_ref{dadj_e4}.

\subsection{Construction of Diffeomorphism}
\label{strthm_subs}

\noindent
We continue with the notation of the previous subsection.
For each $\ups\!=\!(b,v)\!\in\!\wt{\cal FT}_{\de}^{\eset}$,
we define the homomorphism 
$$\pi^{0,1}_{\ups;-}\!: \Ga^{0,1}(\ups;J)\lra\Ga^{0,1}_-(b;J)
\qquad\hbox{by}\qquad 
\pi^{0,1}_{\ups;-}\eta=\sum_r\lr{\eta,R_{\ups}\eta_r}\eta_r\in\Ga^{0,1}_-(b;J),$$
where $\{\eta_r\}$ is an orthonormal basis for $\Ga_-^{0,1}(b;J)$
as in (5) of the proof of Corollary~\ref{reg1_crl2}.
We denote the kernel of $\pi^{0,1}_{\ups;-}$ by $\Ga^{0,1}_+(\ups;J)$.
By Lemma~\ref{str_lmm}, 
\begin{equation*}\begin{split}
U_{\T}'
&=\big\{[\exp_{b_0(\ups)}\!\ze]\!:
\ups\!=\!(b,v)\!\in\!\ti{\cal F}\T_{\de},~
\ze\!\in\!\Ga_+(\ups),~\|\ze\|_{\ups,p,1}\!<\!\de(b)\big\}\\
&\subset\big\{[\exp_{b_0(\ups)}\!\ze]\!:
\ups\!=\!(b,v)\!\in\!\wt{\cal FT}_{\de},~
\ze\!\in\!\Ga_+(\ups),~\|\ze\|_{\ups,p,1}\!<\!\de(b)\big\}
\end{split}\end{equation*}
is an open neighborhood of ${\cal U}_{\T}(X;J)$ in $\X_{1,M}(X,A)$.
Thus, we need to solve the equation
\begin{equation}\label{str1_crl1e1}
\bar{\partial}_J\exp_{u_{\ups,0}}\!\ze\!=\!0  ~\Llra~
\begin{cases}
\pi^{0,1}_{\ups;-}
\big(\bar{\partial}_Ju_{\ups,0}\!+\!D_{J,b_0(\ups)}\ze\!+\!N_{\ups,0}\ze\big)\!=\!0
\in\Ga^{0,1}_-(b;J),\\
\bar{\partial}_Ju_{\ups,0}\!+\!D_{J,b_0(\ups)}\ze\!+\!N_{\ups,0}\ze\!=\!0
\in\Ga^{0,1}_+(\ups;J),
\end{cases}
\end{equation}
where $N_{\ups,0}$ is the quadratic term satisfying \e_ref{reg1_crl2e2},
for $\ups\!=\!(b,v)\!\in\!\wt{\cal FT}_{\de}^{\eset}$ and
$\ze\!\in\!\Ga_+(\ups)$ such that $\|\ze\|_{\ups,p,1}\!<\!\de(b)$.
By the proof of~(1) of Corollary~\ref{reg1_crl3},
there exist $\ti\de,\ti{C}\!\in\!({\cal U}_{\T}(X;J);\R^+)$ such that
$\ti\de\!<\!\de$ and every solution $(\ups,\ze)$ of~\e_ref{str1_crl1e1}
\begin{equation}\label{str1_crl1e2}
|\ups|<\ti{\de}(b),~~ \|\ze\|_{\ups,p,1}<\ti{\de}(b)
\qquad\Lra\qquad
\|\ze\|_{\ups,p,1}<\ti{C}(b)\cdot\big|\rho(\ups)\big|.
\end{equation}
On the other hand, by (4) and (6) of Lemma~\ref{str_lmm} and \e_ref{reg1_crl2e2},
\begin{equation}\label{str1_crl1e3a}
\pi^{0,1}_-b_0(\ups,\ze)\equiv \pi^{0,1}_{\ups;-}
\big(\bar{\partial}_Ju_{\ups,0}\!+\!D_{J,b_0(\ups)}\ze\!+\!N_{\ups,0}\ze\big)=
\cD_{\T}(\ups)+\ve(\ups,\ze),
\end{equation}
where $\cD_{\T}$ is as defined in Subsection~\ref{notation1_subs} and
\begin{gather}\label{str1_crl1e3c}
\big\|\ve(\ups,\ze)\big\|
\le C(b)(\big(|\ups|^{1/p}\!+\!|\ups|^{(p-2)/p}\!+\!\|\ze\|_{\ups,p,1}\big)\cdot
\big(|\rho(\ups)|\!+\!\|\ze\|_{\ups,p,1}\big)\\
\hbox{if}\quad \ze\!\in\!\Ga(\ups),~~\|\ze\|_{\ups,p,1}\!\le\!\de(b). \notag
\end{gather}
We will first solve the top equation in~\e_ref{str1_crl1e1} for $b'\!=\!\mu_0(\ups,\ze)$
and then use the Contraction Principle to show that the resulting bottom equation
has a unique solution in~$\ze$ for each $\ups\!\in\!\ti{\cal F}^1\T_{\ti{\de}}^{\eset}$.\\

\noindent
For each $s\!\in\![\lr{\T}\!+\!1]$ and $h\!\in\!I_1$, let
\begin{equation*}\begin{split}
{\cal U}_{\T_h}^{(s)}(X;J)=\Psi_{\T}^{-1}\big(
\big\{(0,r_i)_{i\in I_h}\!\in\!(\C\!\times\!\R)^{I_h}\!:
~& r_i\!=\!\frac{1}{2}~\forall i\!\in\!I\!-\!\chi({\T}),\\
&r_i\!\in\!\Big(\frac{1}{2}\!-\!\frac{s}{4(\lr{\T}\!+\!1)},
\frac{1}{2}\!+\!\frac{s}{4(\lr{\T}\!+\!1)}\Big) ~\forall i\!\in\!\chi_h({\T})\big\};
\end{split}\end{equation*}
see the end of Subsection~\ref{notation0_subs}.
We put
\begin{equation*}\begin{split}
{\cal U}_{\T}^{(s)}(X;J)=
\big\{\big(b_0,(b_h)_{h\in I_1}\big)\!\in\!{\cal U}_{{\T}_0}(X;J)\!\times\!\!
\prod_{h\in I_1}\!\!{\cal U}_{\T_h}^{(s)}(X;J)\!: \qquad\qquad&\\
\ev_0(b_h)\!=\!\ev_{\io_h}(b_0)~\forall h\!\in\!I_1\big\}&
\subset \ti{\cal U}_{\T}^{(0)}(X;J);
\end{split}\end{equation*}
see Subsection~\ref{notation1_subs} for notation.\\

\noindent
For each $s\!\in\![\lr{\T}]$ and $h\!\in\!\cI_{s-1}^*$, let
$$\pi_{s;h}\!: \bP\wt{\F_h\T}\big|_{{\cal U}_{\T}^{(s+1)}(X;J)} \lra 
{\cal U}_{\T}^{(s+1)}(X;J)
\qquad\hbox{and}\qquad
\ga_{s;h}\lra  \bP\wt{\F_h\T}\big|_{{\cal U}_{\T}^{(s+1)}(X;J)}$$
be the natural projection map and the tautological line bundle.
With
$$V_{s;h}=\pi_{s;h}^*\ev_P^*TX \lra \bP\wt{\F_h\T}\big|_{{\cal U}_{\T}^{(s+1)}(X;J)},$$
we define
\begin{gather*}
\al_{s;h}\in\Ga\big(\bP\wt{\F_h\T}\big|_{{\cal U}_{\T}^{(s+1)}(X;J)};
\ga_{s;h}^*\!\otimes\!V_{s;h}\big) \qquad\hbox{by}\\
\al_{s;h}\big(b,(v_i)_{i\in\chi_h(\T)}\big)
=\!\sum_{i\in\chi_h(\T)}\!\!\!\!\cD_{J,i}\big(b,v_i) \in T_{\ev_P(b)}X
\qquad\hbox{if}\qquad
\big(b,(v_i)_{i\in\chi_h(\T)}\big)\!\in\!\ga_{s;h}.
\end{gather*}
We denote by
$$\pi_s\!: \bP_s\T \equiv \prod_{h\in\cI_{s-1}^*}\!\!\!\big(\bP\wt{\F_h\T},\pi_{s;h}\big)
\lra {\cal U}_{\T}^{(s+1)}(X;J)$$
the fiber product of the bundles $\bP\wt{\F_h\T}$ over ${\cal U}_{\T}^{(s+1)}(X;J)$.
Let 
$$V_s=\pi_s^*\ev_P^*TX\lra\bP_s\T \qquad\hbox{and}\qquad
\ga_s=\bigoplus_{{h\in\cI_{s-1}^*}}\!\!\!\ti\pi_{s;h}^*\ga_{s;h},$$
where $\ti\pi_{s;h}\!:\bP_s\T\!\lra\!\bP\wt{\F_h\T}$
is the natural projection map.
We denote by
$$\al_s\in\Ga\big(\bP_s\T;\ga_s^*\!\otimes\!V_s\big)$$
the section induced by the sections $\ti\pi_{s;h}^*\al_{s;h}$ with $s\!\in\!\cI_{s-1}^*$.
Similarly, let
$$\pi_0\!:  \bP_0\T=\bP\wt{\F\T}\big|_{{\cal U}_{\T}^{(1)}(X;J)}\lra{\cal U}_{\T}^{(1)}(X;J)
\qquad\hbox{and}\qquad
V_0=\pi_0^*\big(\ti\pi_P^*\E^*\!\otimes\!\ev_P^*TX\big) \lra \bP_0\T.$$
Let $\ga_0\!\lra\!\bP_0\T$ be the tautological line bundle.
We define $\al_0\!\in\!\Ga\big(\bP_0\T;\ga_0^*\!\otimes\!V_0\big)$ by
\begin{gather*}
\big\{\al_0\big(b,(v_i)_{i\in\chi(b)}\big)\big\}(b,\psi)
=-2\pi\I\!\sum_{i\in\chi(\T)}\!\!\cD_{J,i}(b,\psi_{x_{h(i)}(b)}v_i\big) \in T_{\ev_P(b)}X\\
\hbox{if}\qquad  \big(b,(v_i)_{i\in\chi(b)}\big)\in\ga_0, 
\quad (b,\psi)\in\E_{\ti\pi_P(b)}.
\end{gather*}\\

\noindent
Let
$$\vph\!: \ti\Ga_-(\cdot)\!\equiv\!
\equiv\!\! \bigcup_{b\in{\cal U}_{\T}^{(\lr{\T})}(X;J)}\!\!\!\!\!\!\!\!\!\! \ti{\Ga}_-(b)
~\lra\ti{\cal U}_{\T}^{(0)}(X;J)$$
be the smooth map induced from the maps~$\vph$ of the remark following
Lemma~\ref{reg0_lmm2} via the decomposition~\e_ref{g1decomp_e3a}.
In particular, 
$$\Si_{\vph(b;\xi)}=\Si_b \quad\hbox{and}\quad
\ev_P\big(\vph(b;\xi)\big)=\ev_P(b) \qquad\forall~ 
b\!\in\!{\cal U}_{\T}^{(\lr{\T})}(X;J), ~\xi\!\in\!\ti\Ga_-(b).$$
Thus, the fibers of the vector bundles
$$\wt{\cal FT},\qquad
\ti\pi_P^*\E^*, \qquad \ev_P^*TX,  \qquad \wt{\F\T}, \quad\hbox{and}\quad 
\wt{\F_h\T} ~~\hbox{for}~~h\!\in\!\cI_{s-1}^*,~s\!\in\![\lr{\T}],$$
at $b$ and at $\vph(b,\xi)$ are canonically isomorphic. 
If $\xi\!\in\!\ti\Ga_-(b)$ is sufficiently small, 
$$\ups\!\equiv\!(b,v)\in\wt{\cal FT}_{\de}^{\eset} \qquad\hbox{and}\qquad 
\ups'\!=\!(\vph(b,\xi),v)\in\wt{\cal FT}_{\de}^{\eset}$$
are corresponding elements of the fibers of $\wt{\cal FT}$ at $b$ and at $\vph(b,\xi)$,
let $\ze_{\ups,0}(\xi)\!\in\!\Ga(\ups)$ be given~by
$$\exp_{u_{\ups,0}}\!\ze_{\ups,0}(\xi)=u_{\ups',0}
\qquad\hbox{and}\qquad \|u_{\ups,0}\|_{C^0}<r_J.$$
We identify $\Ga_+(\ups)$ and $\Ga^{0,1}_+(\ups;J)$ with $\Ga_+(\ups')$ 
and $\Ga^{0,1}_+(\ups';J)$ by composing the $\na^J$-parallel transports
$$\Ga_+(\ups)\lra\Ga(\ups') \qquad\hbox{and}\qquad 
\Ga_+^{0,1}(\ups;J)\lra\Ga^{0,1}(\ups';J)$$
along the geodesics corresponding to $\ze_{\ups,0}(\xi)$ with $L^2$-projection maps
$$\Ga(\ups')\lra\Ga_+(\ups')  \qquad\hbox{and}\qquad
\Ga^{0,1}(\ups';J)\lra\Ga_+^{0,1}(\ups';J)$$
corresponding to the metric $g_{\ups}\!=\!g_{\ups'}$ on $\Si_{\ups}$.\\

\noindent
For each $s\!\in\!\{0\}\!\cup\![\lr{\T}]$,  the map $\vph$ induces a smooth~map
$$\vph_s\!: {\frak G}_s\!\equiv\!\pi_s^*\ti\Ga_-(\cdot)\big|_{{\cal U}_{\T}^{(s)}(X;J)}
\lra \bP_s\T.$$
Similarly to the previous paragraph, the fibers of vector bundles
$$\pi_s^*\wt{\cal FT}, \qquad  \ga_s, \quad\hbox{and}\quad V_s$$
at $\ti{b}$ and $\vph_s(\ti{b},\xi)$ are canonically isomorphic,
if $(b,\xi)\!\in\!{\frak G}_s$ is sufficiently small.
By the regularity assumptions (a) and (b-i) of Definition~\ref{g1reg_dfn}, 
the differential
$$\na^{\T}\!\al_s\!:{\frak G}_s\lra\ga_s^*\!\otimes\!V_s$$
of $\al_s$, defined via the above isomorphisms, is surjective.
Let ${\frak G}_s^{\perp}$ be the $L^2$-orthogonal complement
of $\ker\na^{\T}\!\al_s$ in~${\frak G}_s$.
By the surjectivity of~$\na^{\T}\!\al_s$, the Contraction Principle, 
the precompactness of the fibers of 
$$\pi\circ\pi_s\!:\bP_s\T \lra {\cal U}_{\T}^{(s+1)}(X;J) \lra {\cal U}_{\T}(X;J),$$
there exists $\ep,C\!\in\!C({\cal U}_{\T}(X;J);\R^+)$ with the following property.
If $b\!\in\!{\cal U}_{\T}(X;J)$ and
$$\ka\in \Ga\big(\bP_s\T|_{\pi^{-1}(b)};\ga_s^*\!\otimes\!V_s\big)$$ 
is a smooth section such~that
$$\|\ka(\ti{b})\|,\, \big\|\na^{\T}\ka(\ti{b})\big\| \le \ep(b)
\qquad\forall ~\ti{b}\in\bP_s\T\big|_{\pi^{-1}(b)},$$
then for every 
$\ti{b}^*\!\in\!\bP_s\T|_{\pi^{-1}(b)\cap{\cal U}_{\T}^{(s)}(X;J)}$
the equation
$$\al_s\big(\vph_s(\ti{b}^*,\xi)\big)
+\ka\big(\vph_s(\ti{b}^*,\xi)\big)=\al_s(b^*)\in\ga_s^*\!\otimes\!V_s,
\qquad\xi\!\in\!{\frak G}_s\big|_{\ti{b}^*},~~|\xi|<2C(b)\ep(b),$$
has a unique solution $\xi_{\ka}(\ti{b}^*)$.
Furthermore,
$$\big|\xi_{\ka}(\ti{b}^*)\big|\le 
2C(b)\max\big\{\|\ka(\ti{b})\|,\,\big\|\na^{\T}\ka(\ti{b})\big\|\!: 
\ti{b}\in\bP_s\T\big|_{\pi^{-1}(b)}\big\};$$
see Subsection~3.6 in~\cite{Z4}, for example.\\

\noindent
We are now ready to return to the gluing construction of the previous subsection.
For every element $\ups\!=(b,v)$ of $\wt{\cal FT}_{\de}^{\eset}$, let
$$\mu_{\lr{\T}+1}(\ups)=b\in\ti{\cal U}_{\T}^{(0)}(X;J)  
\qquad\hbox{and}\qquad 
\ti\mu_{\lr{\T}+1}(\ups)=(\mu_{\lr{\T}+1}(\ups),v)\in\wt{\cal FT}_{\de}^{\eset}.$$
Suppose $s\!\in\![\lr{\T}]$ and for all $t\!\in\![\lr{\T}]$ such that $t\!>\!s$
and $\ups\!\in\!\wt{\cal FT}_{\de}^{\eset}|_{{\cal U}_{\T}^{(t)}(X;J)}$ as above
we have constructed
$$\mu_t(\ups)\in\ti{\cal U}_{\T}^{(0)}(X;J)  
\qquad\hbox{and}\qquad 
\ti\mu_t(\ups)=\big(\mu_t(\ups),v\big)\in\wt{\cal FT}_{\de}^{\eset}$$
such that
\begin{gather}\label{indgluing_e1}
\cD_h\ti{b}_t\big(\ti\mu_t(\ups)\big)
=\big\{\al_{t;h}(b)\big\}\big(\rho_{t;h}(\ups)\big)
 \qquad\forall~h\in\cI_{t-1}^*   \qquad\hbox{and}\\
\label{indgluing_e2}\begin{split}
&\mu_t(\ups)=\mu_{t+1}\big(\vph\big(b,\xi_{\ups,t}\big),v\big)
\qquad \hbox{for some}\quad \xi_{\ups,t}\in\ti\Ga_-(b)\\ 
&\quad\hbox{s.t.}\quad
\big|\xi_{\ups,t}\big|,\big|\na^{\T}\xi_{\ups,t}\big| \le C(b)|\ups|^{1/p},~~
\big([\rho_t(\ups)],\xi_{\ups,t}\big)\in{\frak G}_t^{\perp}\big|_{[\rho_t(\ups)]},
\end{split}\end{gather}
where $[\rho_t(\ups)]\!\in\!\bP_t\T$ denotes the image of 
$\rho_t(\ups)\!\in\!\wt{\F_t\T}^{\eset}$ under the projection map 
$\wt{\F_t\T}^{\eset}\!\lra\!\bP_t\T$ and $\na^{\T}\xi_{\ups,t}$ is
the covariant derivative of $\xi_{\ups,t}$ along the directions in $\ti\Ga_-(b)$
as before.\\

\noindent
By~(5) of Lemma~\ref{str_lmm}, its proof, and~\e_ref{indgluing_e2},
\begin{equation}\label{indgluing_e3}
\cD_h\ti{b}_s\big(\ti\mu_{s+1}(\ups)\big)
=\big\{\al_{s;h}(b)+\ve_{s;h}(\ups)\big\}\big(\rho_{s;h}(\ups)\big)
\qquad\forall~h\in\cI_{s-1}^*,
\end{equation}
where $\ve_{s;h}\!\in\!T_{\ev_P(b)}X$ satisfies
\begin{equation}\label{indgluing_e4}
\big|\ve_{s;h}(\ups)\big|,\big|\na^{\T}\ve_{s;h}(\ups)\big| \le C(b)|\ups|^{1/p}.
\end{equation}
The estimate~\e_ref{indgluing_e3} can be restated as 
$$\big(\cD_h\ti{b}_s(\ti\mu_{s+1}(\ups))\big)_{h\in\cI_{s-1}^*}
=\big\{\al_s\big([\rho_s(\ups)]\big)\!+\!\ve_s(\ups)\big\}\big(\rho_s(\ups)\big),$$
where $\ve_s(\ups)\!\in\!\ga_s^*\!\otimes\!V_s|_{[\rho_s(\ups)]}$ satisfies the analogue 
of~\e_ref{indgluing_e4}.
Thus, by the previous paragraph, \e_ref{indgluing_e3}, and~\e_ref{indgluing_e4},
there exists a unique small element $\xi_{\ups,s}\!\in\!\ti\Ga_-(b)$ such that
\begin{gather*}
\big([\rho_s(\ups)],\xi_{\ups,s}\big)\in{\frak G}_s^{\perp}\big|_{[\rho_s(\ups)]} 
\qquad\hbox{and}\\
\cD_h\ti{b}_s\big(\ti\mu_{s+1}\big(\vph(b,\xi_{\ups,s}),v\big)\big)
=\big\{\al_{s;h}(b)\big\}\big(\rho_{s;h}(\ups)\big) \qquad\forall~h\in\cI_{s-1}^*. 
\end{gather*}
Furthermore, $\xi_{\ups,s}$ satisfies the first estimate in~\e_ref{indgluing_e2},
with $t\!=\!s$, for some $C\!\in\!C({\cal U}_{\T}(X;J);\R^+)$.
The second estimate is obtained by differentiating~\e_ref{indgluing_e3}.
Thus, we~take 
$$\mu_s(\ups)=\mu_{s+1}\big(\vph\big(b,\xi_{\ups,s}\big),v\big).$$\\

\noindent
Suppose we have defined $\mu_s(\ups)$ for all $s\!\in\![\lr{\T}]$.
By~\e_ref{str1_crl1e3a}, \e_ref{str1_crl1e3c}, their proof, and~\e_ref{indgluing_e2},
\begin{equation}\label{indgluing_e6a}
\pi^{0,1}_- b_0\big(\ti\mu_1(\ups),\ze\big) = 
\big\{\al_0(b)+\ve_0(\ups,\ze)\big\}\big(\rho(\ups)\big),
\end{equation}
where $\ve_0(\ups,\ze)\!\in\!\ga_0^*\!\otimes\!V_0|_{[\rho(\ups)]}$ satisfies
\begin{gather}\label{indgluing_e6b}
\big\|\ve_0(\ups,\ze)\big\|,\big\|\na^{\T}\ve_0(\ups,\ze)\big\|
\le C(b)\big(|\ups|^{1/p}\!+\!|\ups|^{(p-2)/p}\!+\!\|\ze\|_{\ups,p,1}\big)\cdot
\big(1\!+\!|\rho(\ups)|^{-1}\|\ze\|_{\ti\mu_1(\ups),p,1}\big)\\
\hbox{if}\qquad\ze\!\in\!\Ga\big(\ti\mu_1(\ups)\big),~
\|\ze\|_{\ti\mu_1(\ups),p,1}\!\le\!\de(b).\notag
\end{gather}
Thus, for every 
$$\ups\!=\!(b,v)\in\ti{\cal F}^1{\cal T}^{\eset}_{\de}  \qquad\hbox{and}\qquad
\ze\!\in\!\Ga_+\big(\ti\mu_1(\ups)\big) ~~~\hbox{s.t.}~~~
\|\ze\|_{\ti\mu_1(\ups),p,1}\le2\ti{C}(b)|\rho(\ups)|,$$
the equation 
$$\pi^{0,1}_- b_0\big(\ti\mu_1\big(\vph(b,\xi),v\big),\ze\big) =
\big\{\al_0(b)\big\}\big(\rho(\ups)\big) \equiv 0$$
has a unique small solution $\xi_{\ups,0}(\ze)\!\in\!\ti\Ga_-(b)$ such that
$$\big([\rho(\ups)],\xi_{\ups,0}(\ze)\big)\in{\frak G}_0^{\perp}\big|_{[\rho(\ups)]},$$
provided that
\begin{equation}\label{str1_crl1e21}
C(b)\big(|\ups|^{1/p}\!+\!|\ups|^{(p-2)/p}\!+\!\ti{C}(b)|\rho(\ups)|\big)\cdot
\big(1\!+\!2\ti{C}(b)\big) \le C(b)\ve(b).
\end{equation}
Furthermore, this solution satisfies
\begin{equation}\label{indgluing_e8}
\big|\xi_{\ups,0}(\ze)\big|\le 
C(b)\big(|\ups|^{1/p}\!+\!|\ups|^{(p-2)/p}\big)\cdot
\big(1\!+\!|\rho(\ups)|^{-1}\|\ze\|_{\ti\mu_1(\ups),p,1}\big).
\end{equation}
We~put 
$$\mu_0(\ups,\ze)=\mu_1\big(\vph\big(b,\xi_{\ups,0}(\ze)\big),v\big)
\qquad\hbox{and}\qquad
\ti\mu_0(\ups,\ze)=\big(\mu_0(\ups,\ze),v\big).$$
Let $\mu_0(\ups)\!=\!\mu_0(\ups,0)$ and $\ti\mu_0(\ups,0)\!=\!\ti\mu_0(\ups,0)$.\\

\noindent
For every $\ups\!=\!(b,v)$ in  $\tilde{\cal F}^1{\cal T}^{\eset}$
sufficiently small, we define the~map
\begin{gather*}
\Psi_{\ups}\!:
\big\{\ze\!\in\!\Ga_+(\ti\mu_0(\ups))\!:
\|\ze\|_{\ti\mu_0(\ups),p,1}\!\le\!2\ti{C}(b)|\rho(\ups)|\big\}
\lra \Ga^{0,1}_+\big(\ti\mu_0(\ups);J\big)               \qquad\hbox{by}\\
\Psi_{\ups}(\ze)= \bar{\partial}_Ju_{\ti\mu_0(\ups,\ze)}
+D_{J,b_0(\mu_0(\ups,\ze))}\ze+N_{J,\mu_0(\ups,\ze)}\ze.
\end{gather*}
Let 
$$D_{J,\ups}^+\!:\Ga_+\big(\ti\mu_0(\ups)\big) \lra 
\Ga^{0,1}_+\big(\ti\mu_0(\ups);J\big)$$ 
be the derivative of $\Psi_{\ups}$ at $\ze\!=\!0$ and let
$N_{J,\ups}^+\ze\!\in\!\Ga^{0,1}_+(\ti\mu_0(\ups);J)$ be given~by 
$$\Psi_{\ups}(\ze)=\Psi_{\ups}(0)+D_{J,\ups}^+\ze+N_{J,\ups}^+\ze.$$
By the construction of $\Psi_{\ups}$, (3) of Lemma~\ref{str_lmm}, \e_ref{reg1_crl2e2},
and~\e_ref{indgluing_e8}
\begin{gather}\label{str1_crl1e24a}
\|\Psi_{\ups}(0)\|_{\ti\mu_0(\ups),p,1}\le 2C(b)\big|\rho(\ups)\big|, \\
\label{str1_crl1e24b}
\big(2C(b)\big)^{-1}\|\ze\|_{\ti\mu_0(\ups),p,1}\le 
\|D_{J,\ups}^+\ze\|_{\ti\mu_0(\ups),p}
\le 2C(b)\|\ze\|_{\ti\mu_0(\ups),p,1} ~~\forall\, \ze\!\in\!\Ga_+\big(\ti\mu_0(\ups)\big),\\
\label{str1_crl1e24c}
\big\|N_{J,\ups}^+\ze\!-\!N_{J,\ups}^+\ze'\big\|_{\ti\mu_0((\ups),p}
\le 2C(b)\big(\|\ze\|_{\ti{\mu}_0(\ups),p,1}\!+\!\|\ze'\|_{\ti\mu_0(\ups),p,1}\big)
\|\ze\!-\!\ze'\|_{\ti{\mu}_0(\ups),p,1}
~\forall\ze,\ze'\!\in\!\hbox{Dom}\,\Psi_{\ups},
\end{gather}
provided that $\ups$ is sufficiently small.
Since the index of $D_{J,\ups}^+$ is zero, if $\ti{C}\!\in\!C({\cal U}_{\T}(X;J);\R^+)$
is sufficiently large and $\ups\!\in\!{\cal F}^1{\T}_{\de}^1$ is sufficiently small, 
by \e_ref{str1_crl1e24a}-\e_ref{str1_crl1e24c} and  the Contraction Principle, the equation 
$$\Psi_{\ups}(\ze)=0$$
has a unique solution $\ze_{\ups,0}\!\in\!\Ga_+(\ti\mu_0(\ups))$.\\

\noindent
If $\ti\de\!\in\!C({\cal U}_{\T}(X;J);\R^+)$ is sufficiently small, 
we define the map
\begin{equation}\label{str1_crl1e29}
\phi\!:\ti{\cal F}^1\T^{\eset}_{\ti\de}\lra \M_{1,M}^0(X,A;J)
\qquad\hbox{by}\qquad
\phi(\ups)=\big[\exp_{b_0(\ti\mu_0(\ups,\ze_{\ups,0}))}\!\ze_{\ups,0}\big].
\end{equation}
By construction, the map $\phi$ is $\Aut(\T)\!\propto\!(S^1)^{\hat{I}}$-invariant
and smooth, and thus descends to a smooth~map
$$\phi\!:{\cal F}^1\T^{\eset}_{\ti{\de}}\lra\M_{1,M}^0(X,A;J).$$
By an argument analogous to that in Subsection~4.2 of~\cite{Z4},
the map $\phi$ is an immersion into $\X_{1,M}^0(X,A)$, if $\ti{\de}$ is sufficiently small.
By the proof of Corollary~\ref{reg1_crl3} and the construction of the map~$\phi$,
the image of~$\phi$ contains $\M_{1,M}^0(X,A;J)\cap U_{\T}$
for a neighborhood $U_{\cal T}$ of ${\cal U}_{\T}(X;J)$ in $\X_{1,M}(X,A)$.
Thus, the map
$$\phi\!:{\cal F}^1\T^{\eset}_{\ti{\de}} \lra\M_{1,M}^0(X,A;J)\cap U_{\T}$$
is a diffeomorphism.
It can be seen to be orientation-preserving by an argument similar to
that of Subsection~3.9 in~\cite{Z4}.

\subsection{Extension to Homeomorphism}
\label{strthmcont_subs}

\noindent
In the rest of this section, we show that the map $\phi$ extends continuously 
over $\ti{\cal F}^1{\T}_{\ti{\de}}$.
This will be achieved by combining the approach of Subsections~3.9 and~4.1 in~\cite{Z4}
with the conditions~\e_ref{indgluing_e1} and~\e_ref{indgluing_e6a} on the corrections
$\xi_s(\ups)$ to the maps $\ti{b}_s(\ups)$.\\\

\noindent
For every $b\!\in\!\ti{\cal U}_{\T}^{(0)}(X;J)$, $s\!\in\![\lr{\T}]$, and
$h\!\in\!\cI_{s-1}^*$, let
\begin{gather*}
\Si_b^h=\bigcup_{h\le i}\!\Si_{b,i} \subset\Si_b, \qquad
\Ga_h^{0,1}(b;J)=\big\{\eta\!\in\!\Ga_B^{0,1}(b;J)\!:\eta|_{\Si_b-\Si_b^h}\!=\!0\big\},\\
\Ga_h(b)=\big\{\xi\!\in\!\Ga_B(b)\!: \xi|_{\Si_b-\Si_b^h}\!=\!0\big\}, 
\quad\hbox{and}\quad \ti\Ga_{h;-}(b)=\Ga_h(b)\cap\ti\Ga_-(b).
\end{gather*}
We note that by (a) of Definition~\ref{g0reg_dfn}, the operator 
$$D_{J,b}^h\!:\Ga_h(b)\lra\Ga_h^{0,1}(b;J)$$
induced by $D_{J,b}$ is surjective.
By the regularity assumptions (a) and (b-i) of Definition~\ref{g1reg_dfn}, 
the differential
$$\na^{\T}\!\al_{s;h}\!:{\frak G}_{s;h}\!\!\equiv\!
\pi_{s;h}^*\!\!\! \bigcup_{b\in\wt{\cal U}_{\T}^{(0)}(X;J)} \!\!\!\!\!\!\!\!\!\!
\ti\Ga_{h;-}(b)     \lra\ga_{s;h}^*\!\otimes\!V_{s;h}$$
of $\al_{s;h}$ is surjective. 
We denote by ${\frak G}_{s;h}^{\perp}$ be the $L^2$-orthogonal complement
of $\ker\na^{\T}\!\al_{s;h}$ in~${\frak G}_{s;h}$.
If $w_h\!\in\!\bP\wt{\F_h\T}|_b$ and $w\!\in\!\bP_s\T|_b$, 
let
$$\ti\Ga_{h;-}(b;w_h)=\big\{\xi\!\in\!\ti\Ga_{h;-}(b)\!: 
(w_h,\xi)\!\in\!{\frak G}_{s,h}^{\perp}\big\}
\quad\hbox{and}\quad
\ti\Ga_-(b;w)=\big\{\xi\!\in\!\ti\Ga_-(b)\!: 
(w,\xi)\!\in\!{\frak G}_s^{\perp}\big\}.$$
We note that
\begin{equation}\label{pertsplit_e1}
\ti\Ga_-(b;w)=\bigoplus_{h\in\cI_{s-1}^*}\!\!\!\ti\Ga_{h;-}(b;w_h)
\qquad\hbox{if}\quad  w=(w_h)_{h\in\cI_{s-1}^*}.
\end{equation}
If $\ups\!=\!(b,v)\!\in\!\wt{\cal FT}$ and $s$ and $h$ are as above, let
\begin{gather*}
\cI_{s-1}^0(\ups)=\big\{h\!\in\!\cI_{s-1}\!:\rho_{s;h}(\ups)\!=\!0\big\}
\qquad\hbox{and}\qquad
\Si_{\ups_{\lr{s}}}^h=q_{\ups_{\lr{s}}}^{-1}\big(\Si_b^h\big).
\end{gather*}
We note that $\Si_{\ups_{\lr{s}}}^h$ is a union of components of $\Si_{\ups_{\lr{s}}}$. \\

\noindent
The multi-step gluing construction of Subsection~\ref{g1gluing2_subs} extends 
continuously over $\wt{\cal FT}_{\de}$. 
This extension is formally described in exactly the same way as the construction itself;
see Subsection~3.9 in~\cite{Z4} for a description of the continuous extension for 
a similar gluing construction and Subsection~4.1 in~\cite{Z4} for a proof of its continuity.
Using the notation of Subsection~\ref{g1gluing2_subs}, 
we now make an observation regarding this extension.
For each $s\!\in\![\lr{\T}]$ and $h\!\in\!\cI_{s-1}^*$, let
\begin{gather*}
\Ga_h(\ups_{\lr{s}})=
\big\{\xi\!\in\!\Ga_B(\ups_{\lr{s}})\!:
\xi|_{\Si_{\ups_{\lr{s}}}-\Si_{\ups_{\lr{s}}}^h}\!=\!0\big\}, \\
\Ga_h^{0,1}(\ups_{\lr{s}};J)=
\big\{\eta\!\in\!\Ga_B^{0,1}(\ups_{\lr{s}};J)\!:
\eta|_{\Si_{\ups_{\lr{s}}}-\Si_{\ups_{\lr{s}}}^h}\!=\!0\big\}.
\end{gather*}
By the surjectivity of the operators $D_{J,b}^h$, with $h\!\in\!\cI_{s-1}^*$,
for all $h\!\in\!\cI_{s-1}^*$ the operators 
$$D_{J,b_s(\ups)}^h\!: \Ga_h(\ups_{\lr{s}})\lra\Ga_h^{0,1}(\ups_{\lr{s}};J)$$
induced by $D_{J,b_s(\ups)}$  are surjective, provided $\ups$ is sufficiently small.
Since $\bar\partial_Ju_{\ups,s}$ vanishes on $\Si_{\ups_{\lr{s}}}^h$ for all
$h\!\in\!\cI_{s-1}^0(\ups)$ and $\ze_{\ups,s}$ is 
the unique small solution of~\e_ref{indpert_e2c}, it follows that
\begin{equation}\label{indcont_e1}
\ze_{\ups,s}\big|_{\Si_{\ups_{\lr{s}}}^h}=0~~~\forall~ h\!\in\!\cI_{s-1}^0(\ups).
\end{equation}\\

\noindent
We next extend the construction of perturbations $\xi_{\ups,s}$ for
$\ups\!\in\!\wt{\cal FT}_{\de}^{\eset}|_{{\cal U}_{\T}^{(s)}(X;J)}$
in Subsection~\ref{strthm_subs} to $\wt{\cal FT}_{\de}|_{{\cal U}_{\T}^{(s)}(X;J)}$.
Suppose $s\!\in\![\lr{\T}]$ and for all $t\!\in\![\lr{\T}]$ such that $t\!>\!s$
and for all elements $\ups\!=\!(b,v)$ in $\wt{\cal FT}_{\de}|_{{\cal U}_{\T}^{(t)}(X;J)}$
we have constructed
$$\mu_t(\ups)\in\ti{\cal U}_{\T}^{(0)}(X;J)  
\qquad\hbox{and}\qquad 
\ti\mu_t(\ups)=\big(\mu_t(\ups),v\big)\in\wt{\cal FT}_{\de}$$
such that
\begin{gather}\label{indcon2_e2}
\cD_h\ti{b}_t\big(\ti\mu_t(\ups)\big)
=\big\{\al_{t;h}(b)\big\}\big(\rho_h(\ups)\big)
 \qquad\forall~h\in\cI_{t-1}^*   \qquad\hbox{and}\\
\label{indcont_e3}\begin{split}
\mu_t(\ups)=\mu_{t+1}\big(\vph\big(b,\xi_{\ups,t}\big),v\big)
&\qquad \hbox{for some}\quad 
\xi_{\ups,t}\in\bigoplus_{h\in\cI_{t-1}^*-\cI_{t-1}^0(\ups)} \!\!\!\!\!\!\!\!\!\!\!\!\!
\ti\Ga_{h;-}\big(b;[\rho_h(\ups)\big)\\ 
&\hbox{s.t.}\quad
\big|\xi_{\ups,t}\big|,\big|\na^{\T}\xi_{\ups,t}\big| \le C(b)|\ups|^{1/p}.
\end{split}\end{gather}
We note that $\cD_h\ti{b}_t(\ti\mu_t(\ups))\!=\!0$ for any $h\!\in\!\cI_{t-1}^0(\ups)$
and $\mu_t(\ups)\in\ti{\cal U}_{\T}^{(0)}(X;J)$.
Thus, \e_ref{indcon2_e2} is a nontrivial condition only for 
$h\!\in\!\cI_{t-1}^*\!-\!\cI_{t-1}^0(\ups)$.
In particular, if $\cI_{s-1}^0(\ups)\!=\!\cI_{s-1}^*$, for the inductive step in
the construction of the previous subsection we simply take~$\xi_{\ups,s}\!=\!0$.
On the other hand, if $\cI_{s-1}^0(\ups)\!\neq\!\cI_{s-1}^*$,  the inductive step
is nearly the same as the before.
The only difference is that instead of working with the section 
$$\al_s\equiv  \{\al_{s;h}\}_{h\in\cI_{s-1}^*}$$
over $\bP_s\T$, we work with the section
$\{\al_{s;h}\}_{h\in\cI_{s-1}^*-\cI_{s-1}^0(\ups)}$ over the fiber product of the bundles 
$$\{\bP\wt{\F_h\T}\}_{h\in\cI_{s-1}^*-\cI_{s-1}^0(\ups)}
\lra \ti{\cal U}_{\T}^{(s+1)}(X;J).$$
We note that in this case the orthogonal complement of the kernel of 
$$\na^{\T}\{\al_{s;h}\}_{h\in\cI_{s-1}^*-\cI_{s-1}^0(\ups)}$$ 
is given by~\e_ref{pertsplit_e1}
with $\cI_{s-1}^*$ replaced by $\cI_{s-1}^*\!-\!\cI_{s-1}^0(\ups)$.
In summary, 
$$\xi_{\ups,s}\in \!\!
\bigoplus_{h\in\cI_{s-1}^*-\cI_{s-1}^0(\ups)}\!\!\!\!\!\!\!\!\!\!\!\!
\ti\Ga_{h;-}\big(b;[\rho_{s;h}(\ups)\big)$$
is the unique small solution to the system of equations
$$\al_{s;h}\big(\vph(b,\xi_{\ups,s}),w_h\big)+
\ve_{s;h}\big(\vph(b,\xi_{\ups,s}),v\big)
=\al_{s;h}\big(b,w_h\big), \qquad h\in\cI_{s-1}^*\!-\!\cI_{s-1}^0(\ups).$$
The final, $s\!=\!0$, step splits into two cases as well.
If $\rho(\ups)\!=\!0$, then we take 
$$\ze_{\ups,0}=0  \qquad\hbox{and}\qquad 
\xi_{\ups,0}\equiv \xi_{\ups,0}(\ze_{\ups,0})=0.$$
Otherwise, the argument of the previous subsection still applies.\\

\noindent
We will show that the above extension of the construction described in 
Subsection~\ref{strthm_subs} is continuous at every step. 
First, note that by definition, for every $s\!\in\![\lr{\T}\!-\!1]$ and 
$h\!\in\!\cI_{s-1}^*\!\cap\!\chi^-(\T)$,
\begin{equation}\label{contrho_e}
\rho_{s;h}(\ups)=\big(v_{h'}\rho_{s+1;h'}(\ups)\big)_{\io_{h'}=h}
\qquad\forall~ 
\ups\!\equiv\!\big(b,(v_i)_{i\in\aleph\cup\hat{I}}\big)\in\wt{\cal FT}.
\end{equation}
In this case, by the proof of Lemma~\ref{reg0_lmm4},
\begin{equation}\label{contest_e1}
\cD_h\ti{b}_s(\ups) =  
\sum_{\io_{h'}=h}\!\!v_{h'}\cD_{h'}\ti{b}_{s+1}(\ups)+
\sum_{i\in\chi_h(\T)}\!\!\!\!\ve_{h;i}(\ups)\rho_{h;i}(\ups)
\qquad\forall~\ups\!\equiv\!(b,v)\in\wt{\cal FT}_{\de}^{\eset},
\end{equation}
where $\ve_{h;i}(\ups)\!\in\!T_{\ev_P(b)}X$ is given by the right-hand side 
of~\e_ref{reg0_lmm4e7b}, with 
$$k=1,\quad \ti\ze_{\ti{J},\ups}\!=\!\ze_{\ups,s}, \quad
\vt_b=\vt_{\ti{b}_{s+1}(\ups)}, \quad \Phi_b=\Phi_{J,\ups_{\lr{s+1}}},
\quad \Phi_{\ti{J},\ups}=\Phi_{J,\ups_{\lr{s}}},$$
and with $\de_K$ replaced by a function $\de\!=\!\de(b)$.
Since $\partial^-A_i^-(\de(b))\!\subset\!\Si_{\ups_{\lr{s}}}^h$ 
for all $i\!\in\!\chi_h(\T)$, it follows that
\begin{gather*}
\big|\ve_{h;i}(\ups)\big|,\,\big|\na^{\T}\ve_{h;i}(\ups)\big|
\le C(b)\big(\big\|\ze_{\ups,s}|_{\Si_{\ups_{\lr{s}}}^h}\big\|_{\ups,p,1}
\!+\!\big\|\na^{\T}\ze_{\ups,s}|_{\Si_{\ups_{\lr{s}}}^h}\big\|_{\ups,p,1}\big)\\
\forall\quad\ups\!\equiv\!(b,v)\in\wt{\cal FT}_{\de}^{\eset},~
h\!\in\!\cI_{s-1}^*\!\cap\!\chi^-(\T),~i\!\in\!\chi_h(\T).
\end{gather*}
Let $\ve_{s;h}\!=\!\ve_{s;h}(\ups)$ be the bundle map as in~\e_ref{indgluing_e3}.
Combining the last estimate with \e_ref{contest_e1}, \e_ref{contrho_e}, 
and~\e_ref{indgluing_e1} with $t\!=\!s\!+\!1$, we find that
\begin{gather}\label{contest_e3}
\big|\ve_{s;h}(\ups)\big|,~\big|\na^{\T}\ve_{s;h}(\ups)\big|
\le C(b)\big(\big\|\ze_{\ups,s}|_{\Si_{\ups_{\lr{s}}}^h}\big\|_{\ups,p,1}
\!+\!\big\|\na^{\T}\ze_{\ups,s}|_{\Si_{\ups_{\lr{s}}}^h}\big\|_{\ups,p,1}\big)\\
\forall\quad\ups\!\equiv\!(b,v)\in
 \wt{\cal FT}_{\de}^{\eset}\big|_{{\cal U}_{\T}^{(s+1)}(X;J)},~
 h\!\in\!\cI_{s-1}^*.\notag
\end{gather}
We note that if $h\!\in\!\cI_{s-1}^*\!-\!\chi^-(\T)$, \e_ref{contest_e3} follows immediately
from~\e_ref{indgluing_e1} with $t\!=\!s\!+\!1$.\\

\noindent
We next observe that by the proof of the estimate in (6) of Lemma~\ref{reg1_lmm3},
\begin{equation}\label{contest_e5}
\pi_{\ups;-}^{0,1}\bar\partial_J u_{\ups,0} = \!\! 
\sum_{h\in(\chi(\T)\cup\chi^-(\T))\cap I_1}\!\!\!\!\!\!\!\!\!\!\!\!\!\!\!\!\!
\big(-2\pi\I\,v_h\cD_h\ti{b}_1(\ups)\!+\!\ve_h(\ups)\!+\!\ti\ve_h(\ups)\big)
\qquad\forall~\ups\!\equiv\!(b,v)\in\wt{\cal FT}_{\de}^{\eset}.
\end{equation}
The error term $\ve_h(\ups)$ is described by the left-hand side of~\e_ref{dadj_e11a}.
This term and its $\na^{\T}$-derivative are bounded by the last
expression in~\e_ref{dadj_e14}.
The other term is given~by
\begin{gather*}
\ti\ve_h(\ups)=v_h \Big(2\pi\I\, \cD_h\ti{b}_1(\ups)-
\oint_{\partial A_{\ups_{\lr{1}},h}^-(|v_h|^2/\de(b))}\!
\ti\ze_{b;\ups} \frac{dw_h}{w_h^2}\Big),  \qquad\hbox{where}\\
\ti\ze_{b;\ups}\!: A_{\ups_{\lr{1}},h}^-(\de(b))\lra T_{\ev_P(b)}X, \qquad
\exp_{\ev_P(b)}\!\ti\ze_{b;\ups}=\ti{u}_{\ups,1}, \qquad
\big\|\ti\ze_{b;\ups}\big\|_{C^0}< r_J; 
\end{gather*}
see the proof of~(6) of Lemma~\ref{reg1_lmm3}.
Applying the approach of Lemma~\ref{reg0_lmm2} and Cauchy's formula, 
we find~that
\begin{equation}\label{contest_e8}
\ti\ve_h(\ups)=v_h \oint_{\partial A_{\ups_{\lr{1}},h}^-(|v_h|^2/\de(b))}\!
\big(\ti\Phi_{\ups,1}\!-\!\id\big)\ti\vt_{\ups,1} \frac{dw_h}{w_h^2},
\end{equation}
where
\begin{gather}\label{contest_e9a}
\ti\Phi_{\ups,1}\big|_{w=0}=\id, \qquad
\big\|\ti\Phi_{\ups,1}\!-\!\id\|_{\ups_1,p,1}\le C(b)
\big\|du_{\ups,1}|_{A_{\ups_{\lr{1}},h}^-(\de(b))}\big\|_{\ups_1,p,1} \le 
C(b)\big|\ti{\rho}_h(\ups)\big|,\\
\label{contest_e9b}
\big|\ti\vt_{\ups,1}\big|_{w_h} \le C(b)|w_h|\big|\ti{\rho}_h(\ups)\big|,
\end{gather}
by (2b) and the first estimate in (2c) of Corollary~\ref{reg0_crl4}.
By \e_ref{contest_e8}-\e_ref{contest_e9b} and Holder's inequality
\begin{equation}\label{contest_e10}\begin{split}
\big|\ti\ve_h(\ups)\big| 
&\le C(b)|v_h|\big|\ti\rho_h(\ups)\big|
\int_0^{2\pi}\big|\ti\Phi_{\ups,1}(|v_h|^2/\de(b),\th)\!-\!\id\big|\,d\th\\
&\le  C(b)|\rho_h(\ups)|
\int_0^{2\pi}\!\!\!\int_0^{|v_h|^2/\de(b)} \big|d\ti\Phi_{\ups,1}\big|_{r,\th}\,drd\th\\
&\le C(b) |\rho(\ups)|
\big\|\ti\Phi_{\ups,1}\!-\!\id\big\|_{\ups_1,p,1}
\Big(\int_{A_{\ups_{\lr{1}},h}^-(|v_h|^2/\de(b))}
\!|w_h|^{-\frac{p}{p-1}}\Big)^{\frac{p-1}{p}}\\
& \le C(b)|\rho(\ups)|\big|\ti\rho_h(\ups)\big|\cdot|v_h|^{\frac{p-2}{p}}
\le C(b)|\rho(\ups)|^{\frac{p-2}{p}}|\rho(\ups)|.
\end{split}\end{equation}
Let $\ve_0\!=\!\ve_0(\ups,\ze)$ be the bundle map as in~\e_ref{indgluing_e6a}.
Combining the last estimate with \e_ref{contest_e1}, \e_ref{contest_e5}, 
\e_ref{indgluing_e1} with $t\!=\!1$, and~(4) of Lemma~\ref{str_lmm} we find that
\begin{gather}\label{contest_e15}
\big|\ve_0(\ups,\ze)\big|,~\big|\na^{\T}\ve_0(\ups,\ze)\big|
\le C(b)|\rho(\ups)|^{\frac{p-2}{p}}\\
\forall\quad\ups\!\equiv\!(b,v)\in
\wt{\cal FT}_{\de}^{\eset}\big|_{{\cal U}_{\T}^{(1)}(X;J)},~~\ze\!\in\!\Ga(\ups)~\hbox{s.t.}~
\|\ze\|_{\ups,p,1}\le \ti{C}(b)|\rho(\ups)|.\notag
\end{gather}
The same holds for the derivatives of $\ve_0(\ups,\ze)$ and $\na^{\T}\ve_0(\ups,\ze)$
with respect to~$\ze$.\\

\noindent
Suppose $s\!\in\![\lr{\T}]$ and for every $t\!\in\![\lr{\T}]$ such that $t\!>\!s$
the bundle map
$$\wt{\cal FT}_{\de} \lra \ti\Ga_-(\cdot), \qquad \ups\lra\xi_{\ups,t},$$ 
and its $\na^{\T}$-derivative are continuous over ${\cal U}_{\T}^{(t)}(X;J)$.
We will show that this must also be the case for $t\!=\!s$.
Since the maps 
\begin{alignat*}{2}
&\wt{\cal FT}_{\de} \lra \bigcup_{\ups\in\wt{\cal FT}_{\de}}\!\!\!\!\Ga_{B;+}(\ups_{\lr{s}}),
&\qquad& \ups\!\lra\!\ze_{\ups,s}, \qquad\hbox{and}\\
&\wt{\cal FT}_{\de}\big|_{{\cal U}_{\T}^{(s+1)}(X;J)} \lra \wt{\cal FT}_{\de},
&\qquad& \ups\lra\ti\mu_{s+1}(\ups),
\end{alignat*}
are continuous, so are the maps 
$$\wt{\cal FT}_{\de}\big|_{{\cal U}_{\T}^{(s+1)}(X;J)} \lra \ga_{s;h}^*\!\otimes\!V_{s;h}, 
\qquad \ups\!\lra\!\ve_{s;h}(\ups),$$ 
for all $h\!\in\!\cI_{s-1}^*$.
Suppose $\ups_r\!\equiv\!(b_r,v_r)$ is a sequence of elements in 
$\wt{\cal FT}_{\de}^{\eset}\big|_{{\cal U}_{\T}^{(s)}(X;J)}$ such that
$$\lim_{r\lra\i}\!b_r=b'\in{\cal U}_{\T}^{(s)}(X;J) 
\qquad\hbox{and}\qquad
\lim_{r\lra\i}\!\ups_r=\ups'\!\equiv\!(b',v')\in\wt{\cal FT}_{\de}.$$
Let 
\begin{gather*}
w_h=\lim_{r\lra\i}\!\big[\rho_{s;h}(\ups_r)\big] \in \bP\wt{\F_h\T}\big|_{b'}
\quad\hbox{if}~~ h\!\in\!\cI_{s-1}^*; \qquad
w=(w_h)_{h\in\cI_{s-1}^*}\in\bP_s\T\big|_{b'};\\
\xi_{\ups',s}'\!\equiv\!(\xi_{\ups',s;h}')_{h\in\cI_{s-1}^*}
=\lim_{r\lra\i}\!\xi_{\ups_r,s}\in \ti\Ga_-(b';w)
=\bigoplus_{h\in\cI_{s-1}^*}\!\!\!\ti\Ga_{h;-}(b;w_h).
\end{gather*}
We recall that
$$\xi_{\ups_r,s}\!\equiv\!(\xi_{\ups_r,s;h})_{h\in\cI_{s-1}^*}
\in \ti\Ga_-\big(b_r;[\rho_s(\ups_r)]\big)
=\bigoplus_{h\in\cI_{s-1}^*}\!\!\!\ti\Ga_{h;-}\big(b;[\rho_{s;h}(\ups_r)]\big)$$
is the unique small solution to the system of equations
$$\al_{s;h}\big(\vph(b,\xi_{\ups_r,s}),[\rho_{s;h}(\ups_r)]\big)+
\ve_{s;h}\big(\vph(b,\xi_{\ups_r,s}),v_r\big)
=\al_{s;h}\big(b,[\rho_{s;h}(\ups_r)]\big), \qquad h\in\cI_{s-1}^*.$$
Thus,  $\xi_{\ups',s}'\!\in\!\ti\Ga_-(b';w)$ 
is the unique small solution to the system of equations
$$\al_{s;h}\big(\vph(b,\xi_{\ups',s}'),w_h\big)+
\ve_{s;h}\big(\vph(b,\xi_{\ups',s}'),v'\big)
=\al_{s;h}\big(b,w_h\big), \qquad h\in\cI_{s-1}^*.$$
Since 
$$\ve_{s;h}\big(\vph(b,\xi_{\ups',s}'),v'\big)=0  \qquad\forall~ h\in\cI_{s-1}^0(\ups')$$
by~\e_ref{contest_e3} and~\e_ref{indcont_e1}, 
$\xi_{\ups',s;h}'\!=\!0$  for all $h\!\in\!\cI_{s-1}^0(\ups')$ and
$\xi_{\ups',s}'\!=\!\xi_{\ups',s}$, as needed.\\

\noindent
We finally show that the map 
\begin{equation}\label{indcont_e25}
\ups\lra \ti{b}_0(\ups)\!\equiv\!
\big(\Si_{\ups};\exp_{b_0(\ti\mu_0(\ups,\ze_{\ups,0}))}\!\ze_{\ups,0}\big)
\end{equation}
is continuous over $\wt{\cal F}^1\T_{\de}^{\eset}$.
First, the map
$$\bigcup_{\ups\in\wt{\cal FT}_{\de},\rho(\ups)\neq0} \!\!\!\!\!\!\!\!\!\!\!\!
\big\{\ze\!\in\!\Ga_+(\ups)\!:\|\ze\|_{\ups,p,1}\!<\!\de(b)\big\}
\lra \bigcup_{b\in\ti{\cal U}_{\T}^{(0)}(X;J)}\!\!\!\!\!\!\!\!\Ga_-^{0,1}(b;J),
\qquad (\ups,\ze) \lra \ve(\ups,\ze),$$
of~\e_ref{str1_crl1e3a} is continuous.
Since so is the map $\ups\!\lra\!\xi_{\ups,1}$, the map
$$\bigcup_{\ups\in\wt{\cal FT}_{\de},\rho(\ups)\neq0} \!\!\!\!\!\!\!\!\!\!\!\!
\big\{\ze\!\in\!\Ga_+(\ups)\!:\|\ze\|_{\ups,p,1}\!<\!\ti{C}(b)|\rho(\ups)|\big\}
\lra \ga_0^*\!\otimes\!V_0, 
\qquad (\ups,\ze) \lra \ve_0(\ups,\ze),$$
is also continuous.
It then follows immediately from the construction that the maps
$$(\ups,\ze)\lra\xi_{\ups,0}(\ze) \qquad\hbox{and}\qquad \ups\lra\ze_{\ups,0}$$ 
are continuous over $\wt{\cal F}^1\T_{\de}\!-\!\rho^{-1}(0)$.
On the other hand, suppose $\ups_r\!\equiv\!(b_r,v_r)$ is a sequence of elements in 
$\wt{\cal F}^1\T_{\de}^{\eset}$ such that
$$\lim_{r\lra\i}\!b_r=b'\in{\cal U}_{\T}^{(0)}(X;J), \qquad
\lim_{r\lra\i}\!\ups_r=\ups'\!\equiv\!(b',v')\in\wt{\cal F}^1\T_{\de},
\quad\hbox{and}\quad \rho(\ups')=0.$$
Let 
\begin{gather*}
w=\lim_{r\lra\i}\!\big[\rho(\ups_r)\big] \in \bP\wt{\F\T}\big|_{b'}
\qquad\hbox{and}\qquad
\xi_{\ups',0}'=\lim_{r\lra\i}\!\xi_{\ups_r,0}(\ze_{\ups_r,0})\in \ti\Ga_-(b';w).
\end{gather*}
Since $\rho(\ups_r)\!\lra\!0$, \e_ref{contest_e15} implies that 
$\xi_{\ups',0}'\!\in\!\ti\Ga_-(b';w)$ is the unique small solution~of
$$\al_0\big(\vph(b',\xi_{\ups',0}'),w\big)=0=\al_0(b',w).$$
Thus, $\xi_{\ups',0}'\!=\!0\!=\!\xi_{\ups',0}$.
Furthermore, by~\e_ref{str1_crl1e2}
$$\lim_{r\lra\i}\!\ze_{\ups_r,0}=0= \ze_{\ups',0}.$$
It follows that the map~\e_ref{indcont_e25} is continuous.\\

\noindent
We have thus constructed a continuous map
$$\phi\!:{\cal F}^1\T_{\de} \lra\ov\M_{1,M}^0(X,A;J)\cap U_{\T},$$
where $U_{\T}$ is a neighborhood of ${\cal U}_{\T;1}(X;J)$ in $\X_{1,M}(X,A)$.
By the same argument as in Subsections~4.2 and~4.5,
this map is injective if $\de\!\in\!C({\cal U}_{\T}(X;J);\R^+)$
and surjective if $U_{\T}$ is sufficiently small.
Since the space $\ov\M_{1,M}^0(X,A;J)$ is Hausdorff and $\phi|_{{\cal U}_{\T}(X;J)}$
is the identity map, it follows that $\phi$ is a homeomorphism for $\de$ sufficiently
small.\\

\vspace{.2in}

{\it 
\begin{tabbing}
${}\qquad$
\= Department of Mathematics, Stanford University, Stanford,
CA 94305-2125\\ 
\> \textnormal{Current Address:} Department of Mathematics, SUNY, Stony Brook, NY 11794-3651\\
\> azinger@math.sunysb.edu
\end{tabbing}}


\vspace{.2in}



\begin{thebibliography}{99}

\bibitem[BCOV]{BCOV} M.~Bershadsky, S.~Cecotti, H.~Ooguri, and C.~Vafa,
{\it Holomorphic Anomalies in Topological Field Theories}, 
Nucl.~Phys.~B405 (1993), 279--304. 

\bibitem[BeFa]{BeFa} K.~Behrend and B.~Fantechi, 
{\it The Intrinsic Normal Cone}, Invent.~Math.~128 (1997), no.~1, 45--88.

\bibitem[FlHSa]{FlHSa} A.~Floer, H.~Hofer, and D.~Salamon, 
{\it Transversality in Elliptic Morse Theory for the Symplectic Action},
Duke Math.~J.~80 (1996), no.~1, 251-292.

\bibitem[FuOn]{FuOn} K.~Fukaya and K.~Ono,
{\it Arnold Conjecture and Gromov-Witten Invariant},
Topology 38 (1999), no.~5, 933--1048.

\bibitem[LT]{LT}  J.~Li and G.~Tian, 
{\it Virtual Moduli Cycles and Gromov-Witten Invariants of General Symplectic Manifolds}, 
Topics in Symplectic \hbox{$4$-Manifolds},
47-83, First Int.~Press Lect.~Ser., I, Internat.~Press, 1998.

\bibitem[LZ]{LZ} J.~Li and A.~Zinger, 
{\it On the Genus-One Gromov-Witten Invariants 
of Complete Intersection Threefolds}, math/0507104.

\bibitem[McSa]{McSa} D.~McDuff and D.~Salamon, 
{\it Introduction to $J$-Holomorphic Curves},
American Mathematical Society, 1994.

\bibitem[RT1]{RT1} Y.~Ruan and G.~Tian, {\it A Mathematical Theory
of Quantum Cohomology}, J.~Diff.~Geom.~42 (1995), no.~2, 259--367. 

\bibitem[RT2]{RT2} Y.~Ruan and G.~Tian, 
{\it Higher Genus Symplectic Invariants and Sigma Models Coupled with Gravity},
Invent.~Math.~130 (1997), no.~2, 455--516.

\bibitem[Va]{Va} R.~Vakil, 
{\it A Tool for Stable Reduction of Curves on Surfaces},
Advances in Algebraic Geometry Motivated by Physics,
145--154, Amer.~Math.~Soc., 2001.

\bibitem[VaZ]{VaZ} R.~Vakil and A.~Zinger, 
{\it A Desingularization of the Main Component of
the Moduli Space of Genus-One Stable Maps into $\bP^n$}, 
math.AG/0603353.

\bibitem[Z1]{Z1}
A.~Zinger, {\it Basic Estimates of Riemannian Geometry
Used in Gluing Pseudoholomorphic Maps}, notes.

\bibitem[Z2]{Z2} A.~Zinger, 
{\it Completion of Katz-Qin-Ruan's Enumeration of Genus-Two Plane Curves}, 
J.~Algebraic Geom.~13 (2004), no.~3, 547-561.

\bibitem[Z3]{Z3} A.~Zinger,
{\it Enumeration of Genus-Two Curves with a Fixed Complex Structure
in $\PP$ and $\PPP$}, J.~Diff.~Geom.~65 (2003), no.~3, 341-467. 

\bibitem[Z4]{Z4} A.~Zinger, 
{\it Enumerative vs.~Symplectic Invariants and Obstruction Bundles},
J.~Symplectic Geom.~2 (2004), no.~4, 445--543.

\bibitem[Z5]{Z5} A.~Zinger,
{\it Counting Rational Curves of Arbitrary Shape in Projective Spaces}, 
Geom.~Top.~9 (2005), 571-697.

\bibitem[Z6]{Z6} A.~Zinger,
{\it Reduced Genus-One Gromov-Witten Invariants}, math/0507103.

\bibitem[Z7]{Z7} A.~Zinger,
{\it On the Structure of Certain Natural Cones 
over Moduli Spaces of Genus-One Holomorphic Maps},
Adv.~Math.~214 (2007) 878–933.

\bibitem[Z8]{Z8} A.~Zinger,
{\it The Reduced Genus-One Gromov-Witten Invariants of Calabi-Yau Hypersurfaces},
math/0705.2397.

\bibitem[Z9]{Z9} A.~Zinger,
{\it  Standard vs.~Reduced Genus-One Gromov-Witten Invariants},
math/0706.0715.





\end{thebibliography}
\end{document}